\documentclass[reqno,11pt]{amsart}

\usepackage{texdraw}
\usepackage{vmargin}
\setpapersize{A4}

\setmarginsrb{2.5cm}{1.2cm}{2.5cm}{1.9cm}{1.2cm}{1.2cm}{1.9cm}
{1.9cm}

\usepackage{amssymb}

\title[]{Analytic asymptotic expansions of the 
Reshetikhin--Turaev invariants of Seifert $3$--manifolds for 
$\SU(2)$}

\author{S\o ren Kold Hansen}

\address{Department of Mathematics, Kansas State University, 
Manhattan, KS 66506, USA}

\email{hansen@math.ksu.edu}

\newcommand{\ra}{\rightarrow}
\newcommand{\sm}{\setminus}

\newcommand{\mA}{\mathcal A}
\newcommand{\mC}{\mathcal C}
\newcommand{\mI}{\mathcal I}
\newcommand{\mM}{\mathcal M}
\newcommand{\mQ}{\mathcal Q}
\newcommand{\mR}{\mathcal R}
\newcommand{\mS}{\mathcal S}

\newcommand{\frsl}{\mathfrak sl}
\newcommand{\frg}{\mathfrak g}

\newcommand{\frD}{\mathfrak D}

\newcommand{\C}{\mathbb C}

\newcommand{\N}{\mathbb N}
\newcommand{\Q}{\mathbb Q}
\newcommand{\R}{\mathbb R}
\newcommand{\Sec}{\mathbb S}
\newcommand{\Z}{\mathbb Z}
\newcommand{\B}{ \{ \pm 1 \}^{n}}
\newcommand{\Bzeta}{\frac{B}{2\pi} + \frac{1}{2r}\zeta}

\newcommand{\ep}{\epsilon}
\newcommand{\vep}{\varepsilon}
\newcommand{\sdel}{\sqrt{\delta}}
\newcommand{\svep}{\sqrt{\varepsilon}}
\newcommand{\livep}{\lim_{\varepsilon \rightarrow 0_{+}}}
\newcommand{\lixi}{\lim_{\xi \rightarrow 0_{+}}}
\newcommand{\lidel}{\lim_{\delta \rightarrow 0_{+}}}
\newcommand{\lieta}{\lim_{\eta \rightarrow 0_{+}}}
\newcommand{\iny}{\int_{-\infty}^{\infty}}

\newcommand{\Arg}{\text{\rm Arg}}
\newcommand{\CS}{\text{\rm CS}}
\newcommand{\dS}{\text{\rm S}}
\newcommand{\dte}{\text{\rm d}}
\newcommand{\e}{\text{\rm e}}

\newcommand{\F}{\text{\rm F}}

\newcommand{\Hom}{\text{\rm Hom}}
\newcommand{\id}{\text{\rm id}}
\newcommand{\ima}{\text{\rm Im}}
\newcommand{\inte}{\text{\rm int}}
\newcommand{\intpolar}{\text{\rm int,polar}}
\newcommand{\K}{\text{\rm K}}
\newcommand{\lcm}{\text{\rm lcm}}
\newcommand{\ns}{\text{\rm n}}
\newcommand{\os}{\text{\rm o}}
\newcommand{\tO}{\text{\rm O}}
\newcommand{\pol}{\text{\rm polar}}
\newcommand{\rea}{\text{\rm Re}}
\newcommand{\Res}{\text{\rm Res}}
\newcommand{\s}{\text{\rm s}}
\newcommand{\sd}{\text{\rm sd}}
\newcommand{\sign}{\text{\rm sign}}
\newcommand{\SL}{\text{\rm SL}}
\newcommand{\Sp}{\text{\rm Sp}}
\newcommand{\spec}{\text{\rm spec}}
\newcommand{\sqr}{\text{\rm sqr}}
\newcommand{\st}{\text{\rm st}}
\newcommand{\SU}{\text{\rm SU}}
\newcommand{\Sym}{\text{\rm Sym}}
\newcommand{\tr}{\text{\rm tr}}
\newcommand{\tC}{\text{\rm C}}
\newcommand{\tCW}{\text{\rm CW}}
\newcommand{\tCLW}{\text{\rm CLW}}
\newcommand{\tu}{\text{\rm u}}
\newcommand{\tU}{\text{\rm U}}

\newcommand{\ualpha}{\underline{\alpha}}
\newcommand{\ubeta}{\underline{\beta}}
\newcommand{\uep}{\underline{\epsilon}}
\newcommand{\um}{\underline{m}}
\newcommand{\umu}{\underline{\mu}}
\newcommand{\un}{\underline{n}}
\newcommand{\uone}{\underline{1}}
\newcommand{\uzero}{\underline{0}}

\newcommand{\musum}{\sum_{\underline{\mu} \in \{ \pm 1 \}^{n} }}

\newcommand{\mumsum}{\sum_{\underline{\mu}' \in \{ \pm 1 \}^{n} }}

\newcommand{\pmualphasum}{\left( \sum_{j=1}^{n} \frac{\mu_{j}}{\alpha_{j}} \right)}

\newcommand{\mualphasum}{\sum_{j=1}^{n} \frac{\mu_{j}}{\alpha_{j}}}

\newcommand{\rhoalphasum}{\sum_{j=1}^{n} \frac{\rho_{j}}{\alpha_{j}}}

\newcommand{\muprod}{\left( \prod_{j=1}^{n} \mu_{j} \right)}

\newcommand{\msumc}{\sum_{\stackrel{m \in \Z}{m \geq \sum_{j=1}^{n} \frac{\mu_{j}' n_{j}'}{\alpha_{j}}}}}

\newcommand{\msumd}{\sum_{\stackrel{m \in \Z}{z_{\st} \in \R \sm \Z}}}

\newcommand{\msume}{\sum_{\stackrel{m \in \Z}{m \geq \sum_{j=1}^{n} \frac{n_{j}'}{\alpha_{j}}}}}

\newcommand{\msumi}{\sum_{\stackrel{m \in \Z}{z_{\st}(m,\lambda) \in \Z}}}

\newcommand{\msumj}{\sum_{\stackrel{m \in \Z}{m +\frac{B}{2\pi}+\frac{A}{\pi}l > \frac{|A|}{\pi}\eta}}}

\newcommand{\msumk}{\sum_{\stackrel{m \in \Z}{m +\frac{B}{2\pi}+\frac{A}{\pi}l \geq 0}}}

\newcommand{\msuml}{\sum_{\stackrel{m \in \Z}{m - V(x) \geq 0}}}

\newcommand{\msumml}{\sum_{\stackrel{m \in \Z}{m - V(\nu) \geq 0}}}

\newcommand{\mualphammsum}{\sum_{j=1}^{n} \frac{\mu_{j}' n_{j}'}{\alpha_{j}}}

\newcommand{\pmualphammsum}{\left( \sum_{j=1}^{n} \frac{\mu_{j}' n_{j}'}{\alpha_{j}} \right)}

\newcommand{\nmualphammsum}{\left| \sum_{j=1}^{n} \frac{\mu_{j}' n_{j}'}{\alpha_{j}} \right|}

\newcommand{\sgmualphammsum}{{\text{\rm sign}} \left( \sum_{j=1}^{n} \frac{\mu_{j}' n_{j}'}{\alpha_{j}} \right) }

\newcommand{\mualphammsuma}{\sum_{j=1}^{n} \frac{\mu_{j}' n_{j}'}{\alpha_{j}} + \frac{1}{2r} \sum_{j=1}^{n} \frac{\mu_{j}}{\alpha_{j}} }

\newcommand{\nsum}{\sum_{\underline{n}=\underline{0}}^{\underline{\alpha}-\underline{1}}}

\newcommand{\nalphasuma}{\sum_{j=1}^{n} \frac{n_{j}}{\alpha_{j}}}

\newcommand{\nmalphasum}{\sum_{j=1}^{n} \frac{n_{j}'}{\alpha_{j}}}

\newcommand{\pnmalphasum}{\left( \sum_{j=1}^{n} \frac{n_{j}'}{\alpha_{j}} \right)}

\newcommand{\nnmalphasum}{\left| \sum_{j=1}^{n} \frac{n_{j}'}{\alpha_{j}} \right|}

\newcommand{\sgnmalphasum}{{\text{\rm sign}} \left( \sum_{j=1}^{n} \frac{n_{j}'}{\alpha_{j}} \right) }

\newcommand{\lsuma}{\sum_{\stackrel{l \in \Z}{\frac{B}{2\pi}+\frac{A}{\pi}l \in \Z}}}

\newcommand{\lambdasuma}{\sum_{\stackrel{\lambda \in \Lambda}{\frac{B}{2\pi}+\frac{A}{\pi}l \in \Z}}}

\newcommand{\msyma}{\frac{1}{\prod_{j=1}^{n} {\text{\rm Sym}}_{\Z_{\pm}} \left( \frac{n_{j}'}{\alpha_{j}} \right) }}

\newcommand{\sgE}{{\text{\rm sign}}(E)}
\newcommand{\sgA}{{\text{\rm sign}}(A)}
\newcommand{\sgV}{{\text{\rm sign}}(V)}

\newcommand{\eksgA}{e^{\frac{i \pi}{4} {\text{\rm sign}} (A)}}

\newcommand{\gex}{\exp \left( 2 \pi i \sum_{j=1}^{n} \frac{\rho_{j}}{\alpha_{j}} [ r n_{j}^{2} + \mu_{j} n_{j} ] \right)}

\newcommand{\mgex}{\exp \left[ 2 \pi i \sum_{j=1}^{n} \mu_{j} \left( \frac{\rho_{j}}{\alpha_{j}} n_{j}' - \frac{1}{2} \sigma_{j}l \right) \right]}
\newcommand{\mmugex}{\exp \left[ 2 \pi i \sum_{j=1}^{n} \mu_{j} \left( \frac{\rho_{j}}{\alpha_{j}} \mu_{j}' n_{j}' - \frac{1}{2} \sigma_{j}l \right) \right]}

\newtheorem{thm}{Theorem}[section]
\newtheorem{conj}[thm]{Conjecture}
\newtheorem{cor}[thm]{Corollary}
\newtheorem{lem}[thm]{Lemma}
\newtheorem{prop}[thm]{Proposition}

\theoremstyle{definition}

\newtheorem{defi}[thm]{Definition}

\newtheorem{rem}[thm]{Remark}

\newcommand{\refthm}[1]{Theorem~\ref{#1}}
    
    \newcommand{\refconj}[1]{Conjecture~\ref{#1}}
    \newcommand{\refcor}[1]{Corollary~\ref{#1}}
    \newcommand{\refdefi}[1]{Definition~\ref{#1}}
    
    \newcommand{\reflem}[1]{Lemma~\ref{#1}}
    \newcommand{\refprop}[1]{Proposition~\ref{#1}}
    \newcommand{\refrem}[1]{Remark~\ref{#1}}

\newcommand{\HS}{\noindent \hfill{$\square$}}

\newenvironment{prf}{\mbox{ } \newline {\bf Proof. }}{\noindent \hfill{$\square$} \newline}
\newenvironment{prfa}{\mbox{ } \newline {\bf Proof }}{\noindent \hfill{$\square$} \newline}

\begin{document}

\maketitle

\begin{abstract}
We calculate the large quantum level asymptotic expansion of 
the RT--invariants associated to $\SU(2)$ of all oriented 
Seifert $3$--manifolds $X$ with orientable base or 
non-orientable base with even genus. Moreover, we identify the 
Chern--Simons invariants of flat $\SU(2)$--connections on $X$ 
in the asymptotic formula thereby proving the so-called 
asymptotic expansion conjecture (AEC) due to J.\ E.\ Andersen 
\cite{Andersen1}, \cite{Andersen2} for these manifolds. For
the case of Seifert manifolds with base $S^2$ we actually prove
a little weaker result, namely that the asymptotic formula has
a form as predicted by the AEC but contains some extra terms
which should be zero according to the AEC. We prove that these
`extra' terms are indeed zero if the number of exceptional 
fibers $n$ is less than 4 and conjecture that this is also the
case if $n \geq 4$.
For the case of Seifert fibered rational homology spheres we 
identify the Casson--Walker invariant in the asymptotic formula. 

Our calculations demonstrate a general method for calculating 
the large $r$ asymptotics of a finite sum 
$\Sigma_{k=1}^{r} f(k)$, where $f$ is a meromorphic function
depending on the integer parameter $r$ and satisfying certain 
symmetries. Basically the method, which is due to Rozansky
\cite{Rozansky1}, \cite{Rozansky3}, is based on a limiting 
version of the Poisson summation formula together with an 
application of the steepest descent method from asymptotic 
analysis.
\end{abstract}

\tableofcontents

\section{Introduction}\label{sec-Introduction}

\noindent In this paper we investigate the large quantum level 
asymptotics of the Reshetikhin--Turaev invariants of the 
Seifert manifolds. Here and elsewhere a $3$--manifold is 
closed and oriented unless otherwise stated. In particular, a 
Seifert manifold is an oriented Seifert manifold.

It is now a well established fact that the Reshetikhin--Turaev 
approach \cite{ReshetikhinTuraev}, \cite{Turaev} leads to a 
family of $2+1$--dimensional topological quantum field theories
(TQFTs) indexed by a simply connected compact simple Lie group 
$G$ and an integer $r$ (the shifted quantum level) bigger than 
or equal to the dual Coxeter number $h^{\vee}$ of $G$. In 
particular, we have invariants $\tau_{r}^{G}$ of the 
$3$--manifolds called the quantum $G$--invariants or the 
RT--invariants associated to $G$.

There are several other approaches to these quantum invariants 
and their underlying TQFTs, see e.g.\ \cite{Turaev}, 
\cite{BakalovKirillov} and references therein. For the 
following discussion we note that J.\ E.\ Andersen and K.\ Ueno
\cite{AndersenUeno1} have recently constructed a family of 
TQFTs using ideas from conformal field theory, notably the 
works \cite{TsuchiyaUenoYamada}, \cite{Ueno},
\cite{KawamotoNamikawaTsuchiyaYamada}. They have so far proved 
that these TQFTs for the groups $\SU(n)$ coincide with the 
TQFTs of Reshetikhin and Turaev associated to these Lie groups,
cf.\ \cite{AndersenUeno2}. Via work of Laszlo \cite{Laszlo}
one can combine the work of Andersen and Ueno with work of 
Axelrod, Della Pietra and Witten \cite{AxelrodPietraWitten}, 
Hitchin \cite{Hitchin}, Faltings \cite{Faltings} and others on
the geometric quantization of the moduli space of flat 
connections on surfaces to obtain an alternative gauge 
theoretic approach to the TQFTs of Andersen and Ueno. We refer 
to \cite{Atiyah} and \cite{Andersen3} and references therein 
for details on this gauge theoretic approach which was 
originally outlined by Witten in \cite{Witten}.

The main aim of this paper is to prove the following conjecture
for Seifert manifolds for $G=\SU(2)$.

\begin{conj}[J.\ E.\ Andersen \cite{Andersen1}, 
\cite{Andersen2}. The asymptotic expansion conjecture (AEC)]
\label{conj:AEC'}
\footnote{We give the conjecture in a slightly different form 
than in \cite{Andersen2}, where $\theta_{l}=l$.}
For any $3$--manifold $X$ there exist constants $D_{j} \in \Q$ 
and $b_{j} \in  \C$ for $j=0,1, \ldots, n$ and 
$a_{l}^{j} \in  \C$ for $j=0,1, \ldots, n$ and $l=1,2,\ldots$ 
such that the asymptotic expansion of $\tau_{r}^{G}(X)$ in the 
limit $r \ra \infty$ is given by
\begin{equation}\label{eq:AEC'}
\tau_{r}^{G}(X) \sim \sum_{j=0}^{n} e^{2\pi i r q_{j}}r^{D_{j}}
b_{j}
\left( 1 + \sum_{l=1}^\infty a_{l}^{j} r^{-\theta_{l}}\right),
\end{equation}
where $q_{0} = 0, q_{1}, \ldots, q_{n}$ are the finitely many 
different values of the Chern--Simons functional on the moduli 
space of flat $G$--connections on $X$, and $\{\theta_{l}\}$ is
a strictly increasing sequence in $\frac{1}{b}\Z_{>0}$, where 
$b$ is the least positive number such that $bD_{j} \in \Z$ for
$j=0,1,\ldots,n$.

Here {\bf $\sim$} means that there for all non-negative $N$ is 
a $C_{N} \in \R$ such that
$$
\left| \tau_{r}^{G}(X) - \sum_{j=0}^{n} e^{2\pi i r q_{j}} 
r^{D_{j}}b_{j}\left( 1 
+ \sum_{l=1}^{N} a_{l}^{j} r^{-\theta_{l}}\right) \right| 
\leq C_{N} r^{D-\theta_{N+1}}
$$
for all $r \in \Z_{\geq h^{\vee}}$, where 
$D = \max\{D_{0},\ldots,D_{n}\}$.
\end{conj}

The formulation above seems a little awkward from the point of 
view of asymptotic analysis. We shall later on see a more
compact formulation of the conjecture, cf.\ 
\refconj{conj:AEC}. The main reason for the above formulation 
is to expose the quantities $D_{j}$ and $b_{j}$ which are 
believed to have interpretations in terms of certain 
topological/geometric quantities, cf.\ \refconj{conj:d'} 
below and \cite{Andersen2}.

Let us give a review of some of the results leading to this 
conjecture. We want by this review to stress what is known 
today from a rigorous mathematical point of view and what is
based on arguments using techniques from theoretical
physics not yet well understood from a mathematical point of 
view. In 1988 E.\ Witten \cite{Witten} proposed new invariants
$Z_{k}^{G}(X,L)$ of a $3$--manifold $X$ with an embedded 
(colored) link $L$, which generalize the famous Jones 
polynomial \cite{Jones} of knots in the $3$--sphere. Here $G$ 
is a nice Lie group as above and $k$ is a non-negative integer 
called the (quantum) level. In case there is no link present, 
$Z_{k}^{G}(X)$ is given (formally) by the following path 
integral over the infinite dimensional space of gauge 
equivalence classes of $G$--connections on $X$
\begin{equation}\label{eq:Witten-invariant}
Z^{G}_{k}(X) = \int e^{2\pi \sqrt{-1} k\CS (A)} \frD A,
\end{equation}
where
\begin{equation}\label{eq:CSfunctional}
\CS (A) = \frac{1}{8\pi^{2}} \int_{X}
\tr\left(A \wedge dA + \frac{2}{3} A \wedge A \wedge A \right)
\end{equation}
is the Chern--Simons functional associated to $G$. It is 
well-known that $\CS$ considered as a function into $\R/\Z$ is 
invariant under gauge transformations if the Ad-invariant inner
product $\tr$ on $\frg$ (the Lie algebra of $G$) is normalized 
properly, see e.g.\ \cite[Sect.~2]{Freed}. Although no rigorous
mathematical definition of the path integral
(\ref{eq:Witten-invariant}) has been given yet, see 
\cite[Sect.~20.2.A]{JohnsonLapidus} for some comments, Witten 
was able using path integral techniques to give a recipe for 
calculating the invariant $Z_{k}^{G}(X)$ via a surgery 
description of $X$. These ideas enabled people to obtain 
explicit (and rigorous) formulas for the invariants 
$Z_{k}^{G}(X)$ for some classes of manifolds $X$ such as lens 
spaces \cite{FreedGompf}, \cite{Garoufalidis1}, \cite{Jeffrey1}
and more generally of Seifert manifolds with orientable base 
\cite{FreedGompf}, \cite{Garoufalidis1}, \cite{Rozansky1},
\cite{Rozansky3}. 

In quantum field theory it is natural to study Feynman path 
integrals as in (\ref{eq:Witten-invariant}) and a main approach
to understand them is via their behaviour in the large $k$ 
limit. By using stationary phase approximation techniques 
together with path integral arguments, Witten was able 
\cite{Witten} to express the leading large $k$ asymptotics 
of $Z_{k}^{G}(X)$ (also called the semiclassical approximation)
as a sum over the set $\mM$ of stationary points of the 
Chern--Simons functional $\CS$, the summand being given by such
invariants as Chern--Simons invariants, spectral flows and 
Reidemeister torsions. Here $\mM$ is the moduli space of gauge 
equivalence classes of flat connections in the trivial 
$G$--bundle over $X$ (and $\mM$ is assumed discrete hence 
finite in Witten's derivation). Freed and Gompf 
\cite{FreedGompf} and Jeffrey \cite{Jeffrey1} suggested 
different refinements to Witten's semiclassical approximation. 
In particular Jeffrey allowed $\mM$ to have smooth components 
of non-zero dimension, and in these cases the sum over $\mM$
should be replaced by some integral over $\mM$. Let us give 
some details. Let $A$ be a flat $G$--connection on $X$ and 
consider the elliptic complex
$d_{A}: \Omega^{*}(X,\frg) \to \Omega^{*+1}(X,\frg)$, where 
$d_{A}f=df+[A \wedge f]$ is the covariant derivative in the 
adjoint representation. Let $h_{A}^{i}$ be the dimension of the
$i$th cohomology group $H^{i}(X,d_{A})$ of this complex. If 
$\mM$ is discrete and the covariant derivative complex is 
acyclic for all $A$ (like in Witten's derivation) then the 
conjectured formula for the semiclassical approximation of 
$Z_{k}^{G}(X)$ is
\begin{eqnarray}\label{eq:semiclassical-G'}
Z^{G}_{k}(X) &\sim_{k \ra \infty}& \frac{1}{|Z(G)|}
e^{-\pi \sqrt{-1}\dim(G)(1+b_{1}(X))/4} \\
&& \times \sum_{A \in \mM} 
e^{2\pi \sqrt{-1} r CS(A)}r^{D_{A}}
e^{-2\pi\sqrt{-1}(I_{A}/4+B_{A})} \sqrt{\tau_{X}(A)}, \nonumber
\end{eqnarray}
where $r=k+h^{\vee}$, $Z(G)$ is the center of $G$ and 
$b_{1}(X)$ is the first betti number of $X$. Moreover, 
$\tau_{X}(A) \in \R_{+}$ is the Reidemeister torsion of the 
complex $(\Omega^{*}(X,\frg),d_{A})$ and $I_{A} \in \Z/8\Z$ is 
the spectral flow of the family of operators
$
\left( \begin{array}{cc}
\ast d_{A_{t}} & -d_{A_{t}} \ast \\
d_{A_{t}} \ast & 0
\end{array}
\right)
$ 
on $\Omega^{1}(X;\frg) \otimes \Omega^{3}(X;\frg)$, where
$A_{t}$ is a path in the affine space of $G$--connections on 
$X$ from the trivial connection to the flat connection $A$. 
Under the above assumptions $B_{A}=D_{A}=0$ for all $A$, but 
the factor $r^{D_A}e^{-2\pi\sqrt{-1}B_A}$ is included for the 
following discussion. (Here and elsewhere we follow the usual 
practice by denoting a connection and its gauge equivalence 
class by the same symbol.)

If we allow some of the $A$'s to be reducible (implying that 
$h_{A}^{0}>0$) (still assuming that $\mM$ is finite) then it is
conjectured that the above formula holds with 
$D_{A}=-h_{A}^{0}/2$ and $B_{A}=h_{A}^{0}/8$. In particular, 
these reducible connections do not contribute to the leading 
large $k$ asymptotics of $Z^{G}_{k}(X)$ (which is given by the 
sum of terms with $D_{A}=\max_{A \in \mM}\{D_{A}\}$). If the 
gauge equivalence class of a flat connection $A$ is 
non-isolated in $\mM$ then the situation becomes more 
complicated. It is conjectured that we in that case have
$D_{A}=(h_{A}^{1}-h_{A}^{0})/2$ and 
$B_{A}=(h_{A}^{0}+h_{A}^{1})/8$, and that the above sum over 
$\mM$ should be replaced by some integral over $\mM$. In cases,
where $A$ is a non-acyclic connection, the interpretation of 
the factor $\sqrt{\tau_{X}(A)}$ is not clear, $A$ being called 
acyclic if $(\Omega^{*}(X,\frg),d_{A})$ is acyclic. (One can 
still define a torsion $\in \R_{+}$ of the covariant derivative
complex once a basis of $H^{*}(X,d_{A})$ has been fixed, but 
there seem to be some normalization problems even in the case 
of isolated reducible points in $\mM$.) In the case of 
non-isolated points in $\mM$ the factor $\sqrt{\tau_{X}(A)}$ 
should somehow be interpreted as a density function giving the 
needed measure on $\mM$. We refer to \cite{Jeffrey1} and 
\cite{Andersen2} for some comments on this and also to
\cite{Rozansky1}, \cite{Rozansky2} for some futher comments on 
the semiclassical approximation of $Z_k^G(X)$. We note that the
Chern--Simons functional is constant on the connected 
components of $\mM$, so in all cases one gets a (conjectured)
formula for the leading large $k$ asymptotics of $Z_k^G(X)$ of 
the form
\begin{equation}\label{eq:semiclassical-G''}
Z^{G}_{k}(X) \sim_{k \ra \infty} \sum_{j=1}^{m} 
e^{2\pi \sqrt{-1} r \CS_{j}} \int_{\mM_j} 
r^{(h_{A}^{1}-h_{A}^{0})/2}
e^{-2\pi\sqrt{-1}(I_{A}/4+(h_{A}^{0}+h_{A}^{1})/8)} f(A),
\end{equation}
where the sum is over the connected components of $\mM$, 
$\CS_j$ is the Chern--Simons invariant of the elements in the 
$j$th component $\mM_j$ of $\mM$ and the integral over $\mM_j$ 
is w.r.t.\ some density function $f(A)$ generalizing the factor
$\frac{1}{|Z(G)|}e^{-\pi \sqrt{-1}\dim(G)(1+b_{1}(X))/4}
\sqrt{\tau_{X}(A)}$
in (\ref{eq:semiclassical-G'}).

A main first observation is that the expression in the 
right-hand side of (\ref{eq:semiclassical-G'}) is mathematical
rigorous in case $\mM$ is finite (although it is not clear what
the meaning of the factor $\sqrt{\tau_X(A)}$ should be in the
non-acyclic case). Moreover, as stated above, one can in
principal obtain a rigorous expression for the invariants 
$Z_k^G(X)$ via a surgery description of $X$. Thus a kind of 
verification of (\ref{eq:semiclassical-G'}) should be possible.
This program was first carried out partly by Freed and Gompf 
\cite{FreedGompf} presenting a large amount of computer 
calculations for the $\SU(2)$--invariants of lens spaces and 
some $3$--fibered Seifert manifolds, and about the same time 
by Jeffrey \cite{Jeffrey1} and Garoufalidis 
\cite{Garoufalidis1} who independently gave exact calculations 
of the semiclassical approximation of the $\SU(2)$--invariants 
of lens spaces, starting from explicit expressions for the 
invariants. Jeffrey also verified parts of the conjecture 
(\ref{eq:semiclassical-G'}) for $G$ arbitrary and $X$ belonging
to a class of mapping tori of the torus. Garoufalidis verified 
parts of the conjecture (\ref{eq:semiclassical-G'}) for 
$G=\SU(2)$ and $X$ any Seifert fibered integral homology 
sphere. 

The perturbative (or asymptotic) expansion of the Chern--Simons
path integral $Z_{k}^{G}(X)$ is given by the semiclassical 
approximation and its so-called higher loop correction terms.
This leads together with the semiclassical approximation 
(\ref{eq:semiclassical-G''}) to a conjecture for the (full) 
perturbative expansion of $Z^{G}_{k}(X)$ in the large $k$ 
limit, see \cite[Conjecture 7.6]{Andersen2} for a version of 
this conjecture. Roughly speaking one should obtain an 
expansion given by the right-hand side of
(\ref{eq:semiclassical-G''}) with each integrand being 
multiplied by an asymptotic series of the form 
$1+\sum_{l=1}^{\infty} c_{l}^{A} r^{-l}$ with the $c_{l}^{A}$'s
given in terms of certain contributions of Feynman diagrams 
determined by the Feynman rules of the 
Chern--Simons theory.

S.\ Axelrod and I.\ M.\ Singer have considered the higher loop
contributions, cf.\ \cite{AxelrodSinger1},
\cite{AxelrodSinger2}, \cite{Axelrod}. Actually they were able 
to deduce rigorously defined formulas for the coefficients 
$c_{l}^{A}$ in case $A$ is an acyclic point in $\mM$ or an 
element in a smooth component of $\mM$. Moreover, they showed 
that their $c_{l}^{A}$ define topological invariants of 
$(X,A)$. There are, however, major difficulties which have not 
yet been worked out. We refer to the above works of Axelrod and
Singer and \cite{Andersen2} for more comments and details.

It is generally believed that the TQFTs of Reshetikhin and 
Turaev are a mathematical realization of Witten's TQFTs, and
one can say that the AEC, \refconj{conj:AEC'}, offers in a 
sense a converse point of view to the above works on the 
perturbative expansion of the Witten invariants, where one 
seeks to derive the final output of perturbation theory after 
all cancellations have been made (i.e.\ collect all terms with 
the same Chern--Simons value). We stress that, although the 
ideas leading to the AEC stem from Witten's approach, the AEC 
is a conjecture for the rigorously defined RT--invariants, and 
it is completely independent of Witten's Chern--Simons path 
integral approach. Seen in the above light we should also 
expect the following slightly different version of the AEC,
where we have seperated terms coming from different components
of $\mM$ (even if there are different components with the same 
Chern--Simons value).

\begin{conj}[J.\ E.\ Andersen \cite{Andersen1}, 
\cite{Andersen2}. The asymptotic expansion conjecture (AEC)]
\label{conj:AEC''}

For any $3$--manifold $X$ there exist constants 
$\tilde{D}_{j} \in \Q$, $I_{j} \in \Q \bmod{8}$, 
$v_{j} \in  \R_{+}$ for $j=1,2,\ldots, m$ and 
$\tilde{a}_{j}^{l} \in  \C$ for $j=1,2,\ldots, m$ and
$l=1,2,\ldots$ such that the asymptotic expansion of 
$\tau_{r}^{G}(X)$ in the limit $r \ra \infty$ is given by
$$
\tau_{r}^{G}(X) \sim \sum_{j=1}^{m} e^{2\pi i r q_{j}}
r^{\tilde{D}_{j}} v_{j}e^{\frac{\pi i}{4}I_{j}}
\left( 1 + \sum_{l=1}^\infty \tilde{a}_{l}^{j} r^{-l}\right),
$$
where the $j$-sum is over the connected components of the 
moduli space of flat $G$--connections on $X$ and $q_{j}$ is the
value of the Chern--Simons functional on the component indexed 
by $j$.

Here {\bf $\sim$} means that there for all non-negative $N$ is 
a $C_{N} \in \R$ such that
$$
\left| \tau_{r}^{G}(X) - \sum_{j=1}^{m} e^{2\pi i r q_{j}} 
r^{\tilde{D}_{j}} v_{j}e^{\frac{\pi i}{4}I_{j}}
\left( 1 + \sum_{l=1}^{N} \tilde{a}_{l}^{j} r^{-l} \right) 
\right| \leq C_{N} r^{D-N-1}
$$
for all $r \in \Z_{\geq h^{\vee}}$, where 
$D = \max\{\tilde{D}_{1},\ldots,\tilde{D}_{m}\}$.
\end{conj}

We note that \refconj{conj:AEC'} follows from 
\refconj{conj:AEC''}, and that $D$ is equal in the two versions
of the conjecture. The AEC concerns (analytic) asymptotic 
expansions of (a class of) complex functions defined on the 
positive integers. As noted by Andersen \cite{Andersen1}, 
\cite{Andersen2}, such a function has at most one asymptotic 
expansion of the form (\ref{eq:AEC'}) (up to some trivialities 
such as permutations of the terms in the asymptotic formula). 
(Here it is of course important that the $q_{j}$'s are mutually
different.) This means that if the AEC is true, then the 
$D_{j}$'s, $b_{j}$'s and the $a_{l}^{j}$'s in (\ref{eq:AEC'}) 
are all uniquely determined by the function
$r \mapsto \tau_{r}^{G}(X)$, hence they are also topological 
invariants of $X$. In particular Andersen has proposed the 
following conjecture (compare with 
(\ref{eq:semiclassical-G'}) and (\ref{eq:semiclassical-G''})):

\begin{conj}[J.\ E.\ Andersen \cite{Andersen1}, 
\cite{Andersen2}. Topological interpretation of the 
$D_{j}$'s]\label{conj:d'}
Let $\mM_{j}$ be the union of the components of the moduli 
space of flat $G$--connections on $X$ having Chern--Simons 
value $q_{j}$. Then
\begin{equation}\label{eq:dconj}
D_{j} = \frac{1}{2} 
\max_{A \in \mM_{j}} \left( h_{A}^{1} - h_{A}^{0} \right),
\end{equation}
where $\max$ here means the maximum value of 
$h_{A}^{1} - h_{A}^{0}$ on a Zariski open subset of $\mM_{j}$.
\end{conj}

The results in \cite{Andersen1} showed the importance of being 
careful about the interpretation of the $\max$ in this 
conjecture. The $\tilde{a}_{l}^{j}$'s, 
$v_{j}$'s, $\tilde{I}_{j}$'s and $\tilde{D}_{j}$'s in 
\refconj{conj:AEC''} are also expected to be topological 
invariants. In particular, $\tilde{D}_{j}$ should be given by
the right-hand side of (\ref{eq:dconj}) with $\mM_{j}$ now 
being the $j$th component of the moduli space.

\refconj{conj:d'} implies that the sequence $\{\theta_{l}\}$ in
(\ref{eq:AEC'}) is the sequence of positive half-integers 
$\frac{1}{2}\Z_{>0}$. This conjecture gives (in combination 
with \refconj{conj:AEC'}) a strong connection between the 
fundamental group of $X$ and the quantum invariants of $X$, 
recalling that the moduli space of flat $G$--connections on $X$
is in bijection with $\Hom(\pi_{1}(X),G)/G$ (for $X$ 
connected). Garoufalidis \cite{Garoufalidis2} and Andersen 
\cite{Andersen4} have e.g.\ pointed out that the above growth 
rate in $r$ for the invariants $\tau^{\SU(n)}_{r}$ leads to the
fact that the colored Jones polynomials associated to $\SU(n)$
(for all colorings and all $n$) can distinguish the unknot 
from all other knots. As also pointed out by Andersen 
\cite{Andersen2}, \refconj{conj:d'} implies together with a 
recent result of Kronheimer and Mrowka \cite{KronheimerMrowka}
that the colored Jones polynomials associated to $\SU(2)$ 
alone can distinguish the unknot from all other knots.
 
Andersen was the first to give a complete proof of the AEC for 
a class of $3$--manifolds. In fact he proved in 
\cite{Andersen1} the AEC (together with \refconj{conj:d'}) for 
the mapping tori of finite order diffeomorphisms of orientable 
surfaces of genus at least $2$ using the gauge theoretic 
approach to the quantum invariants. These manifolds are 
Seifert manifolds with orientable base and Seifert Euler number
equal to zero. Andersen proved the conjectures for $G$ an 
arbitrary simply connected compact simple Lie group.

As already stated Jeffrey \cite{Jeffrey1} and Garoufalidis 
\cite{Garoufalidis1} made completely rigorous calculations of 
the semiclassical approximation of the $\SU(2)$--invariants of 
lens spaces. Actually these calculations contain a complete
verification of the AEC for the lens spaces for $\SU(2)$. 
Moreover, they confirm \refconj{conj:d'} in this case. The 
paper \cite{Jeffrey1} also contains a proof of the AEC for a 
certain class of mapping tori over the torus (for the 
invariants associated to an arbitrary simply connected compact 
simple $G$). Except for the identification of Chern--Simons 
invariants, the AEC has recently been proved in 
\cite{HansenTakata} for all lens spaces and $G$ an arbitrary
simply connected compact simple Lie group. 

In \cite{Rozansky1}, \cite{Rozansky3} Rozansky calculated the 
Witten $\SU(2)$--invariants of all Seifert manifolds with 
orientable base and carried through a rather technical analysis
leading to a candidate for the full asymptotic expansion of 
these invariants. However, the error estimates required to 
prove that his calculations lead to the asymptotics of the 
invariants are missing. In this paper we establish those error
etsimates thereby proving that Rozansky's calculations in
\cite{Rozansky3} really leads to the large $r$ asymptotic
expansion of the quantum $\SU(2)$--invariants of the Seifert
manifolds. In \cite[Sect.~8]{Hansen2} it was proved that the 
Witten $\SU(2)$--invariants of the Seifert manifolds calculated
by Rozansky are equal to the RT--invariants associated to 
$\SU(2)$.

For a certain subclass of the Seifert manifolds Lawrence and 
Rozansky \cite{LawrenceRozansky} have obtained a much shorter 
calculation of the asymptotic expansion of the 
$\SU(2)$--invariants. In this calculation one does not have to
apply the method of this paper to establish the required error
estimates. These estimates directly follow from the method 
without futher work. The method (as it stands) only works for 
the Seifert manifolds satisfying that the first coordinates of 
the oriented Seifert invariants are mutually coprime (the 
oriented Seifert invariants being pairs of coprime integers, 
one for each singular fiber). In \cite{LawrenceRozansky} the 
result is only stated for the Seifert manifolds with base 
$S^{2}$, but their method works for all Seifert manifolds with 
orientable base or nonorientable base with even genus 
satisfying the above coprime condition. However, it does not 
seem to be easy to extend this method to cases where the 
coprime condition is not satisfied. 

In \cite{Yoshida2} Yoshida defines a family of invariants
$Y_{k}^{\SU(2)}$, $k \in \Z_{\geq 0}$, of $3$--manifolds via 
the abelianization of the $\SU(2)$ WZW model obtained in
\cite{Yoshida1}. Moreover, he calculates the semiclassical 
approximation of his invariants for the integral homology 
spheres satisfying that the moduli space $\mM$ of flat 
$\SU(2)$--connections is finite set of acyclic points. This 
calculation includes a complete description of the 
semiclassical formula in terms of Chern--Simons invariants, 
spectral flows and Reidemeister torsions. In fact, he finds 
that the semiclassical approximation of his invariants for 
these manifolds is given by the right-hand side of 
(\ref{eq:semiclassical-G'}) with $G=\SU(2)$ (and $B_A =D_A=0$ 
for all $A$). It is expected (but not yet proved) that the 
invariant $Y_{k}^{\SU(2)}$ is equal to the invariant 
$\tau_{k+2}^{\SU(2)}$. 

Let us finally mention that we via recent private communication
have learned that J.\ E.\ Andersen for the groups $G=\SU(n)$ 
has proved Conjectures~\ref{conj:AEC'} and \ref{conj:d'} for 
all $3$--manifolds via the gauge theoretic approach, see
\cite{Andersen3} for some comments. \footnote{It should thus be
a theorem that the family of colored Jones polynomials 
associated to $\SU(2)$ can distinguish the unknot from all 
other knots.} The proof involves asymptotics of Hitchin's 
connection over Teichm{\"u}ller space, approximations to all 
orders, of the boundary states of handle bodies and techniques 
similar to the ones presented in \cite{Andersen3}. Where 
Andersen works with the gauge theoretic definition of the 
quantum invariants, we work with the definition of Reshetikhin 
and Turaev and our proof of the AEC for the Seifert manifolds
is very different from Andersen's general proof.

Next let us give a review of the content of this paper. The 
main result is

\begin{thm}\label{thm:AEC-Seifert-manifolds}
The {\em AEC, \refconj{conj:AEC'}}, is true for $G=\SU(2)$ and 
$X$ any Seifert manifold with orientable base of positive genus
or non-orientable base with even genus. 
If $X$ is a Seifert manifold with base $S^{2}$ and $n$ 
exceptional fibers, then the AEC holds if $n \leq 3$ and in
case $n>3$ the large $r$ asymptotic expansion of 
$\tau_{r}(X)$ has the form {\em (\ref{eq:AEC'})} with a set 
$\{q_{j}\}_{j=0}^{n}$ containing the image set of 
the Chern--Simons functional on the moduli space 
of flat $\SU(2)$--connections on $X$. In all cases the sequence
$\{ \theta_{l} \} \subseteq \frac{1}{2}\Z$
in accordance with {\em \refconj{conj:d'}}.
\end{thm}

A part of this theorem was proved in the authors thesis 
\cite{Hansen1}, where the case of Seifert manifolds with 
orientable base was handled. In this paper we calculate the 
asymptotics of a class of functions generalizing the 
invariants of the Seifert manifolds with orientable base or 
non-orientable base with even genus. Unfortunately certain 
parts of the proof do not go through as they stand for the 
remaining Seifert manfolds, i.e.\ the ones with non-orientable 
base of odd genus. If $X$ is a Seifert manifold with base 
$S^{2}$ the asymptotic expansion of $\tau_{r}(X)$ has
the form (\ref{eq:AEC'}) as stated above, but some of the
$q_{j}$'s are not Chern--Simons values. If the number of 
exceptional fibers is $n \leq 2$, i.e.\ if $X$ is a lens
space, it is easy to prove that all the terms in our asymptotic
formula, corresponding to $q_{j}$'s which are not
Chern--Simons values, are zero, thus proving the AEC for those
manifolds. We also prove that this is the case for the 
Seifert manifolds with $3$ exceptional fibers, cf.\ 
\refthm{thm:AECSeifertn3}. We do not in this paper carry 
through the analysis needed to prove the AEC for the 
cases of $n \geq 4$ exceptional fibers, see 
Sec.~\ref{secgenus0} for more details. Let us finally
mention that we also prove the AEC for some of the Seifert
manifolds with non-orientable base with odd genus,
namely for the Seifert manifold with base
$\R \text{P}^{2}$ and zero or one excepotional fiber, cf.\
\refcor{cor:oddcase}, this indicating that nothing special
happens in the case of a non-orientable base of odd genus.
 
The proof of \refthm{thm:AEC-Seifert-manifolds} is long and 
technical, but each step in the proof only uses elementary 
mathematics. Let us for the convenience of the reader give an 
outline of the proof here. The method used follows closely 
ideas of Rozansky \cite{Rozansky3} and this paper owes much to
Rozansky's work. As stated earlier we supplement the work of
Rozansky by establishing certain required error estimates. 
Unfortunately the calculations needed to establish those
estimates take up a big space explaining to a great extend the
length of this paper.

Let $\tau_{r}=\tau_{r}^{\SU(2)}$ and let $X$ be a Seifert 
manifold. By \cite[Theorem 8.4]{Hansen2} (see 
\refthm{thm:RT-invariants-Seifert-manifolds}) we have
\begin{equation}\label{eq:Zfunction-Seifert'}
\tau_{r}(X) = a(r)\sum_{\gamma =1}^{r-1} 
\frac{(-1)^{\gamma dg} h(\gamma)}
{\sin ^{n+dg-2} \left( \frac{\pi}{r} \gamma \right) },
\end{equation}
where $a(r)=br^{dg/2-1}e^{ic/r}$ for some constants $c \in \R$ 
and $b \in \C \sm \{0\}$ (depending on $X$). Here $n$ is the 
number of singular fibers, $g$ is the genus of the base of $X$ 
and $d=2$ if this base is orientable and $d=1$ if it is 
non-orientable. (Here we say that the non-orientable connected 
sum $\# ^{k} \R \text{P}^{2}$ has genus $k$.) The function $h$
is extendable to an entire function with some nice symmetry
properties (see (\ref{eq:hfunction-symmetry})).

Because of the simple form of $a(r)$ we can concentrate on
$Z(X;r)=\tau_{r}(X)/a(r)$. For technical reasons (see
Remarks~\ref{rem:P1} and \ref{rem:P2}) we will assume that 
$dg$ is even, so from now on $X$ is a Seifert manifold with
orientable base or nonorientable base of even genus. 

The first task is to get rid of the dependence of $r$ in the 
summation limit in (\ref{eq:Zfunction-Seifert'}). To this end 
we use that a function $f:\Z \ra \C$ periodic with a period of
$P$ satisfies
\begin{equation}\label{eq:periodicity-one-dimensional'}
\sum_{k=0}^{P-1} f(k) = P \livep \svep \sum_{k \in \Z}
e^{-\pi \vep k^{2}} f(k).
\end{equation}
This leads together with some other symmetry considerations to 
the formula
\begin{equation}\label{eq:Zfunction-formula2'}
Z(X;r)=r \livep \svep \lixi 
\sum_{\gamma \in \Z} e^{-\pi \vep \gamma^{2}}\frac{h(\gamma)}
{{\sin^{n+dg-2}\left(\frac{\pi}{r} (\gamma-i\xi) \right)}},
\end{equation}
where we have introduced an extra parameter $\xi$ to avoid 
singularities. (If there are no singular fibers, i.e.\ if 
$n=0$, then we have to include the extra factor
$$
\frac{\sin\left( \frac{\pi}{r}\gamma \right)}
{\sin\left( \frac{\pi}{r}(\gamma -i\xi)\right)}
$$
in the summand.)
 
To proceed we change (\ref{eq:Zfunction-formula2'}) to a sum of
integrals by using the Poisson summation formula
\begin{equation}\label{eq:Poisson'}
\sum_{k \in \Z} \varphi(k) 
= \sum_{m \in \Z} \iny e^{2\pi i mx} \varphi(x) \dte x.
\end{equation}
It is well-known that this formula is valid for functions 
$\varphi$ in the Schwartz space $\mS(\R)$ of smooth functions, 
that, together with their derivatives, are rapidly decreasing 
at infinity, see e.g.\ \cite{Hormander}. We use Poisson's 
formula with 
$\varphi(\gamma)=e^{-\pi \vep \gamma^{2}}\frac{h(\gamma)}
{\sin^{n+dg-2}\left(\pi(\gamma-i\xi)/r \right)}$
(with fixed $\xi$ and $\vep$). This is in $\mS(\R)$ because 
$h$ is periodic and (putting $z=\gamma/r$, 
$\delta=\vep r^{2}$ and $\eta=\xi/r$) we arrive at an 
expression of the form
\begin{equation}\label{eq:Zfunction-formula3'}
Z(X;r)=r \lidel \sdel \lieta \sum_{m \in \Z} 
\sum_{\lambda \in \Lambda} g(\lambda) 
\iny e^{-\pi\delta z^{2}} \tilde{\kappa}(z)
e^{r (Q(z)+2\pi imz)} \dte z,
\end{equation}
where $\Lambda$ is a finite set (only depending on $X$), $Q$ 
is a polynomial in $z$ with coefficients only depending on 
$\lambda$ and $X$, and
$$
\tilde{\kappa}(z) = 
\frac{e^{i b z}}{\sin^{k} \left( \pi (z -i \eta) \right) },
$$
where $b \in \R$ only depends on $\lambda$ and $X$ and 
$k=n+dg-2$. Put 
$\kappa(z)=e^{-\pi \delta z^{2}} \tilde{\kappa}(z)$ and
$Q_{m}(z)=Q(z)+2\pi imz$ in the following. We note that 
$Q_{m}$ is of degree $1$ if the Seifert Euler number $E$ of 
$X$ is equal to zero and of degree $2$ if $E \neq 0$. This 
makes the asymptotic analysis in the case $E=0$ different from
(and also simpler than) that in the case $E \neq 0$. To make 
the notation in the following compatible with the notation 
used in Sec.~\ref{sec-Proof-of} and \ref{sec-Proof-of-E0}, we 
put $A=\pi \frac{E}{2}$.

Let us first assume that $E \neq 0$. In this case the integrals
$I(r)=\iny \kappa(z) e^{r Q_{m}(z)} \dte z$ in 
(\ref{eq:Zfunction-formula3'}) are approximated by using the 
steepest descent method. Thus, if $z_{\st}$ denotes the
stationary point of $Q_{m}$, i.e.\ $Q_{m}'(z_{\st})=0$, we 
deform the integration contour (here the real axis) to a new 
contour $C$ passing through $z_{\st}$ in such a way that the 
large $r$ asymptotics of 
$\int_{C} \kappa(z) e^{r Q_{m}(z)} \dte z$ is easy to obtain 
by standard methods (such as the Laplace method), see e.g.\ 
\cite{BleisteinHandelsman}, \cite{deBruijn}, 
\cite{Wong}. It can be shown that a good such contour is a 
contour $C_{\sd}$ along which the imaginary part of $Q_{m}$ is
constant. In the case of $I(r)$ this so-called steepest 
descend contour is the straight line
\begin{equation}\label{eq:sdcontour'}
\ima(z)=\sgA ( \rea(z)-z_{\st}).
\end{equation}
(We note that the stationary point $z_{\st}$ is non-degenerate,
i.e.\ $Q_{m}''(z_{\st}) \neq 0$, and that this point and 
thereby $C_{\sd}$ depend on $(m,\lambda)$.) To determine the 
asymptotic expansion of $I(r)$ in the limit $r \ra \infty$ we 
need to calculate the large $r$ asymptotics of the new contour
integral $J(r)=\int_{C_{\sd}} \kappa(z) e^{r Q_{m}(z)} \dte z$,
and of the difference $I(r)-J(r)$. This difference is zero for 
$k \leq 0$ and is for $k>0$ given by the infinite sum of 
residues in the poles of $\kappa(z)$ crossed when deforming the
real axes to $C_{\sd}$, i.e.\ the poles $z_{l}=l+i\eta$, 
$l \in \Z$, for which
\begin{equation}\label{eq:condition-poles'}
\sgA (l-z_{\st}) > \eta.
\end{equation}
(We refer to \reflem{lem:Lemma8.7}.) By futher showing that 
certain infinite sums are absolutely convergent (see 
\reflem{lem:Lemma8.12}), we therefore end up with the 
following partition of $Z(X;r)$:
\begin{equation}\label{eq:Zfunction-splitting'}
Z(X;r)=r\lidel \sqrt{\delta} \lieta 
\left( Z_{\inte,1}(\delta,\eta)
+ Z_{\intpolar}(\delta,\eta)+Z_{\pol,1}(\delta,\eta) \right),
\end{equation}
where
\begin{eqnarray}\label{eq:Zterms'}
Z_{\inte,1}(\delta,\eta) &=& \sum_{(m,\lambda) \in W} 
g(\lambda) 
\int_{C_{\sd}(m,\lambda)} \kappa(z)e^{rQ_{m}(z)} \dte z, \\ 
Z_{\intpolar}(\delta,\eta) &=& \sum_{\lambda \in \Lambda} 
g(\lambda) \msumi \int_{C_{\sd}(m,\lambda)} 
\kappa(z)e^{rQ_{m}(z)} \dte z, \nonumber \\
Z_{\pol,1}(\delta,\eta) &=& 2\pi i \sum_{\lambda \in \Lambda} 
g(\lambda) \sum_{m \in \Z} 
\sum_{\stackrel{l \in \Z}{\sgA (l-z_{\st}(m,\lambda)) > \eta}}
\Res _{z=z_{l}(\eta)} \left\{ \kappa(z)e^{rQ_{m}(z)} \right\}, 
\nonumber
\end{eqnarray}
where
$W=\{ (m,\lambda) \in \Lambda \times \Z\;|\;
z_{\st}(m,\lambda) \notin \Z\}$.
(Here we suppress the dependency on $r$ and also to some 
extend the dependency on $\lambda$ in our notation for not 
making it too clumsy.) The above sums are all absolutely 
convergent. By rearranging the terms in the sum 
$Z_{\pol,1}(\delta,\eta)$ we obtain an expression of the form
(see \reflem{lem:Lemma8.12}) 
$$
Z_{\pol,1}(\delta,\eta) = 2\pi i \sum_{\lambda \in \Lambda} 
g(\lambda)\sum_{l \in \Z} \beta(l) \msumj \Res_{z=i\eta} 
\left\{ \tilde{\phi}(z;l,\lambda)e^{2\pi irmz} \right\},
$$
where $B/\pi$ is a certain rational number only depending on 
$\lambda$ and $X$, and $\tilde{\phi}$ is a meromorphic 
function independent of $m$. Along the same lines we obtain 
the expression
$$
Z_{\intpolar}(\delta,\eta) = \sum_{\lambda \in \Lambda} 
g(\lambda)\sum_{l \in \Z} \beta(l)
\sum_{\stackrel{m \in \Z}{m+\frac{B}{2\pi}+\frac{A}{\pi}l}=0}
\int_{C(\lambda,0)} \tilde{\phi}(z;l,\lambda)e^{2\pi irmz} 
\dte z,
$$
where $C(\lambda,\rho)$, $\rho \in \R$, is the line
$$
\ima(z) = \sgA\rea(z)-\rho.
$$
We have a decomposition (see \reflem{lem:Lemma8.13} for 
details) 
$$
Z_{\intpolar}(\delta,\eta) = Z_{\pol,2}(\delta,\eta)
+ Z_{\inte,2}(\delta,\eta),
$$
where $Z_{\pol,2}(\delta,\eta)$ is given by a sum of residues 
like $Z_{\pol,1}(\delta,\eta)$, in fact
$$
Z_{\pol,2}(\delta,\eta) = \pi i \sum_{\lambda \in \Lambda} 
g(\lambda) \sum_{l \in \Z} \beta(l)
\sum_{\stackrel{m \in \Z}{m+\frac{B}{2\pi}+\frac{A}{\pi}l}=0}
\Res_{z=i\eta} \left\{ \tilde{\phi}(z;l,\lambda)e^{2\pi irmz} 
\right\},
$$
and
\begin{equation}\label{eq:Zint2-deltaeta'}
Z_{\inte,2}(\delta,\eta) = \sum_{\lambda \in \Lambda} 
g(\lambda) \lsuma \beta(l) 
\int_{C(\lambda,\eta)} \tilde{F}(y;l,\lambda) \dte y,
\end{equation}
where $\tilde{F}(y)$ is a meromorphic function with poles in 
$\Z$ (if $k>0$) (explaining the name $Z_{\intpolar}$). The 
next observation is that there exists a $\eta_{0}>0$ such that
$$
Z_{\pol,1}(\delta,\eta) = 2\pi i \sum_{\lambda \in \Lambda} 
g(\lambda) \sum_{l \in \Z} \beta(l)
\sum_{\stackrel{m \in \Z}{m+\frac{B}{2\pi}+\frac{A}{\pi}l}>0}
\Res_{z=i\eta} \left\{ \tilde{\phi}(z;l,\lambda)e^{2\pi irmz} 
\right\},
$$
for all $\eta \in ]0,\eta_{0}]$. By letting
$Z_{\pol}(\delta,\eta) = Z_{\pol,1}(\delta,\eta)
+ Z_{\pol,2}(\delta,\eta)$ 
we therefore obtain
\begin{equation}\label{eq:Zfunction-splitting1'}
Z(X;r) = r\lidel \sqrt{\delta} \lieta \left( 
Z_{\inte,1}(\delta,\eta) + Z_{\inte,2}(\delta,\eta)
+ Z_{\pol}(\delta,\eta) \right),
\end{equation}
where $Z_{\inte,1}(\delta,\eta)$ and 
$Z_{\inte,2}(\delta,\eta)$ are given respectively by 
(\ref{eq:Zterms'}) and (\ref{eq:Zint2-deltaeta'}) and
\begin{equation}\label{eq:Zpolar-intro}
Z_{\pol}(\delta,\eta) = 2\pi i \sum_{\lambda \in \Lambda} 
g(\lambda) \sum_{l \in \Z} \beta(l) \msumk 
\Res_{z=i\eta} \left\{ 
\frac{\tilde{\phi}(z;l,\lambda)e^{2\pi irmz}}
{\Sym_{\pm}\left(m+\frac{B}{2\pi}+\frac{A}{\pi}l \right)} 
\right\}.
\end{equation}
Here we have introduced the function 
$\Sym_{\pm} : \R \to \{1,2\}$ (borrowing notation from
\cite{Rozansky3}) which is $1$ in all points 
except in zero where it is $2$.

The next natural step would be to try to calculate the limit 
$\lidel \sdel \lieta f(\delta,\eta)$ (or at least the large
$r$ asymptotics of that limit) for $f$ equal to each of 
the functions $Z_{\inte,1}$, $Z_{\inte,2}$, and $Z_{\pol}$ 
(assuming these limits exist). Thus we start by showing that 
\begin{equation}\label{eq:Zpolar1'}
Z_{\pol}(r):=r\lidel \sdel \lieta Z_{\pol}(\delta,\eta)
\end{equation}
exists by calculating this limit explicitly.

To handle the two other limits
$\lidel \sdel \lieta Z_{\inte,\nu}(\delta,\eta)$, $\nu=1,2$, 
more care need to be taken. Since the limits on the right-hand
sides of (\ref{eq:Zfunction-splitting1'}) and 
(\ref{eq:Zpolar1'}) exist we have that
\begin{equation}\label{eq:Zinte'}
Z_{\inte}(r) := r\lidel \sdel \lieta Z_{\inte}(\delta,\eta)
\end{equation}
exists, where 
$Z_{\inte}(\delta,\eta) = Z_{\inte,1}(\delta,\eta)
+ Z_{\inte,2}(\delta,\eta)$.
However, we do not calculate this limit explicitly, but only 
give it an asymptotic description. This is of course sufficient
for the proof of \refthm{thm:AEC-Seifert-manifolds}. To find 
the asymptotic expansion of $Z_{\inte}(r)$ we show that we for
each $N \in \Z_{\geq 0}$ have a decomposition
\begin{equation}\label{eq:Zintedecomp}
Z_{\inte,\nu}(\delta,\eta) = Z_{\inte,\nu}(N;\delta,\eta)
+ R_{\nu}(N;\delta,\eta)
\end{equation}
such that 
$$
Z_{\inte,\nu}(r;N) = r\lidel \sdel \lieta 
Z_{\inte,\nu}(N;\delta,\eta)
$$
exists and can be calculated explicitly, $\nu =1,2$. This 
implies that
$$
\lidel \sdel \lieta
\left( R_{1}(N;\delta,\eta) + R_{2}(N;\delta,\eta) \right)
$$
exists. We will not show that
$\lidel \sdel \lieta R_{\nu}(N;\delta,\eta)$ exists for each
$\nu=1,2$ separately. Instead we show that 
$R_{\nu}(N;\delta,\eta)$ is bounded from above by a certain 
function $A_{\nu}(N;\delta,\eta)$ for which
$$
r\lidel \sdel \lieta A_{\nu}(N;\delta,\eta)
$$
can easily be calculated and shown to be small in the correct 
asymptotic sense compared to $Z_{\inte,\nu}(r;N)$ in the limit
of large $r$, $\nu=1,2$. This proves that 
$Z_{\inte,1}(r;N) + Z_{\inte,2}(r,N)$ is the large $r$ 
asymptotics of $Z_{\inte}(r)$ to an order depending on $N$.

We observe that all the sums $Z_{\pol}(\delta,\eta)$,
$Z_{\inte,1}(\delta,\eta)$, and $Z_{\inte,2}(\delta,\eta)$
are infinite. (For the sum $Z_{\inte,2}(\delta,\eta)$
note here that $A/\pi \in \Q$ and that
$B(\lambda)/(2\pi) \in \Z$ for some $\lambda$.)
According to the AEC we should end up with a finite sum. 
Basically these infinite sums are changed to finite sums by 
using the periodicity result 
(\ref{eq:periodicity-one-dimensional'}) `backwards'. The sum 
$Z_{\pol}(\delta,\eta)$ in (\ref{eq:Zpolar-intro}) actually 
contains two infinite sums. We calculate one of these, namely
the sum over $m$, explicitly. Let us give some futher details.
First we change $Z_{\pol}(\delta,\eta)$ to an expression of 
the form
\begin{equation}\label{eq:Zpol-intro}
Z_{\pol}(\delta,\eta) =\sum_{l \in \Z}
\sum_{\nu \in \Gamma(l)} 
\sum_{\stackrel{m \in \Z}{m-V(\nu) \geq 0}}
\frac{\Res_{z = i\eta} \left\{ e^{-\pi\delta(z+l)^{2}}
G(z;l,\nu,\eta)e^{2\pi i rmz}\right\}}
{\Sym_{\pm}(m-V(\nu))},
\end{equation}
where $\Gamma(l)$ is a finite index set depending on $l$ in a 
periodic way, and $G$ is periodic w.r.t.\ $l$. We use certain 
symmetries to establish that. Now, if $a \in \R$ and 
$f: \Z \to \C$ then
$$
\sum_{\stackrel{m \in \Z}{m \geq a}} 
\frac{f(m)}{\Sym_{\pm}(m-a)}
= \sum_{m=0}^{\infty} \frac{f(m)}{\Sym_{\pm}(m)} 
- \sum_{\stackrel{m \in \Z}{0 \leq m \leq |a|}} 
\frac{\sign (a) f(\sign (a) m) }
{ \Sym_{\pm}(m) \Sym_{\pm} (m-|a|)}.
$$
in the sense that if the left-hand side is convergent, then 
the infinite sum in the right-hand side is also convergent and
the identity holds. Using this identity on 
(\ref{eq:Zpol-intro}) together with
$$
\sum_{m=0}^{\infty} \frac{e^{2\pi i rmz}}{\Sym_{\pm}(m)}
= \frac{i}{2} \cot(\pi r z),
$$
valid for all $z$ in the open disk in the complex plane with 
center $i\eta$ and radius $\eta$, we thus obtain
\begin{eqnarray}\label{eq:Zpol-intro1}
Z_{\pol}(\delta,\eta) &=& \sum_{l \in \Z} 
\sum_{\nu \in \Gamma(l)} \left( \frac{i}{2}
\Res_{z=i\eta} \left\{ e^{-\pi\delta (z+l)^{2}}
G(z;l,\nu,\eta)\cot(\pi rz)\right\} \right. \\
&& \left. - \sum_{\stackrel{m \in \Z}{0 \leq m \leq |V|}}
\frac{\sgV \Res_{z=i\eta}\left\{ e^{-\pi\delta (z+l)^{2}} 
G(z;l,\nu,\eta)e^{2\pi ir\sgV mz}\right\}}
{\Sym_{\pm}(m)\Sym_{\pm}(m-|V|)}\right). \nonumber
\end{eqnarray}
From this expression we obtain that
$Z_{\pol}(\delta) := \lieta Z_{\pol}(\delta,\eta)$ is equal to 
the right-hand side of (\ref{eq:Zpol-intro1}) with $\eta=0$.
Finally $\lidel \sdel Z_{\pol}(\delta)$ is calculated by using 
(\ref{eq:periodicity-one-dimensional'}). We thus obtain an 
exact expression for $Z_{\pol}(r)$ equal to some Laurent 
polynomial in $r$.

Let us also give a few more details concerning the calculation
of the asymptotic expansion of $Z_{\inte}(r)$ in 
(\ref{eq:Zinte'}). First we observe that
$$
Z_{\inte,1}(\delta) := \lieta Z_{\inte,1}(\delta,\eta)
= \sum_{(m,\lambda) \in W} g(\lambda) 
\int_{C_{\sd}(m,\lambda)} \kappa_{0}(z)e^{rQ_{m}(z)} \dte z,
$$
where $\kappa_{0}(z)$ is equal to $\kappa_{\eta}(z)$ with 
$\eta=0$. (For $\nu=1$ we thus calculate this limit before 
making the decomposition (\ref{eq:Zintedecomp}).) For an 
arbitrary $(m,\lambda) \in W$ we can expand $\kappa_{0}$ as a 
power series in a small disk around $z_{\st}(m,\lambda)$ (with
a radius independent of $(m,\lambda)$). Let $N$ be an 
arbitrary but fixed non-negative integer, and let 
$\kappa_{0}(z;N)$ be equal to the first $N+1$ terms in this 
power series. Then we let
$$
Z_{\inte,1}(N;\delta) = \sum_{(m,\lambda) \in W} g(\lambda) 
\int_{C_{\sd}(m,\lambda)} \kappa_{0}(z;N) e^{rQ_{m}(z)} \dte z.
$$
The remainder term (to order $N$) is then simply given by
$R_{1}(N;\delta) = Z_{\inte,1}(\delta) 
- Z_{\inte,1}(N;\delta)$.
Since the power series expansion of $\kappa_{0}(z)$ is only 
valid in a small neigborhood of $z_{\st}(m,\lambda)$ this 
remainder term is not simply given by replacing 
$\kappa_{0}(z)$ by the remaining terms in the power series 
expansion in the expression for $Z_{\inte,1}(\delta)$, see 
\reflem{lem:Lemma8.21} for details.

The decomposition (\ref{eq:Zintedecomp}) for $\nu=2$ is 
obtained in a similar way by representing a preexponential 
factor of the integrand $\tilde{F}$ as a power series around 
$y=0$. However, in this case the positive parameter $\eta$ 
plays a crucial role and the limit $\lieta$ can only be 
calculated in the final step.

Let us finally say a few words about the proof of
\refthm{thm:AEC-Seifert-manifolds} in the case $E=0$. Again we 
obtain the expression (\ref{eq:Zfunction-formula3'}). However, 
in this case $Q(z)=2\pi i az$, where $a$ is a constant 
depending on $\lambda$ but not on $m$, and the steepest 
descent method can not be used to calculate the large $r$ 
asymptotics of the integrals 
$\iny \kappa(z)e^{rQ_{m}(z)} \dte z$ since
$Q_{m}(z)=Q(z)+2\pi imz$ has no stationary points in this 
case. Instead we show that there exist certain contours 
$C(m,\lambda)$ so that
$$
Z(X;r) = r\lidel\sdel\lieta \left( Z_{\inte}(\delta,\eta)
+ Z_{\pol}(\delta,\eta)\right),
$$
where
{\allowdisplaybreaks
\begin{eqnarray*}
Z_{\inte}(\delta,\eta) &=& \sum_{\lambda \in \Lambda} 
g(\lambda) \sum_{m \in \Z} 
\int_{C(m,\lambda)} \kappa(z)e^{rQ_{m}(z)} \dte z, \\ 
Z_{\pol}(\delta,\eta) &=& 2\pi i \sum_{\lambda \in \Lambda} 
g(\lambda) \sum_{\stackrel{m \in \Z}{m-a(\lambda) \geq 0}} 
\sum_{l \in \Z}
\Res _{z=z_{l}(\eta)} \left\{ \kappa(z)e^{rQ_{m}(z)} \right\}.
\end{eqnarray*}}\noindent
By estimating the integrals inside the sum 
$Z_{\inte}(\delta,\eta)$ it can be proved that
$$
\lidel\sdel\lieta Z_{\inte}(\delta,\eta)=0.
$$
Hence
$$
Z(X;r)=r\lidel\sdel\lieta Z_{\pol}(\delta,\eta),
$$
and this double limit can be calculated explicitly along the 
same lines as in the case $E \neq 0$.

From a general point of view the main problem when calculating 
the large $r$ asymptotics of the quantum invariants as defined 
by Reshetikhin and Turaev, is that these invariants are given 
by a finite sum in which both the number of terms and the 
terms themselves depend on the level $r$. It has been 
speculated that there is some kind of number theoretic 
principle which can be used to calculate these asymptotics. In
\cite{Jeffrey1} Jeffrey revealed that one can calculate the 
asymptotics of the $\SU(2)$--invariants of lens spaces by an 
inductive argument using a reciprocity formula for Gauss sums.
This reciprocity formula follows from an application of the 
Poisson summation formula, which is a special version of the 
Fourier transform (namely the Fourier transform of the sum of 
Dirac measures on some discrete subgroup of $\R^{n}$). More 
generally Jeffrey speculated that the Fourier transform may be
the basic mechanism behind identities needed for calculating 
the asymptotics.
 
Rozansky \cite{Rozansky1}, \cite{Rozansky3} took up these 
ideas and showed how to use the Poisson summation formula to 
calculate the asymptotics for the broader class of Seifert 
manifolds with orientable base space. For a general Seifert 
manifold the main idea is as illustrated above to change the 
finite sum formula for the invariant to an infinite sum using 
certain periodicity results. This involves introducing certain
limits with respect to certain parameters. One can then use 
the Poisson summation formula on the infinite sum. This leads 
to a formula for the invariant given by an infinite sum of 
integrals. Each of these integrals can be given an asymptotic 
description via the steepest descent method. As proved in this
paper this leads ultimatively to an asymptotic formula like 
(\ref{eq:AEC'}).

For the calculations of the asymptotics of the 
$\SU(2)$--invariants of Seifert manifolds one only needs a 
one-dimensional Poisson summation formula (involving only a 
sum over $\Z$), see (\ref{eq:Poisson'}). In \cite{Jeffrey1} 
Jeffrey also obtained a multi-dimensional reciprocity formula 
for Gauss sums using a multi-dimensional version of the 
Poisson formula (see also \cite{Jeffrey2} and 
\cite[Appendix A]{HansenTakata}). She then used this to 
calculate the asymptotics of the $G$--invariants of some 
mapping tori over the torus, $G$ being any compact simply 
connected simple Lie group. In \cite{HansenTakata} the author 
and Takata use the multi-dimensional reciprocity formula to 
calculate the asymptotics of the $G$--invariants of the lens 
spaces by a similar approach as in \cite{Jeffrey1}.

In \cite{AndersenHansen}, Andersen and the author investigate
the asymptotics of the $\SU(2)$--invariants of the 
$3$--manifolds obtained by rational surgeries on the 
$3$--sphere along the figure $8$ knot. In these calculations 
the Fourier transform is not used. Instead the asymptotics is 
calculated by first obtaining a formula for the invariant 
given by a (finite) sum of certain contour integrals and then 
use the saddle point approximation method on each of the 
integrals.

It seems that a general principle behind calculating the large 
quantum level asymptotics of the RT--invariants of 
$3$--manifolds is to obtain in some way or the other a formula
of the invariants given by some sum of integrals, where the 
number of terms in the sum (allowed to be infinite) is 
independent of the level $r$, and where the integrals can be 
asymptotically described by using some well-known 
approximation method such as the method of steepest descent or
more generally the saddle point method. 
The involved integrals will be of a 
form $\int_{C} \psi(z)e^{r\phi(z)} \dte z$, where $C$ is some 
(multi-dimensional) contour in $\C$ ($\C^{n}$) and the phase
function $\phi$ and the function $\psi$ are some complex 
functions independent of the level $r$. Note that even for 
the case of lens spaces and mapping tori over the torus, where
a reciprocity formula for Gauss sums is used, one actually 
uses this approach, namely, the involved reciprocity formula 
can be proved by this approach, see e.g.\ 
\cite[Appendix]{HansenTakata}. The calculations in 
\cite{AndersenHansen} actually first lead to an expression for
the invariants given by a finite sum of integrals, where the
integrands are not of the nice form 
$\int_{C} \psi(z)e^{r\phi(z)} \dte z$. One has to make
certain initial approximations to obtain an (approximative)
expression for the invariants given by a finite sum of 
integrals of the nice form. This phenomenon should be expected
in general. 

The paper is organized as follows. In 
Sec.~\ref{sec-Asymptotic-expansions} we introduce the notions 
from asymptotic analysis needed. In this section we study from
a general perspective asymptotic expansions of the form 
(\ref{eq:AEC'}). These expansions constitute a certain 
subclass of a class of asymptotic expansions which we denote 
asymptotic expansions of generalized Poincar\'{e} type, since 
they are finite sums with each term being an asymptotic 
expansion of (ordinary) Poincar\'{e} type. A small but 
important result is \refthm{thm:uniqueness}, which extend the 
well-known uniqueness property \reflem{lem:Poincare} for 
asymptotics of Poincar\'{e} type to asymptotics of generalized
Poincar\'{e} type. In Sec.~\ref{sec-The-Reshetikhin} we 
introduce notation for the Seifert manifolds and state from 
\cite{Hansen2} formulas for the RT--invariants associated to 
$\SU(2)$ of all Seifert manifolds. In 
Sec.~\ref{sec-The-analytic} we present detailed formulas for 
the large $r$ asymptotic expansion of $\tau_{r}^{\SU(2)}(X)$,
$X$ a Seifert manifold as in 
\refthm{thm:AEC-Seifert-manifolds}. Moreover, we review results
of Auckly \cite{Auckly1} on the classical $\SU(2)$ 
Chern--Simons theory on Seifert manifolds. (We correct a small
error in Auckly's result, see \refthm{thm:Auckly} and below.) 
We demonstrate that Auckly's results combine with our formulas 
for the asymptotic expansions to give a proof of 
\refthm{thm:AEC-Seifert-manifolds}. Finally, we identify for 
the Seifert fibered rational homology spheres the 
Casson--Walker invariant in the asymptotic formula using 
results of Lescop \cite{Lescop}. In Sec.~\ref{sec-Proof-of} 
and \ref{sec-Proof-of-E0} the proofs of
Theorems~\ref{thm:asymptotics-Enot0} and 
\ref{thm:asymptotics-E0} are given. We focus on the main ideas 
deferring technicalities to appendices.

We end this introduction by a side remark. The part of the 
theory of quantum invariants concerned with asymptotic 
expansions is often denoted the perturbative theory, and 
(parts of) the asymptotic expansions of the RT--invariants (or
Witten's invariants) are called perturbative quantum 
invariants. One should note, however, that the notion of 
perturbative invariants has different meanings in the 
mathematical literature. There are the analytic asymptotics 
studied in this paper, and there are the algebraically defined
perturbative invariants of Ohtsuki \cite{Ohtsuki1}, 
\cite{Ohtsuki2}. Note that Ohtsuki's perturbative invariant of
a $3$--manifold $X$ is trivial unless $X$ is a rational 
homology sphere \cite[Remark 1.3]{Ohtsuki2}. This is certainly
not the case for the analytic asymptotics (as the results of 
this paper reveals). It is, however, generally believed 
that these two notions of perturbative invariants are closely 
related. It is speculated that the coefficients in the formal 
Ohtsuki series determine (some of) the coefficients in the 
analytic asymptotic expansion of the RT--invariants of 
rational homology spheres, see 
\cite[Example 1.6 and Sect.~5]{Ohtsuki2}, \cite{Lawrence},
\cite{Rozansky4}, \cite{Rozansky5} and 
\cite[Sect.~7]{Ohtsuki3} for some results and conjectures. The
results in this paper should be useful to investigate this 
question for the Seifert fibered rational homology spheres. We
will however not pursue this issue futher in this paper. 
Recently K.\ Habiro has proposed new invariants for integral 
homology spheres which according to a conjecture of Habiro and
Le should determine both Ohtsuki's perturbative invariants and
the RT--invariants for these manifolds, see 
\cite[Conjecture 7.29]{Ohtsuki3}. (As stated in 
\cite[p.~492]{Ohtsuki3} this should now be a theorem of Habiro
and Le.) Let us finally mention that there seem to be some very
interesting number theoretic aspects of the asymptotics of the
Seifert fibered integral homology spheres. We will not consider
this issue in this paper but refer to papers of
Lawrence and Zagier \cite{LawrenceZagier} and
Hikami \cite{Hikami1}--\cite{Hikami3}.

\bigskip

\noindent {\bf Acknowledgements} The author thanks Institut de 
Recherche Math\'{e}matique Avanc\'{e}e (IRMA) at Universit\'{e}
Louis Pasteur and C.N.R.S., Strasbourg, the School of 
Mathematics at the University of Edinburgh, the 
Max--Planck--Institut f\"{u}r Mathematik in Bonn, and the
Department of Mathematical Sciences, University of Aarhus, for
their hospitality during this work. He was supported by the 
European Commission, the Danish Natural Science Research 
Council and the Max--Planck--Institut f\"{u}r Mathematik.
In addition this article is based upon work
supported by the National Science Foundation under Grant No.\ 
EPS-0236913 and matching support from the State of Kansas
through Kansas Technology Enterprise Corporation. 
The paper extends results
obtained in the authors thesis \cite{Hansen1}, and he would 
like to thank his thesis advisor J.\ E.\ Andersen for having 
given valuable comments to earlier versions of this paper, and
also for the many helpful conversations about quantum 
invariants in general and asymptotic expansions of these in 
particular.

\section{Asymptotic expansions og quantum type}
\label{sec-Asymptotic-expansions}

\noindent In this section we study asymptotics of the form 
(\ref{eq:AEC'}) from a general point of view. We shall call 
such asymptotics for asymptotic expansions of quantum type, 
see \refdefi{defi:quantum-asymptotics}. In this section we 
shall mention the few facts from asymptotic analysis needed to
study them. For basic introductions to asymptotic analysis in 
general, see e.g.\ \cite{BleisteinHandelsman}, 
\cite{deBruijn}, \cite{Olver}, and \cite{Wong}.

Let in the following $A =\Z_{\geq h}=\{ m \in \Z | m \geq h\}$
for some fixed positive integer $h$. Recall the usual 
order-notation denoted by the $O$- and $o$-symbols: If 
$f,g: A \to \C$ are two functions, then $g=O(f)$ if there 
exists a constant $C$ such that $|f(k)| \leq C|g(k)|$ for all 
$k \in A$. Of course such a condition only puts limits on the 
large $k$ behaviour of $f(k)$. We say that $g=o(f)$ in the 
limit $k \ra \infty$ if for all $\vep >0$ there exists a 
$k_{\vep} \in A$ such that $| g(k) | \leq \vep |f(k)|$ for all 
$k \geq k_{\vep}$. If $f(k) \neq 0$ for all $k$ (or for all 
sufficiently large $k$) this is equivalent to 
$g(k)/f(k) \ra 0$ as $k \ra \infty$. As examples we have 
$k^{-(n+1)}=o(k^{-n})$ in the limit 
$k \ra \infty$ and $k^{n} + \log(k)k^{n-1}=O(k^{n})$ for every
$n \in \Z$. A sequence of functions 
$\{ \varphi_{n}: A \ra \C \}_{n=0}^{\infty}$ is called an 
{\em asymptotic sequence} as $k \ra \infty$ if for all 
$n \geq 0$
$$
\varphi_{n+1}(k)=o(\varphi_{n}(k))
$$
as $k \ra \infty$ and if we for all $n$ have that
$\varphi_{n}(k) \neq 0$ for all sufficiently large $k$. For 
functions $f:A \ra \C$ and $f_{n} : A \ra \C$, the formal 
series $\sum_{n=0}^{\infty} f_{n}(k)$ is called a 
{\em generalized asymptotic expansion} of $f(k)$ with respect 
to the asymptotic sequence $\{ \varphi_{n}(k) \}$, as 
$k \ra \infty$, if 
\begin{equation}\label{eq:generalized-asymptotics}
f(k) = \sum_{n=0}^{N} f_{n}(k) + o(\varphi_{N}(k))
\end{equation}
in the limit $k \ra \infty$ for every $N \geq 0$. In this case
we write
\begin{equation}\label{eq:asymptotic-notation}
f(k) \sim \sum_{n=0}^{\infty} f_{n}(k); \hspace{.2in} 
\{ \varphi_{n} \} , \hspace{.2in} k \ra \infty.
\end{equation}
When $f_{n}(k) = a_{n}\varphi_{n}(k)$, $a_{n}$ a fixed complex 
number, for every $n$, then the above expansion is said to be 
of {\em Poincar\'{e} type}. Furthermore, if the expansion is of
this type, and $\varphi_{n}(k)=(\xi(k))^{\lambda_{n}}$, 
$\lambda_{n}$ a fixed complex number, the expansion is said to
be of {\em power series type}. We have the following immediate
uniqueness result for asymptotic expansions of Poincar\'{e} 
type (see e.g.\ \cite[pp.~16-17]{BleisteinHandelsman}).

\begin{lem}\label{lem:Poincare}
Let $\{ \varphi_{n} \}_{n \in \N}$ be an asymptotic sequence 
and let $f:A \ra \C$. Then there is at most one sequence of 
complex numbers $a_{0},a_1,a_2,\ldots$ such that
$$
f(k) \sim \sum_{n=0}^{\infty} a_{n} \varphi_{n}(k);
\hspace{.2in} \{ \varphi_{n} \} , \hspace{.2in} k \ra \infty.
$$\HS
\end{lem}

The following definition gives a generalization of asymptotic 
expansions of Poincar\'{e} type and power series type suitable 
for our needs.

\begin{defi}\label{defi:generalized-asymptotics}
A function $f:A \ra \C$ is said to have an {\it asymptotic 
expansion of generalized Poincar\'{e} type} (w.r.t.\ the 
asymptotic sequence $\{ \varphi_{n}(k) \}$) if we for all 
nonnegative integers $N$ have
\begin{equation}\label{eq:generalized-Poincare-type}
f(k) = \left[ \sum_{j=1}^{M} e^{2\pi ik q_{j}} \sum_{m=0}^{N} 
c_{m}^{j} \varphi_{m}(k) \right] + o(\varphi_{N}(k))
\end{equation}
in the limit $k \ra \infty$, where $M$ is a non-negative 
integer, $q_{j} \in [0,1[$ and $c_{m}^{j} \in \C$. If, 
moreover, $\varphi_{m}(k)=(\xi(k))^{\lambda_{m}}$, 
$\lambda_{m}$ fixed complex numbers, then we will say that $f$
has an {\it asymptotic expansion of generalized power series 
type}. The asymptotic expansion in 
(\ref{eq:generalized-Poincare-type}) is called {\it trivial} if
$c_{m}^{j}=0$ for all $j$ and $m$ (this case being equivalent 
to the case $M=0$, i.e.\ the case where the sum
$\sum_{j=1}^{M}$ is put equal to zero).  
\end{defi}

By the notation (\ref{eq:asymptotic-notation}), the identities
(\ref{eq:generalized-Poincare-type}) can also be written
$$
f(k) \sim \sum_{m=0}^{\infty} 
\left( \sum_{j=1}^{M} e^{2\pi ik q_{j}} c_{m}^{j} \right)
\varphi_{m}(k);\hspace{.2in} \{ \varphi_{m} \},
\hspace{.2in} k \ra \infty.
$$
We immediately get that the set of functions $f : A \to \C$ 
having an asymptotic expansion of generalized Poincar\'{e} type
w.r.t.\ a fixed asymptotic sequence is a subspace of the 
$\C$--vector space of all functions $f : A \to \C$. 
(The zero function $0:\N \ra \C$ has trivially a trivial 
asymptotic expansion of generalized Poincar\'{e} type (w.r.t.\
any asymptotic sequence).)

It is natural not to destinguish between asymptotic expansions
of generalized Poincar\'{e} type which can be obtain from each
other by trivial means. We will thus say that two such 
asymptotics are equivalent if the one can be obtained from the 
other by a finite sequence of the following operations: i) 
interchange two terms in the sum $\sum_{j=1}^{M}$ in 
(\ref{eq:generalized-Poincare-type}), ii) remove or add a term
in $\sum_{j=1}^{M}$ with $c_m^j=0$ for all $m$ and iii) collect
two terms with the same $q_j$--value. These operations also 
lead to a kind of minimal representation for each equivalence 
class of asymptotic expansions of generalized Poincar\'{e} 
type, namely if we have such an asymptotics we will say that it
is on minimal form if it has no terms with $c_m^j=0$ for all 
$m$ and if the $q_j$'s are mutually different. Thus, up to the
summation order in $\sum_{j=1}^M$ each equivalence class has a
unique representative on minimal form. In particular, a trivial
expansion is on minimal form if and only if $M=0$ in 
(\ref{eq:generalized-Poincare-type}). We shall not destinguish
between an asymptotic expansion of generalized Poincar\'{e} 
type and its equivalence class.  

The following theorem extends the uniqueness property for 
asymptotic expansions of Poincar\'{e} type, 
\reflem{lem:Poincare}, to asymptotic expansions of generalized 
Poincar\'{e} type. The proof is very short due to an idea of 
Pieter Moree.

\begin{thm}\label{thm:uniqueness}
Let $f : A \to \C$ be an arbitrary function. Then $f$ has at 
most one asymptotic expansion of generalized Poincar\'{e} type 
w.r.t.\ a given asymptotic sequence. The zero function has only
the trivial asymptotic expansions of generalized Poincar\'{e} 
type.
\end{thm}

\noindent Said in another way: If $f : A \to \C$ has an 
asymptotic expansion of generalized Poincar\'{e} type as in
(\ref{eq:generalized-Poincare-type}) and if this expanion is on
minimal form, then $M$ and the $q_{j}$'s and $c_{m}^{j}$'s are 
uniquely determined by $f$.

\begin{prf}
Since a linear combination of asymptotic expansions of
generalized Poincar\'{e} type w.r.t.\ a fixed asymptotic 
sequence is again an asymptotic expansion of generalized 
Poincar\'{e} type w.r.t.\ this asymptotic sequence it is enough
to consider the case where $f=0$. Assume therefore that we for
all $N \geq 0$ have
\begin{equation}\label{eq:UA}
0= \left[ \sum_{j=1}^{M} e^{2\pi ik q_{j}} \sum_{m=0}^{N}
c_{m}^{j} \varphi_{m}(k) \right] +o(\varphi_{N}(k))
\end{equation}
in the limit $k \ra \infty$, where the $q_{j}$'s are pairwise 
different numbers in $[0,1[$ and $c_{m}^{j} \in \C$. By 
multiplying by $e^{-2\pi i k q_{1}}$ we can assume that 
$q_{1}=0$ and $q_{j} \in ]0,1[$, $j=2,3,\ldots,M$. We have to 
show that all the $c_{m}^{j}$'s are zero. If $M=1$ this 
trivially follows from \reflem{lem:Poincare}, so we can assume
that $M>1$. By putting $N=0$ we see that 
$g(k):=c_{1}+\sum_{j=2}^{M} c_{j}e^{2\pi ik q_{j}}$ converges 
to zero as $k \ra \infty$, where we have put 
$c_{j}=c_{0}^{j}$. Thus $|g(k)|^{2} \ra 0$ as $k \ra \infty$ 
implying that $\frac{1}{L}\sum_{k=0}^{L-1} |g(k)|^{2}$ 
converges to zero as $L \ra \infty$. Here
$$
|g(k)|^{2} = |c_{1}|^{2} + |c_{2}|^{2} + \ldots + |c_{M}|^{2}
+ \sum_{\stackrel{i,j \in \{1,\ldots,M\}}{i \neq j}} 
c_{i}\bar{c}_{j}e^{2\pi i k (q_{i}-q_{j})},
$$
and therefore
$$
\frac{1}{L}\sum_{k=0}^{L-1} |g(k)|^{2} = 
|c_{1}|^{2} + |c_{2}|^{2} + \ldots + |c_{M}|^{2} + \frac{1}{L}
\sum_{\stackrel{i,j \in \{1,\ldots,M\}}{i \neq j}} 
c_{i}\bar{c}_{j}
\frac{1-e^{2\pi i L (q_{i}-q_{j})}}{1-e^{2\pi i (q_{i}-q_{j})}}
$$
so $\frac{1}{L}\sum_{k=0}^{L-1} |g(k)|^{2}$ converges to
$|c_{1}|^{2} + |c_{2}|^{2} + \ldots + |c_{M}|^{2}$ as 
$L \ra \infty$, hence $c_{j}=0$ for all $j$.

The proof is finalized by induction: Assume that $c_{m}^{j}=0$ 
for all $j$ and all $m=0,1,\ldots,N-1$ and get from 
(\ref{eq:UA}) that 
$c_{N}^{1} + \sum_{j=2}^{M} c_{N}^{j}e^{2\pi ik q_{j}}$ 
converges to zero as $k \ra \infty$. As before we then have
$c_{N}^{j}=0$ for all $j$.
\end{prf}

Let us see how the above fits into the theory of quantum 
invariants. Therefore, let $f : A \ra \C$ and assume that there
exist constants $q_{j} \in [0,1[$, $D_{j} \in \Q$ and
$b_{j} \in \C$ for $j=1,2,\ldots,n$, and $a_{m}^{j} \in \C$ for
$j=1,2,\ldots,n$ and $m = 1,2,\ldots$ such that we for all 
$N \in \{0,1,2,\ldots\}$ have
\begin{equation}\label{eq:asymptotics-quantum1}
f(k) = \sum_{j=1}^{M} b_{j} e^{2\pi ikq_{j}} k^{D_{j}}
\left( 1 + \sum_{m=1}^{N} a_{m}^{j} k^{-\theta_{m}} \right) 
+ o(k^{D-\theta_{N}})
\end{equation}
in the limit $k \ra \infty$, where 
$D = \max\{D_{1},\ldots,D_{M}\}$ and $\{\theta_{m}\}$ is a 
strictly increasing sequence in $\frac{1}{b}\Z_{>0}$, where $b$
is the least positive integer such that $bD_{j} \in \Z$ for all
$j$, compare with \refconj{conj:AEC'}. By adding terms of the 
form $0k^{-l/b}$ we can assume that 
$\{\theta_{m}\}_{m=1}^{\infty} = \frac{1}{b}\Z_{>0}$. It is 
straightforward to see then that we for any 
$N \in \{0,1,2,\ldots\}$ have 
\begin{equation}\label{eq:asymptotics-quantum2}
f(k) = \sum_{j=1}^{M} e^{2\pi ikq_{j}}
\sum_{m=0}^{N} c_{m}^{j} k^{\lambda_{m}} + o(k^{\lambda_{N}})
\end{equation}
in the limit $k \ra \infty$, where $c_{m}^{j} = 0$ for 
$m=0,1,\ldots,e_{j}-1$, $e_{j}=b(D-D_{j})$, and 
$c_{m}^{j}=b_{j}a_{m-e_{j}}^{j}$ for $m=e_{j},e_{j}+1,\ldots$, 
and $\lambda_{m}=D-\theta_{m}$, $m=0,1,\ldots$. Here we have 
put $a_{0}^{j}=1$ and $\theta_{0}=0$. Thus $f$ has an 
asymptotic expansion of generalized power series type with 
$\xi=\id$ and $\{\lambda_{m}\}$ a strictly decreasing sequence
equal to $\frac{1}{b}\Z_{\leq bD}$ as a set.

If oppositely $f: A \ra \C$ has an asymptotic expansion like in
(\ref{eq:asymptotics-quantum2}) with $\{\lambda_{m}\}$ a 
strictly decreasing sequence equal to $\frac{1}{b}\Z_{\leq a}$
as a set, $a \in \Z$ and $b$ some positive integer, then it is 
straightforward to bring the expansion on a form like in 
(\ref{eq:asymptotics-quantum1}). Simply put 
$D_{j} = \lambda_{e_{j}}$, where
$$
e_{j} = \min\{ m \in \Z_{\geq 0} | c_{m}^{j} \neq 0 \}
$$
and $b_{j}=c_{e_{j}}^{j}$, $a_{m}^{j} = c_{m+e_{j}}^{j}/b_{j}$ 
and $\theta_{m} = D - \lambda_{m+e}$, where
$e=\min\{e_1,\ldots,e_M\}$. (We assume here that there are no 
$j$'s with $c_{m}^{j} =0$ for all $m$. If $c_{m}^{j}=0$ for all
$m$, simply ignore that term.)

Because of the above, we will pay special attention to the 
following asymptotic expansions of generalized power series 
type.

\begin{defi}\label{defi:quantum-asymptotics} 
A function $f: A \to \C$ is said to have an {\it asymptotic 
expansion of quantum type} w.r.t.\ $b \in \N$ or a quantum 
$b$--expansion for short if $f$ has an asymptotic expansion of 
generalized power series type as in 
\refdefi{defi:generalized-asymptotics} with $\xi=\id$ and
$\{\lambda_{m}\} \subseteq \frac{1}{b}\Z$. The expansion is 
said to be on minimal form if the expansion is on minimal form
as an asymptotic expansion of generalized power series type.
\end{defi}

We can assume w.l.o.g.\ that
$\{\lambda_{m}\}_{m=0}^{\infty}= \frac{1}{b}\Z_{\leq a}$, where
$a=\lambda_{0}b$. Moreover, it is obvious that a quantum 
$b$--expansion is also a quantum $b'$--expansion for any $b'$ 
divisible by $b$. We note that the set $\mQ$ of functions 
$f : A \ra \C$ having an asymptotic expansion of quantum type 
is a subalgebra of the $\C$--algebra of all complex functions 
on $A$: It is obvious that a scalar times a $f \in \mQ$ is 
again in $\mQ$. If $f_{i} \in \mQ$ has a quantum 
$b_{i}$--expansion, $i=1,2$, then both $f_{1}$ and $f_{2}$ have
a quantum $b$--expansion, where $b=\lcm(b_{1},b_{2})$. Now the
sum and product of quantum  $b$--expansions are again quantum 
$b$--expansions, so $f_{1}+f_{2}$ and $f_{1}f_{2}$ are both in 
$\mQ$.

We can now restate the AEC, \refconj{conj:AEC'}, in the 
following form:

\begin{conj}\label{conj:AEC}
For any $3$--manifold $X$ the function 
$r \mapsto \tau_{r}^{G}(X)$, $\Z_{\geq h^{\vee}} \to \C$, has 
an asymptotic expansion of quantum type as in 
{\em (\ref{eq:asymptotics-quantum2})} with 
$\{q_{j}\}_{j=1}^{M}$ being the image set of the Chern--Simons
functional on the moduli space of flat $G$--connections on $X$.
\end{conj}

By combining the AEC and \refconj{conj:d'} we get the stronger
conjecture that $r \mapsto \tau_{r}^{G}(X)$ has a quantum 
$2$--expansion. Note that the fact that $\mQ$ is multiplicative
stable is well in accordance with the AEC and the fact that the
RT--invariants behave multiplicatively w.r.t.\ connected sums, 
see e.g.\ \cite[Chap.~II]{Turaev}.

We end this section by giving miscellaneous facts about 
asymptotic expansions of quantum type. First we make the 
trivial remark that not every function $f: A \to \C$ has 
an asymptotic expansion of quantum type. If namely $f$ has such
an expansion there exists a positive constant $C$ and a 
rational number $\lambda$ such that
$$
|f(k)| \leq C k^{\lambda}
$$
for all $k \in A$. Thus it is easy to give examples of classes 
of functions which do not have an asymptotic expansion of 
quantum type. One such class is given by $f(k)=P(k)e^{Q(k)}$, 
where $P(k)$ is a Laurent polynomial in $k$ (different from 
zero), and $Q(k)$ is a polynomial in $k$ of degree $\geq 1$ 
with positive leading coefficient.

\refthm{thm:uniqueness} shows that a function having an 
asymptotic expansion of quantum type determines the quantities 
building up the expansion (if the expansion is on minimal 
form). The opposite is obviously not the case. That is, two 
different functions $f,g: A \to \C$ can have the same 
asymptotic expansion of quantum type. If namely $f$ has an 
expansion like in (\ref{eq:asymptotics-quantum2}), if $g=f+h$,
and if $h = o(k^{\lambda_{m}})$ as $k \ra \infty$ for all $m$,
then $g$ has the same asymptotic expansion of quantum type as 
$f$. As an example one can let $h(k)=\exp(\beta k)$ with 
$\beta$ a negative real number or more generally 
$h(k)=P(k)\exp(Q(k))$, where $P(k)$ is a (non-zero) Laurent 
polynomial and $Q(k)$ is a polynomial of degree $\geq 1$ with 
negative leading coefficient. In particular, such a function 
$h$ has a trivial asymptotic expansion of quantum type 
(without being zero).

Consider a function $g: A \to \C$ given by an expression of 
the form
\begin{equation}\label{eq:generic-type-alternative}
g(k)=\sum_{j=1}^{M} e^{2\pi ik q_{j}}k^{D_{j}} b_{j} g_{j}(k),
\end{equation}
where $b_{j} \in \C \sm \{0\}$, $D_{j} \in \Q$ and 
$q_{j} \in [0,1[$, and $g_{j}(k)$ are polynomials in 
$k^{-1/b}$ with constant term $1$, where $b$ is the least 
positive integer such that $bD_{j} \in \Z$ for all $j$. Then 
$f$ has obviously an asymptotic expansion of quantum type. 
More generally this is true if the functions $g_{j}$ satisfy 
asymptotic identities
\begin{equation}\label{eq:gfunction-asymptotics}
g_{j}(k) \sim \sum_{m=0}^{\infty} a_{m}^{j}k^{-m/b};
\hspace{.2in} \{ k^{-m/b} \} , \hspace{.2in} k \ra \infty
\end{equation}
with $a_{0}^{j}=1$ for all $j$. A particular example of the 
above is given by the case where
$$
g_{j}(k) \sim \sum_{m=0}^{\infty} \tilde{a}_{m}^{j}k^{-m};
\hspace{.2in} \{ k^{-m} \} , \hspace{.2in} k \ra \infty.
$$
If we put $a_{bm}^{j}=\tilde{a}_{m}^{j}$ for 
$m \in \Z_{\geq 0}$ and $a_{m}^{j}=0$ for $m$ not divisible by
$b$ we get (\ref{eq:gfunction-asymptotics}).

For any function $f: A \to \C$ having an asymptotic expansion 
of quantum type there exists a function $g: A \to \C$ given as
in (\ref{eq:generic-type-alternative}) and satisfying that $f$
and $g$ have the same asymptotic expansion of quantum type. In
fact, assume that the expansion of $f$ is non-trivial and 
assume it is on minimal form (to avoid terms with 
$c_{m}^{j}=0$ for all $m$). By our comments in connection to 
(\ref{eq:asymptotics-quantum2}) we can find constants $b_{j}$,
$D_{j}$ and $a_{m}^{j}$ and a sequence $\{\theta_{m}\}$ such 
that (\ref{eq:asymptotics-quantum1}) is satisfied.  According 
to \cite[Sec.~1.9]{Olver} we can find functions 
$g_{j}: A \to \C$ such that (\ref{eq:gfunction-asymptotics}) is
satisfied. (In fact, if $0<a<1$, $0< \sigma < \pi$ and
$\Sec=\{ z \in \C : |z| \geq a, |\Arg(z)| \leq \sigma \}$
then we can find analytic functions $f_{j} : \Sec \to \C$ such 
that
$$
f_{j}(z) \sim \sum_{m=0}^{\infty} a_{m}^{j}z^{-m};
\hspace{.2in} \{ z^{-m} \} , \hspace{.2in} 
z \ra \infty \hspace{.05in} \text{in} \hspace{.05in} \Sec.
$$
By letting $g_{j}(k)=f_{j}(k^{1/b})$, $k \in \N$, we get the
desired functions.) Now, if we let $g: A \to \C$ be given by
(\ref{eq:generic-type-alternative}) with the $b_{j}$'s, 
$D_{j}$'s and the $q_{j}$'s coming from the expansion 
(\ref{eq:asymptotics-quantum1}) of $f$, then $g$ has an 
asymptotic expansion of quantum type equal to that of $f$. The 
functions $g_{j}$ are of course not uniquely determined by the
numbers $a_{m}^{j}$. We see that the function $g$ can even be 
chosen to be analytic in any sector shaped region in the 
complex plane containing the positive integers. However, from
an asymptotic point of view we do of course not gain anything.
The asymptotic behaviour is solely determined by the 
coefficients $c_{m}^{j}$, the numbers $q_{j}$ and the sequence 
$\{\lambda_{m}\}$.

Assume $f: A \to \C$ has a non-trivial asymptotic expansion of
quantum type as in (\ref{eq:asymptotics-quantum1}), and let
$J=\{\;j \in \{1,\ldots,M\}\;|\; d_{j}=D\;\}$. Then
\begin{equation}\label{eq:leading}
\sum_{j \in J} b_{j} e^{2\pi ik q_{j}} k^{D}
\end{equation}
is called the leading term of the asymptotic expansion of $f$.
For a $3$--manifold $X$ with discrete moduli space of flat
$G$--connections the semiclassical approximation of 
$Z_{k}^{G}(X)$, i.e.\ the right-hand side of 
(\ref{eq:semiclassical-G'}), is (conjecturally) the leading 
order large $k$ asymptotics of $Z_k^G(X)$ if and only if all 
the $D_{A}$'s are equal. In general the leading asymptotics is
given by the sum of terms with maximal $D_A$.

\section{The Reshetikhin--Turaev invariants of Seifert
manifolds for $\SU(2)$}\label{sec-The-Reshetikhin}

\noindent We use the notation introduced by Seifert in his 
classification results for the Seifert manifolds (or rather 
fibrations), see \cite{Seifert1}, \cite{Seifert2}, 
\cite[Sect.~2]{Hansen2}. That is, 
$(\ep;g\;|\;b;(\alpha_{1},\beta_{1}),\ldots,
(\alpha_{n},\beta_{n}))$ 
is the Seifert manifold with orientable base of genus 
$g \geq 0$ if $\ep=\os$ and non-orientable base of genus $g>0$
if $\ep=\ns$ (where the genus of the non-orientable connected 
sum $\# ^{k} \R \text{P}^{2}$ is $k$). (In \cite{Seifert1}, 
\cite{Seifert2}
$(\ep;g\;|\;b;(\alpha_{1},\beta_{1}),\ldots,
(\alpha_{n},\beta_{n}))$ is
denoted $(\tO,\ep;g\;|\;b;\alpha_{1},\beta_{1};\ldots;
\alpha_{n},\beta_{n})$
but we leave out the $\tO$, since we are only dealing with
oriented Seifert manifolds.) The pair $(\alpha_{j},\beta_{j})$
of coprime integers is the oriented Seifert invariant of the 
$j$'th exceptional (or singular) fiber. We have 
$0< \beta_{j}<\alpha_{j}$. The integer $-b$ is equal to the
Euler number of the Seifert fibration $(\ep;g\;|\;b)$ (which is
a locally trivial $S^{1}$--bundle). The sign is chosen so that
the Euler number of the spherical (or unit) tangent bundle over
an orientable surface $\Sigma$ is equal to the Euler
characteristic of $\Sigma$, see 
\cite[Chap.~1 and 4]{Montesinos}, \cite[Sec.~3]{Scott}. More 
generally, the Seifert Euler number of
$(\ep;g\;|\;b;(\alpha_{1},\beta_{1}),\ldots,
(\alpha_{n},\beta_{n}))$
is $E=-\left(b+\sum_{j=1}^{n} \beta_{j}/\alpha_{j} \right)$.
We note that except for a small class of Seifert manifolds,
the Seifert invariants are actually topological invariants, 
i.e.\ they classify the Seifert manifolds up to orientation 
preserving homeomorphisms, see 
e.g.\ \cite[Chap.~5 Theorem 6 p.~97 and Sec.~5.4]{Orlik} or
\cite[Theorem 5.1 p.~32]{JankinsNeumann}. The exceptions are
the lens spaces, the prism manifolds 
$(\os;0\;|\;b;(2,1),(2,1),(\alpha,\beta))$, and the manifolds
$(\ns;2\;|0)$, $(\os;0\;|\;-2;(2,1),(2,1),(2,1),(2,1))$, 
$(\ns;1\;|\;b)$ with $b \neq 0$ and 
$(\ns;1\;|\;b;(\alpha,\beta))$. These exceptions constitute a
proper subclass of the 'small Seifert manifolds' listed in 
\cite[pp~91-92]{Orlik}. The lens spaces are the Seifert 
manifolds with $\ep=\os$, $g=0$
and $n \leq 2$. A Seifert manifold which is not small is 
sometimes called a large Seifert manifold, see 
\cite[p.~92]{Orlik}.

For the convenience of the reader we will state all results in
the following both in terms of the Seifert invariants and in
terms of the so-called non-normalized Seifert invariants due to
W.\ D.\ Neumann. For a Seifert manifold $X$ with non-normalized
Seifert invariants
$\{\ep;g;(\alpha_{1},\beta_{1}),\ldots,$ 
$(\alpha_{n},\beta_{n})\}$
the invariants $\ep$ and $g$ are as above. The
$(\alpha_{j},\beta_{j})$'s are here pairs of coprime integers 
with $\alpha_{j} >0$ but not necessarily with 
$0< \beta_{j}<\alpha_{j}$. These pairs are actually not 
invariants of the fibration $X$. Thus 
$\{\ep;g;(\alpha_{1},\beta_{1}),\ldots, 
(\alpha_{n},\beta_{n})\}$ 
and
$\{\ep;g;(\alpha_{1}',\beta_{1}'),\ldots, 
(\alpha_{m}',\beta_{m}')\}$
are two sets of non-normalized Seifert invariants of $X$ if and
only if the set of pairs $(\alpha_{j}',\beta_{j}')$ can be 
obtained from the set of pairs $(\alpha_i,\beta_i)$ by a finite
number of the following two moves: i) add or delete a pair 
$(1,0)$, ii) replace each $(\alpha_i,\beta_i)$ by
$(\alpha_i,\beta_i + K_i\alpha_i)$ provided $\sum K_i =0$. For
details, see \cite{JankinsNeumann}. The Seifert Euler number of
$X$ is $-\sum_{i=1}^{n} \beta_{i}/\alpha_{i}$ (which is an 
invariant of the Seifert fibration $X$). The main advantage 
obtained by working with the non-normalized invariants is to 
make more `symmetric' expressions, since the invariant $b$ is 
treated formally as a Seifert invariant $(1,b)$ of an 
exceptional fiber. Thus the Seifert manifold 
$(\ep;g\;|\;b;(\alpha_{1},\beta_{1}),\ldots,
(\alpha_{n},\beta_{n}))$ is equal to
$\{\ep;g;(1,b),(\alpha_{1},\beta_{1}),\ldots,
(\alpha_{n},\beta_{n})\}$.
Of course the number of exceptional fibers is a constant for a
Seifert manifold. This number can be read off directly from the
(normalized) Seifert invariants, but can also be read off from
a set of non-normalized invariants as the number of pairs 
$(\alpha_{j},\beta_{j})$ with $\alpha_{j} > 1$.

In \cite{JankinsNeumann} the authors work with a more general 
class of oriented fibered spaces, denoted Generalized Seifert 
fibrations. According to \cite[Theorem 5.1]{JankinsNeumann} 
these spaces are (up to orientation preserving homeomorphism)
the classical oriented Seifert manifolds as considered here and
connected sums of lens spaces (considering $S^{1} \times S^{2}$
as a lens space). Since the RT--invariant of a connected sum
of $3$--manifolds is the product of the RT--invariants of
the $3$--manifolds in that connected sum (up to a normalization
factor), the results of \cite{Jeffrey1} actually show the
AEC for these connected sums of lens spaces and $G=\SU(2)$.

To state the next theorem we need some notation. For a pair of
coprime integers $c,d$ with $|c| \geq 1$ the Dedekind sum is 
given by
\begin{equation}\label{eq:Dedekind-sum}
\s (d,c)= \frac{1}{4|c|} \sum_{j=1}^{|c|-1} \cot\frac{\pi j}{c}
\cot \frac{\pi d j}{c}
\end{equation}
for $|c|>1$ and $\s (d,\pm 1)=0$. We refer to 
\cite{RademacherGrosswald} for a comprehensive description of 
this function. The Dedekind symbol is given by
\begin{equation}\label{eq:Dedekind-symbol}
\dS (d/c) = 12\sign(c)\s (d,c).
\end{equation}
Multi-indices are denoted by an underline (e.g.\ $\un$). For
$\underline{k}=(k_{1},\ldots,k_{n}), 
\underline{l} = (l_{1},\ldots,l_{n}) \in \Z^{n}$, 
$\underline{k} < \underline{l}$ if and only if $k_{j} < l_{j}$ 
for all $j=1,\ldots,n$. We let $\uone =(1,\ldots,1)$. For
$\underline{k}=(k_{1},\ldots,k_{n}) \in \Z_{+}^{n}$ we write
$\sum_{\um =\uzero}^{\underline{k}}$ for
$\sum_{m_{1}=0}^{k_{1}}\ldots\sum_{m_{n}=0}^{k_{n}}$ etc.
Let $a_{\os}=2$ and $a_{\ns}=1$. Given pairs of coprime 
integers $\alpha_{j},\beta_{j}$ we choose integers 
$\rho_{j},\sigma_{j}$ such that 
$\alpha_{j}\sigma_{j}-\beta_{j}\rho_{j}=1$.

\begin{thm}[{\cite[Theorem 8.4]{Hansen2}}]
\label{thm:RT-invariants-Seifert-manifolds}
The RT--invariant at level $r-2$ of the Seifert manifold $X$ 
with {\em (}normalized\,{\em )} Seifert invariants
$(\ep;g\;|\;b;(\alpha_{1},\beta_{1}),\ldots,
(\alpha_{n},\beta_{n}))$ is
\begin{eqnarray}\label{eq:RT-invariants-Seifert-manifolds}
\tau_{r}(X) &=& (-1)^{a_{\ep}g} 
\frac{r^{a_{\ep}g/2-1}}{2^{a_{\ep}g/2-1}} 
\frac{1}{ \sqrt{\mA}} e^{i\frac{3\pi}{4}(1-a_{\ep})\sgE } \\  
&& \times
\exp \left( \frac{i \pi}{2r} \left[ 3(a_{\ep}-1)\sgE - E - 
\sum_{j=1}^{n} \dS \left(\frac{\beta_{j}}{\alpha_{j}}\right) 
\right] \right) Z(X;r), \nonumber
\end{eqnarray}
where $\mA=\prod_{j=1}^{n} \alpha_{j}$ and
\begin{equation}\label{eq:Zfunction-Seifert}
Z(X;r) = \left(\frac{i}{2}\right)^{n} \sum_{\gamma =1}^{r-1} 
\frac{(-1)^{\gamma a_{\ep}g} h(\gamma)}
{\sin ^{n+a_{\ep}g-2} \left( \frac{\pi}{r} \gamma \right) },
\end{equation}
where 
\begin{eqnarray}\label{eq:hfunction}
h(\gamma) &=& \exp \left( \frac{i \pi}{2r} E \gamma^{2} \right)
\musum \nsum \muprod \\
&& \times \gex 
\exp \left( - \frac{i \pi}{r} \gamma \sum_{j=1}^{n} 
\frac{ 2rn_{j} + \mu_{j}}{\alpha_{j}} \right). \nonumber
\end{eqnarray}
The RT--invariant at level $r-2$ of the Seifert manifold $X$
with non-normalized Seifert invariants
$\{\ep;g;(\alpha_{1},\beta_{1}),\ldots,
(\alpha_{n},\beta_{n})\}$
is given by the same expression. 
\end{thm}

The theorem is also valid in case $n=0$, where one has to put 
$\sum_{j=1}^{n} \dS(\beta_{j}/\alpha_{j})=0$ and $\mA=1$ in 
(\ref{eq:RT-invariants-Seifert-manifolds}) and let 
$h(\gamma) = \exp \left( \frac{i \pi}{2r} E \gamma^{2} \right)$
in (\ref{eq:Zfunction-Seifert}). The reason for including the
factor $(i/2)^{n}$ in the expression for $Z(X;r)$ is that 
then $Z(X;r)$ does not depend on the choice of non-normalized 
Seifert invariants for $X$. In fact, $Z(X;r)$ only depends on 
the $\beta_j$'s through the Seifert Euler number $E$, and if we
add a pair $(\alpha,\beta)=(1,0)$ each of the terms in the 
sum-expression for $Z(X;r)$ changes by a factor $-2i$, hence
$Z(X;r)$ does not change.

\section{The analytic asymptotic expansions of the 
Reshetikin-Turaev invariants of Seifert manifolds}
\label{sec-The-analytic}

\noindent Let $X$ be a Seifert manifold with (normalized) 
Seifert invariants
$(\ep;g\;|\;b;(\alpha_{1},\beta_{1}),\ldots,
(\alpha_{n},\beta_{n}))$
or non-normalized Seifert invariants
$\{\ep;g;(\alpha_{1},\beta_{1}),\ldots,
(\alpha_{n},\beta_{n})\}$ 
with $\ep =\os$ or $g$ even. Let $\beta_{0}=b$ if $X$ is given
by normalized Seifert invariants and $\beta_0=0$ otherwise. 
We will below present detailed expressions for the asymptotic 
expansions of the RT--invariants of $X$. It will follow that 
all parts of the asymptotics are expressible by the Seifert 
invariants.  

By the Theorems~\ref{thm:asymptotics-Enot0} and 
\ref{thm:asymptotics-E0} below, $Z(X;r)$ has an asymptotic 
expansion of the form (\ref{eq:AEC'}), and according to 
\refthm{thm:RT-invariants-Seifert-manifolds}
$$
\tau_{r}(X)=b r^{a_{\ep}g/2-1} \exp\left( \frac{ia}{r} \right) 
Z(X;r),
$$
where $a \in \R$ and $b \in \C \sm \{0\}$, so $\tau_{r}(X)$
has also an asymptotic expansion of the form (\ref{eq:AEC'}).

Before stating Theorems~\ref{thm:asymptotics-Enot0} and
\ref{thm:asymptotics-E0} we need some preliminaries. 
We treat the case $n=0$ as the case $n=1$ with 
$(\alpha_{1},\beta_{1})=(1,0)$.
For $E \neq 0$ we let $z_{\st}:\Z \times \Z^{n} \ra \R$ be 
given by
\begin{equation}\label{eq:SeifertStationary}
z_{\st}(m,\un) = -\frac{2}{E} \left(m-\nalphasuma \right).
\end{equation}
Moreover, let
\begin{eqnarray}
q_{(m,\un)} &=& \rhoalphasum n_{j}^{2} - \frac{1}{4} E 
z_{\st}(m,\un)^{2} \pmod{\Z}, 
\hspace{.2in} (m,\un) \in \Z \times \Z^{n}, \label{eq:CSa} \\
q_{(l,\un')} &=& -\frac{1}{4}\beta_{0}l^{2} + \sum_{j=1}^{n} 
\left( \frac{\rho_{j}}{\alpha_{j}} {n_{j}'}^{2} 
- \frac{1}{4} \sigma_{j}\beta_{j}l^{2} \right) \pmod{\Z},
\hspace{.2in} (l,\un') \in M, \label{eq:CSb}
\end{eqnarray}
where
$
M = \left\{\;\left.(l,\un') \in \Z \times 
\left( \frac{1}{2} \Z \right)^{n} \; \right| \; 
\un' \in \Z^{n} + \frac{1}{2}l\ubeta,\;\right\}.
$
We let
\begin{equation}\label{eq:S-set}
\mS = \{ \;\un \in \Z^{n} \; | 
\; \uzero \leq \un \leq \ualpha -\uone \;\}
\end{equation}
and
\begin{eqnarray}
\mI_{1} &=& \left\{\;(m,\un) \in \Z \times \mS \; | \; 
 0 < z_{\st}(m,\un) < 1 \;\right\}, \label{eq:I1} \\
\mI_{2} &=& \left\{\;(l,\un') \in M  \; \left| \; 
l \in \{0,1\},\,\uzero \leq \un' \leq \frac{1}{2} \ualpha 
\;\right.\right\}. \label{eq:I2}
\end{eqnarray}
Above $q_{(m,n)}$ and $\mI_{1}$ are of course only defined for
$E \neq 0$. We need a partition
$$
\mI_{2}^{a} = \left\{ \;(l,\un') \in \mI_{2} \; \left| \; 
\exists \umu' \in \B : \sum_{j=1}^{n} 
\frac{\mu_{j}' n_{j}'}{\alpha_{j}} \in \Z \;\right.\right\}
$$
and $\mI_{2}^{b}=\mI_{2} \sm \mI_{2}^{a}$ of $\mI_2$. We will 
need to write $1$ and $2$ in a sophisticated way, namely let
\begin{equation}\label{eq:Sym}
\Sym_{\pm}(x) = \left\{ \begin{array}{ll}
1 \hspace{.2in} & ,x \neq 0, \\
2 \hspace{.2in} & ,x=0, \end{array} \right.
\end{equation}
and
\begin{equation}\label{eq:SymZ}
\Sym_{\Z_{\pm}}(x) = \left\{ \begin{array}{ll}
1 \hspace{.2in} & ,x \in \R \sm \frac{1}{2} \Z, \\
2 \hspace{.2in} & ,x \in \frac{1}{2} \Z. \end{array} \right.
\end{equation}
These functions have a group theoretical explanation, see
\cite[p.~36]{Rozansky3}. However, this can be neglected here. 
Basically these functions are used to keep track of how many 
times certain terms contribute to certain sums being part of 
the asymptotic expansion of $Z(X;r)$.

As already indicated above we have to consider the cases $E=0$
and $E \neq 0$ separately. \refthm{thm:asymptotics-Enot0} 
concerns the case $E \neq 0$ while \refthm{thm:asymptotics-E0}
takes care of the case $E=0$.

\begin{thm}\label{thm:asymptotics-Enot0}
Assume the Seifert Euler number $E \neq 0$.
Then
$$
Z(X,r)=Z_{\pol}(X;r)+Z_{\inte}(X;r),
$$
where $Z_{\pol}(X;r)$ is a sum of residues, while the other
term $Z_{\inte}(X;r)$ is given by certain limits of sums of
certain integrals along contours in the complex plane. If 
$n+a_{\ep}g-2 \leq 0$ then $Z_{\pol}(X;r)=0$. For 
$n+a_{\ep}g-2 >0$ we have
\begin{equation}\label{eq:Zpolar}
Z_{\pol}(X;r) = \sum_{(l,\un') \in \mI_{2}} b_{(l,\un')} 
\exp(2\pi i r q_{(l,\un')}) r \left( Z^{(l,\un')}_{0}(r) 
+ Z^{(l,\un')}_{1}(r) \right),
\end{equation}
where
\begin{equation}\label{eq:bfactor}
b_{(l,\un')} = 
\frac{2^{n}(-1)^{\left(n+\sum_{j=1}^{n}\sigma_{j}\right)l}}
{\prod_{j=1}^{n} \Sym_{\Z_{\pm}} 
\left( \frac{n_{j}'}{\alpha_{j}} \right) },
\end{equation}
$q_{(l,\un')}$ is given by {\em (\ref{eq:CSb})}, and the 
functions $Z^{(l,\un')}_{0}$ and $Z^{(l,\un')}_{1}$ are
given by
\begin{eqnarray}\label{eq:Zpolarreg}
Z^{(l,\un')}_{0}(r) &=& -\frac{\pi}{2} 
\Res_{z=0} \left\{ 
\frac{ \exp \left( \frac{i \pi r}{2} E z^{2} \right) }
{\sin^{n+a_{\ep}g-2}(\pi z)} \cot(\pi r z) \right. \\
&& \times \prod_{j=1}^{n} \left[ i 
\sin \left( 2 \pi \frac{\rho_{j}}{\alpha_{j}} n_{j}'\right) 
\cos \left( \pi \frac{z}{\alpha_{j}} \right) 
\sin \left( 2\pi r z \frac{n_{j}'}{\alpha_{j}} \right) 
\right. \nonumber \\
&& \left. \left. + 
\cos \left( 2 \pi \frac{\rho_{j}}{\alpha_{j}} n_{j}' \right) 
\sin \left( \pi \frac{z}{\alpha_{j}} \right) 
\cos \left( 2\pi r z \frac{n_{j}'}{\alpha_{j}} \right) 
\right] \right\} \nonumber
\end{eqnarray}
and
\begin{eqnarray}\label{eq:Zpolarsing}
Z^{(l,\un')}_{1}(r) &=& -\frac{\pi i}{(-2)^{n}}\mumsum 
\sum_{\stackrel{m \in \Z}{0 \leq m \leq |a|} }
\frac{\sign(a)}{\Sym_{\pm} (m) \Sym_{\pm} \left(m-|a| \right)} 
\nonumber \\
&& \times\Res_{z=0}
\left\{ \frac{ \exp \left[ \frac{i \pi r}{2} E z^{2} + 
2 \pi i r \sign(a) \left( m-|a| \right) z \right] }
{\sin^{n+a_{\ep}g-2}(\pi z)} \right. \nonumber \\
&& \times \left. \prod_{j=1}^{n} \mu_{j}'
\sin \left[ 2 \pi  \left( \frac{\rho_{j}}{\alpha_{j}} n_{j}' -
\frac{z}{2} \frac{\mu_{j}'}{\alpha_{j}} \right) \right] 
\right\}, \nonumber
\end{eqnarray}
where $a=a(\umu',\un')=\mualphammsum$. In particular, 
$Z_{0}^{(l,\un')}(r)=0$ if $g=0$. 

Let $\mI = \mI_{1} \cup \mI_{2}^{a}$, and let $k_{0}=0$ if $n$
is even and $k_{0}=1$ otherwise. Moreover, let 
$k_{2}=\min\{0,k_{1}\}$, where $k_{1}=(k_{0}-n-a_{\ep}g+2)/2$.
Then
$$
r^{-1/2} Z_{\inte}(X;r) \sim \sum_{k =k_{2}}^{\infty}
\left( \sum_{\gamma \in \mI} \exp( 2\pi i r  q_{\gamma} )
c_{k}^{\gamma} \right) r^{-k}; \hspace{.2in} \{ r^{-k} \},
\hspace{.2in} r \ra \infty,
$$
where the quantities $q_{\gamma}$ and $c_{k}^{\gamma}$ are 
given as follows: For $\gamma=(m,\un) \in \mI_{1}$, 
$q_{\gamma}$ is given by {\em (\ref{eq:CSa})}. Moreover,
\begin{equation}\label{eq:A3}
c_{k}^{\gamma} = (-1)^{n} e^{\frac{i\pi}{4}\sgE}
\sqrt{\frac{2}{|E|}}\frac{1}{k!}
\left( \frac{i}{2\pi E} \right)^{k}
\left. \partial_{z}^{(2k)} \frac{ \prod_{j=1}^{n}
\sin \left( \frac{\pi}{\alpha_{j}} (2\rho_{j}n_{j}-z) \right) }
{ \sin^{n+a_{\ep}g-2} (\pi z) } \right|_{z = z_{\st}(m,\un)}
\end{equation}
for $k \geq 0$ and $c_{k}^{\gamma}=0$ for $k<0$
{\em (}in case $k_{2}<0${\em )}. Here $z_{\st}(m,\un)$ is given
by {\em (\ref{eq:SeifertStationary})}. For 
$\gamma=(l,\un') \in \mI_{2}^{a}$, $q_{\gamma}$ is given by 
{\em (\ref{eq:CSb})}. For $k \in \{k_{1},k_{1}+1,\ldots\}$ we 
put $k'=2k+n+a_{\ep}g-2$ and have
\begin{equation}\label{eq:A3'}
c_{k}^{\gamma}= \frac{b_{(l,\un')}}{(-2)^{n}\pi^{n+a_{\ep}g-2}}
\frac{e^{\frac{i\pi}{4}\sgE}}{\sqrt{2\pi|E|}}
\frac{\Gamma\left(k+\frac{1}{2}\right)}{k'!}
\left( \frac{2i}{\pi E} \right)^{k} 
\sum_{ \umu' \in \B \; : \; \mualphammsum \in \Z }
\left. \partial_{z}^{(k')} f(z;l,\un',\umu') \right|_{z=0},
\end{equation}
where $b_{(l,\un')}$ is given by {\em (\ref{eq:bfactor})} and
$$
f(z;l,\un',\umu') = 
\left( \frac{\pi z}{\sin(\pi z)} \right)^{n+a_{\ep}g-2}
\prod_{j=1}^{n} \sin \left( \pi 
\frac{2\rho_{j} \mu_{j}' n_{j}' -z}{\alpha_{j}} \right).
$$
We note that $k_{1} \leq 1$, and $k_{1}=1$ if and only if 
$n=k_{0}=1$, $\ep=\os$ and $g=0$. In that case 
$c_{0}^{\gamma}=0$.
\end{thm}

We have in the above theorem focused on the asymptotic 
expansion of $Z(X;r)$ and have therefore chosen not to give an
explicit expression for the term $Z_{\inte}(X;r)$. Such an 
explicit expression can be found in the proof of 
\refthm{thm:asymptotics-Enot0}, see (\ref{eq:Zinte}) together 
with (\ref{eq:Zterms}) and (\ref{eq:Zint2-deltaeta}), but is 
not relevant here and is not very informative.

\begin{rem}
The above theorem is also true in case $n=0$ if we as usual 
put all products $\prod_{j=1}^{n}$ equal to $1$ and all sums 
$\sum_{j=1}^{n}$ equal to $0$ and if we make the following 
natural adjustments: The function $z_{\st}$ in
(\ref{eq:SeifertStationary}) and the 
$q$--function in (\ref{eq:CSa}) only depend on $m \in \Z$,
namely
\begin{eqnarray*}
z_{\st}(m) &=& -\frac{2m}{E}, \\
q_{m} &=& -\frac{m^{2}}{E} \pmod{\Z}.
\end{eqnarray*}
The index sets $\mI_{1}$, $\mI_{2}$ and $\mI_{2}^{a}$ become
\begin{eqnarray*}
\mI_{1} &=& 
\left\{ \;m \in \Z \; | \; 0 < z_{\st}(m) < 1 \;\right\}, \\
\mI_{2} &=& \mI_{2}^{a} = \{0,1\}.
\end{eqnarray*}
The $q$--function in (\ref{eq:CSb}) degenerates to
$q_{l} = -\frac{1}{4}\beta_{0}l^{2} \pmod{\Z}$, $l=0,1$. The 
sum
$$
\sum_{ \umu' \in \B \; : \; \mualphammsum \in \Z }
\left. \partial_{z}^{(k')} f(z;l,\un',\umu') \right|_{z=0}
$$
in (\ref{eq:A3'}) simply becomes
$$
\left. \partial_{z}^{(k')} 
\left( \frac{\pi z}{\sin(\pi z)} \right)^{a_{\ep}g-2}
 \right|_{z=0}.
$$
The coefficients $b_{(l,\un')}$ in (\ref{eq:bfactor}) become
$b_{l}=1$, $l=0,1$.
Finally the functions $Z^{(l,\un')}_{\nu}$, $\nu=0,1$,
degenerate to $Z_{1}^{l}(r)=0$ and
$$
Z^{l}_{0}(r) = -\frac{\pi}{2} 
\Res_{z=0} \left\{ 
\frac{ \exp \left( \frac{i \pi r}{2} E z^{2} \right) }
{\sin^{a_{\ep}g-2}(\pi z)} \cot(\pi r z)\right\}
$$
for $l \in \mI_{2}$, both actually independent of $l$. 
Thus $$
Z_{\pol}(X;r) = -\pi r 
\Res_{z=0} \left\{ 
\frac{ \exp \left( \frac{i \pi r}{2} E z^{2} \right) }
{\sin^{a_{\ep}g-2}(\pi z)} \cot(\pi r z)\right\}.
$$
To verify the above claims recall that we in 
\refthm{thm:asymptotics-Enot0} have treated the case $n=0$
formally as the case $n=1$ with $(\alpha_1,\beta_1)=(1,0)$. 
Putting $(\rho_1,\sigma_1)=(0,1)$
we immediately get the claims above, see also the proof of
\refprop{prop:dependence}.  
The case $n=0$ and $E \neq 0$ only occurs if we work
with normalized Seifert invariants, namely for the spaces
$X=(\ep;g\;|\;b)$, $b \neq 0$. If we work with non-normalized
Seifert invariants there will always be a
$(\alpha,\beta)$--pair of the form $(1,\beta)$ coming from 
(modifications of) the pair $(1,b)$.
\end{rem}

There is a peculiar thing in the above theorem. Namely
if $X$ is described by non-normalized Seifert invariants we
can always add as many $(\alpha,\beta)$--pairs $(1,0)$ as we
want. In that way we can always obtain that the criterion
$n+a_{\ep}g-2>0$ is satisfied. We will see that this does not
change anything. In fact both terms $Z_{\pol}(X;r)$ and 
$Z_{\inte}(X;r)$ are independent of the choice of 
non-normalized Seifert invariants for $X$, see
\refprop{prop:dependence} for the precise statement.

Before considering the case $E=0$ let us take a more careful 
look at $Z_{\pol}(X;r)$.

\begin{prop}\label{prop:residue}
Assume that $n+a_{\ep}g-2>0$. For each 
$\gamma=(l,\un') \in \mI_{2}$ and $\nu \in \{0,1\}$ the 
function $Z_{\nu}^{\gamma}(r)$ is zero or is a Laurent 
polynomial in $r$ of the form
$\sum_{k=-\infty}^{n+a_{\ep}g-3} a_{k}^{\gamma} r^{k}$, where
the coefficients $a_{k}^{\gamma}=0$ for all but finitely 
many $k$.
\end{prop}

\begin{prf}
The proof is straightforward. Let $(l,\un') \in \mI_{2}$ be 
fixed and put 
$\lambda_{j} = 
\sin\left(\pi \frac{2\rho_{j}n_{j}'}{\alpha_{j}}\right)$ 
and
$\kappa_{j} =
\cos\left(\pi \frac{2\rho_{j}n_{j}'}{\alpha_{j}}\right)$. 
We start by analyzing $Z_{0}^{(l,\un')}$. Put
\begin{eqnarray*}
\tau(z) &=& \prod_{j=1}^{n} \left[ i 
\sin \left( 2 \pi \frac{\rho_{j}}{\alpha_{j}} n_{j}'\right) 
\cos \left( \pi \frac{z}{\alpha_{j}} \right) 
\sin \left( 2\pi r z \frac{n_{j}'}{\alpha_{j}} \right) 
\right. \\
&& \left. + 
\cos \left( 2 \pi \frac{\rho_{j}}{\alpha_{j}} n_{j}' \right) 
\sin \left( \pi \frac{z}{\alpha_{j}} \right) 
\cos \left( 2\pi r z \frac{n_{j}'}{\alpha_{j}} \right) \right]
\end{eqnarray*}
and write
$\tau(z) = \sum_{\nu =0}^{n} i^{\nu} 
\sum_{1\leq j_{1} < \ldots < j_{\nu} \leq n} 
g_{(j_{1},\ldots,j_{\nu})}(z)$,
where
\begin{eqnarray*}
g_{(j_{1},\ldots,j_{\nu})}(z) &=&
\prod_{j \in \{ j_{1},\ldots,j_{\nu}\}}
\lambda_{j}\cos\left( \pi \frac{z}{\alpha_{j}}\right)
\sin\left( 2\pi rz \frac{n_{j}'}{\alpha_{j}}\right) \\
&& \times
\prod_{j \in \{1,\ldots,n\} \sm \{j_{1},\ldots,j_{\nu}\}}
\kappa_{j}\sin\left( \pi \frac{z}{\alpha_{j}}\right)
\cos\left( 2\pi rz \frac{n_{j}'}{\alpha_{j}}\right).
\end{eqnarray*}
Note that $g_{(j_{1},\ldots,j_{\nu})}$ has a zero in $0$ of 
order $n$ or is constantly zero. We have
$$
g_{(j_{1},\ldots,j_{\nu})}(z) = r^{\nu} \sum_{k=n}^{\infty}
e_{k}^{(j_{1},\ldots,j_{\nu})}(r) z^{k},
$$
where $e_{k}^{(j_{1},\ldots,j_{\nu})}(r)$ is a real polynomial
in $r$ of order at most $k-n$. This description gives the 
result for $Z_{0}^{(l,\un')}(r)$. For the 
analysis of $Z_{1}^{(l,\un')}(r)$ we fix 
$\umu' \in \{\pm 1\}^{n}$ and put
$$
\sigma(z) = \prod_{j=1}^{n} \mu_{j}'
\sin \left[ 2 \pi  \left( \frac{\rho_{j}}{\alpha_{j}} n_{j}' 
- \frac{z}{2} \frac{\mu_{j}'}{\alpha_{j}} \right) \right] 
= \prod_{j=1}^{n} \left[ \mu_{j}' \lambda_{j}
\cos\left(\pi \frac{z}{\alpha_{j}}\right) -
\kappa_{j}\sin\left(\pi \frac{z}{\alpha_{j}}\right)\right]
$$
and write
$\sigma(z) = \sum_{\nu =0}^{n} 
\sum_{1\leq j_{1} < \ldots < j_{\nu} \leq n} 
h_{(j_{1},\ldots,j_{\nu})}(z)$,
where
$$
h_{(j_{1},\ldots,j_{\nu})}(z) =
(-1)^{n-\nu}\prod_{j \in \{ j_{1},\ldots,j_{\nu}\}}
\mu_{j}'\lambda_{j}\cos\left( \pi \frac{z}{\alpha_{j}}\right)
\prod_{j \in \{1,\ldots,n\} \sm \{j_{1},\ldots,j_{\nu}\}}
\kappa_{j}\sin\left( \pi \frac{z}{\alpha_{j}}\right).
$$
We note that $h_{(j_{1},\ldots,j_{\nu})}(z)$ is independent of
$r$ and has a zero in $0$ of order $n-\nu$ or is constantly 
zero.
\end{prf}

Let $\lambda_{j}$ and $\kappa_{j}$ be as in the above proof. We
note that $\lambda_{j}=0$ if and only if 
$n_{j}' \in \{0,\alpha_{j}/2\}$ and $\kappa_{j} =0$ if and only
if $n_{j}' = \alpha_{j}/4$ (so this happens only in case 
$\alpha_{j}$ is even and $l$ is odd). If $\lambda_{j}=0$ for a
$j \in \{1,2,\ldots,n\}$ then the highest degree term in both 
$Z_{0}^{(l,\un')}(r)$ and $Z_{1}^{(l,\un')}(r)$ is of an order
strictly less than $n+a_{\ep}g-3$. If $\alpha_{j}=1$ for a 
$j \in \{1,2,\ldots,n\}$ we have $\lambda_{j}=0$ for all
$(l,\un') \in \mI_{2}$. If 
$\lambda_{j}=0$ for $m$ indices $j \in \{1,2,\ldots,n\}$ then 
the highest degree term in both $Z_{0}^{(l,\un')}(r)$ and
$Z_{1}^{(l,\un')}(r)$ is of an order less than or equal to
$n+a_{\ep}g-3-m$.

In the following theorem the case $n=0$ is treated as the
case $n=1$ with $(\alpha_1,\beta_1)=(1,0)$. For that case,
see also the remarks following the theorem.

\begin{thm}\label{thm:asymptotics-E0}
Assume the Seifert Euler number $E=0$. If $n+a_{\ep}g-2 \leq 0$
then $X = S^{2} \times S^{1}$ {\em (so $\tau_{r}(X)=1$)}.
For $n+a_{\ep}g-2 >0$ we have
$$
Z(X;r) = \sum_{(l,\un') \in \mI_{2}} b_{(l,\un')}
\exp(2\pi i r q_{(l,\un')}) r \left(Z^{(l,\un')}_{0}(r) +
\tilde{Z}^{(l,\un')}_{1}(r)\right),
$$
where $b_{(l,\un')}$ is as in {\em (\ref{eq:Zpolar})}, 
$q_{(l,\un')}$ is given by {\em (\ref{eq:CSb})}, 
$Z^{(l,\un')}_{0}$ is given by {\em (\ref{eq:Zpolarreg})} with
$E=0$, and
\begin{eqnarray*}
\tilde{Z}^{(l,\un')}_{1}(r)
&=& -\pi i \left(\frac{i}{4}\right)^{n} \musum \mumsum 
\sum_{\stackrel{m \in \Z}{0 \leq m \leq |a|} } \muprod \\
&& \times 
\frac{\sign (a)}{\Sym_{\pm} (m) \Sym_{\pm}( m - |a|)}
\exp\left(2\pi i \sum_{j=1}^{n} \mu_{j} 
\frac{\rho_{j}}{\alpha_{j}} \mu_{j}'n_{j}' \right) \\ 
&& \times \Res_{z=0} \left\{ 
\frac{ \exp \left[ 2 \pi i r \sign (a)) ( m-|a|) z \right] }
{\sin^{n+a_{\ep}g-2}(\pi z)} \right\},
\end{eqnarray*}
where $a=a(\umu,\umu',\un')=\mualphammsuma$.
Moreover, if 
$\sum_{j=1}^{n} \frac{1}{\alpha_{j}} < n+a_{\ep}g-2$ we have
that $Z(X;r) = Z_{\pol}(X;r)$ for all levels $r$, where
$Z_{\pol}(X;r)$ is given by {\em (\ref{eq:Zpolar})} with 
$E=0$. Finally, in all cases there exists a positive number 
$r(\ualpha,\ubeta)$ such that
$$
Z(X;r) = Z_{\pol}(X;r) + rZ_{\spec}(X;r)
$$
for all $r \geq r(\ualpha,\ubeta)$, where $Z_{\pol}(X;r)$ is 
given by {\em (\ref{eq:Zpolar})} with $E=0$ and
\begin{eqnarray*}
Z_{\spec}(X;r) &=& \pi i \left(\frac{i}{4}\right)^{n}
\sum_{(l,\un') \in \mI_{2}^{a}}
b_{(l,\un')}\exp(2\pi i r q_{(l,\un')})
\sum_{\stackrel{\umu' \in \B}{\mualphammsum \in \Z}} 
\sum_{\stackrel{\umu \in \B}{\mualphasum < 0}} \muprod \\
&& \times \exp\left( 2\pi i\sum_{j=1}^{n} 
\mu_{j}\frac{\rho_{j}}{\alpha_{j}}\mu_{j}'n_{j}'\right)
\Res_{z=0} \left\{ \frac{f(z;\umu)}
{\sin^{n+a_{\ep}g-2}(\pi z)} \right\},
\end{eqnarray*}
where
$$
f(z;\umu) = \left\{ \begin{array}{cc}
\cos\left( \pi \pmualphasum z \right) 
&,n \hspace{.2in}\text{odd},\\
-i\sin\left( \pi \pmualphasum z \right) 
&,n \hspace{.2in}\text{even}.
\end{array}\right.
$$
In particular, $Z_{\spec}(X;r)=0$ for all 
$r \geq r(\ualpha,\ubeta)$ if
$\sum_{j=1}^{n} \frac{1}{\alpha_{j}} < n+a_{\ep}g-2$. 
\end{thm}

Let us look at the special case $n=0$, i.e.\ $X$ is given
by the Seifert invariants $(\ep;g\;|\;0)$ or non-normalized
Seifert invariants $\{\ep;g\;|\}$. As stated above this case 
is treated as the case $n=1$ with $(\alpha_1,\beta_1)=(1,0)$. 
If $a_{\ep}g \geq 4$ we get from the above theorem that
\begin{equation}\label{eq:Zpoln0}
Z(X;r)= Z_{\pol}(X;r) = -\pi r 
\Res_{z=0} \left\{ 
\frac{ \cot(\pi r z) }{\sin^{a_{\ep}g-2}(\pi z)} \right\}.
\end{equation}
If $a_{\ep}g=2$ we have $X=(\os;1\;|\;0)$ or $X=(\ns;2\;|\;0)$.
For this case we use that $r(\ualpha,\ubeta)$ in 
\refthm{thm:asymptotics-E0} can be put equal to $2$ since 
$n=1$. Thus we get for all levels $r$ that
$$
Z(X;r)=Z_{\pol}(X;r) + rZ_{\spec}(X;r).
$$
Here $Z_{\pol}(X;r)$ is given by (\ref{eq:Zpoln0}), i.e.\
$Z_{\pol}(X;r)=-\pi r \Res_{z=0} (\cot(\pi r z))=-1$. Moreover
$Z_{\spec}(X;r)=1$ so $Z_{\pol}(X;r)=r-1$. Note that this
infact follows directly from (\ref{eq:Zfunction-Seifert}).
Finally, if $a_{\ep}g=0$, i.e.\ $\ep=\os$ and $g=0$, then
$X=S^{2} \times S^1$.

By the above theorems we see that
\begin{equation}\label{eq:asymptotic-Z}
Z(X;r) = \left[ \sum_{j \in \mI} e^{2\pi ir q_{j}} 
\sum_{l=l_{0}}^{N} c_{l}^{j} r^{-l/2} \right] + o(r^{-N/2})
\end{equation}
for every $N \geq l_{0}$ in the limit $r \ra \infty$, where 
$l_{0}$ is some fixed integer (depending on $X$), and
$\mI=\mI_{1} \amalg \mI_{2}$ if $E \neq 0$ and $\mI=\mI_{2}$ 
if $E=0$. By \refprop{prop:residue} the term $Z_{\pol}(X;r)$ 
has a structure like (\ref{eq:generic-type-alternative}). In 
\refthm{thm:asymptotics-Enot0} we have separated the 
`exact part' $Z_{\pol}(X;r)$ from a part $Z_{\inte}(X;r)$, 
which is only given an asymptotic description. However, one 
should note that in the expression (\ref{eq:asymptotic-Z}) 
there are contributions to the terms with $j \in \mI_{2}^{a}$ 
coming both from $Z_{\pol}(X;r)$ and $Z_{\inte}(X;r)$.

If we work with non-normalized Seifert invariants, 
the different quantities involved in the
above Theorems~\ref{thm:asymptotics-Enot0} and 
\ref{thm:asymptotics-E0} potentially depends on the actual
choice of invariants. The precise
situation is described in the following proposition, where we
remind the reader that for the large Seifert manifolds
the Seifert invariants in fact classify the Seifert manifolds
up to an orientation preserving homeomorphism. Thus all
quantities independent of the choice of non-normalized Seifert
invariants are topological invariants for these large Seifert
manifolds. Moreover, these topological invariants can then be
calculated using any set of non-normalized Seifert invariants
for a given large Seifert manifold.

\begin{prop}\label{prop:dependence}
We have
$$
Z(X;r) = Z_{\pol}(X;r) + Z_{\inte}(X;r)
$$
where $Z_{\pol}(X;r)$ and $Z_{\inte}(X;r)$ are as in
{\em \refthm{thm:asymptotics-Enot0}} if $E \neq 0$ and 
$Z_{\inte}(X;r)=0$ if $E=0$. In both cases
$$
Z_{\pol}(X;r) = \sum_{(l,\un') \in \mI_{2}} 
\exp(2\pi irq_{(l,\un')})Q^{(l,\un')}(r)
$$
where $q_{(l,\un')}$ is given by {\em (\ref{eq:CSb})} and the 
$Q^{(l,\un')}(r)$ are Laurent polynomials in $r$. The 
expressions $Z_{\pol}(X;r)$ and $Z_{\inte}(X;r)$ do not 
depend on the choice of non-normalized Seifert invariants
for $X$. More specifically, let 
$\Delta=\{\ep;g;(\alpha_{1},\beta_{1}),\ldots,
(\alpha_{n},\beta_{n})\}$
and
$\tilde{\Delta}=\{\ep;g;(\tilde{\alpha}_{1},\tilde{\beta}_{1}),
\ldots,
(\tilde{\alpha}_{\tilde{n}},\tilde{\beta}_{\tilde{n}})\}$
be two sets of non-normalized Seifert invariants for $X$,
let $\mI_{1}$ and $\mI_{2}$ be the index sets related to
$\Delta$ given by {\em (\ref{eq:I1})} and {\em (\ref{eq:I2})}
respectively and let $\tilde{\mI}_{1}$ and $\tilde{\mI}_{2}$ be
the corresponding index sets related to $\tilde{\Delta}$. Then 
we have a one-one correspondence between $\mI_{j}$ and
$\tilde{\mI}_{j}$ for $j=1,2$, and under these 
correspondences the functions $q$ in {\em (\ref{eq:CSa})} and 
{\em (\ref{eq:CSb})} are preserved. Moreover, the large $r$ 
asymptotic expansion of $Z_{\inte}(X;r)$, i.e.\ the 
coefficients in {\em (\ref{eq:A3})} and {\em (\ref{eq:A3'})} 
are preserved. Finally the Laurent polynomials 
$Q^{(l,\un')}(r)$ are preserved if $E \neq 0$ or $E=0$ and 
$\sum_{j=1}^{n} \frac{1}{\alpha_{j}} <n+a_{\ep}g-2$.
\end{prop}

\begin{prf}
The proof is simply a matter of routine checks. We first 
consider the case 
$E = -\sum_{j=1}^{n} \frac{\beta_j}{\alpha_j} \neq 0$. Let us
write a tilde on each quantity or set relating to 
$\tilde{\Delta}$. Thus we write $\tilde{\mS}$ instead of $\mS$
etc. Let us first assume that $\tilde{\Delta}$ is obtained from
$\Delta$ by adding a trivial $(\alpha,\beta)$--pair 
$(1,0)$. Thus $\tilde{n}=n+1$ and 
$(\tilde{\alpha}_{j},\tilde{\beta}_{j})=(\alpha_{j},\beta_{j})$
for $j=1,2,\ldots,n$ and
$(\tilde{\alpha}_{n+1},\tilde{\beta}_{n+1})=(1,0)$.
We then have 
$\tilde{\mS}=\mS \times \{0\}$ and
$\tilde{z}_{\st}(m,(\un,0))=z_{\st}(m,\un)$ for all 
$(m,\un) \in \Z \times \mS$. Thus we also have
$\tilde{\mI}_{1} = \mI_{1} \times \{0\}$. We let
$(\tilde{\rho}_{j},\tilde{\sigma}_{j})=(\rho_{j},\sigma_{j})$
for $j=1,2,\ldots,n$ and put 
$(\tilde{\rho}_{n+1},\tilde{\sigma}_{n+1})=(0,1)$ and find that
$\tilde{q}_{(m,(\un,0))}=q_{(m,\un)}$ for all 
$(m,\un) \in \mI_{1}$.
For $(m,\un) \in \mI_{1}$ and $k =0,1,2,\ldots$ it follows
immediately from (\ref{eq:A3}) that 
$\tilde{c}_{k}^{(m,(\un,0))}=c_{k}^{(m,\un)}$. Next, let us 
look at the set $\tilde{\mI}_{2}$ and the related quantities. 
We have $\tilde{M}= M \times \Z$ and then 
$\tilde{\mI}_{2}=\mI_{2} \times \{0\}$. Thus we also have that
$\tilde{\mI}_{2}^{a}=\mI_{2}^{a} \times \{0\}$. For 
$(l,\un') \in \mI_{2}$ we immediately get that
$\tilde{q}_{(l,(\un',0))}=q_{(l,\un')}$, 
$\tilde{b}_{(l,(\un',0))}=b_{(l,\un')}$ and
$\tilde{Z}_{\nu}^{(l,(\un',0))}(r) = Z_{\nu}^{(l,\un')}(r)$, 
$\nu=0,1$. Finally, let us check that the coefficients in 
(\ref{eq:A3'}) stay unchanged. First assume that $n$ is even. 
Then $k_0 = 0$, $n+1$ is odd and 
$\tilde{k}_{0}=1$ so 
$\tilde{k}_{1} =(\tilde{k}_{0}-(n+1)-a_{\vep}g+2)/2 = k_{1}$. 
If $k \in \{k_1,k_{1}+1,\ldots\}$, then 
$\tilde{k}'=2k+(n+1)+a_{\vep}g-2=k'+1$.
We note that
$\tilde{f}(z;l,(\un',0),(\umu',\mu_{n+1}')) = 
-\pi z f(z;l,\un',\umu')$, 
and that
$$
\left. \partial_{z}^{(k'+1)} zf(z;l,\un',\umu') \right|_{z=0}
= (k'+1)\left. \partial_{z}^{(k')} 
f(z;l,\un',\umu') \right|_{z=0}.
$$
This shows that $\tilde{c}_{k}^{(l,(\un',0)}=c_{k}^{(l,\un')}$.
Next assume that $n$ is odd. Then $k_0 =1$, $n+1$ is even and
$\tilde{k}_{0}=0$, so $\tilde{k}_{1}=k_{1}-1$. Like before we 
find that $\tilde{c}_{k}^{(l,(\un',0)}=c_{k}^{(l,\un')}$ for
$k \in \{k_{1},k_{1}+1,\ldots\}$. Finally
$\tilde{c}_{\tilde{k}_{1}}^{(l,(\un',0)}=0$ since
$\tilde{f}(0;l,(\un',0),(\umu',\mu_{n+1}'))=0$.

Next, consider the case where $\tilde{n}=n$,
$\tilde{\alpha}_{j}=\alpha_{j}$ and 
$\tilde{\beta}_{j}=\beta_{j}+K_{j}\alpha_{j}$ for some integers
$K_{j}$ with $\sum_{j=1}^{n} K_{j} =0$.
Let $\tilde{\rho}_{j}=\rho_{j}$ and
$\tilde{\sigma}_{j}=\sigma_{j}+K_{j}\rho_{j}$, 
$j=1,2,\ldots,n$. 
We immediately get that $\tilde{\mS}=\mS$ and that 
$\tilde{z}_{\st}(m,\un)=z_{\st}(m,\un)$ and 
$\tilde{q}_{(m,\un)}=q_{(m,\un)}$ for 
$(m,\un) \in \Z \times \mS$.
Moreover, $\tilde{\mI}_{1}=\mI_{1}$ and the coefficients in
(\ref{eq:A3}) are unchanged. Next let us consider 
$\tilde{\mI}_{2}$ and the quantities related to that index set.
First note that
$(0,\un') \in \tilde{\mI}_{2}$ if and only if 
$(0,\un') \in \mI_{2}$ and that all quantities related to such 
a point $(0,\un')$ are the same w.r.t.\ the two sets of 
non-normalized Seifert invariants $\Delta$ and 
$\tilde{\Delta}$. Let us next consider a point 
$(1,\un') \in \tilde{\mI}_{2}$. Thus 
$n_{j}' \in \Z + \frac{1}{2}\beta_{j}$ if $\alpha_{j}$ or
$K_{j}$ is even and $n_{j}' \in \Z + \frac{1}{2}\beta_{j} + 
\frac{1}{2}$ if both $\alpha_{j}$ and $K_{j}$ are odd.
But if $\alpha_{j}$ is odd then 
$\alpha_{j}/2 \in \frac{1}{2} + \Z$ so in all cases there is
a one-one correspondence between the set of points 
$(l,\tilde{\un}') \in \tilde{\mI}_{2}$ with $l=1$ and the set 
of points $(l,\un') \in \mI_{2}$ with $l=1$. 
A correspondence is given as follows: For 
$(1,\un') \in \mI_{2}$, let 
$n_{j}'' = n_{j}'+\frac{1}{2}K_{j}\alpha_{j} \pmod{\alpha_j}$
such that $0 \leq n_{j}'' < \alpha_{j}$
and let $\tilde{n}_{j}' = n_{j}''$ if 
$0 \leq n_{j}'' \leq \alpha_{j}/2$ and 
$\tilde{n}_{j}' = -n_{j}''+\alpha_{j}$ otherwise.
Then $(1,\tilde{\un}') \in \tilde{\mI}_{2}$. To see that
$\tilde{q}_{(1,\tilde{\un}')} = q_{(1,\un')}$ we observe that
the $q$-function in (\ref{eq:CSb}) only depends on 
$n_{j}' \pmod{\alpha_{j}}$ and is unchanged when changing the
sign of $n_{j}'$. Moreover,
$$
\frac{\rho_{j}}{\alpha_{j}}
\left(n_{j}' + \frac{1}{2}K_{j}\alpha_{j}\right)^{2}
- \left(\sigma_{j} + \rho_{j}K_{j}\right)
\left(\beta_{j} + \alpha_{j}K_{j}\right)
= \frac{1}{4}K_{j} \pmod{\Z}.
$$
We have 
$\tilde{b}_{(1,\tilde{\un}')} = 
(-1)^{\sum_{j=1}^{n}K_{j}\rho_{j}} b_{(1,\un')}$
since $\tilde{n}_{j}'/\alpha_{j} \in \frac{1}{2}\Z$ if and only
if $n_{j}'/\alpha_{j} \in \frac{1}{2}\Z$.
Assume that 
$\sum_{j=1}^{n} \frac{\mu_{j}'\tilde{n}_{j}'}{\alpha_{j}}
\in \Z$ and let $\tilde{\mu}_{j}'=\mu_{j}'$ if 
$0 \leq n_{j}'' \leq \alpha_{j}/2$ and 
$\tilde{\mu}_{j}'=-\mu_{j}'$ otherwise, where $n_{j}''$ is as
above. Then
$$
\sum_{j=1}^{n} 
\frac{\tilde{\mu}_{j}'\tilde{n}_{j}'}{\alpha_{j}}
= \sum_{j=1}^{n} \frac{\mu_{j}' n_{j}'}{\alpha_{j}}
+ \frac{1}{2}\sum_{j}^{n}\mu_{j}'K_{j} \pmod{\Z}.
$$
But $\frac{1}{2}K_{j} - \frac{1}{2}\mu_{j}'K_{j} \in \Z$ so
$\frac{1}{2}\sum_{j}^{n}\mu_{j}'K_{j} \in \Z$.
We also have
$$
\tilde{f}(z;1,\tilde{\un}',\tilde{\umu}')
(-1)^{\sum_{j=1}^{n}K_{j}\rho_{j}}f(z;l,\un',\umu').
$$
These observations show that 
$\tilde{c}_{k}^{(1,\tilde{\un}')}=c_{k}^{(1,\un')}$ in
(\ref{eq:A3'}). Finally, let us notice that
$$
\tilde{b}_{(1,\tilde{\un}')}
\left( \tilde{Z}_{0}^{(1,\tilde{\un}')}(r) 
+ \tilde{Z}_{1}^{(1,\tilde{\un}')}(r) \right)
= b_{(1,\un')}\left( Z_{0}^{(1,\un')}(r) 
+ Z_{1}^{(l,\un')}(r) \right).
$$ 
The easiest way to see this is to refer to the actual 
calculation of $Z_{\pol}(X;r)$. Thus, for $l \in \Z$ let 
$J_{l}'$ be the set in (\ref{eq:Jl'}) and let $\tilde{J}_{l}'$ 
be the corresponding set w.r.t.\ $\tilde{\Delta}$. 
Then $\tilde{J}_{l}'=J_{l}'$ for each even $l$ and for $l$
odd the map $\un' \mapsto \tilde{\un}'$ described above gives
a bijection from $J_{l}'$ to $\tilde{J}_{l}'$. Moreover, if 
$Z^{l}(\delta,\eta)$ is given by (\ref{eq:Zl}) and
$\tilde{Z}^{l}(\delta,\eta)$ is the corresponding function 
w.r.t.\ $\tilde{\Delta}$ then
$\tilde{Z}^{l}(\delta,\eta)=Z^{l}(\delta,\eta)$ proving the 
claim. Note that we do not necessarily have
$\tilde{b}_{(1,\tilde{\un}')}
\tilde{Z}_{\nu}^{(1,\tilde{\un}')}(r) =
b_{(1,\un)}Z_{\nu}^{(1,\un')}(r)$
for $\nu=0,1$.

Let us next assume that $E=0$. Exactly as above we get that
there is a one-one correspondence between $\mI_{2}$ and
$\tilde{\mI}_{2}$ and under this correspondence the function
$q$ in (\ref{eq:CSb}) is preserved. Moreover, we know that
$Z(X;r)$ is independent of the choice of non-normalized Seifert
invariants. Assume 
$\sum_{j=1}^{n} \frac{1}{\alpha_{j}} < n +a_{\ep}g-2$. Note
that this condition is independent of the choice of 
non-normalized Seifert invariants. In this case we have
$Z(X;r)=Z_{\pol}(X;r)$, where $Z_{\pol}(X;r)$ is given by
(\ref{eq:Zpolar}) with $E=0$. Thus in this case we get that
$Q^{(l,\un')}(r) = b_{(1,\un')}\left( Z_{0}^{(1,\un')}(r) 
+ Z_{1}^{(l,\un')}(r) \right)$
is independent of the choice of non-normalized Seifert 
invariants exactly as in the case $E=0$.
\end{prf}

In case $E=0$ it follows from \refthm{thm:asymptotics-E0} that
the large $r$ asymptotics of $Z(X;r)$ is given by 
$Z_{\pol}(X;r)$ in (\ref{eq:Zpolar}) (with $E=0$) except for 
the special cases 
$\sum_{j=1}^{n} \frac{1}{\alpha_{j}} \geq n+a_{\ep}g-2$, where
we also have the extra term $Z_{\spec}(X;r)$.

Let us determine which of the Seifert manifolds with $E=0$
that satisfy the condition 
$\sum_{j=1}^{n} \frac{1}{\alpha_{j}} \geq n+a_{\ep}g-2$. 
As noted above the condition is independent of the choice
of non-normalized Seifert invariants if we work with such
invariants. We can therefore safely restrict to consider 
(normalized) Seifert invariants. This also allows us to keep 
track of the types of Seifert fibrations. Therefore, let
$X=(\ep;g\;|\;b;(\alpha_{1},\beta_{1}),\ldots,
(\alpha_{n},\beta_{n}))$. Recall that 
$0 < \beta_{j} < \alpha_{j}$ so $\alpha_{j} \geq 2$ for all
$j=1,\ldots,n$. Thus
$\sum_{j=1}^{n} \frac{1}{\alpha_{j}} \leq \frac{n}{2}$ so 
$\sum_{j=1}^{n} \frac{1}{\alpha_{j}} \geq n + a_{\ep}g-2$
implies that $2-a_{\ep}g \geq \frac{n}{2}$, and this is only 
satisfied in the cases $g=0$, $\ep=\os$, 
$n \in \{0,1,2,3,4\}$ and $g=1$, $\ep=\os$, $n=0$ and $g=2$, 
$\ep=\ns$, $n=0$ (recall that we only consider even genus for 
$\ep=\ns$). In the last two cases the Seifert Euler number 
$E=-b$, so since this number is zero we have 
$X=(\os;1|\;0)= T^{2} \times S^{1}$ (the $3$--torus) or 
$X=(\ns;2|\;0)$, both being small Seifert manifolds. 
In case $\ep=\os$, $g=0$, and 
$n \in \{0,1,2\}$, $X$ is a lens space $L(p,q)$, see 
\cite[Sect.~5.4 (i) pp.~99--100]{Orlik}. We find in all these 
cases that the Seifert Euler number $E=0$ if and only
if $p=0$, hence $X=S^{2} \times S^{1}$ if $E=0$.
The remaining cases to consider are
$\ep=\os$, $g=0$ and $n \in \{3,4\}$. Let us first consider the
easy case $n=4$. Here the $\alpha_{j}$'s have to satisfy 
$\sum_{j=1}^{4} \frac{1}{\alpha_{j}} \geq 2$, hence we have
$\alpha_{j}=2$ for all $j\in\{1,2,3,4\}$. In this case
$E=-b-4\frac{1}{2}=-b-2$, so $b=-2$ when $E=0$. Thus
$X=(\os;0\;|\;-2;(2,1),(2,1),(2,1),(2,1))$. This is the 
small Seifert manifold considered in 
\cite[Sec.~5.4 (iii) p.~101]{Orlik}. In fact, this manifold is
homeomorphic by an orientation preserving homeomorphism to the
Seifert manifold $(\ns;2|\;0)$, cf.\ 
\cite[Theorem 5.1]{JankinsNeumann}. Let us finally consider the
case $\ep=\os$, $g=0$, $n=3$. Here we search for solutions to
$\sum_{j=1}^{3} \frac{1}{\alpha_{j}} \geq 1$. The subset of 
these satisfying the strict inequality 
$\sum_{j=1}^{3} \frac{1}{\alpha_{j}} > 1$ are the small Seifert
manifolds considered in 
\cite[Sec.~5.4 (ii) pp.~100--101]{Orlik}. By a direct 
inspection one finds that none of these manifolds have Seifert 
Euler number $E=0$. Finally we consider the manifolds with 
$\sum_{j=1}^{3} \frac{1}{\alpha_{j}} = 1$. The only solutions 
to this equation with 
$\alpha_{1} \leq \alpha_{2} \leq \alpha_{3}$ are $(2,3,6)$,
$(2,4,4)$ and $(3,3,3)$. A direct inspection reveals that we 
have the following possibilities with $E=0$:
\begin{eqnarray*}
&&X_{1} = (\os;0\;|\;-1;(2,1),(3,1),(6,1)), \\
&&X_{2} = (\os;0\;|\;-2;(2,1),(3,2),(6,5))=-X_{1}, \\
&&X_{3} = (\os;0\;|\;-1;(2,1),(4,1),(4,1)), \\
&&X_{4} = (\os;0\;|\;-2;(2,1),(4,3),(4,3))=-X_{3}, \\
&&X_{5} = (\os;0\;|\;-1;(3,1),(3,1),(3,1)), \\
&&X_{6} = (\os;0\;|\;-2;(3,2),(3,2),(3,2))=-X_{5},
\end{eqnarray*}
where, as usual, $-X$ means $X$ with the opposite orientation.
Thus the Seifert fibrations with 
$\sum_{j=1}^{n} \frac{1}{\alpha_{j}} \geq n+a_{\ep}g-2$, 
$a_{\ep}g$ even and 
$E=0$ are the small Seifert manifolds $S^{2} \times S^{1}$,
$(\os;1\;|\;0)=T^{2} \times S^{1}$, 
$(\ns;2\;|\;0) \cong (\os;0\;|\;-2;(2,1),(2,1),(2,1),(2,1))$
and the six $3$--fibered large Seifert manifolds 
$X_{j}$, $j=1,\ldots,6$. 

We saw above that all but 9 of the Seifert fibrations 
satisfying
$\sum_{j=1}^{n} \frac{1}{\alpha_{j}} \geq n +a_{\ep}g-2$, 
$a_{\ep}g$ even and 
$E=0$ are topologically the space $S^{2} \times S^{1}$. There 
is an infinite number of such fibrations for 
$S^{2} \times S^{1}$. If 
$X=(\os;0\;|\;b,(\alpha_{1},\beta_{1}),
(\alpha_{2},\beta_{2}))$,
then $E=-b-\beta_{1}/\alpha_{1}-\beta_{2}/\alpha_{2}=0$ if and
only if $b=-1$, $\alpha_{1}=\alpha_{2}$ and 
$\beta_{2}=\alpha_{1}-\beta_{1}$. Thus $(\os;0\;|\;0)$ and
$(\os;0\;|\;-1;(\alpha,\beta),(\alpha,\alpha-\beta))$, 
$\alpha \geq 2$, give the possible Seifert fibrations for 
$S^{2} \times S^{1}$. (Note that there are no 
$(\os;0\;|\;b;(\alpha,\beta))$ with $E=-b-\beta/\alpha=0$
and $0 < \beta < \alpha$.)
For all the Seifert fibrations satisfying
$\sum_{j=1}^{n} \frac{1}{\alpha_{j}} \geq n +a_{\ep}g-2$, 
$a_{\ep}g$ even and $E=0$ the term
$Z_{\spec}(X;r) \neq 0$.

\subsection{The moduli space of flat $\SU(2)$--connections on 
Seifert manifolds}\label{sec-The-moduli}

It follows from the results in the previous section that the 
proof of \refthm{thm:AEC-Seifert-manifolds} it finalized once 
we have shown that the values of the $q$--functions defined in
(\ref{eq:CSa}) and (\ref{eq:CSb}) are in fact the Chern--Simons
values of the flat $\SU(2)$--connections on $X$. This is done 
by a simple comparison between the values of the $q$--functions
and the values of the Chern--Simons functional of flat
$\SU(2)$--connections on the Seifert manifolds calculated by 
D.\ Auckly, cf.\ \cite{Auckly1}.

Let $X$ be a Seifert manifold given by
the Seifert invariants 
$(\ep;g\;|\;b;(\alpha_{1},\beta_{1}),\ldots,
(\alpha_{n},\beta_{n}))$ 
or the non-normalized Seifert invariants
$\{\ep;g;(\alpha_{1},\beta_{1}),\ldots,
(\alpha_{n},\beta_{n})\}$. Let $\beta_{0}$ be as in the 
previous section and let ${\mathcal M}_{X}$ denote the moduli 
space of flat $\SU(2)$--connections on $X$. There is a 
classical identification
$$
\mM_{X} = \Hom(\pi_{1}(X),\SU(2))/\SU(2),
$$
where $/\SU(2)$ means moduli conjugation. Let us recall 
Auckly's results. Assume first that $\ep=\os$ (so $g \geq 0$).
It is well-known that
\begin{eqnarray*}
\pi_{1}(X) &=&  \langle h,q_{1},\ldots,q_{n},a_{1},b_{1},
\ldots,a_{g},b_{g}\; | \; q_{k}^{\alpha_{k}}h^{\beta_{k}}=1,\;
h^{-\beta_0}q_{1}\ldots q_{n}\prod_{i=1}^{g} [a_{i},b_{i}]=1,\\
&& [a_{j},h]=[b_{j},h]=[q_{k},h]=1,\; 
j=1,2,\ldots,g,\; k=1,2,\ldots,n\;\rangle,
\end{eqnarray*}
cf.\ \cite{Auckly1}, \cite[Theorem 6.1]{JankinsNeumann}. Let 
$g_{1},\ldots,g_{n} \in \SU(2)$, $\un \in Z^n$ and 
$\vep \in \{0,1/2\}$. Then Auckly defines representations
$$
\omega = \omega(\vep,\un)[g_{1},\ldots,g_{n}]:
\pi_{1}(X) \to \SU(2)
$$
in the following way, cf.\ \cite[Definition p.~56]{Auckly2},
\cite[Definition p.~231]{Auckly1}. Let $\tilde{\omega}$ be the 
homomorphism from the free group generated by
$h,q_{1}\ldots,q_{n},a_{1},b_{1},\ldots,a_{g-1},b_{g-1}$ to 
$\SU(2)$ given by $\tilde{\omega}(h)=e^{2\pi i\vep}$, 
$\tilde{\omega}(a_{j})=\tilde{\omega}(b_{j})=1$, 
$j=1,\ldots,g-1$, and
$\tilde{\omega}(q_{k})=g_{k}\exp\left( 2\pi i 
(n_{k}-\vep\beta_{k})/\alpha_{k} \right)g_{k}^{-1}$,
$k=1,\ldots,n$, where we identify $\SU(2)$ by the unit 
quaternions $\Sp_{1}$ in the usual way. Let
$$
a+jb = \tilde{\omega}(q_{n}^{-1}\ldots q_{1}^{-1}h^{\beta_0}), 
\hspace{.2in}a,b \in \C.
$$
If $a=1$ (so $b=0$) we let $x=y=1$. If not we let $x \in S^1$ 
be a square root of
$\frac{ a(2-2\rea(a)) - |b|^{2}}{2-2\rea(a)-|b|^{2}}$ and
$y = \frac{\bar{b}}{r(x^2 - 1)} + jr$,
where $r \in [-1,1]$ satisfies
$r^2 = 1 - \frac{|b|^{2}}{2-2\rea(a)}$
(noting that $1- \frac{|b|^{2}}{2-2\rea(a)} \geq 0$).
We then define $\omega$ by letting
$\omega(\gamma)=\tilde{\omega}(\gamma)$,
$\gamma \in \{ h,q_{1},\ldots,q_{n},a_{1},b_{1},\ldots,
a_{g-1},b_{g-1} \}$,
and by letting $\omega(a_{g})=x$ and
$\omega(b_{g})=y$. (Thus, in case $a \neq 1$, there are four 
possible choices of $\omega$ corresponding to the two possible 
choices of each of $x$ and $y$. All four choices are needed in
\refprop{prop:modulioriented}.) If $g \geq 1$ then $\omega$ 
extends to a $\SU(2)$ representation of $\pi_{1}(X)$. If $g=0$ 
then $\omega$ extends to a representation 
$\pi_{1}(X) \to \SU(2)$ if and only if
\begin{equation}\label{eq:genus0condition}
e^{-2\pi i\beta_{0}\vep}\prod_{j=1}^{n} g_{j}\exp\left( 2\pi i 
\left( \frac{n_{j}-\vep\beta_{j}}{\alpha_{j}} \right)\right) 
g_{j}^{-1} =1.
\end{equation}
Above we can of course restrict to $\un \in \mS$. If the 
Seifert Euler number 
$E = -\beta_0-\sum_{j=1}^{n} \frac{\beta_{j}}{\alpha_{j}}
\neq 0$ 
then we define a representation
$$
\rho=\rho(m,\un) : \pi_{1}(X) \to \SU(2)
$$
for each integer $m$ and each $\un \in \mS$ by letting
$\rho(a_{j})=\rho(b_{j})=1$, $j=1,\ldots,g$, and by putting
\begin{eqnarray*}
\rho(h) &=& \exp\left( -2\pi i \frac{1}{E} \left( m + 
\sum_{k=1}^{n} \frac{n_{k}}{\alpha_{k}} \right)\right), \\
\rho(q_{j}) &=& \exp\left( 2\pi i \left( 
\frac{n_{j}}{\alpha_{j}} + 
\frac{1}{E}\frac{\beta_{j}}{\alpha_{j}} \left( m + 
\sum_{k=1}^{n} \frac{n_{k}}{\alpha_{k}} \right)\right)\right).
\end{eqnarray*}
Note that the image of $\rho$ is contained in 
$\C \cap \Sp_{1}=\tU(1)=S^{1}$, so $\rho$ is reducible. We have

\begin{prop}[{\cite[Lemma p.~57]{Auckly2}}]
\label{prop:modulioriented}
Let the situation be as above. If $E=0$ then any element of 
$\mM_{X}$ is in the same path component as one of the 
conjugacy classes of the $\omega$'s. If $E \neq 0$ then any 
element of $\mM_{X}$ is in the same path component as one of 
the conjugacy classes of the $\omega$'s and $\rho$'s.\HS
\end{prop}

If $l \in \{0,1\}$, $\un' \in \Z^{n} + \frac{1}{2}l\ubeta$ and
$g_1,\ldots,g_n \in \SU(2)$ then we let $\vep=l/2$ and
$n_j=n_j'+\vep\beta_j$ and let
$\omega'(l,\un')[g_1,\ldots,g_n]$ be equal to 
$\omega(\vep,\un)[g_1,\ldots,g_n]$. We note that
$\omega'(l,-\un')[g_{1}j,\ldots,g_{n}j]
=\omega'(l,\un')[g_1,\ldots,g_n]$ thus we can restrict to
$(l,\un') \in \mI_2$.

Next we consider the Seifert manifolds with non-orientable 
base. Therefore, let $X$ be as above with $\ep=\ns$ (so $g>0$).
We have
\begin{eqnarray*}
\pi_{1}(X) &=&  \langle h,q_{1},\ldots,q_{n},a_{1},\ldots,
a_{g}\; | \; q_{k}^{\alpha_{k}}h^{\beta_{k}}=1,\;
h^{-\beta_{0}}q_{1}\ldots q_{n}\prod_{i=1}^{g} a_{i}^{2}=1, \\
&& a_{j}^{-1}ha_{j}h=[q_{k},h]=1,\; 
j=1,2,\ldots,g,\; k=1,2,\ldots,n\;\rangle,
\end{eqnarray*}
cf.\ \cite{Auckly1}, \cite[Theorem 6.1]{JankinsNeumann}.
Let $g_{1},\ldots,g_{n} \in \SU(2)$, $\un \in \Z^n$ and 
$\vep \in \{0,1/2\}$. Then Auckly defines representations
$$
\nu = \nu(\vep,\un)[g_{1},\ldots,g_{n}]: \pi_{1}(X) \to \SU(2)
$$
in the following way, cf.\ \cite[Definition p.~56]{Auckly2},
\cite[Definition p.~231]{Auckly1}. 
First let $\tilde{\nu}$ be the homomorphism from the free group
generated by
$h,q_{1}\ldots,q_{n},a_{1},\ldots,a_{g-1}$ to $\SU(2)$ given by
$\tilde{\nu}(h)=\exp(2\pi i\vep)$,
$\tilde{\nu}(q_{k})=g_{k}\exp\left( 2\pi i (n_{k} - 
\vep\beta_{k})/\alpha_{k} \right)g_{k}^{-1}$,
and $\tilde{\nu}(a_{j})=1$, $j=1,\ldots,g-1$, $k=1,\ldots,n$. 
We then define $\nu$ by letting 
$\nu(\gamma)=\tilde{\nu}(\gamma)$,
$\gamma \in \{h,q_{1},\ldots,q_{n},a_{1},\ldots,a_{g-1}\}$, and
by putting 
$\nu(a_{g})=\sqr\left( \tilde{\nu}(q_{n}^{-1}\ldots 
q_{1}^{-1}h^{\beta_0})\right)$,
where $\sqr (x)$ is any element of $\SU(2)$ such that
$(\sqr(x))^{2}=x$, $x \in \SU(2)$. (We know that such an 
element exists because there is an element $z \in \SU(2)$ such
that $zxz^{-1} \in S^{1}$. Therefore we have an element 
$y \in S^{1}$ such that $y^{2} =zxz^{-1}$. But then 
$(z^{-1} y z)^{2}=x$.) This defines a representation 
$\nu : \pi_{1}(X) \to \SU(2)$. (We will see below, that the 
Chern--Simons value of $[\nu] \in \mM_{X}$ is independent of 
the choice of this square root, so therefore we don't specify 
this choice here.) As in the oriented base case we can 
parametrize slightly differently by letting 
$\nu'(l,\un')[g_1,\ldots,g_n]$ be equal to 
$\nu(\vep,\un)[g_1,\ldots,g_n]$ for each $l \in \{0,1\}$ and 
$\un' \in \Z^n + \frac{1}{2}l\ubeta$ with $\vep=l/2$ and 
$\un = \un' + \vep\ubeta$. Also here we can restrict to 
$(l,\un') \in \mI_2$.

If the Euler number
$E=-\beta_{0}-\sum_{j=1}^{n} \frac{\beta_{j}}{\alpha_{j}} 
\neq 0$ 
then we define a representation
$$
\sigma=\sigma(m,\un) : \pi_{1}(X) \to \SU(2)
$$
for each integer $m$ and each $\un \in \mS$ by letting
$\sigma(a_{i})=j$, $i=1,\ldots,g$, and by putting
\begin{eqnarray*}
\sigma(h) &=& \exp\left( -2\pi i \frac{1}{E} \left( 
\frac{g}{2} + m + \sum_{k=1}^{n} \frac{n_{k}}{\alpha_{k}} 
\right)\right), \\
\sigma(q_{j}) &=& \exp\left( 2\pi i \left( 
\frac{n_{j}}{\alpha_{j}} + \frac{1}{E}
\frac{\beta_{j}}{\alpha_{j}} \left( 
\frac{g}{2}+ m + \sum_{k=1}^{n} \frac{n_{k}}{\alpha_{k}} 
\right)\right)\right).
\end{eqnarray*}
We then have

\begin{prop}[{\cite[Lemma p.~62]{Auckly2}}]
Let the situation be as above. If $E=0$ then any element of 
$\mM_{X}$ is in the same path component as one of the 
conjugacy classes of the $\nu$'s. If $E \neq 0$ then any 
element of $\mM_{X}$ is in the same path component as one of 
the conjugacy classes of the $\nu$'s and $\sigma$'s.\HS
\end{prop}

\noindent Using \cite[Theorem 4.2]{KirkKlassen1} Auckly proved 
the following result, where $\CS$ denotes the Chern--Simons
functional defined in (\ref{eq:CSfunctional}) (with 
$G=\SU(2)$).

\begin{thm}[{\cite[pp.~232--234]{Auckly1},
\cite[pp.~63--70]{Auckly2}}]\label{thm:Auckly}
Let the situation be as above. Let $\gamma$ be a 
$\SU(2)$--representation of $\pi_{1}(X)$ equal to
$\omega(\vep,\un)[g_{1},\ldots,g_{n}]$ or
$\nu(\vep,\un)[g_{1},\ldots,g_{n}]$. Thinking of the 
conjugacy class $[\gamma]$ as a gauge equivalence class of 
flat connections we have
$$
\CS([\gamma]) = -\vep c - \vep^{2}E - \sum_{j=1}^{n} \left( 
\frac{\rho_{j}n_{j}^{2} + 2n_{j}\vep}{\alpha_{j}} \right) 
\pmod{\Z},
$$
where the $\rho_j$'s are as before and $c=0$ if $\ep=\os$ and 
$c=g$ if $\ep=\ns$.
Next assume that $E \neq 0$ and let $\gamma$ be a 
$\SU(2)$--representation of $\pi_{1}(X)$ equal to 
$\rho(m,\un)$ or $\sigma(m,\un)$. Thinking of the conjugacy
class $[\gamma]$ as a gauge equivalence class of flat 
connections we have
\begin{eqnarray*}
\CS([\gamma]) &=& \frac{1}{E}\left( m + \frac{c}{2} \right)
\left(m + \frac{c}{2} + \sum_{k=1}^{n} \frac{n_{k}}{\alpha_{k}}
\right) \\
&& - \sum_{j=1}^{n} \frac{\rho_{j}n_{j}^{2} 
- \frac{n_{j}}{E}\left(m + \frac{c}{2} + \sum_{k=1}^{n} 
\frac{n_{k}}{\alpha_{k}} \right)}{\alpha_{j}} \pmod{\Z},
\end{eqnarray*}
where $c$ and $\rho_{j}$ are as before.\HS
\end{thm}

Actually Auckly works with (normalized) Seifert invariants all
over. To see that all the above results holds when working
with non-normalized Seifert invarints simply note that all
the calculations of Auckly stay unchanged with the exception
that one should put $b=0$ all over in Auckly's calculations as
indicated above. It is instructive to give a more direct
argument for \refthm{thm:Auckly}. First simply note 
that since this theorem is true for (normalized) 
Seifert invariants it is also true for the non-normalized 
Seifert invariants
$\{\ep;g\;|\;(\alpha_{1},\beta_{1}),\ldots,
(\alpha_{n},\beta_{n}),(1,b)\;\}$, 
with $0<\beta_{j}<\alpha_j$, $j=1,2,\ldots,n$, since $n_{n+1}$
has to be zero. Secondly, we note that two sets of 
non-normalized Seifert invariants for the same Seifert 
fibration give the same
set of Chern--Simons values. Of course this has to be the case
since isomorphic Seifert fibrations are homeomorphic. However,
it is in fact also easy to see directly from the formulas in 
\refthm{thm:Auckly}. Thus adding or deleating 
an $(\alpha,\beta)$--pair $(1,0)$ preserves the set of 
Chern--Simons invariants since if the $j$th pair is such a pair
then $n_j$ has to be zero. Secondly, if we change 
$(\alpha_j,\beta_j)$ to $(\alpha_j,\beta_j +K_{j}\alpha_j)$
with $\sum_{j=1}^{n} K_{j}=0$, then we can keep the $\rho_j$'s
unchanged and thus the set of Chern--Simons invariants don't
change. 

The above \refthm{thm:Auckly} actually corrects a small 
error in Auckly's result. Thus the term $-\vep c$ is missing in
the formula for $\CS(\nu(\vep,\un)[g_{1},\ldots,g_{n}])$ given
in \cite{Auckly1}, \cite{Auckly2}. Let us give a short
account for why this extra term has to be included, following 
Auckly's proof most of the way. As mentioned above Auckly bases
his proof of \refthm{thm:Auckly} on a result of Kirk and
Klassen, namely \cite[Theorem 4.2]{KirkKlassen1}. Let us work
with normalized invariants. Thus assume
$X=(\ns;g\;|\;b;(\alpha_{1},\beta_{1}),\ldots,
(\alpha_{n},\beta_{n}))$ 
and let $X_{j}$ be the space obtained
from $X$ by cutting out the interior of a tubular neighborhood
of the $j$th exceptional fiber, $j=1,\ldots,n$, and let $X_0$ 
be the space obtained from $X$ by cutting out the interior of
a tubular neighborhood of a regular fiber. Thus
\begin{eqnarray*}
\pi_{1}(X_j) &=&  \langle h,q_{1},\ldots,q_{n},a_{1},\ldots,
a_{g}\; | \; 
h^{-b}q_{1}\ldots q_{n}\prod_{i=1}^{g} a_{i}^{2}=1, \;
q_{k}^{\alpha_{k}}h^{\beta_{k}}=1, \;k \neq j, \\
&& a_{m}^{-1}ha_{m}h=[q_{k},h]=1, 
m=1,2,\ldots,g, k=1,2,\ldots,n\;\rangle
\end{eqnarray*}
and
\begin{eqnarray*}
\pi_{1}(X_0) &=&  \langle h,q_{1},\ldots,q_{n},a_{1},\ldots,
a_{g}\; | \; q_{k}^{\alpha_{k}}h^{\beta_{k}}=[q_{k},h]=1, \;
k=1,2,\ldots,n, \\
&& a_{m}^{-1}ha_{m}h=1, \; m=1,2,\ldots,g, \;\rangle.
\end{eqnarray*} 
Let $\rho_j,\sigma_j \in \Z$ such that 
$\alpha_j\sigma_j-\beta_j\rho_j =1$ as usual. We note that
$\mu_j = q_{j}^{\alpha_j}h^{\beta_j}$ and 
$\lambda_j=q_{j}^{\rho_j}h^{\sigma_j}$ are respectively a
meridian and a longitude for the torus neighborhood in $X$
around the $j$th exceptional fiber. Moreover,
$\mu_0 = \left(h^{-b}q_{1}\ldots q_{n}\prod_{k=1}^{g} 
a_{k}^{2} \right)^{-1}$ 
and $\lambda_0=h$ are a meridian and longitude respectively for
the torus neighborhood in $X$ around a regular fiber. 
Let $\nu = \nu(\vep,\un)[g_1,\ldots,g_n]$ and let us define 
curves 
$\nu_j = \nu_{j}(\vep,\un)[g_{1},\ldots,g_n]
: [0,1] \to \Hom(\pi_1(N_j),\SU(2))$
by $\nu_j(t)(x)=\nu(x)$ for $t \in [0,1]$ and 
$x = h$, $x=a_1,\ldots,a_{g-1}$ and $x=q_1,\ldots,q_{j-1}$
(a curve being a continuous curve here and in the following). 
Moreover, we let
$\nu_j(t)(q_i) = g_i 
\exp\left(-2\pi i\vep \frac{\beta_i}{\alpha_i}\right)
g_{i}^{-1}$ 
for $i=j+1,\ldots,n$ and
$\nu_j(t)(q_j) = g_j 
\exp\left(2\pi i\frac{n_{j}t-\vep\beta_j}{\alpha_j}\right)
g_{j}^{-1}$.
Finally we define $\nu_j(t)(A_g)$ in the following way: First
choose curves $h_j : [0,1] \to \SU(2)$ such that
$h_j(t)^{-1}\nu_j(t)(q_{n}^{-1}\ldots q_{1}^{-1}h^{b})h_j(t) 
\in S^{1}$
and choose curves $\varphi_j : [0,1] \to \R$ such that 
$h_j(t)^{-1}\nu_j(t)(q_{n}^{-1}\ldots q_{1}^{-1}h^{b})h_j(t) =
e^{2 i\varphi_j(t)}$.
Since $\nu_j(1)(x)=\nu_{j+1}(0)(x)$ for $j=1,2,\ldots,n$ and
$x=h,q_1,\ldots,q_n,a_1,\ldots,a_{g-1}$, where 
$\nu_{n+1}(0) = \nu$, we can choose the $h_j$ and $\varphi_j$
such that $h_j(1)=h_{j+1}(0)$ and 
$\varphi_j(1)=\varphi_{j+1}(0)$, $j=1,\ldots,n-1$ and such
that $h_n(1)e^{i\varphi_n(1)}h_n(1)^{-1} = \nu(a_g)$. Thus
$u_j : [0,1] \to \SU(2)$ given by 
$u_j(t) = h_j(t)e^{i\varphi_j(t)}h_j(t)^{-1}$ define curves 
such that 
$u_j(t)^{2} = \left(\nu_j(t)\left(h^{-b}q_1 \ldots q_n 
\prod_{k=1}^{g-1} a_{k}^{2} \right)\right)^{-1}$
and such that $u_j(1)=u_{j+1}(0)$, $j=1,\ldots,n-1$ and
$u_n(1)=\nu(a_g)$. We let $\nu_j(t)(a_g) = u_j(t)$. Thus
$\nu_{j}(1)=\nu_{j+1}(0)$ and $\nu_{n}(1)=\nu$. Also we note
that 
$\nu_{j}(0) = \nu(\vep,n_1,\ldots,n_{j-1},0,\ldots,0)
[g_1,\ldots,g_n]$
so $\nu_{j}(0)$ and $\nu_{j}(1)$ both extend to 
$\SU(2)$--representations of $\pi_1(X)$. We note that
$\nu_j(t)(\mu_j) = g_j e^{2\pi i a_j(t)} g_{j}^{-1}$ and
$\nu_j(t)(\lambda_j)=g_j e^{2\pi i b_j(t)} g_{j}^{-1}$,
where $a_j(t)=n_j t$ and 
$b_j(t) =(\rho_{j}n_{j}t+\vep)/\alpha_{j}$.
By \cite[Theorem 4.2]{KirkKlassen1} we thus find that
$$
\CS(\nu) - \CS(\nu_{1}(0)) = 
-2\sum_{j=1}^{n} \int_{0}^{1} b_{j}(t)a_{j}'(t) \dte t
= -\sum_{j=1}^{n} 
\frac{\rho_{j}n_{j}^{2} + 2n_{j}\vep}{\alpha_{j}} \pmod{\Z}.
$$
To calculate $\CS(\nu_1(0))$ we construct a curve
$\gamma : [0,1] \to \Hom(\pi_1(X_0),\SU(2))$ from the trivial
representation to $\nu_{1}(0)$ and apply 
\cite[Theorem 4.2]{KirkKlassen1} once more. To be specific we 
let
$$
\gamma(t)(h) = \left\{ \begin{array}{ll}
1 &, t \in [0,1/4], \\
e^{\pi i\vep(4t-1)} &, t \in [1/4,3/4], \\
e^{2\pi i\vep} &, t \in [3/4,1]
\end{array}\right.
$$
and
$$
\gamma(t)(q_j) = \left\{ \begin{array}{ll}
1 &, t \in [0,1/4], \\
e^{-\pi i\vep\frac{\beta_j}{\alpha_j}(4t-1)} &, t \in 
[1/4,3/4], \\
h_j(4t-3) e^{-2\pi i\vep \frac{\beta_j}{\alpha_j}} 
h_j(4t-3)^{-1} &, t \in [3/4,1]
\end{array}\right.
$$
for $j=1,\ldots,n$, where $h_j : [0,1] \to \SU(2)$ is a curve 
from $1$ to $g_j$. Moreover we let
$$
\gamma(t)(a_k) = \left\{ \begin{array}{ll}
u_1(4t) &, t \in [0,1/4], \\
j &, t \in [1/4,3/4], \\
u_1(4-4t) &, t \in [3/4,1]
\end{array}\right.
$$
for $k=1,\ldots,g-1$, where $u_1 : [0,1] \to \SU(2)$ is a curve
from $1$ to $j$. Finally we let
$$
\gamma(t)(a_g) = \left\{ \begin{array}{ll}
u_1(4t) &, t \in [0,1/4], \\
j &, t \in [1/4,3/4], \\
u_2(4t-3) &, t \in [3/4,1],
\end{array}\right.
$$
where $u_2 : [0,1] \to \SU(2)$ is a curve from $j$ to 
$\nu_1(0)(a_g)$. Thus $\gamma$ indeed defines a curve in
$\Hom(\pi_1(X_0),\SU(2))$ from the trivial representation to
$\nu_1(0)$. We recall that $\lambda_0 = h$ and thus have that
$\gamma(t)(\lambda_0)=e^{2\pi i b(t)}$, where
$$
b(t) = \left\{ \begin{array}{ll}
0 &, t \in [0,1/4], \\
\frac{\vep}{2}(4t-1) &, t \in [1/4,3/4], \\
\vep &, t \in [3/4,1].
\end{array}\right.
$$
Choose piecewise smooth curves $g : [0,1] \to \SU(2)$ and
$a : [0,1] \to \R$ such that 
$\gamma(t)(\mu_0) = g(t)^{-1}e^{2\pi i a(t)}g(t)$. By
\cite[Theorem 4.2]{KirkKlassen1} we find that
$$
\CS(\gamma(1)) - \CS(\gamma(0)) = 
-2\int_{0}^{1} b(t)a'(t) \dte t \pmod{\Z}.
$$
Here
$$
-2 \int_{0}^{1} b(t)a'(t) \dte t =
-4\vep \int_{1/4}^{3/4} t a'(t) \dte t + 
\vep\left(3a(3/4) - a(1/4) \right) - 2\vep a(1).
$$
For $t \in [1/4,3/4]$ we can choose $g(t)=1$ and
$a(t) = -\frac{g}{2} - \frac{\vep}{2} E (4t-1)$, where $E$
is the Seifert Euler number as usual. This gives
$$
-2 \int_{0}^{1} b(t)a'(t) \dte t =
-\vep g - \vep^{2}E - 2\vep a(1).
$$
Finally we note that $a(1) \in \Z$ since 
$\gamma(1)(\mu_0)=1$, since $\gamma(1)=\nu_{1}(0)$ extends to
a representation of $\pi_1(X)$. 

We note that the Chern--Simons invariants only depend on the 
genus of the base in the case of non-orientable base and that 
the set of Chern--Simons invariants in that case only depends 
on the parity of the genus. The 
Chern--Simons value of $\sigma(m,\un)$ for genus $g$ is 
equal to the Chern--Simons value of $\sigma(m+l,\un)$ for
genus $g-2l$. We can of course replace the $\omega$ and $\nu$
above by the $\omega'$ and $\nu'$.

We are now ready to show that the values of the $q$--functions 
in (\ref{eq:CSa}) and (\ref{eq:CSb}) are equal to the 
Chern--Simons invariants of the flat $\SU(2)$--connections on 
the Seifert manifold. In fact we find that the values differ by
a sign. This sign discrepancy is either due to a general sign 
error in \cite{KirkKlassen1} or else the AEC is only true if we
work with the complex conjugated invariants (which are the 
invariants associated to the mirror category of the modular 
category induced by the representation theory of 
$U_{q}(\frsl_{2}(\C))$). A similar phenomenon was observed in
\cite{AndersenHansen}. In any case this sign is not serious 
and is solely related to sign conventions.

Let $X$ and $\beta_0$ be as before and let $\gamma$ be an equal
to $\omega'(l,\un')[g_{1},\ldots,g_{n}]$ or
$\nu'(l,\un')[g_{1},\ldots,g_{n}]$, where $(l,\un') \in \mI_2$
and $g_1,\ldots,g_n \in \SU(2)$. Let $\vep=l/2$ and 
$\un = \un' + \vep\ubeta$. We assume that
(\ref{eq:genus0condition}) is satisfied if $(\ep,g)=(\os,0)$
i.e.\ we assume that
\begin{equation}\label{eq:genus0condition'}
e^{-\pi i\beta_{0}l}\prod_{j=1}^{n} g_{j}
\exp\left(2\pi i \frac{n_{j}'}{\alpha_j}\right)g_{j}^{-1}=1
\end{equation}
in that case. 
The number $q_{(l,\un')}$ in (\ref{eq:CSb}) is given by
\begin{eqnarray*}
q_{(l,\un')} &=& \sum_{j=1}^{n} 
\frac{\rho_{j}}{\alpha_{j}} {n_{j}'}^{2}
- \vep^{2} \sum_{j=1}^{n} \sigma_{j}\beta_{j} \\
&=& \sum_{j=1}^{n} \frac{\rho_{j}}{\alpha_{j}} n_{j}^{2}
+ \vep^{2} \sum_{j=1}^{n} 
\left( \frac{\rho_{j}\beta_{j}^{2}}{\alpha_{j}} 
- \sigma_{j}\beta_{j} \right)
- 2\vep \sum_{j=1}^{n} 
\frac{\rho_{j}\beta_{j}}{\alpha_{j}} n_{j} \pmod{\Z}.
\end{eqnarray*}
By using that $\rho_{j}\beta_{j}=\alpha_{j}\sigma_{j}-1$ we get
$\sum_{j=1}^{n} \frac{\rho_{j}\beta_{j}^{2}}{\alpha_{j}}
= E + \sum_{j=1}^{n} \beta_{j}\sigma_{j}$ so
\begin{eqnarray*}
q_{(l,\un')} &=& \sum_{j=1}^{n} 
\frac{\rho_{j}}{\alpha_{j}} n_{j}^{2}
+ \vep^{2} E - 2\vep \sum_{j=1}^{n} 
\frac{\alpha_{j}\sigma_{j}-1}{\alpha_{j}} n_{j} \\
&=& \vep^{2} E + \sum_{j=1}^{n} \frac{\rho_{j}n_{j}^{2} 
+ 2 n_{j}\vep}{\alpha_{j}} \pmod{\Z},
\end{eqnarray*}
where we use that 
$2\vep \sum_{j=1}^{n} \sigma_{j}n_{j} \in \Z$. 
We have thus shown that
$$
q_{(l,\un')} = - \CS([\gamma]).
$$
Let us next assume that the Seifert Euler number is nonzero and
let $\gamma$ be equal to $\rho(m,\un)$ or $\sigma(m,\un)$. Then
$$
\CS([\gamma]) = - \sum_{j}^{n} 
\frac{\rho_{j}}{\alpha_{j}} n_{j}^{2}
+ \frac{1}{4} E \left( \frac{2}{E}\left[ m + \frac{c}{2}
+ \sum_{k=1}^{n} \frac{n_{k}}{\alpha_{k}} \right] \right)^{2}
\pmod{\Z},
$$
where $c$ is the number of cross caps in the base of $X$. If 
$c$ is even we see that
$$
\CS([\gamma]) = -q_{(m',\un)},
$$
where $m'= -m - \frac{c}{2}$ and $q_{(m',\un)}$ is given by 
(\ref{eq:CSa}).

There is a small appropriate remark to make here. Namely, in 
our asymptotic formula in \refthm{thm:asymptotics-Enot0}, only
the values $q_{(m,\un)}$ for which $z_{\st}(m,\un) \in ]0,1[$ 
are present. By the symmetry considerations above
\refcor{cor:Zint1Seifert} we indeed have that the image of the
$q$--function in (\ref{eq:CSa}) is given by the values
$q_{(m,\un)}$ for which $z_{\st}(m,\un) \in [0,1]$. If
$z_{\st}(m,\un)=l \in \{0,1\}$, put 
$n_{j}'=n_{j}+\frac{1}{2}\beta_{j}l$ for all $j$ and get
that $\sum_{j=1}^{n} \frac{n_{j}'}{\alpha_{j}} = m \in \Z$.
Thus $(l,\un') \in \mI_{2}^{a}$. Moreover, $q_{(m,\un)}$ in
(\ref{eq:CSa}) is equal to $q_{(l,\un')}$ in (\ref{eq:CSb}).

\subsection{The genus 0 case}\label{secgenus0}

In the case where the base is $S^{2}$ it follows from the 
previous section that not all values $q_{(l,\un')}$, 
$(l,\un') \in \mI_2$, in (\ref{eq:CSb}) have to be 
Chern--Simons values. Namely, the point $(l,\un')$ corresponds 
to a representation if and only if there exists
$g_{1},\ldots,g_{n} \in \SU(2)$ such that 
(\ref{eq:genus0condition'}) is satisfied. Of course it could
happen that $q_{(l,\un')}$ is still equal to a Chern--Simons
value for a point $(l,\un') \in \mI_2$ that does not satisfy
(\ref{eq:genus0condition'}) but there are certainly cases where
this is not the case. To give an example consider the
manifold $M_{-1}$ obtained by surgery on the $3$--sphere along
the figure $8$ knot with framing $-1$. This manifold is equal
to the Seifert manifold
$(\os;0\;|\;-1;(2,1),(3,1),(7,1))$. It is easy to see that
$\mI_{1}$ is empty in this case. One finds that 
$\mI_{2}^{a}=\{(l=0,\un'=\uzero)\}$ and that $\mI_{2}^{b}$
contains 15 points $(l,\un')$ with $l=0$ and 8 points 
$(l,\un')$ with $l=1$. However, only the single point in
$\mI_{2}^{a}$ and two points
$(l,\un') \in \{(1,(1/2,1/2,3/2)),(1,(1/2,1/2,5/2))\}$
in $\mI_{2}^{b}$ correspond to $\SU(2)$--representations
of $\pi_{1}(M_{-1})$, namely they are the only three points
in $\mI_{2}$ satisfying (\ref{eq:modulin'}) below. The point
in $\mI_{2}^{a}$ corresponds to the trivial connection and
has Chern--Simons invariant $0$. For the two points in
$\mI_{2}^{b}$ corresponding to connections we find
$q_{(1,(1/2,1/2,3/2))} = 121/168 \pmod{\Z}$ and
$q_{(1,(1/2,1/2,5/2))} = 25/168 \pmod{\Z}$. The other 21 points
in $\mI_{2}^{b}$ have the following 21 different values for 
$q_{(l,\un')}$ (all $\pmod{\Z}$): 1/2,2/3,1/6,3/7,5/7,6/7,5/8,
3/14,5/14,13/14,2/21,8/21,11/21,7/24,1/42,25/42,37/42,3/56,
19/56,27/56,1/168, none of which are equal to one of our three 
Chern--Simons invariants.

Thus, to finalize the proof of the AEC 
for the Seifert manifolds $X$ with base $S^2$ we have to prove 
that the terms in our large $r$ asymptotic formulas for 
$\tau_{r}(X)$ in Theorems~\ref{thm:asymptotics-Enot0} and 
\ref{thm:asymptotics-E0}, corresponding to points 
$(l,\un') \in \mI_{2}$ such that $q_{(l,\un')}$ is not a
Chern--Simons value, are zero. In the genus zero case
the Laurent polynomial $Z_{0}^{(l,\un')}(r)$ is always zero
so left is to prove that $Z_{1}^{(l,\un')}(r)=0$ and that the
coefficients $c_{k}^{(l,\un')}$ in (\ref{eq:A3'}) are zero
for such a 'non-contributing' point $(l,\un') \in \mI_{2}$.

Since $je^{i\theta}j^{-1}=e^{-i\theta}$ we have that there for
all $(l,\un') \in \mI_{2}^{a}$ exists
$g_{1},\ldots,g_{n} \in \SU(2)$ such that 
(\ref{eq:genus0condition'}) is satisfied. Thus the part
$Z_{\inte}(X;r)$ in \refthm{thm:asymptotics-Enot0} and also the
part $Z_{\spec}(X;r)$ in \refthm{thm:asymptotics-E0} have an
asymptotic expansion in accordance with the AEC, so the AEC
will follow if we can prove the following

\begin{conj}
Let $X$ be a Seifert manifold with base $S^{2}$ and let
$(l,\un') \in \mI_{2}^{b}$ be a point such that there do not 
exist $g_{1},\ldots,g_{n} \in \SU(2)$ such that
{\em (\ref{eq:genus0condition'})} is satisfied. Then
$Z_{1}^{(l,\un')}(r)=0$.
\end{conj}

Let us prove this conjecture for the cases with $3$ or less
exceptional fibers.

\begin{thm}\label{thm:AECSeifertn3}
The AEC is true for $G=\SU(2)$ and $X$ any Seifert 
$3$--manifold with base $S^{2}$ and $n \leq 3$ exceptional 
fibers.
\end{thm}

\begin{prf}
The case $n \leq 2$ are the lens spaces (considering
$S^{2} \times S^{1}$ a lens space) and for that case the AEC
immediately follows from our formulas since
$Z_{1}^{(l,\un')}(r)=0$ for all $(l,\un') \in \mI_{2}$ in that
case (since we calculate the residue of an entire function).
(As mentioned in the introduction Jeffrey has already proved 
that the AEC holds for the lens spaces.)

Let us next look at the case $n=3$. Let us work with
non-normalized Seifert invariants so let
$X=\{\os;0\;|\;(\alpha_{1},\beta_{1}),(\alpha_{2},\beta_{2}),
(\alpha_{3},\beta_{3})\}$.
Let us first determine which of the pairs 
$(l,\un') \in \mI_{2}^{b}$ that corresponds to a
$\SU(2)$--representation of $\pi_{1}(X)$. 
Any element $q \in \SU(2)=\Sp_{1}$
can be written $q=r+ai+bj+ck$, where $r,a,b,c \in \R$ with
$r^2 +a^2 + b^2 + c^2=1$. The number $r$ is called the real 
part of $q$ and any two elements of $\SU(2)$ are conjugate if 
and only if they have the same real part. Thus the conjugacy 
class of any $q \in \SU(2) \sm \{\pm 1\}$ is topologically a 
$2$-sphere. Now, for any element $q \in \SU(2)$ we write
$q=\cos(\theta)+\sin(\theta)(ai+bj+ck)$ for a unique
$\theta \in [0,\pi]$ and $a,b,c \in \R$ with 
$a^2 + b^2 + c^2=1$. (If $q=\pm 1$, then of course 
$\sin(\theta)=0$ and $a,b,c$ are redundant.) Recall here 
that any element $q \in \SU(2)$ is conjugate to an 
element of $S^{1} = \{ x+yi \;|\;x,y \in \R,\;x^2 + y^2=1 \;\} 
\subseteq \SU(2)$. Thus, $q=ge^{i\theta}g^{-1}$ for a 
$g \in \SU(2)$, $\theta$ as above. Let $S_{\theta}$ be the 
conjugacy class consisting of the elements of $\SU(2)$ with 
real part $\cos(\theta)$, $\theta \in [0,\pi]$.

Let us find all possible 
$(\theta_1,\theta_2,\theta_3) \in [0,\pi]^{3}$ such that
\begin{equation}\label{eq:modulicondition}
g_{1}e^{i\theta_1}g_{1}^{-1}g_{2}e^{i\theta_{2}}g_{2}^{-1}
g_{3}e^{i\theta_3}g_{3}^{-1} =1
\end{equation}
for some $g_{1},g_{2},g_{3} \in \SU(2)$. By conjugation we
can assume that $g_{1}=1$ and the above equation then becomes
$$
e^{i\theta_1}ge^{i\theta_{2}}g^{-1}
=he^{i\theta_3}h^{-1}
$$
for some $g,h \in \SU(2)$, that is, we search for triples
$(\theta_1,\theta_2,\theta_3) \in [0,\pi]^{3}$ such that there
exists a $u \in S_{\theta_2}$ with 
$e^{i\theta_1}u \in S_{\theta_{3}}$, i.e.\ such that the real
part of $e^{i\theta_1}u$ is $\cos(\theta_{3})$.
Writing $u =\cos(\theta_{2})+\sin(\theta_{2})(ai+bj+ck)$ 
we find that the real part of $e^{i\theta_1}u$ is
$\cos(\theta_1)\cos(\theta_2)-a\sin(\theta_1)\sin(\theta_2)$.
If $\theta_{2} = 0$ then $\theta_{3}=\theta_1$. If 
$\theta_{2} =\pi$, then $\theta_{3}=\pi-\theta_{1}$. In the 
other cases $S_{\theta_{2}}$ is a $2$--sphere and $a$ runs 
through $[-1,1]$ when $u$ runs through $S_{\theta}$. Thus we 
get that (\ref{eq:modulicondition}) is satisfied for some
$g_{1},g_{2},g_{3} \in \SU(2)$ if and only if
\begin{equation}\label{eq:thetamoduli}
\left|\theta_{1}-\theta_{2}\right| \leq \theta_{3}
\leq \min\{\theta_{1}+\theta_{2},2\pi-\theta_{1}-\theta_{2}\}.
\end{equation}
Note that the special cases $\theta_{2} \in \{0,\pi\}$ are
covered by this condition. Thus 
$\omega'(l,\un')[g_{1},g_{2},g_{3}]$ is a 
$\SU(2)$--representation of $\pi_{1}(X)$ if and only if
$(l,\un') \in \mI_{2}$ satisfies the condition
\begin{equation}\label{eq:modulin'}
\left|\frac{n_{1}'}{\alpha_{1}} -
\frac{n_{2}'}{\alpha_{2}}\right| \leq 
\frac{n_{3}'}{\alpha_{3}} \leq
\min\left\{ \frac{n_{1}'}{\alpha_{1}} +
\frac{n_{2}'}{\alpha_{2}}, 1- \frac{n_{1}'}{\alpha_{1}} -
\frac{n_{2}'}{\alpha_{2}} \right\}.
\end{equation}
The points $(l,\un') \in \mI_{2}^{a}$ are the points
satisfying (\ref{eq:modulin'}) with at least one of the two 
$\leq$ being an equality.
We have to prove that if $(l,\un') \in \mI_{2}^{b}$ does not
satisfy (\ref{eq:modulin'}), then $Z_{1}^{(l,\un')}(r)=0$. It
is not hard to see that
\begin{eqnarray}\label{eq:Z1specific}
Z_{1}^{(l,\un')}(r) &=& \frac{i}{4} \left( \prod_{j=1}^{3}
\sin\left(\pi \rho_{j} \frac{2n_{j}'}{\alpha_{j}} \right) 
\right)
\sum_{\stackrel{\umu' \in \{\pm 1\}^{3}}{a(\umu',\un') > 0}}
\prod_{j=1}^{3} \mu_{j}'  \nonumber \\
&& \sum_{\stackrel{m \in \Z}{0 \leq m \leq a(\umu',\un')}} 
\frac{1}{\Sym_{\pm}(m)\Sym_{\pm}(m-a(\umu',\un'))}
\end{eqnarray}
for all $(l,\un') \in \mI_{2}$, where 
$a(\umu',\un')=\mualphammsum$.
Now let $(l,\un') \in \mI_{2}^{b}$. Then
$$
Z_{1}^{(l,\un')}(r) = \frac{i}{4} \left( \prod_{j=1}^{3}
\sin\left(\pi \rho_{j} \frac{2n_{j}'}{\alpha_{j}} \right) 
\right)
\sum_{\stackrel{\umu' \in \{\pm 1\}^{3}}{a(\umu',\un') > 0}}
\prod_{j=1}^{3} \mu_{j}'
\sum_{\stackrel{m \in \Z}{0 \leq m < a(\umu',\un')}} 
\frac{1}{\Sym_{\pm}(m)}.
$$
The double inequality (\ref{eq:modulin'}) is 
completely symmetric in $n_{1}'$, $n_{2}'$ and $n_{3}'$. 
This should be so since we in the above argument 
have numbered the exceptional fibers aribitrary. To see the
symmetry directly from (\ref{eq:modulin'}), note that 
$(\theta_{1},\theta_{2},\theta_{3})$ satisfies 
(\ref{eq:thetamoduli}) if and only if there exist
$u_{j} \in S_{\theta_{j}}$ such that $u_{1}u_{2}u_{3} =1$. But
this is equivalent to $u_{2}u_{3}u_{1}^{-1}=1$ etc.\ (where of 
course $u_{1}^{-1} \in S_{\theta_{1}}$).
Thus it is enough to consider the case 
$\frac{n_{1}'}{\alpha_{1}} \leq \frac{n_{2}'}{\alpha_{2}}
\leq \frac{n_{3}'}{\alpha_{3}}$. Assume this and assume that
$\un'$ does not satisfy (\ref{eq:modulin'}). Since
the first inequality in (\ref{eq:modulin'}) is
satisfied we have that the second is not satisfied. 
Assume first that 
$\sum_{j=1}^{3} \frac{n_{j}'}{\alpha_j} >1$ (in which case
we have that the right-hand side of (\ref{eq:modulin'}) is
$1- \frac{n_{1}'}{\alpha_1} - \frac{n_{2}'}{\alpha_2}$
and that the second inequality in (\ref{eq:modulin'}) is not
satisfied). In this case we have
$$
Z_{1}^{(l,\un')}(r) = \frac{i}{4} \left( \prod_{j=1}^{3}
\sin\left(\pi \rho_{j} \frac{2n_{j}'}{\alpha_{j}} \right) 
\right) \left( \frac{3}{2} +
\frac{1}{2}
\sum_{\stackrel{\umu' \in \{\pm 1\}^{3}}{0<a(\umu',\un') < 1}}
\prod_{j=1}^{3} \mu_{j}' \right).
$$
By going through all the $8$ possibilities for 
$\umu' \in \{\pm 1\}^{3}$ one finds that 
$a(\umu',\un') \in ]0,1[$ if and only if
$\umu' \in \{ (1,1,-1), (1,-1,1), (-1,1,1) \}$, thus
$Z_{1}^{(l,\un')}(r)=0$. 

Finally assume that
$\sum_{j=1}^{3} \frac{n_{j}'}{\alpha_j} < 1$.
In that case
$$
Z_{1}^{(l,\un')}(r) = \frac{i}{8} \left( \prod_{j=1}^{3}
\sin\left(\pi \rho_{j} \frac{2n_{j}'}{\alpha_{j}} \right)
\right) 
\sum_{\stackrel{\umu' \in \{\pm 1\}^{3}}{0<a(\umu',\un') < 1}}
\prod_{j=1}^{3} \mu_{j}'.
$$
We still restrict to the case
$\frac{n_{1}'}{\alpha_{1}} \leq \frac{n_{2}'}{\alpha_{2}}
\leq \frac{n_{3}'}{\alpha_{3}}$
and have then that the second inequality in 
(\ref{eq:modulin'}) is not satisfied if and only if
$\frac{n_{3}'}{\alpha_{3}} > \frac{n_{1}'}{\alpha_1} +
\frac{n_{2}'}{\alpha_2}$. Assuming this we find that
$a(\umu',\un') \in ]0,1[$ if and only if
$\umu' \in \{(1,1,1),(1,-1,1),(-1,1,1),(-1,-1,1)\}$, so
$Z_{1}^{(l,\un')}(r) = 0$.
\end{prf}

Let us for completeness also calculate $Z_{1}^{(l,\un')}(r)$
for $(l,\un') \in \mI_{2}^{a}$. Note first that
$\rho_{j}\frac{2n_{j}'}{\alpha_{j}} \in \Z$ if and only if
$n_{j}' \in \{0,\alpha_{j}/2\}$ and in these cases 
$Z_{1}^{(l,\un')}(r)=0$ by (\ref{eq:Z1specific}).
For $\umu' \in \{\pm 1\}^{3}$,
let $a(\umu,\un')=\mualphammsum$ as above. By assumption there
exists a $\umu' \in \{\pm 1\}^{3}$ such that
$a(\umu',\un') \in \{0,1\}$. 
Assume first that $a(\umu',\un')=1$ for a 
$\umu' \in \{\pm 1\}^{3}$. Then
either $\mu_{j}'=1$ for $j=1,2,3$ or there is a 
$j \in \{1,2,3\}$ with $n_{j}'=0$. In the last case
$Z_{1}^{(l,\un)}(r)=0$. In both cases 
$\sum_{j=1}^{3} \frac{n_{j}'}{\alpha_{j}}=1$.
By (\ref{eq:Z1specific}) we get that
$$
Z_{1}^{(l,\un')}(r) = \frac{i}{4} \left( \prod_{j=1}^{3}
\sin\left(\pi \rho_{j} \frac{2n_{j}'}{\alpha_{j}} \right) 
\right)
\left(1+\frac{1}{2}
\sum_{\stackrel{\umu' \in \{\pm 1\}^{3}}{0<a(\umu',\un') <1}}
\prod_{j=1}^{3} \mu_{j}' \right).
$$
By symmetry we can assume that
$\frac{n_{1}'}{\alpha_{1}} \leq \frac{n_{2}'}{\alpha_{2}}
\leq \frac{n_{3}'}{\alpha_{3}}$
as usual. If $\frac{n_{3}'}{\alpha_{3}}=\frac{1}{2}$, then
$Z_{1}^{(l,\un')}=0$. If that is not the case then
$a(\umu',\un') \in ]0,1[$ if and
only if $\umu' \in \{(1,1,-1),(1,-1,1),(-1,1,1)\}$ and
\begin{equation}\label{eq:Z1n3Ia}
Z_{1}^{(l,\un')}(r) = -\frac{i}{8} \prod_{j=1}^{3}
\sin\left(\pi \rho_{j} \frac{2n_{j}'}{\alpha_{j}} \right) 
\end{equation}
which is also valid in case 
$\frac{n_{3}'}{\alpha_{3}}=\frac{1}{2}$.

Next assume that $a(\umu',\un')=0$ for a 
$\umu' \in \{\pm 1\}^{3}$. If 
$\sum_{j=1}^{3} \frac{n_{j}'}{\alpha_{j}}=1$ we are in the
above case, so assume that 
$\sum_{j=1}^{3} \frac{n_{j}'}{\alpha_{j}} <1$. Assume also
that 
$\frac{n_{1}'}{\alpha_{1}} \leq \frac{n_{2}'}{\alpha_{2}}
\leq \frac{n_{3}'}{\alpha_{3}}$
as always. Then
$$
Z_{1}^{(l,\un')}(r) = \frac{i}{8} \left( \prod_{j=1}^{3}
\sin\left(\pi \rho_{j} \frac{2n_{j}'}{\alpha_{j}} \right) 
\right)
\sum_{\stackrel{\umu' \in \{\pm 1\}^{3}}{0<a(\umu',\un') <1}}
\prod_{j=1}^{3} \mu_{j}'.
$$   
In this case we find that $a(\umu',\un') \in ]0,1[$ if and only
if $\umu' \in \{(1,1,1),(1,-1,1),(-1,1,1)\}$ so again we find
that $Z_{1}^{(l,\un')}(r)$ is given by (\ref{eq:Z1n3Ia}). Thus
we have proved that
$Z_{1}^{(l,\un')}(r)$ is given by (\ref{eq:Z1n3Ia})
for all $(l,\un') \in \mI_{2}^{a}$.

To prove the AEC for Seifert manifolds with base $S^{2}$ and
$n \geq 4$ exceptional fibers, one can proceed as in the above
proof, but the combinatorics is of course harder.
For a given $(\theta_{1},\ldots,\theta_{n}) \in [0,\pi]^{n}$
there exist $g_{1},\ldots,g_{n} \in \SU(2)$ such that
$$
\prod_{j=1}^{n} g_{j}e^{i\theta_{j}}g_{j}^{-1} =1
$$
if and only if there exist $u_{j} \in S_{\theta_{j}}$, 
$j=2,3,\ldots,n-1$ such that the real part of
$e^{i\theta_{1}}u_{2}\ldots u_{n-1}$ is $\cos(\theta_{n})$.
Note that the function that maps $(u_{2},\ldots,u_{n-1})$ to
the real part of $e^{i\theta_{1}}u_{2}\ldots u_{n-1}$ is a
continuous map $S_{\theta_{2}} \times \ldots 
S_{\theta_{n-1}} \to [-1,1]$, so the image is an interval.
Thus this analysis leads to a condition of the form
$$
a \leq \theta_{3} \leq b
$$
where $a,b \in [0,\pi]$ are continuous functions in
$(\theta_{1},\theta_{2},\ldots,\theta_{n-1})$. The condition
(\ref{eq:thetamoduli}) is the special case where $n=3$.

Let $X$ be an integral homology sphere with $n \geq 3$ 
exceptional fibers and let $\mM_{\rho}$ be the component in
the moduli space of irreducible flat $\SU(2)$--connections on 
$X$ containing the class represented by the connection 
corresponding to a representation $\rho:\pi_{1}(X) \to \SU(2)$.
Then Fintushel and Stern \cite{FintushelStern} have proved that
$\mM_{\rho}$ is a closed manifold of dimension $2m-6$ if
$\rho(q_{j}) \notin \{\pm 1\}$ for exactly
$m \geq 3$ of the $j \in \{1,2,\ldots,n\}$.
(Note that the only reducible connection is the trivial
connection in this case.) 
We also refer to \cite{KirkKlassen2} regarding results on the 
components of the moduli space of flat $\SU(2)$--connections
on a Seifert fibered integral homology sphere.

\subsection{Small Seifert manifolds with base $\R\text{P}^{2}$}\label{Secgenus1}

Let $X$ be a Seifert manifold with base $\R\text{P}^{2}$ and
zero or one exceptional fiber. Except for a few cases $X$ is
then homeomorphic to a prism manifold, i.e.\ a Seifert manifold
of the form $(\os;0\;|\;b;(2,1),(2,1),(\alpha,\beta))$. The
few exceptions are homeomorphic to lens spaces except
$(\ns;0\;|\;0)$ which is homeomorphic to 
$\R\text{P}^{3} \# \R\text{P}^{3}$, see 
\cite[Theorem 5.1]{JankinsNeumann} and
\cite[Sec.~5.4 (vii) p.~102]{Orlik} (there seem to be some
small mistakes in \cite{Orlik}. Thus $(\ns;1\;|\;\pm 1)$ both
have cyclic fundamental group, namely $\Z_{4}$, and
$(\ns;1\;|\;0;(\alpha,1)\}$ also has cyclic fundamental group,
namely $\Z_{4\alpha}$. 
All three cases are homeomorphic to lens spaces
according to \cite[Theorem 5.1]{JankinsNeumann}.)
Since $\R\text{P}^{3}$ is homeomorphic to $L(2,1)$ and since
the RT--invariants are multiplicative under connected sums
and since the AEC behaves nicely w.r.t.\ products of 
RT--invariants we obtain by \refthm{thm:AECSeifertn3} the 
following

\begin{cor}\label{cor:oddcase}
The AEC is true for $G=\SU(2)$ and $X$ any Seifert manifold 
with base $\R\text{P}^{2}$ and zero or one exceptional fiber.
\end{cor}

The AEC can actually be proved directly for
$X=(\ns;1\;|\;0)$ by a small calculation: By 
(\ref{eq:Zfunction-Seifert}) we namely have
$$
Z(X;r)=\sum_{\gamma=1}^{r-1} (-1)^{\gamma}
\sin\left(\frac{\pi}{r}\gamma\right),
$$
and by an elementary calculation we find that
$$
Z(X;r) = -f\left(\frac{\pi}{r}\right) 
\left(1+e^{2\pi i r \frac{1}{2}} \right),
$$
where 
$f(x) = \frac{\sin(x)}{2\left( 1 + \cos(x)\right)}$. Thus we
find that
\begin{equation}\label{eq:invariantng1}
\tau_{r}(X) = \frac{1}{\sqrt{2r}}f\left(\frac{\pi}{r}\right)
\left(1+e^{2\pi i r \frac{1}{2}} \right).
\end{equation}
By \refthm{thm:Auckly} the Chern--Simons invariants of flat 
$\SU(2)$--connections on $X$ are $0 \pmod{\Z}$ and
$1/2 \pmod{\Z}$.

\subsection{The contribution from the trivial connection and 
the Casson--Walker invariant}

From a topological point of view it is not necessarily a good 
idea to bring the asymptotics of the quantum invariants on 
minimal form; it is actually more natural to sum over the 
connected components of the moduli space instead of over the 
image set of the Chern--Simons functional. Let us consider the
case of lens spaces. The lens spaces are identical with the 
Seifert manifolds with base $S^{2}$ and zero, one or two 
exceptional fibers. Let $(p,q)$ be a pair of coprime integers 
and let $L(p,q)$ be the associated lens space given by surgery
on $S^{3}$ along the unknot with surgery coefficient $-p/q$. 
Thus we can and will assume in the following that $p \geq 0$. 
We have $\pi_{1}(L(p,q))=\Z/p\Z$ for $p> 0$, and
$L(0,1)=S^{2} \times S^{1}$ with fundamental group $\Z$. Hence
$L(p,q)$ is a rational homology sphere iff $p \neq 0$ and an 
integral homology sphere iff $p= 1$, and $L(p,q)$ and 
$L(p',q')$ are nonhomeomorphic if $p' \neq p$. It is a 
classical result that $L(p,q)$ and $L(p,q')$ are homeomorphic 
if and only if $q' = \pm q \bmod{p}$ or $qq'= \pm 1 \bmod{p}$.
A homeomorphism is orientation preserving if and only if the 
relevant sign is $+$. In particular, we can always assume that
$0<q<p$ if $|p|>1$.

Let us compare $\tau_r(L(p,q))$ with the Casson--Walker 
invariant of $L(p,q)$, $p \neq 0$. First we note that 
$L(1,q) = S^{3}$ for all $q \in \Z$. In our normalization we 
have
$$
\tau_{r}(S^{3}) = \sqrt{\frac{2}{r}} 
\sin\left(\frac{\pi}{r}\right) = 
\sqrt{\frac{2}{r}}\frac{\pi}{r} 
\left( 1 + c_{1} r^{-1} + c_{2}r^{-2} + \ldots\right),
$$
where $c_{1}=0$. Note also that $\lambda_{\tC}(S^{3})=0$, where
$\lambda_{\tC}$ is the Casson--invariant of integral homology 
spheres, cf.\ \cite{AkbulutMcCarthy}. We note that there is 
only one flat $SU(2)$--connection on $S^{3}$, namely the 
trivial one.

Next let us investigate the lens spaces with $|p|>1$. Assume 
that $p > q > 0$ and let $q^{*}$ denote the inverse of 
$q \bmod{p}$. By \cite[Theorem 5]{HansenTakata} we have
$$
\tau_{r}(L(p,q)) = \sqrt{\frac{2}{pr}} 
\exp\left(\frac{\pi i}{2r} \dS \left(\frac{q}{p} \right)\right)
\sum_{n=0}^{p-1} c_{n}(r)
\exp\left(2\pi ir \frac{q^{*}}{p}n^{2} \right),
$$
where
$$
c_{n}(r) = i\cos\left(\frac{\pi}{pr}\right)
\sin\left(2\pi \frac{q^{*}}{p}n\right)
\sin\left(\frac{2\pi n}{p} \right)
+ \sin\left(\frac{\pi}{pr}\right)
\cos\left(2\pi \frac{q^{*}}{p}n\right)
\cos\left(\frac{2\pi n}{p} \right).
$$
(Use that $L(p,q) \cong L(p,q^{*})$ and that 
$\dS(q^{*}/p)=\dS(q/p)$. The shift from $q$ to $q^{*}$ is done
to facilitate a comparison to Chern--Simons invariants, see 
below.) The moduli space $\mM_{L(p,q)}$ of flat 
$\SU(2)$--connections on $L(p,q)$ can be identified with the 
set $\mM = \{0,1,\ldots,[p/2]\}$ via the identification of 
$\mM_{L(p,q)}$ with the set of representations 
$\Z/p\Z \to \SU(2)$ moduli conjugation by $\SU(2)$. Here $[x]$
is the integer part of $x \in [0,\infty[$. To be precise the
integer $n \in \{0,1,\ldots,[p/2]\}$ corresponds to the
representation which maps the generator $e^{2\pi i/p}$ to 
$e^{2\pi i n/p}$. By \cite[Theorem 5.1]{KirkKlassen1} the 
Chern--Simons invariant of the flat $\SU(2)$--connection 
corresponding to $n \in \mM$ is equal to 
$q_{n} = -\frac{q^{*}}{p}n^{2}$. (Note that Kirk and Klassen 
use another convention for the lens spaces. Thus $L(p,q)$ in 
\cite{KirkKlassen1} is equal to $L(p,-q)$ here.) Using a 
symmetry under the change of $n$ to $p-n$ we immediately get a
formula of the form
$$
\tau_{r}(L(p,q)) = \sqrt{\frac{2}{pr}} 
\exp\left(\frac{\pi i}{2r} \dS \left(\frac{q}{p} \right)\right)
\sum_{n=0}^{[p/2]} a_{n}(r)
\exp\left(2\pi ir \frac{q^{*}}{p}n^{2} \right).
$$
Here, for $p$ even we have $a_{n}(r) = 2c_{n}(r)$ for 
$n=1,2,\ldots,p/2-1$ and $a_{0}(r) = \sin(\pi/pr)$ and
$a_{p/2}(r) = (-1)^{1+q^{*}}\sin(\pi/pr) = \sin(\pi/pr)$. For 
$p$ odd we have $a_{n}(r) = 2c_{n}(r)$ for 
$n=1,2,\ldots,(p-1)/2$ and $a_{0}(r) = \sin(\pi/pr)$. The 
trivial connection corresponds to $n=0$. Let us write
$$
\tau_{r}^{n}(p,q) = \sqrt{\frac{2}{pr}} 
\exp\left(\frac{\pi i}{2r} \dS \left(\frac{q}{p} \right)\right)
a_{n}(r)\exp\left(2\pi ir \frac{q^{*}}{p}n^{2} \right)
$$
for $n \in \mM$. Then
$$
\tau_{r}^{0}(p,q) = \sqrt{\frac{2}{pr}} 
\sin\left(\frac{\pi}{pr}\right)
\exp\left(\frac{\pi i}{2r} \dS \left(\frac{q}{p} \right)\right)
= \sqrt{\frac{2}{pr}} \frac{\pi}{pr} 
\left( 1 + c_{1}r^{-1}+c_{2}r^{-2} +\ldots \right),
$$
where
$$
c_{1} = \frac{\pi i}{2} \dS \left(\frac{q}{p} \right)
= 6\pi i \lambda_{\tCW}(L(p,q)).
$$
Here $\lambda_{\tCW}$ is Walker's extension of the 
Casson--invariant to rational homology spheres, cf.\ 
\cite{Walker}.

Let us next investigate the sum of $\tau_{r}^{n}(p,q)$ for 
which $q_{n} = 0 \bmod{\Z}$. We have $q_{n}=0 \bmod{\Z}$ if and
only if $p$ divides $n^{2}$. Let us consider $L(9,q)$, 
$q \in \{1,2,4,5,7,8\}$. We have $q_{n} = 0 \bmod{\Z}$ if and 
only if $n \in \{0,3\}$. Moreover
$$
\tau_{r}^{3}(p,q) = 2\sqrt{\frac{2}{pr}} 
\exp\left(\frac{\pi i}{2r} \dS 
\left(\frac{q}{p} \right)\right)c_{3}(r),
$$
where
\begin{eqnarray*}
c_{3}(r) &=& i\cos\left(\frac{\pi}{pr}\right)
\sin\left(2\pi \frac{q^{*}}{3}\right)
\sin\left(\frac{2\pi}{3} \right)
+ \sin\left(\frac{\pi}{pr}\right)
\cos\left(2\pi \frac{q^{*}}{3}\right)
\cos\left(\frac{2\pi}{3} \right) \\
&=& \frac{i\sqrt{3}}{2}\cos\left(\frac{\pi}{pr}\right)
\sin\left(2\pi \frac{q^{*}}{3}\right)
- \frac{1}{2}\sin\left(\frac{\pi}{pr}\right)
\cos\left(2\pi \frac{q^{*}}{3}\right).
\end{eqnarray*}
Thus
\begin{eqnarray*}
&& \tau_{r}^{0}(p,q)+\tau_{r}^{3}(p,q) \\
&=& \sqrt{\frac{2}{pr}} 
\exp\left(\frac{\pi i}{2r} \dS 
\left(\frac{q}{p} \right)\right) \\
&& \times \left( \sin\left(\frac{\pi}{pr}\right)
\left(1-\cos\left(2\pi \frac{q^{*}}{3}\right)\right)
+ i\sqrt{3}\cos\left(\frac{\pi}{pr}\right)
\sin\left(2\pi \frac{q^{*}}{3}\right)\right) \\
&=& i\sqrt{\frac{6}{pr}}
\sin\left(2\pi \frac{q^{*}}{3}\right)
\left( 1 + c_{1}r^{-1}+c_{2}r^{-2}+\ldots\right),
\end{eqnarray*}
where
$$
c_{1} = 
\frac{\pi\left(1-\cos\left(2\pi \frac{q^{*}}{3}\right)\right)}
{ip\sqrt{3}\sin\left(2\pi \frac{q^{*}}{3}\right)}
+ \frac{\pi i}{2} \dS\left(\frac{q}{p} \right)
= 6\pi i \left( s(q,p) - 
\frac{\left(1-\cos\left(2\pi \frac{q^{*}}{3}\right)\right)}
{6p\sqrt{3}\sin\left(2\pi \frac{q^{*}}{3}\right)}\right).
$$
Note that $\s(q^{*},p)=\s(q,p)$ and put
$$
a(q) = \frac{\left(1-\cos\left(2\pi \frac{q}{3}\right)\right)}
{54\sqrt{3}\sin\left(2\pi \frac{q}{3}\right)}
$$
and $b(q)=s(q,p)-a(q)$. We find that $\s(1,9)=14/27$ and 
$\s(2,9)=4/27$, while $a(2)=-a(1)=-1/54$, so
$b(1)=\frac{27}{28}\s(1,9)$ while
$b(2)=\frac{9}{8}\s(2,9)$. We conclude that the Casson--Walker
invariant of $L(p,q)$ is naturally associated with the trivial
connection's contribution to the asymptotic expansion of 
$\tau_{r}(L(p,q))$ and not with the full contribution to this 
asymptotics coming from the zero set of the Chern--Simons 
functional.

Next let us turn to the general Seifert fibered rational 
homology spheres. A Seifert manifold
$X=(\ep;g\;|\;b;(\alpha_{1},\beta_{1}),\ldots,
(\alpha_{n},\beta_{n}))$ 
is a rational homology sphere iff $\ep=\os$, $g=0$ and 
$E \neq 0$ or $\ep=\ns$ and $g=1$. The Seifert manifold $X$ is
an integral homology sphere if and only if $\ep=\os$, $g=0$ 
and $|E|\mA=1$. Here 
$E=-b-\sum_{j=1}^{n}\frac{\beta_{j}}{\alpha_{j}}$ is the 
Seifert Euler number and $\mA=\prod_{j=1}^{n}\alpha_{j}$ as 
usual. In this paper we have not calculated the asymptotics 
for the cases where $\ep=\ns$ and $g=1$, except for the ones
with zero or one exceptional fiber. Let us first consider
the cases with base $S^{2}$ and nonzero Seifert Euler number. 
For these cases we have by 
\cite[Proposition 6.1.1]{Lescop} that
\begin{equation}\label{eq:Casson}
\lambda_{\tCW}(X) = \frac{1}{12}\left(3\sgE -E 
- \sum_{i=1}^{n} \dS\left(\frac{\beta_{i}}{\alpha_{i}}\right) 
- \frac{1}{E}\left(2-n+\sum_{i=1}^{n} \frac{1}{\alpha_{i}^{2}}
\right)\right),
\end{equation}
where we have used the facts that $|H_{1}(X;Z)|=|E|\mA$ and 
that we for any rational homology sphere $N$ have
$\lambda_{\tCLW}(N) = \frac{1}{2}|H_{1}(N;\Z)|
\lambda_{\tCW}(N)$, 
see \cite[p.~13]{Lescop}.

Under the identification of the moduli space $\mM_{X}$ of flat 
$\SU(2)$--connections on $X$ by the set of 
$\SU(2)$--representations of $\pi_{1}(X)$ moduli conjugation
by $\SU(2)$, the trivial connection is identified by the 
trivial representation. Using notation from the previous 
section we note that $\omega'(l,\un')[g_{1},\ldots,g_{n}]$ is 
the trivial representattion if and only if $l=0$ and 
$\un'=\uzero$ and in this case $\omega'$ is independent of the
$g_{j}$'s. Next assume that $E \neq 0$ and let 
$\rho=\rho(m,\un)$. We note that 
$\rho(h) = \exp\left(\pi i z_{\st}(-m,\un)\right)$, where 
$z_{\st}$ is given by (\ref{eq:SeifertStationary}). Hence 
$\rho$ is trivial only if $z_{\st}(-m,\un) \in 2\Z$ and in
this case $(-m,\un) \notin \mI_{1}$. It
should be noted that for $n \in \{0,1,2\}$, $X$ is a lens 
space, so in that case $\mM_{X}$ is a finite discrete space. 
This is also the case for $n=3$, but not for $n \geq 4$,
see also our remarks at the end of Sec.~\ref{secgenus0}. For 
$n \geq 4$ it is not known whether the connected 
components of $\mM_{X}$ are parametrized by 
$\mI_{1} \cup \mI_{2}^{a} \cup \mI_{2}^{c}$, $\mI_{2}^{c}$
being the set of points in $\mI_{2}^{b}$ satisfying
(\ref{eq:genus0condition'}).

In any case, let us calculate the contribution to the 
asymptotics of $\tau_{r}(X)$ coming from the point 
$(l,\un')=(0,\uzero)$. First note that there is no polar 
contribution, that is, $Z_{0}^{(0,\uzero)}(r)$ and 
$Z_{1}^{(0,\uzero)}(r)$ are both zero. The part of the
large $r$ asymptotics of $\tau_{r}(X)$ associated to the zero
point in $\mI_2$, denoted $\tau_{r}^{0}(X)$ in the
following, is therefore given by 
$$
\frac{i^{n+1}}{2^{n}\sqrt{r}} 
\frac{e^{-i\frac{3\pi}{4}\sgE}}{\sqrt{\mA}}
\exp \left( \frac{i \pi}{2r} \left[ 3\sgE -E - \sum_{j=1}^{n}
\dS \left(\frac{\beta_{j}}{\alpha_{j}}\right) \right] \right)
Z^{0}(X;r),
$$
where $Z^{0}(X;r)$ is the (not necessarily convergent) power 
series given as follows: Let $k_{0}=0$ if $n$ is odd and 
$k_{0}=1$ otherwise. Moreover, let $k_{1}=(k_{0}-n+1)/2$. Then
$Z^{0}(X;r) = \sum_{k=k_{1}}^{\infty} c_{k} r^{-k}$, where
$$
c_{k}=\frac{(-2i)^{n+1}}{\pi^{n-1}}
\frac{e^{i\frac{\pi}{4}\sgE}}{\sqrt{2\pi|E|}}
\frac{\Gamma\left( k+\frac{1}{2}\right)}{k'!}
\left(\frac{2i}{\pi E}\right)^{k}
\left.\partial_{z}^{(k')} 
\left\{ \left( \frac{\pi z}{\sin(\pi z)} \right)^{n-1}
\prod_{j=1}^{n+1} 
\sin\left(\frac{\pi z}{\alpha_{j}}\right)\right\}\right|_{z=0},
$$ 
where $k'=2k+n-1$ and $\alpha_{n+1}=1$. Therefore
$$
c_{k}=\frac{(-2i)^{n+1}}{\pi^{n-2}}
\frac{e^{i\frac{\pi}{4}\sgE}}{\sqrt{2\pi|E|}}
\frac{\Gamma\left( k+\frac{1}{2}\right)}{k'!}
\left(\frac{2i}{\pi E}\right)^{k}
\left.\partial_{z}^{(k'-1)} 
\left\{ \left( \frac{\pi z}{\sin(\pi z)} \right)^{n-2}
\prod_{j=1}^{n} 
\sin\left(\frac{\pi z}{\alpha_{j}}\right)\right\}\right|_{z=0}.
$$
In general we have for $m \in \{0,1,\ldots,n\}$ that
\begin{eqnarray*}
&&\left.\partial_{z}^{(l)} 
\left\{ \left( \frac{\pi z}{\sin(\pi z)} \right)^{m}
\prod_{j=1}^{n} 
\sin\left(\frac{\pi z}{\alpha_{j}}\right)\right\}
\right|_{z=0} \\
&& \hspace{1.2in} = \left\{ \begin{array}{lc}
0 & , l=0,1,\ldots,m-1 \\
\frac{l!\pi^{m}}{(l-m)!} 
\left.\partial_{z}^{(l-m)} \left\{ 
\frac{\prod_{j=1}^{n} 
\sin\left(\frac{\pi z}{\alpha_{j}}\right)}
{\sin^{m}(\pi z)} \right\}\right|_{z=0} & , l=m,m+1,\ldots.
\end{array} \right.
\end{eqnarray*}
Since we have already treated the lens spaces, we can assume 
that $n>2$. In that case $k_{1}<0$ and $k'-1=n-2$ if and only 
if $k=0$, so $c_{k}=0$ for $k=k_{1},\ldots,-1$ and
\begin{eqnarray*}
c_{k}  &=& (-2i)^{n+1}
\frac{e^{i\frac{\pi}{4}\sgE}}{\sqrt{2\pi|E|}}
\frac{\Gamma\left( k+\frac{1}{2}\right)}{(2k)!}
\left(\frac{2i}{\pi E}\right)^{k}
\left.\partial_{z}^{(2k)} 
\left\{ \frac{\prod_{j=1}^{n} 
\sin\left(\frac{\pi z}{\alpha_{j}}\right)}
{\sin^{n-2}(\pi z)} \right\} \right|_{z=0} \\
&=& (-2i)^{n+1} \frac{e^{i\frac{\pi}{4}\sgE}}{\sqrt{2|E|}}
\frac{1}{k!}\left(\frac{i}{2\pi E}\right)^{k}
\left.\partial_{z}^{(2k)} \left\{ \frac{\prod_{j=1}^{n} 
\sin\left(\frac{\pi z}{\alpha_{j}}\right)}
{\sin^{n-2}(\pi z)} \right\} \right|_{z=0}
\end{eqnarray*}
for $k=0,1,\ldots$.
Therefore
$$
\tau^{0}_{r}(X) = -i\sgE\sqrt{\frac{2}{r|E|\mA}} 
\exp \left( \frac{i \pi}{2r} \left[ 3\sgE -E - \sum_{j=1}^{n}
\dS \left(\frac{\beta_{j}}{\alpha_{j}}\right) \right] 
\right) \tilde{Z}^{0}(X;r),
$$
where $\tilde{Z}^{0}(X;r) = \sum_{k=0}^{\infty} d_{k} r^{-k}$
with
$$
d_{k}= \frac{1}{k!}
\left(\frac{i}{2\pi E}\right)^{k}
\left.\partial_{z}^{(2k)} 
\left\{ \frac{\prod_{j=1}^{n} 
\sin\left(\frac{\pi z}{\alpha_{j}}\right)}
{\sin^{n-2}(\pi z)} \right\} \right|_{z=0}.
$$
We have
$$
f(z):=\frac{\prod_{j=1}^{n} 
\sin\left(\frac{\pi z}{\alpha_{j}}\right)}
{\sin^{n-2}(\pi z)} = \frac{(\pi z)^{2}}{\mA}
\left(1 + \frac{1}{6} \left(n-2-\sum_{j=1}^{n} 
\frac{1}{\alpha_{j}^{2}}\right)
(\pi z)^{2} + c(\pi z)^{4}+\ldots\right).
$$
Thus $f(0)=0$, $f^{(2)}(0)=2\frac{\pi^{2}}{\mA}$ and
$$
f^{(4)}(0) = \frac{24\pi^{4}}{\mA}\frac{1}{6} 
\left(n-2-\sum_{j=1}^{n} \frac{1}{\alpha_{j}^{2}}\right)
= 12\pi^{2}\frac{1}{6} 
\left(n-2-\sum_{j=1}^{n} \frac{1}{\alpha_{j}^{2}}\right)
f^{(2)}(0).
$$
Hence
$$
\sum_{k=0}^{\infty} d_{k} r^{-k} = 
rd_{1}\left(1+\frac{\pi i}{2E}
\left(n-2-\sum_{j=1}^{n} \frac{1}{\alpha_{j}^{2}}\right)r^{-1}
+ a_{2}r^{-2} + \ldots\right)
$$
and finally by (\ref{eq:Casson})
$$
\tau_{r}^{0}(X) = \frac{\pi}{|E|\mA}\sqrt{\frac{2r}{|E|\mA}}
\left(1+6\pi i\lambda_{\tCW}(X) r^{-1} + b_{2}r^{-2}
+\ldots \right).
$$

Let us next consider the Seifert manifolds with base 
$\R\text{P}^{2}$ and zero or one exceptional fiber. Except for
the case $(\ns;1\;|\;0)$ all these manifolds are covered by
the above calculation, since they are all homeomorphic to
Seifert manifolds with base $S^{2}$ and nonzero Seifert Euler 
number, see Sec.~\ref{Secgenus1}. Left is to consider the case
$X=(\ns;1\;|\;0)$. By (\ref{eq:invariantng1}) we have to find
the first terms in the Taylor expansion of 
$f(x)=\frac{\sin(x)}{2(1+\cos(x))}$. But $f$ is an odd
function so we have 
$f(x)=cx\sum_{j=0}^{\infty} c_{2j}x^{2j}$ for 
$x \in ]-\pi,\pi[$. Thus
$f\left(\frac{\pi}{r}\right) = 
c\frac{\pi}{r}(1+c_{1}r^{-1}+ \pi^{2} c_{2} r^{-2} +
\ldots....)$
for $r \geq 4$, where $c_{1}=0$. But 
$\lambda_{\tCW}(\R\text{P}^{3})=0$,
thus $\lambda_{\tCW}(X)=0$, so again we find that the part
of the expansion of $\tau_{r}(X)$ related to the trivial
connection has a form 
$ar^{d}\left(1+6\pi i\lambda_{\tCW}(X)r^{-1}
+b_{2}r^{-2}+\ldots\right)$. 

The Casson--Walker invariant $\lambda_{\tCW}$ has been extended
by Lescop to an invariant $\lambda_{\tCLW}$ of all closed 
oriented $3$--manifolds, cf.\ \cite{Lescop}. We note that 
$\tau_{r}(S^{2} \times S^{1}) = 1 = 1 + c_{1}r^{-1} + 
c_{2}r^{-2}+\ldots \ldots$,
where $c_{l}=0$ for all $l=1,2,\ldots$. We have previously seen
that $S^{2} \times S^{1}$ has an infinite number of Seifert 
fibered structures. Thus 
$S^{2} \times S^{1} = (\os;0\;|\;-1;(\alpha,\beta),
(\alpha,\alpha-\beta))$ for any $\alpha \in \Z_{\geq 2}$ and 
any $\beta \in \{1,2,\ldots,\alpha-1\}$ coprime to $\alpha$.
Using this together with \cite[Proposition 6.1.1]{Lescop}, we 
find that $\lambda_{\tCLW}(S^{2} \times S^{1})=-1/12$. The 
moduli space $\mM_{S^{2} \times S^{1}}$ is topologically a 
closed interval since $\pi_{1}(S^{2} \times S^{1})=\Z$. (An 
element $\rho \in \Hom(\Z,\SU(2))$ is determined by $\rho(1)$. 
By conjugation we can assume that $\rho(1) =e^{i\theta}$,
$\theta \in [0,\pi]$.) In particular, the moduli space 
$\mM_{S^{2} \times S^{1}}$ is connected. Thus it seems that
Lescop's extension of the Casson--Walker invariant to 
$3$--manifolds not being rational homology spheres is not
part of the quantum $\SU(2)$--invariants in the same
way as the Casson--Walker invariants of rational homology 
$3$--spheres. At least this extension does not seem to be 
proportional to any of the coefficients in the asymptotic 
expansion of the quantm invariant as the simple example
$S^2 \times S^1$ reveals.

\section{Proof of \refthm{thm:asymptotics-Enot0}}
\label{sec-Proof-of}

\noindent The proof of \refthm{thm:asymptotics-Enot0} is rather
long and technical but the single steps in the proof use 
elementary analysis. To streamline the proof we will in this
section emphasize the main ideas in the arguments and defer 
technical details to a number of appendices. The proof of 
\refthm{thm:asymptotics-E0}, which uses the same ideas as the 
proof of \refthm{thm:asymptotics-Enot0}, is much shorter and is
given in Sect.~\ref{sec-Proof-of-E0}.

In the introduction we gave an outline of the proof. Let us 
recall the basic ideas. Let $X$ be a Seifert manifold described
by (normalized) Seifert invariants 
$(\ep;g\;|\;b;(\alpha_{1},\beta_{1}),\ldots,
(\alpha_{n},\beta_{n}))$
or non-normalized Seifert invariants
$\{\ep;g;(\alpha_{1},\beta_{1}),\ldots,
(\alpha_{n},\beta_{n})\}$.
The first problem we encounter when calculating the large $r$ 
asymptotics of $r \mapsto Z(X;r)$ is that this function is 
given by a sum in which both the terms and the summation range 
depend on $r$, cf.\ 
\refthm{thm:RT-invariants-Seifert-manifolds}. This problem is 
solved by using certain symmetries to replace 
$\sum_{\gamma=1}^{r-1}$ by another sum of the form
$\sum_{\gamma \in \Z}$. After this manoeuvre we use Poisson's 
summation formula to change the expression for $Z(X;r)$ 
to an infinite sum of integrals. Apparently all this seems to 
complicate things. However, we obtain the big advantage that 
the large $r$ asymptotics of the integrals in the thus obtained
sum can be calculated with a strong method called the steepest
descent method. Roughly speaking this will express $Z(X;r)$ as
an infinite sum of main contributions plus an infinite sum of 
remainder terms. The final step will be to rewrite the infinite
sum of main contributions as a finite sum in which the 
summation range does not depend on $r$ and to show that the 
infinite sum of remainder terms is small compared to this 
finite sum in the large $r$ limit.

We want to generalize the situation. Therefore, let us first 
consider the function $Z(X;r)$ in (\ref{eq:Zfunction-Seifert})
in greater detail. As mentioned above the first step is to 
rewrite this expression so that $r$ only appears in the 
summands and not in the summation range. To this end we need to
extend the function $h$ in (\ref{eq:hfunction}) to $\Z$. There 
is of course no problem in doing this by using the expression
(\ref{eq:hfunction}). However, we will need that the extended
function $h$ satisfies the following symmetries
\begin{eqnarray}\label{eq:hfunction-symmetry}
h(r\gamma) &=& 0, \\
h(-\gamma) &=& (-1)^{n} h(\gamma), \nonumber \\
h(\gamma +2r) &=& h(\gamma) \nonumber
\end{eqnarray}
valid for all $\gamma \in \Z$. To establish these symmetries
we need to examine the proof of \cite[Theorem 8.4]{Hansen2}.
Let us note here that the first symmetry $h(r\gamma)=0$,
$\gamma \in \Z$, is not satisfied if $n=0$. When this identity
is needed we will force it to be true also in case $n=0$, see
below (\ref{eq:Zfunction-formula1}). The two other symmetries
are trivially satisfied when $n=0$.

Assume $n>0$. The Rademacher Phi function $\Phi$ is defined on 
$P\SL(2,\Z)=\SL(2,\Z)/\{\pm 1\}$ by
\begin{equation}\label{eq:Rademacher}
\Phi \left[ \begin{array}{ll}
a & b \\
c & d
\end{array}
\right] = \left\{ \begin{array}{ll}
\frac{a+d}{c} - \dS \left(\frac{d}{c}\right) 
& ,c \neq 0, \\
\frac{b}{d} & ,c=0,
\end{array}
\right.
\end{equation}
where the Dedekind symbol $\dS (d/c)$ is given by 
(\ref{eq:Dedekind-symbol}). By the proof of 
\cite[Theorem 8.4]{Hansen2} the function $h$ is given by
\begin{equation}\label{eq:hfunction-compact}
h(\gamma)=\kappa \prod_{i=1}^{n} (\tilde{M}_{i})_{\gamma,1}
\end{equation}
with
$$
\kappa = \sqrt{\mA}(2r)^{n/2}i^{-n}
\exp\left( \frac{i\pi}{4} \sum_{j=1}^{n} \Phi(M_{j}) \right) 
\exp \left( -\frac{i\pi}{2r} \rhoalphasum \right),
$$
where 
$M_{i}=\left( \begin{array}{ll}
-\beta_{i} & -\sigma_{i} \\
\alpha_{i} & \rho_{i}
\end{array} \right)$, 
$i=1,\ldots,n$. Here $A \mapsto \tilde{A}$ is a unitary 
$(r-1)$--dimensional representation of $P\SL(2,\Z)$. By 
unitarity we have
$\tilde{A}_{j,k} = \overline{(\widetilde{A^{-1}})_{k,j}}$ 
for $j,k \in \{1,2,\ldots,r-1\}$, where $\bar{\cdot}$ is 
complex conjugation. By using this with $A=M_{i}$ together with
the formula for $(\widetilde{M_{i}^{-1}})_{1,\gamma}$ given in 
\cite[Propositions 2.7 (a) and 2.8]{Jeffrey1} we get
\begin{eqnarray}\label{eq:Mi-explicit}
(\tilde{M}_{i})_{\gamma,1} &=& i \frac{1}{\sqrt{2r\alpha_{i}}} 
e^{-\frac{i \pi}{4} \Phi(M_{i})} \\
&& \times \sum_{\mu = \pm 1} 
\sum_{m=0}^{\alpha_{i}-1} \mu
\exp \left( \frac{i \pi}{2r\alpha_{i}} 
[-\beta_{i} \gamma^{2} - 2 \gamma (2rm + \mu) + 
\rho_{i} (2rm + \mu)^{2} ] \right), \nonumber
\end{eqnarray}
and inserting this expression for $(\tilde{M}_{i})_{\gamma,1}$
into (\ref{eq:hfunction-compact}) leads to 
(\ref{eq:hfunction}).

Now, to prove (\ref{eq:hfunction-symmetry}), we only need to 
establish the same symmetries for the functions 
$\gamma \mapsto(\tilde{M}_{i})_{\gamma,1}$ (extended to $\Z$).
To do this we make a little simple observation. The group 
$\SL(2,\Z)$ is generated by the two matrices
\begin{equation}\label{eq:generators-SL(2,Z)}
\Xi = \left( \begin{array}{cc}
0 & -1 \\
1 & 0
\end{array}
\right),\hspace{.2in}
\Theta = \left( \begin{array}{cc}
1 & 1 \\
0 & 1
\end{array}
\right) 
\end{equation}
Any $A \in \SL(2,\Z)$ can be written $A=B\Xi$, where 
$B=A\Xi^{-1}=-A\Xi$. In particular
\begin{equation}\label{eq:Matrix-splitting}
\tilde{A}_{j,k} = \sum_{l=1}^{r-1} \tilde{B}_{j,l}
\tilde{\Xi}_{l,k},
\end{equation}
where
$$
\tilde{\Xi}_{l,k}=\sqrt{\frac{2}{r}}
\sin\left(\frac{\pi lk}{r}\right), 
\hspace{.2in}l,k=1,\ldots,r-1. 
$$
This expression is valid for all $k \in \Z$, so we can use 
(\ref{eq:Matrix-splitting}) to extend the function
$k \mapsto \tilde{A}_{j,k}$ from $\{1,\ldots,r-1\}$ to $\Z$.
We can also use the explicit expression for 
$k \mapsto \tilde{A}_{j,k}$ given by 
\cite[Propositions 2.7 (a) and 2.8]{Jeffrey1} to extend this
function to $\Z$, but it follows by the proof of these
propositions in \cite{Jeffrey1} that these two extensions 
coincide. By (\ref{eq:Matrix-splitting}) we immediately get 
that
\begin{eqnarray}\label{eq:Matrix-symmetry}
\tilde{A}_{j,rk} &=& 0, \\
\tilde{A}_{j ,-k} &=& -\tilde{A}_{j,k}, \nonumber \\
\tilde{A}_{j,k+2r} &=& \tilde{A}_{j,k} \nonumber
\end{eqnarray}
for all $k \in \Z$ and $j \in \{1,\ldots,r-1\}$, giving the 
needed symmetries for the entries $(\tilde{M}_{i})_{\gamma,1}$.
In the above kind of argument we do of course not need to use 
unitarity of the representation $A \mapsto \tilde{A}$. The 
reason for using this unitarity is that it leads to an 
expression for $h(\gamma)$ which is slightly better to work 
with in the proof of \refthm{thm:AEC-Seifert-manifolds}.
  
Assume still that $n>0$. By introducing a parameter $\xi >0$ we
get by (\ref{eq:Zfunction-Seifert}) and 
(\ref{eq:hfunction-symmetry}) that
\begin{equation}\label{eq:Zfunction-formula1}
Z(X;r)=\frac{1}{2} \left( \frac{i}{2} \right)^{n} 
\lixi \sum_{\gamma=-r+1}^{r} P(\gamma)H(\gamma,\xi),
\end{equation}
where
$$
H(\gamma,\xi) = \frac{e^{\pi ia_{\ep} g\gamma} h(\gamma)}
{\sin^{n+a_{\ep}g-2}\left(\frac{\pi}{r} (\gamma-i\xi) \right)}
$$
and where $P(\gamma)=(-1)^{a_{\ep}g}$ for 
$\gamma=-r+1,\ldots,-1$, $P(\gamma)=1$ for $\gamma=1,\ldots,r$.
Moreover we let $P(x)=(P(x-)+P(x+))/2$, $x=0,r$, where
$P(x \pm)=\lim_{t \rightarrow x_{\pm}}P(t)$. 
If we multiply $H(\gamma,\xi)$ by
$$
\frac{\sin\left(\frac{\pi}{r}\gamma\right)}
{\sin\left(\frac{\pi}{r}(\gamma-i\xi)\right)}
$$
in case $n=0$, (\ref{eq:Zfunction-formula1}) also holds in that
case. 
For later use we extend $P$ not only to $\Z$ but to all of $\R$
by letting $P(\gamma)=(-1)^{a_{\ep}g}$ for $-r < \gamma < 0$ 
and $P(\gamma)=1$ for $0 < \gamma < r$ and letting $P$ be 
periodic with a period of $2r$. Note that for $a_{\ep}g$ even 
(that is for $\ep=\os$ or $g$ even) $P$ is identically one. The
parameter $\xi$ can of course be ignored in cases where
$n+a_{\ep}g-2 \leq 0$. To continue we use the following formula
which is a special case of \reflem{lem:periodicity}. If 
$f:\Z \ra \C$ is a function with a period of $N$ then
\begin{equation}\label{eq:periodicity-one-dimensional}
\sum_{k=0}^{N-1} f(k) = N \livep \svep
\sum_{k \in \Z} e^{-\pi \vep k^{2}} f(k).
\end{equation}
By applying this to the sum (\ref{eq:Zfunction-formula1}) we 
obtain
$$
Z(X;r) = r \left( \frac{i}{2} \right)^{n}
\lixi \livep \svep \sum_{\gamma \in \Z} 
e^{-\pi \vep \gamma^{2}} P(\gamma) H(\gamma,\xi).
$$
We can exchange the order of the two limits. To see this, 
simply use (\ref{eq:periodicity-one-dimensional}) with 
$f(k)=P(k) \lixi H(k,\xi)$ to get
$$
Z(X;r)=r \left( \frac{i}{2} \right)^{n}
\livep \svep \sum_{\gamma \in \Z} 
e^{-\pi \vep \gamma^{2}} P(\gamma)\lixi H(\gamma,\xi).
$$
By uniform convergence of the sum 
$\sum_{\gamma \in \Z} e^{-\pi \vep \gamma^{2}} P(\gamma)
H(\gamma,\xi)$
with respect to $\xi$ on an interval of the form $]0,\xi_{0}]$,
$\xi_{0}>0$, we then have that
\begin{equation}\label{eq:Zfunction-formula2}
Z(X;r) = r \left( \frac{i}{2} \right)^{n}
\livep \svep \lixi \sum_{\gamma \in \Z} 
e^{-\pi \vep \gamma^{2}} P(\gamma) H(\gamma,\xi).
\end{equation}
The expression (\ref{eq:hfunction}) allows us to consider $h$
as an entire function. In particular we can consider
$\gamma \mapsto H(\gamma,\xi)$ as a meromorphic function on 
$\C$ with a countable number of isolated poles. To proceed we 
change (\ref{eq:Zfunction-formula2}) to a sum of integrals by 
using the Poisson summation formula
\begin{equation}\label{eq:Poisson}
\sum_{k \in \Z} \varphi(k) = \sum_{m \in \Z}
\iny e^{2\pi i mx} \varphi(x) \dte x.
\end{equation}
By a result of Zygmund \cite[p.~68]{Zygmund} this formula is 
valid for every function $\varphi$ which is absolutely 
integrable over $\R$, of bounded variation, and satisfies
$2\varphi(x)=\varphi(x+)+\varphi(x-)$ for all $x \in \R$. 
We use Poisson's formula with 
$\varphi(\gamma) = e^{-\pi \vep \gamma^{2}}P(\gamma)
H(\gamma,\xi)$
(with fixed $\xi$). By (\ref{eq:hfunction}) and the above we 
get
\begin{equation}\label{eq:Zfunction-formula3}
Z(X;r) = r \left( \frac{i}{2} \right)^{n}
\lidel \sdel \lieta \sum_{m \in \Z} \musum 
\nsum g(\umu,\un) \iny \kappa(z) e^{r Q_{m}(z)} \dte z,
\end{equation}
where
\begin{eqnarray*}
g(\umu,\un) &=& \muprod \gex, \\
\kappa(z) &=& \kappa_{\eta}(z;\umu,\delta) = 
\e^{-\pi \delta z^{2}} P(rz)
\frac{\exp \left[ - i \pi \pmualphasum z \right]} 
{\sin^{n+a_{\ep}g-2} \left( \pi (z -i \eta) \right) }, 
\nonumber \\
Q_{m}(z) &=& Q_{m}(z;\un) = 
2\pi i\left(\frac{a_{\ep}g}{2}+m-\nalphasuma \right) z 
+ \frac{i \pi}{2} E z^{2}. \nonumber
\end{eqnarray*}
(Put $z=\gamma/r$, $\delta=\vep r^{2}$ and $\eta=\xi/r$.) 
In case $n=0$ the function 
$\kappa(z)$ contains an extra factor 
$\sin(\pi z)/\sin(\pi(z-i\eta)$. But
$$
\sin(\pi z) = \frac{i}{2} \sum_{\mu \in \{\pm 1\}}
\mu \exp(-i\pi \mu z),
$$
so formally we can treat the case $n=0$ as the case $n=1$
with $(\alpha_1,\beta_1)=(1,0)$ 

We note that it was necessary in the process of changing the 
expression (\ref{eq:Zfunction-Seifert}) to 
(\ref{eq:Zfunction-formula3}) to involve the sum 
$\musum \nsum$. This was needed for establishing (some of) the
properties (\ref{eq:hfunction-symmetry}). Note that 
$z \mapsto P(rz)$ is independent of $r$.

\begin{rem}\label{rem:P1}
A remark is appropriate here. It is well-known that the Poisson
formula (\ref{eq:Poisson}) is valid for functions in the 
Schwartz space $\mS(\R)$ of smooth functions, that together 
with their derivatives are rapidly decreasing at infinity, see 
e.g.\ \cite{Hormander}. If $a_{\ep}g$ is even then $P$ is 
identically $1$ and 
$\varphi(\gamma) = e^{-\pi \vep \gamma^{2}}H(\gamma,\xi)$ is in
$\mS(\R)$ since $\gamma \mapsto H(\gamma,\xi)$ is periodic by 
(\ref{eq:hfunction-symmetry}). If $a_{\ep}g$ is odd it is
not clear how many times 
$\gamma \mapsto P(\gamma)H(\gamma,\xi)$ is differentiable in a
point $\gamma \in r\Z$. We can, however, choose $P$ to be 
smooth also for $a_{\ep}g$ odd. In this case we can simply let
$P$ be a smooth function equal to $1$ on $[3/4,r-3/4]$, equal 
to $0$ on $[0, 1/4]$ and on $[r-1/4,r]$, and equal to 
$-P(-\gamma)$ on $[-r,0]$. Finally we extend $P$ to all of 
$\R$ by letting it have the period $2r$. The problem with this
approach is that $z \mapsto P(rz)$ so defined depends on $r$.

For the method used in the following we will need to extend $P$
to an entire (or at least meromorphic) function with certain 
properties, see above \refrem{rem:P2}. It is not clear at the 
moment for the author if such a $P$ exists, and therefore we
restrict to the case $a_{\ep}g$ even, i.e.\ $\ep=\os$ or $g$ 
even. The calculations presented in this paper were first done
for the case $\ep=\os$ where $a_{\ep}g$ is always even, and 
first later the author started to consider the case $\ep=\ns$.
It was then realized that the calculations could only be 
carried through as they stand for half of the Seifert manifolds
with nonorientable base, namely the ones with a base with even
genus. We note that one can probably calculate the asymptotics
of $Z(X;r)$ by other means than applied in this paper, see 
e.g.\ \cite{Rozansky1} and \cite[Appendix]{Rozansky3} for some
suggestions. It should be possible to use one of these methods
to handle the case of $a_{\ep}g$ odd, but we will defer that to
another paper.
\end{rem}

We approximate the integrals 
$\iny \kappa(z) e^{r Q_{m}(z)} \dte z$ in 
(\ref{eq:Zfunction-formula3}) by the steepest descent method. 
This is a general method for obtaining asymptotic expansions of
contour integrals of the form
$$
I(t)=\int_{C} g(z) e^{t f(z)} \dte z
$$
in the limit $t \ra \infty$, where $f(z)$ and $g(z)$ are 
analytic functions on some domain $\Omega$ in the complex 
plane, both independent of the real positive parameter $t$, and
$C$ is a contour inside $\Omega$. Actually we allow $g(z)$ to 
have poles in $\Omega$ away from the contours involved. The 
basic idea is to deform $C$ inside $\Omega$ into a new path of 
integration $C'$ so that the following conditions hold:
\newline  
\begin{enumerate}
\item[i)] The path $C'$ passes through one or more of the 
stationary points of $f$ (that is the zeroes of $f'(z)$);
\item[ii)] The imaginary part of $f(z)$ is constant on $C'$.
\end{enumerate}
\bigskip

\noindent We will only need to consider the case of 
non-degenerate stationary points of $f$, i.e.\ stationary 
points $p$ such that $f''(p) \neq 0$. In that case one finds
that there is a unique contour $C'$ through $p$ satisfying
ii) above and satisfying that $\rea(f(p)-f(z))<0$ for all $z$
on $C' \sm \{p\}$ near $p$. This contour is called the steepest
descent contour through $p$. It's argument is 
$(\pi-\arg(f''(p)))/2$, see e.g.\ 
\cite[Theorem 7.1]{BleisteinHandelsman} or
\cite[Sec.~5.4]{deBruijn}. To determine the asymptotic 
expansion of $I(t)$ in the limit $t \ra \infty$ we have to 
calculate the large $t$ asymptotics of the new contour integral
$$
J(t)=\int_{C'} g(z) e^{t f(z)} \dte z.
$$
The difference between $I(t)$ and $J(t)$ can be zero or contain
a sum of residue contributions, coming from the poles of $g(z)$
crossed in the process of deforming the contour from $C$ to 
$C'$. Moreover, this difference can contain contributions given
by integrals along contours connecting the contours $C$ and 
$C'$, but these contributions should be $o(t)$ compared to 
$J(t)$ in the limit $t \ra \infty$. The integral $J(t)$ will 
for quite general $f(z)$ and $g(z)$ be of a form for which the 
Laplace method can be used to calculate the large $t$ 
asymptotics. In our case these integrals will however be so 
simple that a general description of the Laplace method is not
relevant. Moreover, we will not go into details about why the 
steepest descent contour $C'$ is a particularly `good' contour
when it comes to finding the large $t$ asymptotics of $I(t)$. 
From a technical point of view a main reason is, as already 
pointed out, that the Laplace method is applicable to the new 
integral $J(t)$ along $C'$. We refer to \cite[pp.~84-90]{Wong},
\cite[pp.~262-268]{BleisteinHandelsman} and
\cite[Chap.~5]{deBruijn} for futher details.

In case $f$ has more than one stationary point one has normally
to determine the contribution to the large $t$ asymptotics of 
$I(t)$ coming from each stationary point by using a `short' 
steepest descent contour through each such point. One should of
cource also secure that it is possible to deform the original 
contour $C$ into the steepest descent contour within the domain
$\Omega$. 

In our situation $t=r$, $f(z)=Q_{m}(z)$, $g(z)=\kappa(z)$, 
$\Omega=\C$, and the contour $C$ is the real axes. If the 
Seifert Euler number $E=0$ the functions $Q_{m}$ do not have
stationary points so in this case the steepest descent method 
is not applicable. The case $E=0$ will be taken care of in 
Sec.~\ref{sec-Proof-of-E0} as already pointed out earlier. If 
$E \neq 0$ the phase function $Q_{m}$ has exactly one 
stationary point
$$
z_{\st} = -\frac{2}{E} 
\left(\frac{a_{\ep}g}{2}+m-\nalphasuma \right)
$$
which is non-degenerate.

\subsection{Calculation of asymptotics. A general case}
\label{sec-Calculation-of}

To make the calculations useable for other situations and for 
shortening notation we will generalize the situation above 
without making things more difficult. To be more concrete we 
will study functions of the form
\begin{equation}\label{eq:Zfunction-typical}
Z(r) = \sum_{\gamma=1}^{r-1} 
\frac{h(\gamma)}{\sin^{k}\left(\frac{\pi}{r}\gamma\right)},
\end{equation}
where $k\in \Z$ and
\begin{equation}\label{eq:hfunction-typical}
h(\gamma) = \sum_{\lambda \in \Lambda} g(\lambda) p(\gamma)
\exp\left(i\left[\frac{A}{r} \gamma^{2}+B\gamma \right]\right).
\end{equation}
We assume that $A$ and $B$ are real and independent of $r$.
Moreover, we assume that $g$ is independent of $\gamma$ and 
that $p$ is an entire function. The index set $\Lambda$ is 
assumed to be a finite subset of $\R^{e}$ for some 
$e \in \Z_{>0}$. We also assume that $h$ satisfies 
(\ref{eq:hfunction-symmetry}) with $n=k$ for all 
$\gamma \in \Z$. We will allow $g$ to depend on $r$ and $A$, 
$B$, and $p$ to depend on $\lambda$. By 
(\ref{eq:periodicity-one-dimensional}) we get
$$
Z(r) = r \livep \svep \, \lixi \sum_{\gamma \in \Z} 
e^{-\pi \vep \gamma^{2}} \frac{h(\gamma)}
{\sin^{k}\left(\frac{\pi}{r}(\gamma-i\xi)\right)}.
$$
(We can exchange the order of the limits $\livep$ and $\lixi$ 
in (\ref{eq:Zfunction-typical-formula1}). This follows exactly
as in the special case where $Z(r)=Z(X;r)$). The function
$x \mapsto e^{-\pi \vep x^{2}} 
\frac{h(x)}{\sin^{k}\left(\frac{\pi}{r}(x-i\xi)\right)}$
belongs to $\mS(\R)$ since
$x \mapsto 
\frac{h(x)}{\sin^{k}\left(\frac{\pi}{r}(x-i\xi)\right)}$
is periodic. Therefore
\begin{equation}\label{eq:Zfunction-typical-formula1}
Z(r) = r \lidel \sdel \, \lieta \sum_{m \in \Z} 
\sum_{\lambda \in \Lambda} g(\lambda) 
\iny f(z;m,\lambda) \dte z
\end{equation}
by the Poisson formula (\ref{eq:Poisson}), where
\begin{equation}\label{eq:ffunction}
f(z;m,\lambda) = e^{-\pi\delta z^{2}} q(z) 
\frac{e^{rQ_{m}(z)}}{\sin^{k}(\pi(z-i\eta))},
\end{equation}
with $q(z)=q(z;\lambda)=p(rz)$ and $Q_{m}(z)=Q(z)+2\pi i mz$,
where $Q(z)=Q(z;\lambda)=i(Az^{2}+Bz)$. In case $k \leq 0$ we 
can and will replace $1/\sin^{k}(\pi(z-i\eta))$ in 
$f(z;m,\lambda)$ by $\sin^{|k|}(\pi z)$ and remove $\lieta$. 
(Note that we do not loose any generality by not allowing the 
polynomial $Q$ to have a constant term. Such a constant term 
$c$ would just produce an extra factor $e^{irc}$ in $g$.) From
the first identity in (\ref{eq:hfunction-symmetry}) we have
\begin{equation}\label{eq:zero-1}
\sum_{\lambda \in \Lambda} g(\lambda) q(z) e^{rQ(z)}=0
\end{equation}
for all $z \in \Z$.

In case $Z(r)=Z(X;r)$ (with $a_{\ep}g$ even) we have by
(\ref{eq:Zfunction-formula3}) that
{\allowdisplaybreaks
\begin{eqnarray}\label{eq:Seifert-case}
\Lambda &=& \{\pm 1\}^{n} \times \mS, \\ 
A &=& \frac{\pi}{2} E = -\frac{\pi}{2}\left(\beta_0 +
\sum_{j=1}^{n} \frac{\beta_j}{\alpha_j}\right), \nonumber \\ 
B &=& B(\un)= -2\pi \sum_{j=1}^{n}\frac{n_{j}}{\alpha_{j}},
\nonumber \\
k &=& n+a_{\ep}g-2, \nonumber \\
g(\umu,\un) &=& \left(\frac{i}{2}\right)^{n}
\muprod \gex, \nonumber \\
q(z) &=& q(z;\umu) =
\exp\left( -i\pi\left(\mualphasum\right)z\right) \nonumber
\end{eqnarray}}\noindent 
for $(\umu,\un) \in \Lambda$, where $\mS$ is given by 
(\ref{eq:S-set}).

Let us return to the general case. We will need to make a 
series of assumptions on $q$. Firstly, we will assume that $q$
is independent of $r$ (or equivalently, that $z \mapsto p(rz)$ 
is independent of $r$). This assumption is made to make the 
steepest descent method applicable to the asymptotic analysis 
of integrals $\iny f(z;m,\lambda)\dte z$. Secondly, we need to
make some assumptions of a technical nature. The need for these
will only be transparent within the proofs of results to 
follow. In general, they are made to secure convergency of 
certain infinite sums. First of all we assume that there exist 
real positive numbers $a_{n},b_{n}$ such that
\begin{equation}\label{eq:qfunction-symmetry1}
\left|\frac{\dte ^{n} q(z;\lambda)}{{\dte z}^{n}} \right| 
\leq b_{n}e^{a_{n}|z|}
\end{equation}
for all $z \in \C$, all $\lambda \in \Lambda$, and 
$n=0,1,\ldots$. Moreover, we assume that $q$ is periodic along 
each line parallel to the real line. In fact we assume that 
there exists an integer $m_{q}>0$ such that
\begin{equation}\label{eq:qfunction-symmetry2}
q(z+m_{q};\lambda)=q(z;\lambda)
\end{equation}
for all $z \in \C$ and all $\lambda \in \Lambda$.
It will follow that we can relax the assumption
(\ref{eq:qfunction-symmetry1}) slightly by replacing 
$e^{a_{n}|z|}$ by $e^{a_{n}|z|^{\tau}}$ with $\tau \in [1,2[$,
but since the above sufficis for our purpose we will continue
by using (\ref{eq:qfunction-symmetry1}). In the case 
$Z(r)=Z(X;r)$ the above conditions are satisfied, see
(\ref{eq:Seifert-case}).

\begin{rem}\label{rem:P2}
This remark is a continuation of \refrem{rem:P1}. In the case 
where $Z(r)=Z(X;r)$ and $a_{\ep}g$ is odd $q$ contains an
extra factor $P(rz)$, where $P$ is the function in 
(\ref{eq:Zfunction-formula1}). To include this case in the 
following calculations one will need to extend $P$ to an entire
function in such a way that the above assumptions on $q$ are 
not violated. It is not clear to the author if this can be 
done.
\end{rem}

By the formulae $\sin(z)=(e^{iz}-e^{-iz})/(2i)$ and
$\cos(z)=(e^{iz}+e^{-iz})/2$ we see that both $| \sin(z)|$ and
$|\cos(z)|$ are less than or equal to $e^{|z|}$ for all 
$z \in \C$. For $k \leq 0$ we can therefore let 
$1/\sin^{k}(\pi z)$ be included in $q(z)$ without violating the
above assumptions on $q$.

In the following we will suppress $\lambda$ in many of the 
expressions. One should keep in mind that all functions 
depending on $A$, $B$ and/or $q$ may depend on $\lambda$. The 
calculation of the large $r$ asymptotics of $Z(r)$ is divided 
naturally into several parts. To make it possible to follow the
rather long computation we have tried to stress this by 
dividing the calculation into several subsections.

\subsubsection{\bf A partition of $Z(r)$ using the steepest
descent method}

In the remaining part of Sec.~\ref{sec-Proof-of} we assume that
$A \neq 0$. The case $A=0$ will be handled in 
Sec.~\ref{sec-Proof-of-E0}. We use the method of steepest 
descent to calculate the large $r$ asymptotics of the integrals
$\iny f(z;m) dz$ in (\ref{eq:Zfunction-typical-formula1}). The
polynomial $Q_{m}(z;\lambda)$ has only one critical point
\begin{equation}\label{eq:stationarypoint}
z_{\st}(m,\lambda) = -\frac{2\pi m+B}{2A}
\end{equation}
which is non-degenerate. The steepest descent contour, denoted
$C_{\sd}(m,\lambda)$, is the contour
\begin{equation}\label{eq:sdcontour}
\ima(z)=\sgA ( \rea(z)-z_{\st}(m,\lambda)).
\end{equation}
We orient $C_{\sd}$ by starting at $-\infty\eksgA +z_{\st}$. In
the process of deforming the real axes into $C_{\sd}$ the 
integration contour crosses (for $k>0$) those poles
\begin{equation}\label{eq:ffunction-poles}
z_{l}(\eta)=l+i\eta,\hspace{.2in} l \in \Z,
\end{equation}
for which
\begin{equation}\label{eq:condition-poles}
\sgA (l-z_{\st}) > \eta.
\end{equation}

\begin{rem}
It is possible in parts of the calculation of the large $r$ 
asymptotics of $Z(r)$ in (\ref{eq:Zfunction-typical-formula1})
to let $Q$ be of a more general form, in particularly if no 
singularities are present, i.e.\ if $k \leq 0$. One can e.g.\ 
assume that $Q$ is an arbitrary analytic function satisfying 
that $Q_{m}(z)=2\pi imz + Q(z)$ has exactly one stationary 
point $z_{m}$ for each $m$. Moreover, we assume that $z_{m}$ is
nondegenerate. For the stationary point $z_{m}$ we deform the 
integration contour in $\iny f(z;m) \dte z$ from the real axes
to the line $\gamma_{m}(t)=z_{m} + \alpha_{m}t$, $t \in \R$, 
where
$$
\alpha_{m}=\exp\left( \frac{i\pi}{2}
-\frac{i}{2} \Arg (Q_{m}''(z_{m})) \right).
$$
We then have
$$
Q_{m}(\gamma_{m}(t)) = Q_{m}(z_{m})
+ \frac{1}{2} Q_{m}''(z_{m}) \alpha_{m}^{2}t^{2}+\ldots,
$$
where $Q_{m}''(z_{m}) \alpha_{m}^{2}=-|Q_{m}''(z_{m})|<0$.
In case of polar contributions, i.e.\ in case $k>0$ we will
also assume that $\gamma_{m}$ is not parallel with the real 
axes, i.e.\ we assume that $Q_{m}''(z_{m})$ is not real and 
negative. This assumption is imposed to avoid problems with the
sum of residue contributions (see \reflem{lem:Lemma8.12} 
below). The integrals $\int_{\gamma_{m}} f(z;m) \dte z$ can now
be handled by using the Laplace method as illustrated in 
\cite[Sec.~4.4]{deBruijn}. There are several reasons why we 
will not do the calculations in this generality. First of all 
it will put us too far away from the example of main interest 
in this paper, the case $Z(r)=Z(X;r)$. A more serious 
reason is however that to be able to carry out the analysis we
will have to impose too many assumptions. Some of these 
assumptions will then in concrete examples be hard to check 
(or at least these checks will demand long computations), so 
one does not gain very much. We have found that the expression
(\ref{eq:ffunction}) is a good compromise between generality
and ability to carry through the analysis without too many
assumptions.
\end{rem}

We will need to impose a periodicity assumption on the function
$z_{\st} : \Z \times \Lambda \to \R$. To be specific we assume
that there exists an integer $H>0$ such that
\begin{equation}\label{eq:zst-periodicity}
z_{\st}(m+H,\lambda) - z_{\st}(m,\lambda) = 
-\frac{\pi}{A} H \in \Z
\end{equation}
for all $(m,\lambda) \in \Z \times \Lambda$. This implies that
the set
$\{ \; z_{\st}(m,\lambda) \pmod{\Z} \; | \;
(m,\lambda) \in \Z \times \Lambda \;\}$
is finite. Note that (\ref{eq:zst-periodicity}) is satisfied if
and only if
\begin{equation}\label{eq:A}
A=A(\lambda) = \pi \frac{H}{P}
\end{equation}
with $P=P(\lambda) \in \Z \sm \{0\}$ for all 
$\lambda \in \Lambda$. In the following we assume that $A$ is 
given by (\ref{eq:A}) and reserve the symbols $H$ and $P$ for 
this. By (\ref{eq:zst-periodicity}) there exists an 
$\eta_{0} \in ]0,1/2[$ such that
\begin{equation}\label{eq:zstassumption}
z_{\st}(m,\lambda) \not\in \Z + \left([-\eta_{0},\eta_{0}] 
\sm \{0 \}\right)
\end{equation}
for all $(m,\lambda) \in \Z \times \Lambda$. For many of the 
next results we only have to assume (\ref{eq:zstassumption}) 
but later on we will have to assume the stronger 
(\ref{eq:zst-periodicity}). In parts of the calculation we will
also need to assume that
\begin{equation}\label{eq:B}
B =B(\lambda) \in \pi\Q
\end{equation}
for all $\lambda \in \Lambda$. When $Z(r)=Z(X;r)$, 
(\ref{eq:zst-periodicity}) and (\ref{eq:B}) are satisfied by 
(\ref{eq:Seifert-case}).

Let us introduce some notation used in the remaining part of 
this paper. The steepest descent contour $C_{\sd}(m,\lambda)$, 
$(m,\lambda) \in \Z \times \Lambda$, is parametrized by
\begin{equation}\label{eq:sdparametrization}
\gamma_{(m,\lambda)}(t) = t\eksgA + z_{\st}(m,\lambda).
\end{equation}
We let
\begin{equation}\label{eq:kappafunction}
\kappa_{\eta}(z;\lambda,\delta) = e^{-\pi\delta z^{2}} 
\frac{q(z;\lambda)}{\sin^{k}(\pi(z-i\eta))},
\end{equation}
where $\delta$ and $\eta$ are non-negative parameters. We note
that
\begin{equation}\label{eq:sdintegral}
e^{rQ_{m}(\gamma_{(m,\lambda)}(t))} = 
e^{-irA z_{\st}^{2}(m,\lambda)}e^{-r|A| t^{2}}.
\end{equation}
Let
\begin{equation}\label{eq:W-set}
W= \{ \; (m,\lambda) \in \Z \times \Lambda \; | \; 
z_{\st}(m,\lambda) \notin \Z \; \},
\end{equation}
and put
\begin{equation}\label{eq:A0}
A_{0} = \min\{ \; |A(\lambda)| \; | \; \lambda \in \Lambda 
\; \}.
\end{equation}
The open disk in the complex plane with centre $z$ and radius 
$\rho >0$ is denoted $D(z,\rho)$, and the punctured disk 
$D(z,\rho) \sm \{z\}$ is denoted $D'(z,\rho)$. 

The next lemma, \reflem{lem:Lemma8.7}, is necessary because we
deal with infinite contours so we have to check that the 
integrals along the steepest descent contours are convergent. 
Since the real axes and the steepest descent contour are not 
attached to each other in the `end points' we also have to
check that there are no contributions from integrals along
contours connecting the real axes to $C_{\sd}(m)$ at 
respectively minus and plus infinity.

\begin{lem}\label{lem:Lemma8.7}
Let $\eta_{0}>0$ be as in {\em (\ref{eq:zstassumption})} and 
let $(m,\lambda) \in \Z \times \Lambda$ be arbitrary but fixed.
Then
\begin{equation}\label{eq:Lemma8.7}
\iny f(z;m,\lambda) \dte z = 
\int_{C_{\sd}(m,\lambda)} f(z;m,\lambda) \dte z + 2 \pi i
\sum_{\stackrel{l \in \Z}{\sgA (l-z_{\st}(m,\lambda)) > \eta}} 
\Res _{z=z_{l}(\eta)} \left\{ f(z;m,\lambda) \right\}
\end{equation}
for all $\eta \in ]0,\eta_{0}]$ and all 
$\delta \in ]0, \frac{r|A|}{\pi} ]$. The sum of residues is 
zero if $k \leq 0$.
\end{lem}

The proof is given in Appendix B. By 
(\ref{eq:Zfunction-typical-formula1}) and 
\reflem{lem:Lemma8.7} we conclude that
\begin{eqnarray}\label{eq:Zfunction-splitting0}
Z(r) &=& r \lidel \sdel \, \lieta \sum_{m \in \Z} 
\sum_{\lambda \in \Lambda} g(\lambda) \\
&& \times \left( 
\int_{C_{\sd}(m,\lambda)} f(z;m,\lambda) \dte z + 2 \pi i
\sum_{\stackrel{l \in \Z}{\sgA (l-z_{\st}(m,\lambda)) > \eta}}
\Res _{z=z_{l}(\eta)} \left\{ f(z;m,\lambda) \right\} 
\right).\nonumber
\end{eqnarray}
Next we are going to separate the contributions coming from 
respectively poles and the integrals along the steepest descent
contours.

\begin{lem}\label{lem:Lemma8.12}
Let $\lambda \in \Lambda$ be arbitary but fixed. Then the 
infinite series
$$
\Sigma_{1}(\lambda) = \sum_{m \in \Z} 
\int_{C_{\sd}(m,\lambda)} f(z;m,\lambda) \dte z
$$
is absolutely convergent for all $\eta \in ]0,\eta_{0}]$ and 
all $\delta \in ]0,2r|A|/\pi[$, and the infinite series
$$
\Sigma_{2}(\lambda) = \sum_{m \in \Z} 
\sum_{\stackrel{l \in \Z}{\sgA (l-z_{\st}(m,\lambda)) > \eta}} 
\Res _{z=z_{l}(\eta)} \left\{ f(z;m,\lambda) \right\}
$$
is absolutely convergent for all $\eta >0$ and all $\delta >0$.
For all $\eta,\delta \in ]0,\infty[$ we have
$$
\Sigma_{2}(\lambda) = \sum_{l \in Z} \beta(l,\lambda) 
\sum_{\stackrel{m \in \Z}
{m +\frac{B}{2\pi}+\frac{A}{\pi}l > \frac{|A|}{\pi}\eta}}
\Res_{z=i\eta} \left\{ \phi(z;l,m,\lambda) \right\},
$$
where
\begin{eqnarray}\label{eq:Lemma8.12}
\beta(l,\lambda) &=& (-1)^{kl} e^{ir(Al^{2}+Bl)}, \\
\phi(z;l,m,\lambda) &=& e^{-\pi \delta (z+l)^{2}} 
q(z+l;\lambda) e^{irAz^{2}}
\frac{e^{2\pi i r \left( m + \frac{B}{2\pi} +
\frac{A}{\pi}l \right)z}}{\sin^{k}(\pi(z-i\eta))}. \nonumber
\end{eqnarray}
Note that $\Sigma_{2}(\lambda)=0$ for $k \leq 0$.
\end{lem}

For the proof, see Appendix B. Note that the sum over $m$ in 
the expression (\ref{eq:Zfunction-typical-formula1}) is 
absolutely convergent by \reflem{lem:Lemma8.12}. By 
(\ref{eq:zero-1})
\begin{equation}\label{eq:zero-2}
\sum_{\lambda \in \Lambda} g(\lambda) \beta(l,\lambda) q(z+l) 
e^{ir(Az^{2}+Bz+2Alz)} = (-1)^{kl} \sum_{\lambda \in \Lambda} 
g(\lambda) q(z+l) e^{rQ(z+l)} = 0
\end{equation}
for all $z,l \in \Z$. By (\ref{eq:Zfunction-splitting0}) and
\reflem{lem:Lemma8.12} we have
\begin{equation}\label{eq:Zfunction-splitting}
Z(r) = r\lidel \sqrt{\delta} \lieta 
\left( Z_{\inte,1}(\delta,\eta) +  Z_{\intpolar}(\delta,\eta)
+ Z_{\pol,1}(\delta,\eta) \right),
\end{equation}
where
\begin{eqnarray}\label{eq:Zterms}
Z_{\inte,1}(\delta,\eta) &=& \sum_{(m,\lambda) \in W} 
g(\lambda) \int_{C_{\sd}(m,\lambda)} f(z;m,\lambda) \dte z, \\ 
Z_{\intpolar}(\delta,\eta) &=& \sum_{\lambda \in \Lambda} 
g(\lambda) \msumi 
\int_{C_{\sd}(m,\lambda)} f(z;m,\lambda) \dte z, \nonumber \\
Z_{\pol,1}(\delta,\eta) &=& 2\pi i \sum_{\lambda \in \Lambda} 
g(\lambda) \sum_{l \in \Z} \beta(l,\lambda) \msumj
\Res_{z=i\eta} \left\{ \phi(z;l,m,\lambda) \right\}. \nonumber
\end{eqnarray}
Here we suppress the dependency on $r$ and also to some extend 
the dependency on $\lambda$ in our notation for not making it 
too clumsy. We have split the sum $\Sigma_{1}$ in 
\reflem{lem:Lemma8.12} into two parts, since this will be 
convenient for the following calculations. The reason for the 
names $Z_{\intpolar}$ and $Z_{\pol,1}$ instead of $Z_{\inte,2}$
and $Z_{\pol}$ will become clear below.

Let us take a closer look at $Z_{\intpolar}(\delta,\eta)$. 
Assume that $z_{\st}(m,\lambda) = l \in \Z$. By changing 
variable to $w=z-l$, as in the proof of 
\reflem{lem:Lemma8.12}, the contour $C_{\sd}(m,\lambda)$ 
changes to $C(\lambda,0)$, where, for $\rho \in \R$, the 
contour $C(\lambda,\rho)$ is given by
$$
\ima(w) = \sgA \rea(w) - \rho
$$
oriented so that we start at $-\infty \eksgA-i\rho$. It follows
that
$$
Z_{\intpolar}(\delta,\eta) = \sum_{\lambda \in \Lambda} 
g(\lambda) \lsuma \beta(l,\lambda) I(l;\lambda,\delta,\eta),
$$
where 
$I(l;\lambda,\delta,\eta) = 
\int_{C(\lambda,0)} \phi(w;\lambda,m,l) \dte w$,
where
$$
\phi(w;\lambda,m,l) = e^{-\pi \delta (w+l)^{2}} q(w+l) 
\frac{e^{irAw^{2}}}{\sin^{k}(\pi(w-i\eta))}
$$
by (\ref{eq:Lemma8.12}) since 
$m +\frac{B}{2\pi}+\frac{A}{\pi}l=0$. By changing variable once
more to $y=w-i\eta$ we get
$I(l;\lambda,\delta,\eta) = 
\int_{C(\lambda,\eta)} \F(y) \dte y$ 
with
\begin{equation}\label{eq:F}
\F(y) = \F(y;l,\lambda,\delta,\eta) =
e^{-\pi \delta (y+l+i\eta)^{2}} q(y+l+i\eta) 
\frac{e^{irA(y+i\eta)^{2}}}{\sin^{k}(\pi y)}.
\end{equation}
We split the function $\F$ into its even and odd parts, i.e.\ 
$\F(y)=\F^{+}(y)+\F^{-}(y)$, where
$$
\F^{\pm}(y) = \frac{1}{2} (\F(y) \pm \F(-y)).
$$

\begin{lem}\label{lem:Lemma8.13}
Let $\lambda \in \Lambda$ and $l \in \Z$ be arbitrary but 
fixed. Then
$$
\int_{C(\lambda,\eta)} \F^{-}(y;l,\lambda,\delta,\eta) \dte y
= \pi i \Res_{y=0} \F(y;\lambda,l,\delta,\eta)
$$
for any $\eta \in ]0,1[$ and any $\delta \in ]0,2r|A|/\pi[$.
Moreover, the series
$$ 
Z_{\pol,2}(\delta,\eta) := \sum_{\lambda \in \Lambda} 
g(\lambda) \lsuma 
\beta(l,\lambda) \pi i \Res_{y=0} \F(y;l,\lambda,\delta,\eta)
$$
is absolutely convergent for all $\delta,\eta \in ]0,\infty[$.
For $k \leq 0$ this sum is zero.
\end{lem}

The proof is given in Appendix B. By \reflem{lem:Lemma8.12} the
sum $Z_{\intpolar}(\delta,\eta)$ is absolutely convergent for 
$\eta \in ]0,\eta_{0}]$ and $\delta \in ]0,2rA_{0}/\pi[$, and 
by \reflem{lem:Lemma8.13} we then get that
\begin{equation}\label{eq:Zint2-deltaeta}
Z_{\inte,2}(\delta,\eta) := \sum_{\lambda \in \Lambda} 
g(\lambda) \lsuma \beta(l,\lambda) 
\int_{C(\lambda,\eta)} \F^{+}(y;l,\lambda,\delta,\eta) \dte y
\end{equation}
is absolutely convergent and that
$$
Z_{\intpolar}(\delta,\eta) = Z_{\pol,2}(\delta,\eta) +
Z_{\inte,2}(\delta,\eta)
$$
for these $(\delta,\eta)$. We put
$$
Z_{\pol}(\delta,\eta) = Z_{\pol,1}(\delta,\eta) + 
Z_{\pol,2}(\delta,\eta).
$$
To conclude we have so far proven that
\begin{equation}\label{eq:Zfunction-splitting1}
Z(r) = r\lidel \sqrt{\delta} \lieta \left( 
Z_{\inte,1}(\delta,\eta) + Z_{\inte,2}(\delta,\eta) + 
Z_{\pol}(\delta,\eta) \right),
\end{equation}
where $Z_{\inte,1}(\delta,\eta)$ and $Z_{\inte,2}(\delta,\eta)$
are given by respectively (\ref{eq:Zterms}) and 
(\ref{eq:Zint2-deltaeta}), and
\begin{equation}\label{eq:Zpolar-deltaeta}
Z_{\pol}(\delta,\eta) = \sum_{l \in \Z} Z^{l}(\delta,\eta),
\end{equation}
where 
$Z^{l}(\delta,\eta) = Z^{l}_{1}(\delta,\eta) +
Z^{l}_{2}(\delta,\eta)$
with
\begin{eqnarray*}
Z^{l}_{1}(\delta,\eta) &=& 2\pi i \sum_{\lambda \in \Lambda} 
g(\lambda) \beta(l,\lambda) \msumj 
\Res_{z=i\eta} \left\{ \phi(z;l,m,\lambda) \right\}, \\
Z^{l}_{2}(\delta,\eta) &=& \pi i \lambdasuma g(\lambda) 
\beta(l,\lambda)
\Res_{z=i\eta} \left\{ \phi\left(z;\lambda,
-\left(\frac{B}{2\pi}+\frac{A}{\pi}l \right),l \right) 
\right\}
\end{eqnarray*}
(where, as usual, an empty sum is put equal to zero). By 
(\ref{eq:stationarypoint}) we have that
$m+\frac{B}{2\pi}+\frac{A}{\pi}l>\frac{|A|}{\pi}\eta$ is
equivalent to $\sgA (l-z_{\st}(m,\lambda)) >\eta$. By 
(\ref{eq:zstassumption}) this is also equivalent to
$\sgA (l-z_{\st}(m,\lambda)) >0$ for $\eta \in ]0,\eta_{0}]$.
For $\eta \in ]0,\eta_{0}]$ we therefore have
\begin{equation}\label{eq:Zpolar-lterm1}
Z^{l}(\delta,\eta) = 2\pi i \sum_{\lambda \in \Lambda} 
g(\lambda) \beta(l,\lambda) \msumk
\Res_{z=i\eta} \left\{ \frac{\phi(z;l,m,\lambda)}
{\Sym_{\pm}\left(m+\frac{B}{2\pi}+\frac{A}{\pi}l \right)} 
\right\},
\end{equation}
where $\Sym_{\pm}$ is given by (\ref{eq:Sym}).

Let us make some comments on the strategy we will follow from 
here. Idealistically we can calculate the large $r$ asymptotics
of $Z(r)$ in (\ref{eq:Zfunction-splitting1}) by calculating the
large $r$ asymptotics of the limit 
$\lidel \sdel \lieta f(\delta,\eta)$ for $f$ equal to each of 
the functions $Z_{\inte,1}$, $Z_{\inte,2}$ and $Z_{\pol}$ 
(assuming that these limits exist). This is actually what we 
will attempt to do. In fact we start by showing that 
\begin{equation}\label{eq:Zpolar1}
Z_{\pol}(r):=r\lidel \sdel \lieta Z_{\pol}(\delta,\eta)
\end{equation}
exists, and at the same time we calculate this limit. Thus we
not only give an asymptotic description of this limit but we
actually find an exact expression for it. To do this we need to
assume that certain ``symmetry'' conditions are satisfied.

To handle the two other limits
$\lidel \sdel \lieta Z_{\inte,\nu}(\delta,\eta)$, $\nu=1,2$,
more care need be taken. Since the limits on the right-hand 
sides of (\ref{eq:Zfunction-splitting1}) and (\ref{eq:Zpolar1})
exist we have that
\begin{equation}\label{eq:Zinte}
Z_{\inte}(r) := r\lidel \sdel \lieta Z_{\inte}(\delta,\eta)
\end{equation}
exists, where
$Z_{\inte}(\delta,\eta) = Z_{\inte,1}(\delta,\eta) +
Z_{\inte,2}(\delta,\eta)$.
We will, however, not calculate this limit explicitly, but only
give it an asymptotic description. This is of cource sufficient
for the proof of \refthm{thm:asymptotics-Enot0}. Thus we will 
show that we for each $N \in \Z_{\geq 0}$ have a decomposition
$$
Z_{\inte,\nu}(\delta,\eta) = Z_{\inte,\nu}(N;\delta,\eta) +
R_{\nu}(N;\delta,\eta)
$$
such that 
$$
Z_{\inte,\nu}(r;N) = r\lidel \sdel \lieta 
Z_{\inte,\nu}(N;\delta,\eta)
$$
exists, $\nu =1,2$. This implies then that
$$
R(r;N) = r\lidel \sdel \lieta \left( R_{1}(N;\delta,\eta) + 
R_{2}(N;\delta,\eta) \right)
$$
exists. We will not show that
$\lidel \sdel \lieta R_{\nu}(N;\delta,\eta)$ exists for each
$\nu=1,2$ separately. Instead we show that 
$R_{\nu}(N;\delta,\eta)$ is bounded from above by a certain 
function $A_{\nu}(N;\delta,\eta)$ for which
$$
A_{\nu}(r;N) = r\lidel \sdel \lieta A_{\nu}(N;\delta,\eta)
$$
can easily be calculated and shown to be small compared to 
$Z_{\inte,\nu}(r;N)$ in the limit of large $r$, $\nu=1,2$. This
will prove that $Z_{\inte,1}(r;N)+Z_{\inte,2}(r,N)$ is the 
large $r$ asymptotics of $Z_{\inte}(r)$ to an order depending 
on $N$. To carry out this program we have to assume that 
certain functions admit certain periodicities.

One final small remark is appropriate. When dealing with
$Z_{\inte,1}(\delta,\eta)$ we will start by getting rid of the
parameter $\eta$ by showing that
$$
Z_{\inte,1}(\delta) := \lieta Z_{\inte,1}(\delta,\eta)
$$
exists and then make a decomposition
$$
Z_{\inte,1}(\delta) = Z_{\inte,1}(N;\delta) + R_{1}(N;\delta)
$$
for each $N \in \Z_{\geq 0}$. This makes things a little 
easier. When dealing with $Z_{\inte,2}(\delta,\eta)$ this is 
not possible. We will see that the positive parameter $\eta$ 
plays a crucial role in the calculations and can only be 
removed in the final step.

\subsubsection{\bf A calculation of the polar contribution 
$Z_{\pol}(r)$}\label{sec-polar-contribution}

In this section we calculate the limit $Z_{\pol}(r)$ in 
(\ref{eq:Zpolar1}). If $k \leq 0$ this is zero so here we 
assume that $k>0$. By (\ref{eq:Zpolar-deltaeta}) 
$Z_{\pol}(\delta,\eta)$ is given by two infinite sums, the 
$l$-sum and the $m$-sum. We will change the $l$-sum by using 
(\ref{eq:periodicity-one-dimensional}) 'backwards'. The $m$-sum
will be calculated explicitly. Before we can do this we need to
write $Z^{l}(\delta,\eta)$ in (\ref{eq:Zpolar-lterm1}) in a 
different way. A main problem is that the lower bound
in the $m$-sum depends on $l$ in the present expression for
$Z^{l}(\delta,\eta)$.

We begin by a list of assumptions. Recall that the index set 
$\Lambda$ is a subset of $\R^{e}$. We assume that $A,B$, $g$, 
and $q(z)$ (for a fixed but arbitrary $z \in \C$) as functions
of $\lambda$ are defined on all of $\R^{e}$, and that $A$ and 
$B$ are real on all of $\R^{e}$.

Let $V_{l} : \R^{e} \to \R$ and 
$\tU_{l} : \C \times \R^{e} \to \C$ be given by
\begin{eqnarray*}
V_{l}(x) &=& -\frac{B(x)}{2\pi} - \frac{A(x)}{\pi}l, \\
\tU_{l}(z;x) &=& g(x) \beta(l,x) q(z+l;x) e^{i r A(x)z^{2} }.
\end{eqnarray*}
Assume that there for each $l \in \Z$ exists a bijection
$T_{l} : \R^{e} \to \R^{e}$ such that
$$
V(x) := V_{l}\left( T_{l}^{-1}(x) \right)
$$
only depends on $x$ (and not on $l$). We will assume that 
$T_{0}=\id_{\R^{e}}$. By (\ref{eq:Zpolar-lterm1}) and 
(\ref{eq:Lemma8.12}) we then have
\begin{equation}\label{eq:Zpolar-lterm2}
Z^{l}(\delta,\eta) = 2\pi i \sum_{x \in \Lambda(l)}\msuml
\frac{\Res_{z=i\eta} \left\{ e^{-\pi \delta (z+l)^{2}} 
\frac{\K_{l}(z;x)}{\sin^{k}(\pi(z-i\eta))} 
e^{2\pi i r ( m -V(x))z} \right\} }{\Sym_{\pm}(m-V(x))},
\end{equation}
where $\Lambda(l)=T_{l}(\Lambda)$ and 
$\K_{l}(z;x)=\tU_{l}\left(z;T_{l}^{-1}(x)\right)$. In the above
expression the lower bound in the $m$-sum still depends on $l$
since $\Lambda(l)$ depends on $l$. Moreover, in order to be
able to use (\ref{eq:periodicity-one-dimensional}) on the
$l$-sum, the expression needs to be 'sufficiently' periodic in
$l$. We will see below what is meant by that. Let us just
remark here that a main part consists of changing the finite 
index set $\Lambda(l)$ to a new index set $\Gamma(l)$ which is
periodic in $l$. To be able to do this we need to assume that 
there are some additional symmetries present. Firstly, we will
assume that there exists a positive integer $L$ and a subset 
$\mR$ of $\R^{e}$ containing $\cup_{l \in \Z} \Lambda(l)$ such
that $\K_{l+L}(z;x)=\K_{l}(z;x)$ for each $l \in \Z$, 
$z \in \C$ and $x \in \mR$. Secondly, we assume that there 
exists a family of bijections 
$\{S_{j} : \mR \to \mR \}_{j \in I}$ all independent of $l$ 
such that the maps $x \mapsto \K_{l}(z;x)$ and 
$x \mapsto V(x) \pmod{\Z}$ restricted to $\mR$ are invariant 
under these transformations, i.e.\ 
$V(S_{j}(x)) \equiv V(x) \pmod{\Z}$ for all $x \in \mR$
and $\K_{l}(z;S_{j}(x))=\K_{l}(z;x)$ for all $x \in \mR$ and 
all $z \in \C$ and $l \in \Z$.

For a tuple of integers $\mC=(j_{1},\ldots,j_{m}) \in I^{m}$ 
and a sequence of signs $\uep \in \{\pm 1\}^{m}$ we let
$$
S^{\uep}_{\mC}=S_{j_{m}}^{\ep_{m}} \circ \ldots \circ 
S_{j_{1}}^{\ep_{1}},
$$
where $S_{j}^{-1}$ is the inverse of $S_{j}$. We then finally 
make the following ``symmetry'' assumption: There exists a 
family of subsets $\Gamma(l) \subseteq \mR$, $l \in \Z$,
periodic in $l$ such that there for any $l \in \Z$ and any
$x \in \Lambda(l)$ exists a tuple of integers $\mC_{x}^{l}$ and
a sequence of signs $\uep_{x}^{l}$ satisfying
$S_{\uep_{x}^{l}}^{\mC_{x}^{l}}(x) \in \Gamma(l)$. For 
simplicity we will assume that the map 
$\Lambda(l) \to \Gamma(l)$, 
$x \mapsto S_{\uep_{x}^{l}}^{\mC_{x}^{l}}(x)$, is a bijection,
but the following analysis can actually be carried out with 
only minor changes in case this map is only surjective.

Since the product of the periods of two periodic functions is a
period for both of the functions we can assume that the family
of sets $\Gamma(l)$ is periodic in $l$ with a period of $L$, 
i.e.\ $\Gamma(l+L)=\Gamma(l)$ for all $l \in \Z$. Since the 
expression (\ref{eq:Zpolar-lterm2}) only depends on $m$ and 
$V(x)$ through their difference $m-V(x)$ we can always modify 
the different values of $V(x)$ by appropriate integers. In 
particular we have 
\begin{equation}\label{eq:Zpolar-lterm3}
Z^{l}(\delta,\eta) = 2\pi i \sum_{\nu \in \Gamma(l)}\msumml
\frac{\Res_{z=i\eta} \left\{ e^{-\pi \delta (z+l)^{2}} 
\frac{\K_{l}(z;\nu)}{\sin^{k}(\pi(z-i\eta))} 
e^{2\pi i r ( m -V(\nu))z} \right\} }{\Sym_{\pm}(m-V(\nu))}.
\end{equation}
Note that
\begin{equation}\label{eq:periodicity1}
\sum_{\nu \in \Gamma(l)} \K_{l}(z;\nu) e^{-2\pi ir V(\nu) z}
\end{equation}
is periodic in $l$ with a period of $L$ for any fixed 
$z \in \C$. By (\ref{eq:zero-2}) we have
\begin{equation}\label{eq:zero-3}
\sum_{\nu \in \Gamma(l)} \K_{l}(z;\nu) e^{-2\pi ir V(\nu) z}
= (-1)^{kl}\sum_{\lambda \in \Lambda} g(\lambda) q(z+l) 
e^{rQ(z+l)} =0
\end{equation}
for all $z \in \Z$.

Let us consider the case $Z(r)=Z(X;r)$. Here we have
$\Lambda \subseteq \R^{2n}$ by (\ref{eq:Seifert-case}), so 
$e=2n$ here. Let $T_{l}(x,y)=(x,y+l\ubeta/2)$ for 
$x,y \in \R^{n}$. Then 
$\Lambda(l)=T_{l}(\Lambda) \subseteq \{ \pm 1\}^{n} 
\times \left( \frac{1}{2}\Z\right)^{n}$ for all $l \in \Z$.
In particular, we can let
$\mR = \{ \pm 1\}^{n} \times \left( \frac{1}{2}\Z\right)^{n}$. 
We have
$$
V(x,y) = V_{l}\left(T_{l}^{-1}(x,y)\right) = 
\sum_{j=1}^{n} \frac{y_{j} -l\beta_{j}/2}{\alpha_{j}}
+ \frac{l}{2} \sum_{j=1}^{n} \frac{\beta_{j}}{\alpha_{j}} = 
\sum_{j=1}^{n} \frac{y_{j}}{\alpha_{j}}
$$
which is independent of $l$ (and $x$). For $j=1,2,\ldots,n$ we 
let $S_{j} : \mR \to \mR$ be the transformation
$S_{j}(\umu,y)=(\umu,y+\alpha_{j}e_{j})$, where $e_{j}$ is the 
standard $j$th unit vector of $\R^{n}$. Then 
$V(S_{j}(\umu,y)) = V(\umu,y) +1$ for all $j$ and all 
$(\umu,y) \in \mR$. By (\ref{eq:Seifert-case}) and
(\ref{eq:Lemma8.12}) we have
\begin{eqnarray*}
\tU_{l}(z;(x,y)) &=& g(x,y)\beta(l,y)q(z+l;y)
e^{irAz^{2}} \\
&=& (-1)^{nl} \left(\frac{i}{2}\right)^{n}
\left( \prod_{j=1}^{n} x_j \right)
\exp\left(2\pi i \sum_{j=1}^{n} \frac{\rho_j}{\alpha_j}
\left[ry_{j}^{2} + x_{j}y_{j} \right] \right)
\exp\left(\frac{\pi i r}{2} E z^2 \right) \\
&& \times \exp\left(\frac{\pi i r}{2} E l^2 \right)
\exp\left(-2\pi i r \left(\sum_{j=1}^{n} 
\frac{y_j}{\alpha_j} \right)l \right)
\exp\left(-\pi i \left( 
\sum_{j=1}^{n} \frac{x_{j}}{\alpha_{j}} \right) (z+l) \right)
\end{eqnarray*}
for $z \in \C$ and $x,y \in \R^{n}$. But then
\begin{eqnarray*}
\K_{l}(z;(x,y')) &=& (-1)^{nl}
\left(\frac{i}{2}\right)^{n} 
\left( \prod_{j=1}^{n} x_j \right)
\exp\left(2\pi i \sum_{j=1}^{n} x_j \left( 
\frac{\rho_j}{\alpha_j} y_{j}' - \frac{1}{2}\sigma_j l \right)
\right) \\
&& \times \exp\left( 2\pi ir \left( -\frac{1}{4}\beta_0 l^{2} +
\sum_{j=1}^{n} \left[ \frac{\rho_{j}}{\alpha_{j}} {y_{j}'}^{2} 
- \frac{1}{4} \sigma_{j}\beta_{j} l^{2} \right] \right) 
\right) \\
&& \times \exp\left( \frac{\pi i r}{2} E z^{2} \right) 
\exp\left(-\pi i \left(\sum_{j=1}^{n} 
\frac{x_{j}}{\alpha_{j}} \right) z \right) \\
&& \times \exp\left(-2\pi i r \sum_{j=1}^{n} \left(
\sigma_j \left(y_{j}' -\frac{1}{2}\beta_{j}l \right) l \right)
\right)
\end{eqnarray*}
for $z \in \C$ and $x,y' \in \R^{n}$. Thus, for fixed 
$(x,y') \in \mR$ and $z \in \C$, $\K_{l}(z;(x,y'))$
is periodic in $l$ with a period of $2$. Moreover, 
$(x,y') \mapsto \K_{l}(z,(x,y'))$, $\mR \to \C$, is invariant 
under the transformations $S_{j}$ for all $l \in \Z$ and all 
$z \in \C$. Finally, it is easy to see that we can put
$$
\Gamma(l) = \{ \pm 1 \}^{n} \times \{ \un' \in \Z^{n} + 
\frac{1}{2}\ubeta l \; | \; \uzero \leq \un' < \ualpha \; \}.
$$
In particular, the period $L=2$ in this case. Note that the 
last factor in the above expression for $\K_{l}(z;(x,y'))$,
namely
$\exp\left(-2\pi i r \sum_{j=1}^{n} \left(
\sigma_j \left(y_{j}' -\frac{1}{2}\beta_{j}l \right) l \right)
\right)$,
is $1$ for $(x,y') \in \Gamma(l)$.

Let us return to the general case. The next step is to get rid
of of the infinite sum over $m$ in (\ref{eq:Zpolar-lterm3}). A
problem here is that the lower bound in this sum is $V(\nu)$.
To overcome this difficulty we use the following small lemma
which proof is left to the reader. The function $\Sym_{\pm}$ is
given by (\ref{eq:Sym}).

\begin{lem}\label{lem:Lemma8.14}
If $a \in \R$ and $f: \Z \ra \C$ is a function such that
$$
\sum_{\stackrel{m \in \Z}{m \geq a}} 
\frac{f(m)}{\Sym_{\pm}(m-a)}
$$
is convergent, then
$$
\sum_{\stackrel{m \in \Z}{m \geq a}} 
\frac{f(m)}{\Sym_{\pm}(m-a)} = 
\sum_{m=0}^{\infty} \frac{f(m)}{\Sym_{\pm}(m)} 
- \sum_{\stackrel{m \in \Z}{0 \leq m \leq |a|}} 
\frac{\sign (a) f(\sign (a) m) }
{ \Sym_{\pm}(m) \Sym_{\pm} (m-|a|)}.
$$\HS
\end{lem}

By using this result on (\ref{eq:Zpolar-lterm3}) we immediately
get
\begin{eqnarray}\label{eq:Zpolar-lterm4}
Z^{l}(\delta,\eta) &=& 2\pi i \sum_{\nu \in \Gamma(l)}
\left\{ \sum_{m=0}^{\infty}
\frac{\Res_{z=i\eta} \left\{ e^{-\pi \delta (z+l)^{2}} 
\frac{\K_{l}(z;\nu)}{\sin^{k}(\pi(z-i\eta))} 
e^{2\pi i r ( m -V)z} \right\}}{\Sym_{\pm}(m)} \right. \\
&& \left. - \sum_{\stackrel{m\in\Z}{0 \leq m \leq |V|}} 
\frac{\sgV \Res_{z=i\eta} \left\{ e^{-\pi \delta (z+l)^{2}} 
\frac{\K_{l}(z;\nu)}{\sin^{k}(\pi(z-i\eta))} 
e^{2\pi i r \sgV ( m -|V|)z} \right\} }
{ \Sym_{\pm}(m) \Sym_{\pm} (m-|V|)} \right\},\nonumber
\end{eqnarray}
where we write $V$ for $V(\nu)$. We could of course have used 
\reflem{lem:Lemma8.14} directly on (\ref{eq:Zpolar-lterm1}). A
problem with this is however that the resulting expression
contains a sum of the form
$\sum_{0 \leq m \leq 
\left|\frac{B}{2\pi}+\frac{A}{\pi}l\right|}$
depending on $l$ in a bad (=nonperiodic) way. If we use 
\reflem{lem:Lemma8.14} on (\ref{eq:Zpolar-lterm2}) the 
resulting expression contains a sum of the form
$\sum_{0 \leq m \leq |V(x)|}$. The problem here is that $V(x)$
depends on $x \in \Lambda(l)$ and thereby on $l$ in a bad way.

Let $\nu \in \Gamma(l)$ be fixed and let
$$
\Psi(z;m)=\Psi(z;l,m,\nu,\delta) = e^{-\pi \delta (z+l)^{2}} 
\K_{l}(z;\nu) e^{2\pi i r ( m -V)z}.
$$
By using the identity
$$
\sum_{m=0}^{\infty} \frac{e^{2 \pi i r m z}}{\Sym_{\pm}(m)} = 
\frac{i}{2} \cot(\pi r z),
$$
valid for all $z \in D(i\eta,\eta)$, we get
\begin{equation}\label{eq:Residuesum}
\sum_{m=0}^{\infty} \frac{\Res_{z=i\eta} \left\{ 
\frac{\Psi(z;m)}{\sin^{k}(\pi(z-i\eta))} \right\}}
{\Sym_{\pm}(m)} = \frac{i}{2}\Res_{z=i\eta} \left\{ 
\frac{\Psi(z;0)}{\sin^{k}(\pi(z-i\eta))} \cot(\pi r z) 
\right\}.
\end{equation}
To see this, assume that $q(z+l)$ has a zero of order 
$k_{l} \geq 0$ in $z=i\eta$. If $k_{l} \geq k$ then both sides
of (\ref{eq:Residuesum}) are zero. If $k_{l} < k$ we let 
$k'=k-k_{l}$ and put 
$$
\tu (z) = \Psi(z;0) 
\frac{(z-i\eta)^{k'}}{\sin^{k}(\pi(z-i\eta))}.
$$
Then
$$
\sum_{m=0}^{\infty} \frac{\Res_{z=i\eta} \left\{ 
\frac{\Psi(z;m)}{\sin^{k}(\pi(z-i\eta))} \right\}}
{\Sym_{\pm}(m)}
= \left. \sum_{m=0}^{\infty} \frac{1}{(k'-1)!} 
\left( \frac{\dte}{\dte z} \right)^{k'-1} 
\tu_{m}(z) \right|_{z=i\eta},
$$
where $\tu_{m}(z)=\tu(z) \frac{e^{2\pi irmz}}{\Sym_{\pm}(m)}$.
Note that $\tu_{m}$ is analytic on $D(i\eta,\eta)$, where we 
assume that $\eta \leq 1$. Now let $\rho \in ]0,\eta[$. For 
$z \in D(i\eta,\rho)$ we have
$$
\sum_{m=0}^{\infty} |\tu_{m}(z)| \leq |\tu (z)| 
\sum_{m =0}^{\infty} 
\frac{e^{-2\pi mr(\eta-\rho)}}{\Sym_{\pm}(m)}
$$
where the infinite sum on the right-hand side is convergent, so
by Weierstrass' M-test
$
\tU (z) = \sum_{m = 0}^{\infty} \tu_{m}(z)
$
is uniformly convergent on $D(i\eta,\rho)$. It follows that 
$\tU$ is analytic on $D(i\eta,\eta)$ and that all its 
derivatives may be calculated by term-by-term differentiation,
see e.g.\ \cite[p.~95]{Titchmarsh}. Therefore
\begin{eqnarray*}
&&\Res_{z=i \eta} \frac{\Psi(z;0)}{\sin^{k}(\pi(z-i\eta))} 
\frac{i}{2} \cot(\pi r z) = \left. \frac{1}{(k'-1)!} 
\left( \frac{\dte}{\dte z} \right)^{k'-1} \tu (z) \frac{i}{2} 
\cot(\pi r z) \right|_{z=i \eta} \\ 
&& \hspace{1.0in} = \left. \frac{1}{(k'-1)!} 
\left( \frac{\dte}{\dte z} \right)^{k'-1} 
\tU (z) \right|_{z=i \eta} 
= \left. \frac{1}{(k'-1)!} \sum_{m=0}^{\infty} 
\left( \frac{\dte}{\dte z} \right)^{k'-1} 
\tu_{m}(z) \right|_{z=i \eta}.
\end{eqnarray*}
By (\ref{eq:Zpolar-lterm4}) and (\ref{eq:Residuesum}) we have
\begin{eqnarray}\label{eq:Zpolar-lterm5}
Z^{l}(\delta,\eta) &=& 2\pi i \sum_{\nu \in \Gamma(l)}
\left\{ \frac{i}{2}\Res_{z=i\eta} 
\left\{ \frac{\Psi(z;l,0,\nu,\delta)}{\sin^{k}(\pi(z-i\eta))} 
\cot(\pi r z) \right\} \right. \\
&& \left. - \sum_{\stackrel{m\in\Z}{0 \leq m \leq |V|}} 
\frac{\sgV \Res_{z=i\eta} \left\{ 
\frac{\Psi(z;l,\sgV m,\nu,\delta)}{\sin^{k}(\pi(z-i\eta))} 
\right\} }
{ \Sym_{\pm}(m) \Sym_{\pm} (m-|V|)} \right\}.\nonumber
\end{eqnarray}

\begin{lem}\label{lem:Lemma8.16}
Let $\eta_{2} \in ]0,1/r[$. Then the infinite series
$\sum_{l \in \Z} Z^{l}(\delta,\eta)$ is uniformly convergent 
with respect to $\eta$ on $]0,\eta_{2}]$. We therefore have
$$
Z_{\pol}(\delta) := \lieta \sum_{l \in \Z} Z^{l}(\delta,\eta) 
= \sum_{ l \in \Z} Z^{l}(\delta),
$$
where
\begin{eqnarray*}
Z^{l}(\delta) &=& \lieta Z^{l}(\delta,\eta) \\
&=& 2\pi i \sum_{\nu \in \Gamma(l)} \left\{ \frac{i}{2} 
\Res_{z=0} \left\{ 
\frac{\Psi(z;l,0,\nu,\delta)}{\sin^{k}(\pi z)} 
\cot(\pi r z) \right\} \right. \\
&& - \left. \sum_{\stackrel{m\in\Z}{0 \leq m \leq |V|}}
\frac{\sgV \Res_{z=0}\left\{ 
\frac{\Psi(z;l,\sgV m,\nu,\delta)}{\sin^{k}(\pi z)} \right\}}
{ \Sym_{\pm}(m)\Sym_{\pm} (m-|V|)} \right\}.
\end{eqnarray*}
\end{lem}

The proofs of this lemma and the next proposition are given in
Appendix C. By the next proposition we end the calculation of 
the contribution to $Z(r)$ coming from the polar term
$Z_{\pol}(\delta,\eta)$. Note that this contribution is given 
an exact (non-asymptotic) description.

\begin{prop}\label{prop:Lemma8.17}
We have
$$
Z_{\pol}(r) := r\lidel \sdel Z_{\pol}(\delta) = 
\sum_{l=0}^{L-1} Z^{l},
$$
where
\begin{eqnarray*}
Z^{l} &=& \frac{2\pi r}{L} \sum_{\nu \in \Gamma(l)}
\left[ \frac{i}{2} \Res_{z=0} \left\{ 
\frac{\K_{l}(z;\nu)e^{-2\pi i r V z}}{\sin^{k}(\pi z)} 
\cot(\pi r z) \right\} \right. \\
&& - \left. \sum_{\stackrel{m\in\Z}{0 \leq m \leq |V|}} 
\frac{\sgV }{ \Sym_{\pm}(m) \Sym_{\pm} (m-|V|)} 
\Res_{z=0} \left\{  
\frac{\K_{l}(z;\nu)e^{2\pi i r \sgV (m-|V|)z}}
{\sin^{k}(\pi z) } \right\}\right],
\end{eqnarray*}
where $V=V(\nu)$.
\end{prop}

Note that this expression does not in general give a function 
as in (\ref{eq:generic-type-alternative}). The factors 
$e^{2\pi i r q_j}$ can, however, be hidden in the expression. 
Let us see that this is indeed the case when 
$Z(r)=Z(X;r)$.

\begin{cor}\label{cor:polarSeifert}
Let $Z(r)=Z(X;r)$. Then
$$
Z_{\pol}(r) = \sum_{(l,\un') \in J} b_{(l,\un')}
\exp\left(2\pi irq_{(l,\un')}\right) r 
\left(Z_{a}^{(l,\un')}(r) + Z_{2}^{(l,\un')}(r)\right), 
$$
where $q_{(l,\un')}$ is given by {\em (\ref{eq:CSb})},
$$
J = \left\{ \left.(l,\un') \in \{0,1\} \times \left( 
\frac{1}{2}\Z \right)^{2} \; \right| \; \un' \in 
\Z^{n} + \frac{1}{2}l\ubeta,\, \uzero \leq \un' < \ualpha 
\; \right\},
$$
and $b_{(l,\un')} = \pi i (-1)^{nl} (-1)^{n}$, and where
\begin{eqnarray*}
Z_{a}^{(l,\un')}(r) &=& \frac{i}{2} \Res_{z=0} \left\{ 
\frac{\exp\left( \frac{\pi i r}{2} Ez^{2} - 
2\pi ir \pnmalphasum z\right)}{\sin^{n+a_{\ep}g-2}(\pi z)} 
\cot(\pi rz) \right. \\
&& \left. \times \prod_{j=1}^{n} \sin \left( 2\pi \left[ 
\frac{\rho_{j}}{\alpha_{j}} n_{j}' - \frac{z}{2\alpha_{j}} - 
\frac{1}{2}\sigma_{j} l \right]\right)\right\}, \\
Z_{b}^{(l,\un')}(r) &=& - \sum_{\stackrel{m \in \Z}
{0 \leq m \leq \nnmalphasum}} \frac{\sgnmalphasum}
{\Sym_{\pm}(m)\Sym_{\pm}\left(m-\nnmalphasum\right)} \\
&& \times \Res_{z=0} \left\{ 
\frac{\exp\left( \frac{\pi i r}{2} Ez^{2} + 
2\pi ir \sgnmalphasum \left(m-\nnmalphasum \right)z \right)}
{\sin^{n+a_{\ep}g-2}(\pi z)} \right. \\
&& \left. \times \prod_{j=1}^{n} \sin \left( 2\pi \left[ 
\frac{\rho_{j}}{\alpha_{j}} n_{j}' - \frac{z}{2\alpha_{j}} - 
\frac{1}{2}\sigma_{j} l \right]\right)\right\}.
\end{eqnarray*}
\end{cor}

\begin{prf}
The corollary is an immediate consequence of 
\refprop{prop:Lemma8.17} and the remarks above 
\reflem{lem:Lemma8.14} about the Seifert case. We use that
\begin{eqnarray*}
&& \musum \muprod \exp \left(-\pi i \pmualphasum z\right) 
\mgex \\
&=& (2i)^{n} \prod_{j=1}^{n}
\sin \left( 2\pi \left[ 
\frac{\rho_{j}}{\alpha_{j}} n_{j}' - \frac{z}{2\alpha_{j}} - 
\frac{1}{2}\sigma_{j} l \right]\right).
\end{eqnarray*} 
\end{prf}

We observe that this expression is not identical with the one 
stated in \refthm{thm:asymptotics-Enot0}. It has the correct 
form, but the set $J$ contains twice as many elements as the 
set $\mI_{2}$ in (\ref{eq:I2}). By (\ref{eq:CSb}) it follows
that if $(l,\un')=(l,n_{1}',\ldots,n_{j}') \in J$ then
$(l,\un'+(\alpha_{j}-2n_{j}')e_{j}) \in J$ and
$$
q_{(l,\un'+(\alpha_{j}-2n_{j}')e_{j})} = q_{(l,\un')}.
$$
What is even more important is that we expect the connected
components of the moduli space of flat $\SU(2)$--connections on
$X$ to be parametrized by $\mI_{1} \cup \mI_{2}$. Therefore it
is desirable to obtain an expression for $Z_{\pol}(r)$ similar
to the one in the above corollary but with the sum over $J$ 
replaced by a sum over $\mI_{2}$. To obtain the result stated
in \refthm{thm:asymptotics-Enot0} we have to use some extra 
symmetries before using \reflem{lem:Lemma8.14}. Let 
$Z(r)=Z(X;r)$ and let
\begin{equation}\label{eq:Jl}
J_{l} = \left\{ \un' \in \Z^{n} + 
\frac{1}{2}\ubeta l \; \left| \; \right. 
\uzero \leq \un' < \ualpha \; \right\}.
\end{equation}
From (\ref{eq:Zpolar-lterm3}) and the remarks above 
\reflem{lem:Lemma8.14} about the Seifert case we get
{\allowdisplaybreaks
\begin{eqnarray*}
Z^{l}(\delta,\eta) &=& 2\pi i (-1)^{nl}
\left(\frac{i}{2}\right)^{n} 
\sum_{\un' \in J_{l}} 
\exp\left( 2\pi i r q_{(l,\un')}\right) \musum \muprod \\*
&& \times \mgex \msume 
\frac{1}{\Sym_{\pm} \left( m - \nmalphasum \right)} \\ 
&& \times \Res_{z=i\eta} \left\{ e^{-\pi \delta (z+l)^{2}}
\frac{\exp \left( \frac{\pi i r}{2} E z^{2} \right)}
{\sin^{n+a_{\ep}g-2}(\pi(z-i\eta))} \right. \\*
&& \left. \times \exp\left( -\pi i \pmualphasum z \right) 
\exp\left(2\pi i r\left( m - \nmalphasum \right) z \right) 
\right\}.
\end{eqnarray*}}
\noindent Let
\begin{equation}\label{eq:Jl'}
J_{l}' = \left\{ \un' \in \Z^{n} + \frac{1}{2}\ubeta l \; 
\left| \; \right.
\uzero \leq \un' \leq \frac{1}{2} \ualpha \; \right\}.
\end{equation}
In the expression for $Z^{l}(\delta,\eta)$ we can change the 
sum $\sum_{\un' \in J_{l}}$ to
$$
\mumsum \sum_{\un' \in J_{l}'} \msyma
$$
if we at the same time substitute $\mu_{j}'n_{j}'$ for 
$n_{j}'$ everywhere. The function $\Sym_{\Z_{\pm}}$ is given by
(\ref{eq:SymZ}). This leads to the identity
{\allowdisplaybreaks
\begin{eqnarray}\label{eq:Zl}
Z^{l}(\delta,\eta) &=& 2\pi i (-1)^{nl} 
\left(\frac{i}{2}\right)^{n}
\sum_{\un' \in J_{l}'} 
\msyma \exp\left( 2\pi i r q_{(l,\un')} \right) \nonumber \\*
&& \times \mumsum \msumc \frac{1}
{\Sym_{\pm} \left( m - \mualphammsum \right)} \musum 
\muprod \nonumber \\
&& \times \mmugex \nonumber \\
&& \times \Res_{z=i\eta} \left\{ e^{-\pi \delta (z+l)^{2}}
\frac{\exp \left( \frac{\pi i r}{2} E z^{2} \right)}
{\sin^{n+a_{\ep}g-2}(\pi(z-i\eta))} \right. \nonumber \\*
&& \left. \times 
\exp\left( -\pi i \pmualphasum z \right) 
\exp\left(2\pi i r\left( m - \mualphammsum \right) z \right) 
\right\}.
\end{eqnarray}}\noindent
By \reflem{lem:Lemma8.14} and the remarks following that lemma
we get
$$
Z^{l}(\delta,\eta) = 2\pi i (-1)^{nl} 
\left(\frac{i}{2}\right)^{n}
\sum_{\un' \in J_{l}'} 
\msyma \exp\left(2\pi i r q_{(l,\un')} \right)
\left(Z_{0}^{(l,\un')}(\delta,\eta) + 
Z_{1}^{(l,\un')}(\delta,\eta) \right),
$$
where
\begin{eqnarray}\label{eq:Z0lun'}
Z_{0}^{(l,\un')}(\delta,\eta) &=& \mumsum \musum \muprod 
\mmugex \nonumber \\
&& \times \Res_{z=i\eta} \left\{ e^{-\pi \delta (z+l)^{2}}
\frac{\exp \left( \frac{\pi i r}{2} E z^{2} -
2\pi ir \pmualphammsum z \right)}
{\sin^{n+a_{\ep}g-2}(\pi(z-i\eta))} \right. \nonumber \\
&& \left. \times \frac{i}{2} \cot(\pi r z) 
\exp\left( -i\pi \pmualphasum z \right) \right\}
\end{eqnarray}
and
\begin{eqnarray}\label{eq:Z1lun'}
Z_{1}^{(l,\un')}(\delta,\eta) &=& -\mumsum \musum \muprod 
\mmugex \nonumber \\
&& \times \sum_{\stackrel{m \in \Z}
{0 \leq m \leq \nmualphammsum} } 
\frac{\sgmualphammsum}
{\Sym_{\pm} (m) \Sym_{\pm} \left(m-\nmualphammsum \right)} 
\nonumber \\
&& \times \Res_{z=i\eta} \left\{ e^{-\pi \delta (z+l)^{2}} 
\exp\left( -i\pi \pmualphasum z \right) \right. \nonumber \\
&& \left. \times \frac{\exp \left( \frac{\pi i r}{2} E z^{2} +
2\pi ir \sgmualphammsum \left(m - \nmualphammsum \right) z 
\right)}{\sin^{n+a_{\ep}g-2}(\pi(z-i\eta))} \right\}.
\end{eqnarray}
Here
\begin{eqnarray*}
&& \musum \muprod \mmugex 
\exp\left( -i\pi \pmualphasum z \right) \\
&=& (2i)^{n} \prod_{j=1}^{n} \mu_{j}' 
\sin\left( 2\pi \left[ \frac{\rho_{j}}{\alpha_{j}} n_{j}' - 
\frac{z}{2} \frac{\mu_{j}'}{\alpha_{j}} - 
\frac{1}{2}\sigma_{j} l \right] \right),
\end{eqnarray*}
where we use that 
$\sin(x-\pi\sigma_{j}l)=\sin(x+\pi\sigma_{j}l)$ for $x \in \R$
since $\sigma_{j}l \in \Z$. Futhermore, we have
\begin{eqnarray*}
&& \mumsum \left[ \prod_{j=1}^{n} \mu_{j}' 
\sin\left( 2\pi \left[ \frac{\rho_{j}}{\alpha_{j}} n_{j}' - 
\frac{z}{2} \frac{\mu_{j}'}{\alpha_{j}} - 
\frac{1}{2}\sigma_{j} l \right] \right) \right] 
\exp\left( -2\pi i r \pmualphammsum z \right) \\
&=& (-2)^{n} \prod_{j=1}^{n} \left[ i 
\sin\left[ 2\pi \left( \frac{\rho_{j}}{\alpha_{j}} n_{j}' - 
\frac{1}{2} \sigma_{j} l \right)\right] 
\cos\left( \pi \frac{z}{\alpha_{j}} \right)
\sin \left( 2\pi r z \frac{n_{j}'}{\alpha_{j}} \right)\right.\\
&& \left. + 
\cos\left[ 2\pi \left( \frac{\rho_{j}}{\alpha_{j}} n_{j}' - 
\frac{1}{2} \sigma_{j} l \right)\right] 
\sin\left( \pi \frac{z}{\alpha_{j}} \right) 
\cos \left( 2\pi r z \frac{n_{j}'}{\alpha_{j}} \right)\right].
\end{eqnarray*}
By inserting these expressions in the expressions for
$Z_{0}^{(l,\un')}(\delta,\eta)$ and 
$Z_{1}^{(l,\un')}(\delta,\eta)$ and by proceeding like in the 
proofs of \reflem{lem:Lemma8.16} and \refprop{prop:Lemma8.17} 
we finally obtain that $Z_{\pol}(r)=Z_{\pol}(X;r)$,
where $Z_{\pol}(X;r)$ is given by (\ref{eq:Zpolar}).

\subsubsection{\bf An asymptotic description of 
$Z_{\inte,1}(\delta,\eta)$'s contribution to $Z_{\inte}(r)$}

We proceed by analyzing the contribution to $Z_{\inte}(r)$ in 
(\ref{eq:Zinte}) coming from
$$
Z_{\inte,1}(\delta,\eta)=\sum_{(m,\lambda) \in W} g(\lambda) 
\int_{C_{\sd}(m,\lambda)} f(z;m,\lambda) \dte z.
$$
We saw in \reflem{lem:Lemma8.12} that this infinite series is 
absolutely convergent for $\delta \in ]0,2rA_{0}/\pi[$ and 
$\eta \in ]0,\eta_{0}]$, where $\eta_{0}$ is the positive 
constant from (\ref{eq:zstassumption}) and $A_{0}$ is given by
(\ref{eq:A0}).

\begin{lem}\label{lem:Lemma8.18}
For every fixed $\delta \in ]0,2rA_{0}/\pi[$ the infinite 
series $Z_{\inte,1}(\delta,\eta)$ is uniformly convergent 
w.r.t.\ $\eta$ on an interval of the form $[0,\eta_{1}]$, 
where $\eta_{1}>0$.
\end{lem}

The proof is given in Appendix D. By this lemma we conclude 
that
$$
Z_{\inte,1}(\delta) := \lieta Z_{\inte,1}(\delta,\eta) = 
\sum_{(m,\lambda) \in W} g(\lambda)
\lieta \int_{C_{\sd}(m,\lambda)} f(z;m,\lambda) \dte z.
$$
Recall here that $f(z;m,\lambda)$ but not $g(\lambda)$ and 
$C_{\sd}(m,\lambda)$ depend on $\eta$. Let in the following 
$(m,\lambda) \in W$ be arbitrary but fixed and let
$z_{\st}=z_{\st}(m,\lambda)$. By 
(\ref{eq:sdparametrization})--(\ref{eq:sdintegral}) we have 
$$
\int_{C_{\sd}(m,\lambda)} f(z;m,\lambda) \dte z = \eksgA 
e^{-irA z_{\st}^{2}}\iny \psi_{\eta}(t;\delta) 
\exp \left( -r|A| t^{2} \right) \dte t,
$$
where
$\psi_{\eta}(t;\delta) = 
\kappa_{\eta}(\gamma_{(m,\lambda)}(t);\lambda,\delta)$.
By (\ref{eq:qfunction-symmetry1}) and (\ref{eq:sinestimate2})
we have
$$
|\psi_{\eta}(t;\delta)| \leq K e^{Ct+a_{0}|t|}
$$
for all $t \in \R$ and all $\eta \in [0,\eta_{1}]$, where
$C = -\pi\sqrt{2}z_{st}\delta$ and 
$K = \frac{b_{0}}{M^{k}}
e^{-\pi\delta z_{\st}^{2}}e^{a_{0}|z_{\st}|}$.
Since 
$t \mapsto K e^{Ct+a_{0}|t|} e^{-r|A| t^{2}} \in L^{1}(\R)$ and
$\lieta \psi_{\eta}(t;\delta) = \psi_{0}(t;\delta)$ for all $t$
we conclude by Lebesgue's dominated convergence theorem that
$$
\lieta \int_{C_{\sd}(m,\lambda)} f(z;m,\lambda) \dte z = 
\int_{C_{\sd}(m,\lambda)} f_{0}(z;m,\lambda) \dte z,
$$
where $f_{0}(z;m,\lambda)$ is equal to $f(z;m,\lambda)$ in 
(\ref{eq:ffunction}) with $\eta=0$, i.e.\
\begin{equation}\label{eq:typicalf1}
f_{0}(z;m,\lambda) = e^{-\pi\delta z^{2}} q(z;\lambda) 
\frac{e^{rQ(z)}e^{2\pi i r m z}}{\sin^{k}(\pi z)}.
\end{equation}
We have shown

\begin{lem}\label{lem:Lemma8.19}
For every fixed $\delta \in ]0,2rA_{0}/\pi[$ we have
$$
Z_{\inte,1}(\delta) =\lieta Z_{\inte,1}(\delta,\eta) = 
\sum_{(m,\lambda) \in W} g(\lambda) 
\int_{C_{\sd}(m,\lambda)} f_{0}(z;m,\lambda) \dte z,
$$
where the functions $f_{0}(z;m,\lambda)$ are given by 
{\em (\ref{eq:typicalf1})}.\HS
\end{lem}

Next we give an asymptotic description of the integrals 
$I(m,\lambda) = \int_{C_{\sd}(m,\lambda)} 
f_{0}(z;m,\lambda) \dte z$,
$(m,\lambda) \in W$. We follow a standard and completely 
elementary proceedure for calculating such asymptotics. The 
idea is to obtain a partition of $I(m,\lambda)$ into parts in 
such a way that one part constitutes the ``main contribution''
and the other parts are remainder terms, which can be evaluated
easily. The starting point is to expand the preexponential 
factor $\kappa_{0}$ of the integrand 
$f_{0}(z) = \kappa_{0}(z)e^{rQ_{m}(z)}$ around the stationary
point of $Q_{m}$.

Let us first give some preliminary remarks. By 
(\ref{eq:zstassumption}) we can choose a $\rho>0$ such that 
$\kappa_{0}$ is analytic on $D(z_{\st}(m,\lambda),4\rho)$ for 
all $(m,\lambda) \in W$. In fact we can use any 
$\rho \leq \eta_{0}/4$. By (\ref{eq:qfunction-symmetry2}) there
exists a constant $C_{\rho}$ such that
\begin{equation}\label{eq:qestimation}
|q(z+z_{\st}(m,\lambda);\lambda)| \leq C_{\rho}
\end{equation}
for all $(m,\lambda) \in \Z \times \Lambda$ and all 
$z \in D(0,4\rho)$. Let in the following $(m,\lambda) \in W$ be
fixed and write $z_{\st}$ for $z_{\st}(m,\lambda)$. Expand 
$\kappa_{0}$ as a power series around $z_{\st}$
\begin{equation}\label{eq:kappapowerseries}
\kappa_{0}(z;\delta) = \sum_{j=0}^{\infty} 
a_{j}(\delta) (z-z_{\st})^{j}
\end{equation}
for $z \in D(z_{\st},4\rho)$. Moreover, let $N \in \Z_{\geq 0}$
be arbitrary but fixed in the following and write
\begin{equation}\label{eq:kappaR}
\kappa_{0}(z;\delta) = \sum_{j=0}^{N} 
a_{j}(\delta) (z-z_{\st})^{j} + R_{N}(z;\delta).
\end{equation}
We claim that the main contribution to an order of $N$ of the 
integral 
$\int_{C_{\sd}(m,\lambda)} f_{0}(z;m,\lambda) \dte z$ is given
by
$$
I(N;\delta) = I(N;m,\lambda,\delta) = 
\int_{C_{\sd}(m,\lambda)} \kappa_{0}(z;N;\delta) 
e^{rQ_{m}(z)} \dte z,
$$
where 
$\kappa_{0}(z;N;\delta) = \sum_{j=0}^{N} 
a_{j}(\delta)(z-z_{\st})^{j}$.
Accordingly, the remainder term is given my
$$
\ep(N;\delta) = \ep(N;m,\lambda,\delta) = 
I(m,\lambda) - I(N;m,\lambda,\delta).
$$
At a first glance this can seem a little peculiar since the 
power series of $\kappa_{0}$ is only defined in a small 
neighborhood around $z_{\st}$. The point is that the main 
contribution comes from integrating $\kappa_{0}(z;N;\delta)$ 
along the part of $C_{\sd}(m,\lambda)$ lying inside 
$D(z_{\st},4\rho)$. However, by integrating along all of 
$C_{\sd}(m,\lambda)$ we obtain an explicit and nice expression 
for the main contribution. By (\ref{eq:sdintegral}) we have
{\allowdisplaybreaks
\begin{eqnarray*}
I(N;\delta) &=& \sum_{j=0}^{N} a_{j}(\delta) 
\int_{C_{\sd}(m,\lambda)} (z-z_{\st})^{j} 
e^{rQ_{m}(z)} \dte z \\
&=& \sum_{j=0}^{N} a_{j}(\delta) \left(\eksgA\right)^{j+1} 
e^{-irAz_{\st}^{2}}\iny t^{j} e^{-r|A|t^{2}} \dte t \\
&=& 2\sum_{\stackrel{0 \leq j \leq N}{j \in 2\Z}} 
a_{j}(\delta) \left(\eksgA\right)^{j+1} 
e^{-irAz_{\st}^{2}}\int_{0}^{\infty} t^{j} 
e^{-r|A|t^{2}} \dte t \\
&=& \sum_{\stackrel{0 \leq j \leq N}{j \in 2\Z}} 
a_{j}(\delta) \left(\eksgA\right)^{j+1} 
e^{-irAz_{\st}^{2}} \Gamma\left(\frac{j+1}{2}\right)
\left(\frac{1}{r|A|}\right)^{\frac{j+1}{2}}.
\end{eqnarray*}}\noindent
Let us next describe $\ep(N;\delta)$. Note that for any nice 
function $g : \R \to \C$ and any $0 \leq a < b \leq \infty$ we 
have
$\int_{-b}^{-a} g(t) e^{-r|A|t^{2}} \dte t =
\int_{a}^{b} g(-t) e^{-r|A|t^{2}} \dte t$.
Thus
$$
\ep(N;\delta) = \eksgA e^{-irA z_{\st}^{2}}
\sum_{\nu \in \{0,1\}} 
\left( E_{\nu}(N;\delta) + J_{\nu}(\delta)
- \sum_{j=0}^{N} a_{j}(\delta) F_{\nu}^{j} \right),
$$
where
\begin{eqnarray*}
J_{\nu}(\delta) &=& \int_{2\rho}^{\infty} 
\psi((-1)^{\nu}t;\delta)\exp\left(-r|A|t^{2}\right) \dte t, \\
E_{\nu}(N;\delta) &=& \int_{0}^{2\rho} 
R_{N}\left( \gamma_{(m,\lambda)}\left((-1)^{\nu}t\right);
\delta\right) \exp\left(-r|A|t^{2}\right) \dte t,
\end{eqnarray*}
for $\nu=0,1$, and 
$
F_{\nu}^{j} = \int_{2\rho}^{\infty} \left( \gamma_{(m,\lambda)}
\left((-1)^{\nu}t\right) - z_{\st} \right)^{j} 
\exp\left(-r|A|t^{2}\right) \dte t
$
for $\nu=0,1$ and $j=0,1,\ldots,N$. By using the incomplete 
gamma function 
$\Gamma(\alpha,u_{0}) = 
\int_{u_{0}}^{\infty} u^{\alpha-1} e^{-u} \dte u$
we get
$$
F_{\nu}^{j} = \frac{1}{2} \left( (-1)^{\nu} \eksgA \right)^{j} 
\left( \frac{1}{r|A|} \right)^{\frac{j+1}{2}}
\Gamma \left( \frac{j+1}{2},4|A|\rho^{2}r \right).
$$
Using that $F_{0}^{j}+F_{1}^{j}=(1+(-1)^{j})F_{0}^{j}$ and
$\eksgA \left(\eksgA\right)^{2j} 
|A|^{-\frac{2j+1}{2}} = (i/A)^{\frac{2j+1}{2}}$
we get
{\allowdisplaybreaks
\begin{eqnarray*}
I(N;\delta) &=& e^{-irAz_{\st}^{2}} 
\sum_{\stackrel{0 \leq j \leq N}{j \in 2\Z}} 
a_{j}(\delta)\Gamma\left(\frac{j+1}{2}\right)
\left(\frac{i}{rA}\right)^{\frac{j+1}{2}}, \\
\ep(N;\delta) &=& e^{-irAz_{\st}^{2}} 
\left\{ \eksgA \sum_{\nu=0,1} 
\left(E_{\nu}(N;\delta) + J_{\nu}(\delta)\right) \right.\\*
&& \left.-\sum_{\stackrel{0 \leq j \leq N}{j \in 2\Z}} 
a_{j}(\delta)\left(\frac{i}{rA}\right)^{\frac{j+1}{2}}
\Gamma \left( \frac{j+1}{2},4|A|\rho^{2}r \right) \right\}.
\end{eqnarray*}}\noindent
In accordance with the above we put
\begin{eqnarray*}
Z_{\inte,1}(N;\delta) &=& \sum_{(m,\lambda) \in W} 
g(\lambda)I(N;m,\lambda,\delta), \\
R_{1}(N;\delta) &=& Z_{\inte,1}(\delta)- Z_{\inte,1}(N;\delta) 
= \sum_{(m,\lambda) \in W}g(\lambda)\ep(N;m,\lambda,\delta).
\end{eqnarray*}
This decomposition of $Z_{\inte,1}(\delta)$ is justified by the
following lemma.

\begin{lem}\label{lem:Lemma8.21}
For any $N \in \Z_{\geq 0}$ and $\delta >0$ the infinite series
\begin{eqnarray*}
Z_{\inte,1}(N;\delta) &=& \sum_{(m,\lambda) \in W} g(\lambda) 
\exp\left(-irAz_{\st}^{2} \right)
\sum_{\stackrel{0 \leq j \leq N}{j \in 2\Z}}
a_{j}(\delta)\Gamma \left( \frac{j+1}{2} \right)
\left(\frac{i}{rA}\right)^{\frac{j+1}{2}} \\
\Sigma^{1}(N;\delta) &=& \sum_{(m,\lambda) \in W} g(\lambda) 
\exp\left(-irAz_{\st}^{2} \right)
\sum_{\stackrel{0 \leq j \leq N}{j \in 2\Z}}a_{j}(\delta)
\left(\frac{i}{rA}\right)^{\frac{j+1}{2}}
\Gamma \left( \frac{j+1}{2},4|A|\rho^{2}r \right), \\
\Sigma_{\nu}^{2}(N;\delta) &=& \sum_{(m,\lambda) \in W} 
g(\lambda) \exp\left(-irAz_{\st}^{2} \right) 
E_{\nu}(N;\delta),\hspace{.2in}\nu=0,1,
\end{eqnarray*}
are absolutely convergent. Moreover, there exists a 
$r_{0} \in \Z_{\geq 2}$ and a $\delta_{0} >0$ such that the 
infinite series
$$
\Sigma_{\nu}^{3}(\delta) = \sum_{(m,\lambda) \in W} g(\lambda) 
\exp\left(-irAz_{\st}^{2} \right) J_{\nu}(\delta),
\hspace{.2in}\nu=0,1,
$$
are absolutely convergent for all $\delta \in ]0,\delta_{0}]$ 
and $r \in \Z_{\geq r_{0}}$. For all such $\delta$ and $r$ and 
all $N \in \Z_{\geq 0}$ we therefore have
$Z_{\inte,1}(\delta) = Z_{\inte,1}(N;\delta) 
+ R_{1}(N;\delta)$, 
where
$$
R_{1}(N;\delta) = -\Sigma^{1}(N;\delta) + \eksgA \sum_{\nu=0,1}
\left( \Sigma_{\nu}^{2}(N;\delta) + \Sigma_{\nu}^{3}(\delta)
\right).
$$
\end{lem}

The proof is given in Appendix D. By the following lemma we
calculate one of the two main contribution to $Z_{\inte}(r)$ in
(\ref{eq:Zinte}) to arbitrary high order $N \in \Z_{\geq 0}$.
Recall here that in general we allow $g(\lambda)$ to depend on
$r$. By the above lemma we can write
$$
R_{1}(N;\delta) =\sum_{\lambda \in \Lambda} g(\lambda)
\sum_{m \in \Z : (m,\lambda) \in W}
e^{-irA z_{\st}^{2}} D(N,\lambda;m,\delta)
$$
for $r \in \Z_{\geq r_0}$ and $\delta \in ]0,\delta_0]$, where
\begin{eqnarray*}
D(N,\lambda;m,\delta) &=& 
-\sum_{\stackrel{0 \leq j \leq N}{j \in 2\Z}}a_{j}(\delta)
\left(\frac{i}{rA}\right)^{\frac{j+1}{2}}
\Gamma \left( \frac{j+1}{2},4|A|\rho^{2}r \right) \\
&& + \eksgA\sum_{\nu \in \{0,1\}} \left(E_{\nu}(N;\delta) + 
J_{\nu}(\delta) \right).
\end{eqnarray*}

\begin{prop}\label{prop:Lemma8.22}
For each fixed $\lambda \in \Lambda$ and $j \in \Z_{\geq 0}$ 
the functions
\begin{equation}\label{eq:Gjfunctions}
G_{j}(m,\lambda)= \exp\left(-irAz_{\st}(m,\lambda)^{2} \right) 
\partial_{z}^{(j)} \left. 
\frac{q(z)}{\sin^{k}(\pi z)} \right|_{z=z_{\st}(m,\lambda)}
\end{equation}
are periodic in $m$ with a period $M_{\lambda}$ independent of
$j \in \Z_{\geq 0}$. Moreover, for each such set of periods
$\{M_{\lambda}\}_{\lambda \in \Lambda}$ and each 
$N = 0,1,2,\ldots$ we have
\begin{eqnarray}\label{eq:Zintmain1}
Z_{\inte,1}(r;N) &:=& r\lidel\sdel Z_{\inte,1}(N;\delta) \\
&=& \sqrt{\frac{ri}{\pi A}} 
\sum_{\stackrel{0 \leq j \leq N}{j \in 2\Z}} \frac{1}{(j/2)!}
\left(\frac{i}{4r} \right)^{\frac{j}{2}} 
\sum_{\lambda \in \Lambda} \frac{|A|}{M_{\lambda}} 
A^{-\frac{j}{2}} g(\lambda) \nonumber \\
&& \times \sum_{\stackrel{m \in \Z \bmod{M_{\lambda}}}
{z_{\st}(m,\lambda) \in \R \sm \Z}} 
\exp\left(-irAz_{\st}^{2}(m,\lambda) \right)
\partial_{z}^{(j)} \left. 
\frac{q(z)}{\sin^{k}(\pi z)} \right|_{z=z_{\st}(m,\lambda)}.
\nonumber
\end{eqnarray}
Finally, let $\delta_{0}>0$ and $r_0 \in Z_{\geq 2}$ be as in 
{\em \reflem{lem:Lemma8.21}}. For each $N =0,1,2,\ldots$ and
$\lambda \in \Lambda$ there exists a positive function 
$\delta \mapsto \vep_{1}(N,\lambda; \delta)$ on 
$]0,\delta_{0}]$ {\em (}also depending on $r${\em )} such that
$$
\sum_{m \in \Z : (m,\lambda) \in W}
\left| D(N,\lambda;m,\delta) \right| \leq 
\vep_{1}(N,\lambda;\delta),
$$
for all $\delta \in ]0,\delta_{0}]$ and all 
$r \in \Z_{\geq r_0}$ and such that
$r\lidel \sdel \vep_{1}(N,\lambda;\delta)$ exists and is
$O\left(r^{-N/2}\right)$ on $\Z_{\geq r_0}$. If $|g(\lambda)|$
is independent of $r$ it follows that there for each
$N=0,1,2,\ldots$ exists a positive function 
$\delta \mapsto \vep_{1}(N;\delta)$ on $]0,\delta_0]$ 
{\em (}also depending on $r${\em )} such that 
$$
|R_{1}(N;\delta)| \leq \vep_{1}(N;\delta)
$$
for $\delta \in ]0,\delta_{0}]$ and $r \in \Z_{\geq r_0}$ and
such that $r\lidel \sdel \vep_{1}(N;\delta)$ exists and is
$O\left(r^{-N/2}\right)$ on $\Z_{\geq r_0}$.
\end{prop}

The proof is given in Appendix D. By this we have finalized
the description of $Z_{\inte,1}(\delta,\eta)$'s contribution to
the large $r$ asymptotics of $Z(r)$.

Before leaving this section, let us take a closer look at the 
Seifert case, i.e.\ the case where $Z(r)=Z(X;r)$. Let us find 
periods $M_{\lambda}$ as in \refprop{prop:Lemma8.22}. Recall 
that $\mA = \prod_{j=1}^{n} \alpha_{j} >0$. Put $P=2\mA$ and
$$
H = \mA E = -\beta_{0}\mA 
-\sum_{j=1}^{n} \frac{\mA}{\alpha_{j}} \beta_{j} = 
-\beta_{0}\alpha_{1}\ldots\alpha_{n} - \sum_{j=1}^{n} 
\alpha_{1}\ldots\alpha_{j-1}\beta_{j}\alpha_{j+1}
\ldots\alpha_{n}.
$$
Then $A=\pi H/P$. Note that $|H|$ is the order of the torsion 
part of $H_{1}(X;\Z)$ in case of orientable base, i.e\ 
$\ep =\os$, cf.\ \cite[Corollary 6.2]{JankinsNeumann}. For the
following we refer to (\ref{eq:Seifert-case}). By 
(\ref{eq:stationarypoint}) we have that
$$
z_{\st}(m,\un)= -\frac{P}{H}\left( m - \nalphasuma \right)
$$
for $(m,\un) \in \Z \times \Z^{n}$. (Thus this is independent
of $\umu \in \{\pm 1\}^{n}$.) For each fixed $\un \in \Z^{n}$
the function 
$m \mapsto \exp\left( -i r A z_{\st}(m,\un)^{2}\right)$
is periodic with a period of $|H|$. This follows by the facts
that $z_{\st}(m+H,\un) = z_{\st}(m,\un) - P$ and $P \in 2\Z$.
The functions
$$
h_{\umu}(z) = \frac{q(z;\umu)}{\sin^{n+a_{\ep}g-2}(\pi z)}
$$
are periodic with a period of $P$ for any 
$\umu \in \{\pm 1\}^{n}$, and therefore the derivatives of 
$h_{\umu}$ in $z_{\st}(m,\un)$ are periodic in $m$ with a 
period of $|H|$. We can therefore put $M_{(\umu,\un)}=|H|$ for
each $(\umu,\un) \in \Lambda$ and get from (\ref{eq:Zintmain1})
that
\begin{eqnarray*}
&&Z_{\inte,1}( r;N) = \frac{1}{2\mA}
\left(\frac{i}{2}\right)^{n}\sqrt{\frac{2ir}{E}} 
\sum_{\stackrel{0 \leq j \leq N}{j \in 2\Z}} \frac{1}{(j/2)!}
\left( \frac{i}{2\pi E} \right)^{\frac{j}{2}} r^{-\frac{j}{2}} 
\sum_{\un \in \mS} \sum_{\stackrel{m \in \Z \bmod{|H|}}
{z_{\st} \in \R \sm \Z}}
\exp\left( 2\pi i r q_{(m,\un)} \right) \\
&& \hspace{.4in} \times \musum \muprod 
\exp\left(2\pi i \sum_{j=1}^{n} \mu_{j}
\frac{\rho_{j}n_{j}}{\alpha_{j}}\right)
\partial_{z}^{(j)} \left. 
\frac{\exp\left( -i\pi\left(\sum_{j=1}^{n}
\frac{\mu_{j}}{\alpha_{j}} \right) z\right)}
{\sin^{n+a_{\ep}g-2}(\pi z)} \right|_{z=z_{\st}},
\end{eqnarray*}
where $q_{(m,\un)}$ is given by (\ref{eq:CSa}). By calculating
the sum over $\umu$ we get
\begin{eqnarray*}
Z_{\inte,1}(r;N) &=& \frac{(-1)^{n}}{2\mA} \sqrt{\frac{2ir}{E}}
\sum_{\stackrel{0 \leq j \leq N}{j \in 2\Z}} \frac{1}{(j/2)!}
\left( \frac{i}{2\pi E} \right)^{\frac{j}{2}} 
r^{-\frac{j}{2}} \\
&& \times \sum_{\un \in \mS} 
\sum_{\stackrel{m \in \Z \bmod{|H|}}{z_{\st} \in \R \sm \Z}} 
\exp\left(2\pi i r q_{(m,\un)}\right) c_{(m,\un)}^{(j)},
\end{eqnarray*}
where we for $j \in \{0,1,2,\ldots\}$ and 
$(m,\un) \in \Z \times \Z^{n}$ have put
$$
c_{(m,\un)}^{(j)} =  \partial_{z}^{(j)} \left. \left\{ 
\frac{\prod_{j=1}^{n} 
\sin\left( \frac{\pi}{\alpha_{j}}[2\rho_{j}n_{j}-z]\right)}
{\sin^{n+a_{\ep}g-2}(\pi z)} \right\} 
\right|_{z=z_{\st}(m,\un)}.
$$
In particular,
\begin{eqnarray*}
Z_{\inte,1}(r;N) &=& \frac{(-1)^{n}}{2\mA} \sqrt{\frac{2ir}{E}}
\sum_{\stackrel{0 \leq j \leq N}{j \in 2\Z}} \frac{1}{(j/2)!}
\left( \frac{i}{2\pi E} \right)^{\frac{j}{2}} 
r^{-\frac{j}{2}} \\
&& \times \sum_{\un \in \mS} \sum_{\stackrel{m \in \Z}
{z_{\st} \in ]0,2\mA[ \cap (\R \sm \Z)}} 
\exp\left(2\pi i r q_{(m,\un)}\right) c_{(m,\un)}^{(j)}.
\end{eqnarray*}
In this expression certain symmetries are present. Let 
$\underline{e}_{1},\ldots,\underline{e}_{n}$ be the standard 
basis in $\R^{n}$ and let
\begin{eqnarray*}
h_{j}(m,\un) &=& (m+1,\un+\alpha_{j}\underline{e}_{j}), 
\hspace{.2in} j=1,2,\ldots,n \\
f(m,\un) &=& (m,\un-\ubeta), \\
g(m,\un) &=& -(m,\un)
\end{eqnarray*}
for $(m,\un) \in \Z \times \Z^{n}$ and $j=1,2,\ldots,n$. The
function $z_{\st}$ is invariant under the transformations 
$h_{j}$, change sign under $g$ and 
$z_{\st}(f(m,\un))=z_{\st}(m,\un)+2$. Moreover, $q_{(m,\un)}$
and $c_{(m,\un)}^{(j)}$, $j \in 2\Z_{\geq 0}$, are invariant 
under all the above transformations (and their inverses). For 
$a \in \Z$ put
$$
W_{a} = \{ (m,\un) \in \Z \times \mS \; | \; 
z_{\st}(m,\un) \in ]2a,2(a+1)[ \sm \{ 2a+1 \} \;\}.
$$
For any $a \in \Z$ there is a bijection between $W_{a}$ and 
$W_{0}$; namely $(m,\un) \in W_{a}$ is maped to
$h_{1}^{k_{1}} \ldots h_{n}^{k_{n}} f^{-a}(m,\un)$, where
$(k_{1},\ldots,k_{n}) \in \Z^{n}$ is the unique tuple of 
integers such that
$n_{j}+a\beta_{j}+k_{j}\alpha_{j} \in 
\{0,1,\ldots,\alpha_{j}-1\}$.
Finally we partition $W_{0}$ into the two sets
$$
W_{0}^{1} = \{ (m,\un) \in \Z \times \mS \; | \; 
z_{\st}(m,\un) \in ]0,1[ \;\}.
$$
and $W_{0}^{2} = W_{0} \sm W_{0}^{1}$ and observe that there is
a bijection between $W_{0}^{1}$ and $W_{0}^{2}$ constructed 
using $fg$ and the $h_{j}$'s. We have thus shown

\begin{cor}\label{cor:Zint1Seifert}
Let $Z(r)=Z(X;r)$. Then
$$
Z_{\inte,1}(r;N) = \sqrt{r}
\sum_{\stackrel{0 \leq j \leq N}{j \in 2\Z}} 
\sum_{(m,\un) \in \mI_{1}}\exp\left(2\pi irq_{(m,\un)}\right)
c^{(m,\un)}_{j/2}r^{-\frac{j}{2}},
$$
where $q_{(m,\un)}$ is given by {\em (\ref{eq:CSa})} and
$c^{(m,\un)}_{j/2}$ is given by {\em (\ref{eq:A3})}.
\end{cor}

Finally observe that $|g(\umu,\un)| =2^{-n}$ for all 
$(\umu,\un) \in \Lambda$ by (\ref{eq:Seifert-case}). Thus the 
estimate for $R_{1}(N;\delta)$ given at the end of
\refprop{prop:Lemma8.22} is valid.

\subsubsection{\bf An asymptotic description of 
$Z_{\inte,2}(\delta,\eta)$'s contribution to $Z_{\inte}(r)$}

We will follow the same strategy as in the analysis of
$Z_{\inte,1}(\delta,\eta)$, i.e.\ we will partition the 
expression $Z_{\inte,2}(\delta,\eta)$ into simple and 
manageable pieces, one piece giving the main contribution and 
the remaining ones adding up to a remainder part, which is 
small in the large $r$ limit compared to the main contribution
piece.

Let us in the following assume that $q(m;\lambda) \neq 0$ for 
all $m \in \Z$ and all $\lambda \in \Lambda$. By periodicity
and continuity of $q$ it then follows that there exists a 
$\rho>0$ such that $q(z;\lambda) \neq 0$ for all 
$z \in \cup_{m \in \Z} D(m,\rho)$ and all 
$\lambda \in \Lambda$. Assume also that the positive parameter
$\eta < \rho$.

Let $\lambda \in \Lambda$ be arbitrary but fixed and suppress
$\lambda$ from the notation everywhere. By (\ref{eq:F}) we have
$$
\F^{+}(y;l,\delta,\eta) = \exp(-irA\eta^{2}) 
G^{+}(y;l,\delta,\eta) \exp(irAy^{2}),
$$
where
$$
G(y;l,\delta,\eta) = e^{-\pi\delta(y+i\eta+l)^{2}} 
\frac{q(y+i\eta+l) \exp(-2rA\eta y)}{\sin^{k}(\pi y)}
$$
and $G^{+}$ is the even part of $G$, i.e.\
$G^{+}(y;l,\delta,\eta) =  
\left( G(y;l,\delta,\eta) + G(-y;l,\delta,\eta) \right)/2$.
We note that $G^{+}$ is entire if $k \leq 1$ and analytic on 
$D'(0,1)$ with a pole in $0$ of order $m>0$ if $k >1$, where 
$m$ is the biggest even integer less than or equal to $k$. 
Write 
$G(y;l,\delta,\eta) = \sum_{s=-k}^{\infty} 
c_{s}(l,\delta,\eta) y^{s}$
for $y \in D'(0,1)$. Then
$$
G^{+}(y;l,\delta,\eta) = 
\sum_{\stackrel{\nu=0}{\nu -k \in 2 \Z}}^{\infty} 
d_{\nu}(l,\delta,\eta) y^{\nu -k},
$$
where
\begin{equation}\label{eq:dnu}
d_{\nu}(l,\delta,\eta) = \frac{1}{\pi^{k}\nu !} 
\left. \partial_{y}^{(\nu)} \Psi(y;l,\delta,\eta) 
\right|_{y=0},
\end{equation}
where
\begin{equation}\label{eq:Psi}
\Psi(y;l,\delta,\eta) = e^{-\pi\delta(y+i\eta+l)^{2}} 
q(y+i\eta+l) \exp(-2rA\eta y)
\left(\frac{\pi y}{\sin(\pi y)}\right)^{k}.
\end{equation}
By (\ref{eq:Zint2-deltaeta}) the relevant integrals to examine
are
$\int_{C_{\eta}} F^{+}(y;l,\delta,\eta) \dte y = 
\exp(-irA\eta^{2}) I(l,\delta,\eta)$,
where we write $C_{\eta}$ for $C(\lambda,\eta)$ and
$I(l,\delta,\eta) = \int_{C_{\eta}} G^{+}(y;l,\delta,\eta) 
\exp\left(irAy^{2}\right) \dte y$.
Let $N \in \Z_{\geq 0}$ be arbitrary but fixed in the following
and write
$$
G^{+}(y;l,\delta,\eta) = 
\sum_{\stackrel{\nu=0}{\nu -k \in 2 \Z}}^{N} 
d_{\nu}(l,\delta,\eta) y^{\nu -k} + R_{N}(y;l,\delta,\eta)
$$
for $y \in D'(0,1)$. Put
\begin{eqnarray*}
I(N;l,\delta,\eta) &=& 
\sum_{\stackrel{\nu=0}{\nu -k \in 2 \Z}}^{N}
d_{\nu}(l,\delta,\eta)
\int_{C_{\eta}} y^{\nu-k} e^{irAy^{2}} \dte y, \\
\ep(N;l,\delta,\eta) &=& I(l,\delta,\eta) - I(N;l,\delta,\eta).
\end{eqnarray*}
To decribe $\ep(N;l,\delta,\eta)$, let
$C_{\eta}(t) = \eksgA t -i\eta$ be a parametrization of 
$C_{\eta}$, and let $C_{\eta}^{0}$ and $C_{\eta}^{\infty}$ be 
the restrictions of $C_{\eta}(t)$ to $[-1/2,1/2]$ and 
$\R \sm ]-1,2,1/2[$ respectively and put
{\allowdisplaybreaks
\begin{eqnarray*}
I_{\infty}(l,\delta,\eta) &=& \int_{C_{\eta}^{\infty}} 
G^{+}(y;l,\delta,\eta) e^{irAy^{2}} \dte y, \\
J_{1}(N;l,\delta,\eta) &=&  \int_{C_{\eta}^{0}} 
R_{N}(y;l,\delta,\eta) e^{irAy^{2}} \dte y, \\
J_{2}(N;l,\delta,\eta) &=& 
-\sum_{\stackrel{\nu=0}{\nu -k \in 2 \Z}}^{N} 
d_{\nu}(l,\delta,\eta)
\int_{C_{\eta}^{\infty}} y^{\nu -k} 
e^{irAy^{2}} \dte y.
\end{eqnarray*}
Then
$$
\ep(N;l,\delta,\eta) = I_{\infty}(l,\delta,\eta) + 
J_{1}(N;l,\delta,\eta) + J_{2}(N;l,\delta,\eta).
$$
For $\lambda \in \Lambda$ put 
$$
I_{\lambda} = \left\{ l \in \Z \; \left| \; 
\frac{B(\lambda)}{2\pi} + \frac{A(\lambda)}{\pi}l \in \Z \;
\right\}. \right.
$$

\begin{lem}\label{lem:Lemma8.25}
There exists a $\eta_{2}>0$ such that the infinite series
\begin{eqnarray*}
Z_{\inte,2}(N;\delta,\eta) &=& \sum_{\lambda \in \Lambda} 
g(\lambda) e^{-irA\eta^{2}} 
\sum_{l \in I_{\lambda}} \beta(l)I(N;l,\delta,\eta), \\
\Sigma^{4}(N;\delta,\eta) &=& \sum_{\lambda \in \Lambda} 
g(\lambda) e^{-irA\eta^{2}}
\sum_{l \in I_{\lambda}} \beta(l)J_{1}(N;l,\delta,\eta) \\
\Sigma^{5}(N;\delta,\eta) &=& \sum_{\lambda \in \Lambda} 
g(\lambda)e^{-irA\eta^{2}} 
\sum_{l \in I_{\lambda}} \beta(l)J_{2}(N;l,\delta,\eta)
\end{eqnarray*}
are absolutely convergent for all $\delta>0$, all 
$\eta \in ]0,\eta_{2}]$ and all $N \in \Z_{\geq 0}$. 
We can choose $\eta_{2}>0$ such that all the above series are
uniformly convergent w.r.t.\ $\eta$ on $]0,\eta_{2}]$ for each
fixed $\delta>0$ and $N \in \Z_{\geq 0}$.

Moreover, by choosing $\eta_{2}$ sufficiently small, there 
exists a $r_{0} \in \Z_{\geq 2}$ and a 
$\delta_{0} \in ]0,\infty[$ such that
$$
\Sigma^{6}(\delta,\eta) = \sum_{\lambda \in \Lambda} g(\lambda)
e^{-irA\eta^{2}} \sum_{l \in I_{\lambda}} 
\beta(l)I_{\infty}(l,\delta,\eta)
$$
is absolutely convergent for all $\delta \in ]0,\delta_{0}]$,
$\eta \in ]0,\eta_{2}]$ and $r \in \Z_{\geq r_{0}}$, and such
that $\Sigma^{6}(\delta,\eta)$ is uniformly convergent w.r.t.\
$\eta$ on $]0,\eta_{2}]$ for each fixed 
$\delta \in ]0,\delta_{0}]$ and $r \in \Z_{\geq r_{0}}$. 
\end{lem}

The proof is given in Appendix E. By the above lemma 
(and proof of that lemma) we get
\begin{eqnarray*}
Z_{\inte,2}(N;\delta) &:=& \lieta Z_{\inte,2}(N;\delta,\eta) \\
&=& \sum_{\lambda \in \Lambda} g(\lambda) 
\frac{e^{\frac{i\pi}{4}\sgA}}{\sqrt{r|A|}}
\sum_{l \in I_{\lambda}} \beta(l,\lambda)
\sum_{\stackrel{\nu=0}{\nu -k \in 2 \Z}}^{N} 
d_{\nu}(l,\delta,0) \left(\frac{i}{rA}\right)^{\frac{\nu-k}{2}}
\Gamma\left( \frac{\nu -k +1}{2} \right),
\end{eqnarray*}
where we have also used (\ref{eq:Int2}) and where 
$d_{\nu}(l,\delta,0)$ is equal to $d_{\nu}(l,\delta,\eta)$ with
$\eta=0$. Similarly we find that
$$
\Sigma^{\nu}(N;\delta) :=\lieta \Sigma^{\nu}(N;\delta,\eta)
=\sum_{\lambda \in \Lambda} g(\lambda) 
\sum_{l \in I_{\lambda}} \beta(l)J_{\nu-3}(N;l,\delta)
$$
for $\nu=4,5$, where $J_{k}(N;l,\delta)$ is equal to
$J_{k}(N;l,\delta,\eta)$ with $\eta=0$, $k=1,2$. Finally, for 
$\delta \in ]0,\delta_{0}]$ and $r \in \Z_{\geq r_{0}}$ we get
$$
\Sigma^{6}(\delta) :=\lieta \Sigma^{6}(\delta,\eta)
=\sum_{\lambda \in \Lambda} g(\lambda) 
\sum_{l \in I_{\lambda}} \beta(l)I_{\infty}(l,\delta),
$$
where $I_{\infty}(l,\delta)$ is equal to
$I_{\infty}(l,\delta,\eta)$ with $\eta=0$. (That
$\lieta I_{\infty}(l,\delta,\eta)=I_{\infty}(l,\delta,0)$
and $\lieta J_{k}(N;l,\delta,\eta)=J_{k}(N;l,\delta,0)$ follows
by Lebesgue's theorem on dominated convergence.)
Thus, for $\delta \in ]0,\delta_{0}]$ and 
$r \in \Z_{\geq r_{0}}$, we have a decomposition
$$
Z_{\inte,2}(\delta) = Z_{\inte,2}(N;\delta) + 
R_{2}(N;\delta),
$$
where
$$
R_{2}(N;\delta) = \sum_{\lambda \in \Lambda} g(\lambda)
e^{-irA\eta^{2}} \sum_{l \in I_{\lambda}}
\beta(l)\ep(N;l,\delta)
= \Sigma^{4}(N;\delta) + \Sigma^{5}(N;\delta)
+ \Sigma^{6}(\delta),
$$
where
$$
\ep(N;l,\delta)=\lieta \ep(N;l,\delta,\eta)
= I_{\infty}(l,\delta) + J_{1}(N;l,\delta) + J_{2}(N;l,\delta).
$$
The task is now to calculate
$$
Z_{\inte,2}(r;N) := r\lidel\sdel Z_{\inte,2}(N;\delta)
$$
and estimate $R_{2}(N;\delta)$ for small $\delta$ and large 
$r$. For the next proposition, recall that
$\ep(N;l,\delta)=\ep(N,\lambda;l,\delta)$ depends on
$\lambda \in \Lambda$.

\begin{prop}\label{prop:Lemma8.26}
With notation from above we have
\begin{eqnarray}\label{eq:Zintmain2}
Z_{\inte,2}(r;N) &=& \frac{r}{\pi^{k}} 
\sum_{\lambda \in \Lambda}\frac{g(\lambda)}{M_{\lambda}}
\frac{e^{\frac{i\pi}{4}\sgA}}{\sqrt{r|A|}}
\sum_{l \in I_{\lambda} \bmod{M_{\lambda}}}\beta(l,\lambda) \\
&& \times \sum_{\stackrel{\nu=0}{\nu -k \in 2 \Z}}^{N}
\left(\frac{i}{rA}\right)^{\frac{\nu-k}{2}}
\frac{\Gamma\left( \frac{\nu -k +1}{2} \right)}{\nu!}
\left. \partial_{y}^{(\nu)} h(y;l) \right|_{y=0},\nonumber
\end{eqnarray}
where
$$
h(y;l) = q(y+l) \left( \frac{\pi y}{\sin(\pi y)} \right)^{k}
$$
and $M_{\lambda}$ is the least common multiplum of $2$, 
$m_{q}$, $b(\lambda)$ and $P(\lambda)$, where $m_{q}$ and 
$P(\lambda)$ are the integers in 
{\em (\ref{eq:qfunction-symmetry2})} and
{\em (\ref{eq:A})} respectively, and $b(\lambda)$ is an integer
such that $b(\lambda)B(\lambda)/\pi \in \Z$, see 
{\em (\ref{eq:B})}.

Moreover, let $\delta_{0}>0$ and $r_{0} \in \Z_{\geq 2}$ be as
in {\em \reflem{lem:Lemma8.25}}. For each $N = 0,1,2,\ldots$ 
and each $\lambda \in \Lambda$ there exists a function 
$\delta \mapsto \vep_{2}(N,\lambda;\delta)$ on
$]0,\delta_{0}]$ {\em (}also depending on $r${\em )} such that 
$$
\sum_{l \in I_{\lambda}} |\ep(N,\lambda;l,\delta)| \leq 
\vep_{2}(N,\lambda;\delta),
$$
for all $\delta \in ]0,\delta_{0}]$
and $r \in \Z_{\geq r_0}$ and such that 
$r\lidel \sdel \vep_{2}(N,\lambda;\delta)$ exists 
and is $O\left(r^{\frac{k-\nu_{1}+1}{2}}\right)$ on 
$Z_{\geq r_{0}}$, where $\nu_{1} \in \{N+1,N+2\}$ such that 
$v_{1}-k \in 2\Z$. If $|g(\lambda)|$ is independent of $r$ it
follows that there for each $N=0,1,2,\ldots$ exists a function 
$\delta \mapsto \vep_{2}(N;\delta)$ on $]0,\delta_{0}]$
{\em (}also depending on $r${\em )} such that 
$$
|R_{2}(N;\delta)| \leq \vep_{2}(N;\delta),
$$
for all $\delta \in ]0,\delta_{0}]$ and $r \in \Z_{\geq r_0}$ 
and such that $r\lidel \sdel \vep_{2}(N;\delta)$ exists 
and is $O\left(r^{\frac{k-\nu_{1}+1}{2}}\right)$ on 
$Z_{\geq r_{0}}$, where $\nu_{1}$ is as above.
\end{prop}

\refprop{prop:Lemma8.26} is proved by applications of 
\reflem{lem:periodicity}. The estimate on 
$R_{2}(N;\delta,\eta)$ mainly follows from the proof of
\reflem{lem:Lemma8.25}. We give the details in Appendix E. 

Let us finally restrict to the Seifert case, where the relevant
data are given in (\ref{eq:Seifert-case}). Here we have that
$m_{q} = 2\mA = 2\prod_{j=1}^{n} \alpha_{j}$. For all 
$\lambda \in \Lambda$ we can let $b(\lambda)=\mA$ and 
$P(\lambda)=2\mA$. Therefore, we can put $M_{\lambda} = 2\mA$
in \refprop{prop:Lemma8.26} and find that
\begin{eqnarray*}
Z_{\inte,2}(r;N) &=& 
\frac{e^{\frac{i\pi}{4}\sgE}}{\pi^{n+a_{\ep}g-2}}
\sqrt{\frac{2r}{\pi|E|}}\musum \sum_{\un \in \mS} 
\frac{g(\umu,\un)}{2\mA}
\sum_{\stackrel{l \in \Z \bmod{2\mA}}
{\sum_{j=1}^{n} \frac{n_{j} + 
\frac{1}{2}\beta_{j}l}{\alpha_{j}} \in \Z}}
\beta(l,\un) \\
&& \times \sum_{\stackrel{\nu=0}{\nu -n \in 2 \Z}}^{N}
\left( \frac{i2}{r\pi E} \right)^{\frac{\nu-n-a_{\ep}g+2}{2}}
\frac{\Gamma\left( \frac{\nu -n-a_{\ep}g +3}{2} \right)}{\nu!}
\left.\partial_{y}^{(\nu)} h(y;l) \right|_{y=0},
\end{eqnarray*}
where
$$
h(y;l) = \exp\left( -i\pi \pmualphasum (y+l) \right) 
\left( \frac{\pi y}{\sin(\pi y)} \right)^{n+a_{\ep}g-2}.
$$
For $l \in \Z$ arbitrary but fixed we put 
$n_{j}'=n_{j}+\frac{1}{2}\beta_{j}l$ for $\un \in \Z^{n}$ and 
obtain by (\ref{eq:Seifert-case}) and (\ref{eq:Lemma8.12}) that
$$
g(\umu,\un)\beta(l,\un)e^{-i\pi \pmualphasum l} = 
\left(\frac{i}{2}\right)^{n}\muprod k(l,\umu,\un'),
$$
where
$$
k(l,\umu,\un') = (-1)^{nl} e^{2\pi irq_{(l,\un')}}\mgex ,
$$
where $q_{(l,\un')}$ is given by (\ref{eq:CSb}). Therefore
\begin{eqnarray*}
Z_{\inte,2}(r;N) &=& 
\frac{e^{\frac{i\pi}{4}\sgE}}{\pi^{n+a_{\ep}g-2}\mA}
\left(\frac{i}{2}\right)^{n}\sqrt{\frac{r}{2\pi|E|}} 
\musum \muprod \sum_{\un \in \mS}
\sum_{\stackrel{l \in \Z \bmod{2\mA}}
{\un' = \un + \frac{1}{2} \ubeta l \; : \; \nmalphasum \in \Z}}
k(l,\umu,\un') \\
&& \times \sum_{\stackrel{\nu=0}{\nu -n \in 2 \Z}}^{N}
\left( \frac{i2}{r\pi E} \right)^{\frac{\nu-n-a_{\ep}g+2}{2}}
\frac{\Gamma\left( \frac{\nu -n-a_{\ep}g +3}{2} \right)}{\nu!} 
\left. \partial_{y}^{(\nu)} h(y;0) \right|_{y=0}.
\end{eqnarray*}
Since $k(l,\umu,\un')$ and the condition $\nmalphasum \in \Z$
are invariant under the transformations
$\frac{1}{2}\Z^{n} \to \frac{1}{2}\Z^{n}$,
$\un' \mapsto \un' + \alpha_{j}\underline{e}_{j}$,
$j=1,2,\ldots,n$, we get
\begin{eqnarray*}
Z_{\inte,2}(r;N) &=& 
\frac{e^{\frac{i\pi}{4}\sgE}}{\pi^{n+a_{\ep}g-2}\mA}
\left(\frac{i}{2}\right)^{n}\sqrt{\frac{r}{2\pi|E|}}
\sum_{\stackrel{\nu=0}{\nu -n \in 2 \Z}}^{N}
\left( \frac{i2}{r\pi E} \right)^{\frac{\nu-n-a_{\ep}g+2}{2}}
\frac{\Gamma\left( \frac{\nu -n-a_{\ep}g +3}{2} \right)}
{\nu!} \\
&& \times \sum_{l \in \Z \bmod{2\mA}} \musum \muprod
\sum_{ \un' \in J_l \; : \; \nmalphasum \in \Z }
k(l,\umu,\un') \left.\partial_{y}^{(\nu)} h(y;0) \right|_{y=0},
\end{eqnarray*}
where $J_l$ is given by (\ref{eq:Jl}). For fixed $(\umu,\un')$
we have that $k(l,\umu,\un')$ is periodic in $l$ with a period 
of $2$. Moreover, the sum
$\sum_{ \un' \in J_l \; : \; \nmalphasum \in \Z }$ is periodic 
in $l$ with a period of $2$. Therefore
\begin{eqnarray*}
Z_{\inte,2}(r;N) &=& 
\frac{e^{\frac{i\pi}{4}\sgE}}{\pi^{n+a_{\ep}g-2}}
\left(\frac{i}{2}\right)^{n}\sqrt{\frac{r}{2\pi|E|}}
\sum_{\stackrel{\nu=0}{\nu -n \in 2 \Z}}^{N}
\left( \frac{i2}{r\pi E} \right)^{\frac{\nu-n-a_{\ep}g+2}{2}}
\frac{\Gamma\left( \frac{\nu -n-a_{\ep}g +3}{2} \right)}
{\nu!} \\
&& \times \sum_{l \in \Z \bmod{2}} \musum \muprod
\sum_{ \un' \in J_l \; : \; \nmalphasum \in \Z }
k(l,\umu,\un') \left.\partial_{y}^{(\nu)} h(y;0) \right|_{y=0}.
\end{eqnarray*}
Finally, if $J_{l}'$ is the set in (\ref{eq:Jl'}), we can 
change the sum
$\sum_{ \un' \in J_l \; : \; \nmalphasum \in \Z }$ to
$$ 
\mumsum \sum_{ \un' \in J_{l}' \; : \; \mualphammsum \in \Z }
\msyma
$$
if we everywhere change $n_{j}'$ to $\mu_{j}'n_{j}$. By this we
end up with
\begin{eqnarray*}
Z_{\inte,2}(r;N) &=& 
\frac{e^{\frac{i\pi}{4}\sgE}}{\pi^{n+a_{\ep}g-2}}
\left(\frac{i}{2}\right)^{n}\sqrt{\frac{r}{2\pi|E|}}
\sum_{\stackrel{\nu=0}{\nu -n \in 2 \Z}}^{N}
\left( \frac{i2}{r\pi E} \right)^{\frac{\nu-n-a_{\ep}g+2}{2}}
\frac{\Gamma\left( \frac{\nu -n-a_{\ep}g +3}{2} \right)}
{\nu!} \\
&& \times \sum_{l \in \Z \bmod{2}} \musum \mumsum 
\sum_{ \un' \in J_{l}' \; : \; \mualphammsum \in \Z } \msyma \\
&& \times \muprod (-1)^{nl}
\exp\left(2\pi ir q_{(l,\un')}\right) \\
&& \times \mmugex 
\left.\partial_{y}^{(\nu)} h(y;0) \right|_{y=0}.
\end{eqnarray*}
By evaluating the sum over $\umu$ we arrive at

\begin{cor}\label{cor:Zint2Seifert}
Let $Z(r)=Z(X;r)$. Then 
$$
Z_{\inte,2}(r;N) = \sqrt{r} 
\sum_{\stackrel{\nu=0}{\nu -n \in 2 \Z}}^{N}
\sum_{(l,\un') \in \mI_{2}} 
\exp\left(2\pi ir q_{(l,\un')}\right)c^{(l,\un')}_{\nu'}
r^{-\nu'}
$$
where $q_{(l,\un')}$ is given by {\em (\ref{eq:CSb})} and
$c^{(l,\un')}_{\nu'}$ is given by {\em (\ref{eq:A3'})}, where
$\nu'=(\nu-n-a_{\ep}g+2)/2$.\HS
\end{cor}

As already noted after \refcor{cor:Zint1Seifert} we have that
$|g(\umu,\un)|=2^{-n}$ for all $(\umu,\un) \in \Lambda$ so the
estimate for $R_{2}(N;\delta,\eta)$ given at the end of
\refprop{prop:Lemma8.26} is valid.

\section{Proof of \refthm{thm:asymptotics-E0}}
\label{sec-Proof-of-E0}

\noindent In this section we calculate the large $r$ 
asymptotics of 
$Z(r)$ in (\ref{eq:Zfunction-typical-formula1}) in case $A=0$.
Like in the case $A \neq 0$ one can calculate this asymptotics
by finding a good deformation of the integration contour (the 
real axes) of the integral
$$
I(m,\lambda)=\iny f(z;m,\lambda) \dte z
$$
to some other contour (depending on $(m,\lambda)$). However, in
this case we do not have a natural choice of a good such
deformation, since the steepest descent method is not 
applicable. In the case $A \neq 0$ we saw that by deforming the
real axes in $I(m,\lambda)$ to the steepest descent contour we 
obtain an integrand with an exponential decay of the form 
$e^{-ct^{2}}$ with $c$ a positive constant. It turns out that 
we can again obtain integrands with an exponential decay by 
deforming the real axes to a line parallel to the real line and
lying in the upper or lower half plane depending on the sign of
a certain parameter. Like in the case $A \neq 0$ this leads to 
a decomposition of $Z(r)$ into a residue part and a part being 
equal to a limit of a certain infinite sum of integrals. If 
$a_{0}-k\pi<0$ we can show that this limit is equal to zero, 
$a_{0}$ being the constant in (\ref{eq:qfunction-symmetry1}). 
However, if $a_{0}-k\pi \geq 0$, we have not been able to find 
a nice reduction of the expression involving the infinite sum 
of integrals. Let us give the details. Let $q(z)$ be an $r$ 
independent function satisfying (\ref{eq:qfunction-symmetry1}) 
and (\ref{eq:qfunction-symmetry2}). Let us write 
$q(z)=q_{0}(z)q_{1}(z)$, where
\begin{equation}\label{eq:qfunction-simple}
q_{0}(z;\lambda) = e^{\pi i\zeta(\lambda)z}
\end{equation}
with $\zeta(\lambda) \in \Q$ for all $\lambda \in \Lambda$. In
particular $q_{1}$ is independent of $r$ and satisfies
(\ref{eq:qfunction-symmetry1}) and 
(\ref{eq:qfunction-symmetry2}) with the numbers $m_{q}$, 
$a_{n}$ and $b_{n}$ replaced by some other numbers 
$m_{q}'=m_{q_{1}}$, $a_{n}'$ and $b_{n}'$. We will use this 
splitting in a flexible way. In particular we are interested in
the extreme cases, with $q_{0}$ respectively $q_{1}$ being 
constantly one. We have
$$
f(z;m,\lambda) = e^{-\pi \delta z^{2}} q_{1}(z)
\frac{\exp \left( 2\pi i r \left[m + \Bzeta \right]z\right)}
{\sin^{k}\left( \pi(z-i\eta)\right)}
$$
with $B=B(\lambda)$ and $\zeta=\zeta(\lambda)$, so
$$
Z(r) = r\lidel \sdel \lieta \sum_{m \in \Z}
\sum_{\lambda \in \Lambda} g(\lambda) 
\iny e^{-\pi \delta z^{2}}q_{1}(z)
\frac{e^{2\pi i r \nu(m,\lambda)z}}
{\sin^{k}\left( \pi(z-i\eta)\right)}
\dte z,
$$
where $\nu(m,\lambda) =  m + \Bzeta$. For 
$R \in \Z + \frac{1}{2}$, $\nu \in \R$ and 
$\rho \in ]0,\infty[$, let 
$\gamma(t) = \gamma_{(\nu,\rho,R)}(t) = R +i\sign(\nu)t$, 
$t \in [0,\rho]$. Here $\sign(0) \in \{\pm 1\}$ (both options 
will be used). Then
\begin{eqnarray*}
\int_{\gamma} f(z;m,\lambda) \dte z
&=& i\sign(\nu) e^{2\pi ir\nu R} e^{-\pi\delta R^{2}} \\
&& \times 
\int_{0}^{\rho} e^{\pi \delta t^{2}} q_{1}(R+i\sign(\nu)t)
\frac{e^{-i2\pi\delta R\sign(\nu)t}e^{-2\pi t|\nu|t}}
{ \sin^{k}\left( \pi(R+i(\sign(\nu)t-\eta))\right)} \dte t,
\end{eqnarray*}
so
$$
\left| \int_{\gamma} f(z;m,\lambda) \dte z \right| \leq
b_{0}' e^{-\pi\delta R^{2}} 
\int_{0}^{\rho} e^{\pi \delta t^{2}}
\frac{e^{a_{0}'(m_{q}'+ t)}e^{-2\pi t|\nu|t}}
{ \left| \sin\left( \pi(R+i(\sign(\nu)t-\eta))\right) 
\right|^{k}} \dte t.
$$
By (\ref{eq:normsin}) and our choice of $R$ we have
$\left| \sin\left( \pi (R + i(\sign(\nu)t-\eta))\right) 
\right| \geq 1$
so
$$
\left| \int_{\gamma} f(z;m,\lambda) \dte z \right| \leq
b_{0}' e^{a_{0}'m_{q}'} e^{-\pi\delta R^{2}} 
\int_{0}^{\rho} e^{\pi \delta t^{2}} e^{a_{0}'t} 
e^{-2\pi t|\nu|t}  \dte t.
$$
Here the integral on the right-hand side is finite and 
independent of $R$, hence $\int_{\gamma} f(z;m,\lambda) \dte z$
converges to zero as $|R| \ra \infty$. Now put 
$M(m,\lambda)=0$ if $\sign(\nu(m,\lambda)) = -1$ and
$$
M(m,\lambda) = 2\pi i \sum_{l \in \Z} 
\Res_{z=l +i\eta} \left\{ f(z;m,\lambda) \right\}
$$
if $\sign(\nu(m,\lambda)) = +1$. For the following, let 
$\eta_{0}>0$ be arbitrary but fixed, and let
$\kappa_{\nu}^{\rho}(t) = t + i\sign(\nu)\rho$, $t \in \R$, for
$\rho>\eta_{0}$. Moreover, let $J$ be any subset of
$I=\{(m,\lambda) \in \Z \times \Lambda \;|\; 
\nu(m,\lambda) = 0 \}$ 
and let $\sign(\nu(m,\lambda)) = 1$ for $(m,\lambda) \in J$ 
and $\sign(\nu(m,\lambda)) = -1$ for 
$(m,\lambda) \in I \sm J$ and put
\begin{eqnarray*}
Z_{\inte}^{J}(\delta,\eta) &=& \sum_{m \in \Z} 
\sum_{\lambda \in \Lambda} g(\lambda) 
\int_{\kappa_{\nu(m,\lambda)}^{\rho}} f(z;m,\lambda) \dte z, \\
Z_{\pol}^{J}(\delta,\eta) &=& 
\sum_{(m,\lambda) \in \Z \times \Lambda} 
g(\lambda)M(m,\lambda).
\end{eqnarray*}

\begin{lem}\label{lem:A01}
For any $J \subseteq I$ the infinite series 
$Z_{\inte}^{J}(\delta,\eta)$ and $Z_{\pol}^{J}(\delta,\eta)$ 
are absolutely convergent for all $\delta >0$ and 
$\eta \in ]0,\eta_{0}]$. In particular
$$
Z(r) = r\lidel\sdel\lieta\left( Z_{\inte}^{J}(\delta,\eta) 
+ Z_{\pol}^{J}(\delta,\eta) \right)
$$
for any such subset $J$.
\end{lem}

The proof of this lemma is given in Appendix F. By the estimate
(\ref{eq:kappaintegral}) the series 
$Z_{\inte}^{J}(\delta,\eta)$ is uniformly convergent with 
respect to $\eta$ on an interval on the form $]0,\eta_{0}]$. 
Using this and Lebesgue's dominated convergence theorem we get
$$
Z_{\inte}^{J}(\delta) := \lieta 
Z_{\inte}^{J}(\delta,\eta) 
= \sum_{m \in \Z}\sum_{\lambda \in \Lambda} g(\lambda) 
\int_{\kappa_{\nu}^{\rho}} f_{0}(z;m,\lambda) \dte z,
$$
where
$$
f_{0}(z;m,\lambda) = e^{-\pi \delta z^{2}}
\frac{q_{1}(z)e^{2\pi i r \nu(m,\lambda)z}}{\sin^{k}(\pi z)}.
$$
Let $d_{\rho} = \sqrt{ e^{2\pi\rho} + e^{-2\pi\rho} - 2}$. 
Similar to (\ref{eq:kappaintegral}) we have
$$
\left| \int_{\kappa_{\nu}^{\rho}} f_{0}(z;m,\lambda) \dte z 
\right| \leq b_{0}'e^{a_{0}'(m_{q}'+\rho)} 
\frac{1}{\sqrt{\delta}} 
\left( \frac{2}{d_{\rho}} \right)^{k}
e^{\pi\delta\rho^{2}} e^{-2\pi r \rho |\nu|}
$$
and then
$$
\left| \sqrt{\delta}Z_{\inte}^{J}(\delta) \right| \leq
b_{0}'e^{a_{0}'(m_{q}'+\rho)} 
\left( \frac{2}{d_{\rho}} \right)^{k}
e^{\pi\delta\rho^{2}} \sum_{\lambda \in \Lambda} |g(\lambda)| 
\sum_{m \in \Z} e^{-2\pi r \rho \left|m + \Bzeta \right|}
$$
for all $\rho >0$. (Note that this estimate is independent of 
$J$.) Therefore (below we show that 
$\lidel\sdel\lieta Z_{\pol}^{J}(\delta,\eta)$ exists, hence 
that $\lidel \sdel Z_{\inte}^{J}(\sdel)$ exists)
$$
\left| \lidel \sqrt{\delta}Z_{\inte}^{J}(\delta) \right| 
\leq b_{0}'e^{a_{0}'(m_{q}'+\rho)}
\left( \frac{2}{d_{\rho}} \right)^{k}
\sum_{\lambda \in \Lambda} |g(\lambda)| 
\sum_{m \in \Z} e^{-2\pi r \rho \left|m + \Bzeta \right|}
$$
Here
\begin{eqnarray*}
\sum_{m \in \Z} e^{-2\pi r \rho \left|m + \Bzeta \right|}
&\leq& \iny e^{-2\pi r\rho|x|} \dte x 
+ \max_{x \in \R} \left\{ e^{-2\pi r\rho |x|} \right\} \\
&=& 2\int_{0}^{\infty} e^{-2\pi r\rho x} \dte x + 1 
= 1+\frac{1}{\pi r\rho}.
\end{eqnarray*}
In the limit $\rho \ra \infty$ we have 
$d_{\rho} \sim e^{\pi\rho}$ so 
$e^{a_{0}'\rho}/d_{\rho}^{k} \sim e^{(a_{0}'-k\pi)\rho}$, 
and since $Z_{\inte}^{J}(\delta)$ is independent of 
$\rho \in ]0,\infty[$ we find by letting $\rho \ra \infty$ that
$$
\lidel \sdel Z_{\inte}^{J}(\delta) = 0
$$
for any $J \subseteq I$ if $a_{0}'-k\pi < 0$. Thus, in that 
case we have that
$$
Z(r)=r\lidel\sdel\lieta Z_{\pol}^{J}(\delta,\eta)
$$
for any $J$, hence the right-hand side is independent of
$J$ in case $a_{0}'-k\pi <0$. Thus we can conclude that non of
the terms with $(m,\lambda) \in I$ contribute to $Z(r)$ in that
case. For $(m,\lambda) \in I$ we have
$$
\lieta M(m,\lambda) = 2\pi i \sum_{l \in \Z}
(-1)^{kl} \Res_{z=0} \left\{ e^{-\pi \delta (z+l)^{2}}
\frac{q_{1}(z+l)}{\sin^{k}(\pi z)} \right\}.
$$
Here $(-1)^{kl}q_{1}(z+l)$ is periodic in $l$ with a period of 
$2m_{q}'$, thus by the standard argument we obtain
$$
\lidel\sdel\lieta M(m,\lambda) = \frac{\pi i}{m_{q}'} 
\sum_{l=0}^{2m_{q}'-1}
(-1)^{kl} \Res_{z=0} \left\{
\frac{q_{1}(z+l)}{\sin^{k}(\pi z)} \right\}
$$
and by the above argument this limit is zero in case 
$a_{0}' -k\pi <0$. If $q_{1}(z)=1$ (so $m_{q}'=1$ and 
$a_{0}'=0$), this is the trivial statement that
$$
\sum_{l=0}^{1}
(-1)^{kl} \Res_{z=0} \left\{
\frac{1}{\sin^{k}(\pi z)} \right\} =0.
$$
We calculate 
$r\lidel \sdel \lieta Z_{\pol}^{J}(\delta,\eta)$ 
by using some of the same techniques as in the case 
$A \neq 0$, however, this case is much simpler. Indeed we have
\begin{eqnarray*}
Z_{\pol}^{J}(\delta,\eta) &=& A_{+}^{J}(\delta,\eta) -
A_{-}^{J}(\delta,\eta) +
2\pi i\sum_{l \in \Z} \sum_{\lambda \in \Lambda}
g(\lambda) e^{2\pi ir\left( \Bzeta \right)l}(-1)^{kl} \\
&& \times \sum_{m \in \Z : \nu(m,\lambda) \geq 0}
\frac{1}{\Sym_{\pm} (\nu(m,\lambda))}
\Res_{z = i\eta} \left\{ e^{-\pi\delta(z+l)^{2}}
\frac{q_{1}(z+l)e^{2\pi ir\nu(m,\lambda)z}}
{\sin^{k}(\pi(z-i\eta))} \right\},
\end{eqnarray*}
where
$$
A_{j}^{J}(\delta,\eta) = \pi i
\sum_{(m,\lambda) \in J_{j}} g(\lambda)
\sum_{l \in \Z}(-1)^{kl} 
\Res_{z = i\eta} \left\{ e^{-\pi\delta(z+l)^{2}}
\frac{q_{1}(z+l)}{\sin^{k}(\pi(z-i\eta))} \right\}
$$
for $j=+,-$, where $J_{+}=J$ and $J_{-}=I \sm J$. Since the 
summation limit in the sum over $m$ is independent of $l$ we 
can immediately use \reflem{lem:Lemma8.14} to obtain
\begin{eqnarray*}
&&Z_{\pol}^{J}(\delta,\eta) = A_{+}^{J}(\delta,\eta)
- A_{-}^{J}(\delta,\eta)
+ 2\pi i\sum_{l \in \Z} \sum_{\lambda \in \Lambda}
g(\lambda) e^{2\pi ir\left( \Bzeta \right)l} (-1)^{kl} \\
&& \times \left( \frac{i}{2}
\Res_{z=i\eta}\left\{ e^{-\pi \delta(z+l)^{2}} 
\frac{q_{1}(z+l)e^{2\pi ir \left( \Bzeta \right)z}}
{\sin^{k}(\pi(z-i\eta))} \cot(\pi r z) \right\} \right. \\
&& \left. + 
\sum_{\stackrel{m \in \Z}{0 \leq m \leq \left| \Bzeta \right|}}
\frac{\sign\left( \Bzeta \right)
\Res_{z=i\eta} \left\{ e^{-\pi\delta(z+l)^{2}}
\frac{q_{1}(z+l) e^{-2\pi ir\sign\left( \Bzeta \right) 
\left( m - \left| \Bzeta \right| \right)z}}
{\sin^{k}(\pi(z-i\eta))} \right\}}
{\Sym_{\pm}(m)\Sym_{\pm}
\left( m - \left| \Bzeta \right| \right)} \right).
\end{eqnarray*}
Like in the case $A \neq 0$ we get by uniform convergence that
$$
Z_{\pol}^{J}(\delta) := \lieta Z_{\pol}^{J}(\delta,\eta)
= Z_{\pol}^{J}(\delta,0).
$$
Finally, since $B(\lambda)$ and $\zeta(\lambda)$ are rational 
for all $\lambda \in \Lambda$ we arrive via 
\reflem{lem:periodicity} at

\begin{thm}\label{thm:A02}
The limit 
$Z_{\pol}^{J}(r) = r\lidel\sdel Z_{\pol}^{J}(\delta)$ 
exists. In fact, there exists a family of positive integers 
$\{ M_{\lambda} \}_{\lambda \in \Lambda}$ such that
\begin{eqnarray*}
&&Z_{\pol}^{J}(r) = r\left(A_{+}^{J}(r) - A_{-}^{J}(r) \right)
+ 2\pi ir\sum_{\lambda \in \Lambda}
\frac{g(\lambda)}{M_{\lambda}} 
\sum_{l \in \Z \bmod{M_{\lambda}}}
e^{2\pi ir\left( \Bzeta \right)l} (-1)^{kl} \\
&& \hspace{.3in} \times 
\left( \frac{i}{2}\Res_{z=0}\left\{  
\frac{q_{1}(z+l) e^{2\pi ir \left( \Bzeta \right)z}}
{\sin^{k}(\pi z)} \cot(\pi r z) \right\} \right. \\
&& \hspace{.3in} \left.
+ \sum_{\stackrel{m \in \Z}
{0 \leq m \leq \left| \Bzeta \right|}}
\frac{\sign\left( \Bzeta \right)
\Res_{z=0} \left\{ 
\frac{q_{1}(z+l)e^{-2\pi ir\sign\left( \Bzeta \right) 
\left( m - \left| \Bzeta \right| \right)z}}
{\sin^{k}(\pi z)} \right\}}
{\Sym_{\pm}(m)\Sym_{\pm}
\left( m - \left| \Bzeta \right| \right)} \right),
\end{eqnarray*}
where
$$
A_{j}^{J}(r) = \pi i
\sum_{(m,\lambda) \in J_{j}} \frac{g(\lambda)}{M_{\lambda}}
\sum_{l \in \Z \bmod{M_{\lambda}}} (-1)^{kl}
\Res_{z = 0} \left\{ 
\frac{q_{1}(z+l)}{\sin^{k}(\pi z)} \right\}
$$
for $j =+,-$. If $a_{0}'-\pi k<0$, then 
$Z(r) = Z_{\pol}^{J}(r)$ for any $J \subseteq I$. In 
particular, $A_{j}^{J}(r)=0$, $j=+,-$ in that case.
\end{thm}

The disadvantage with the above expression for $Z(r)$ is that 
the summation limit depends on $r$ (except in the extreme case 
where $q_{1}=q$, where $\zeta(\lambda)=0$ for all 
$\lambda \in \Lambda$). Let us derive a slightly different 
expression valid for all large $r$ avoiding this problem. Let 
$\lambda \in \Lambda$. If $B(\lambda) \in \R \sm 2\pi\Z$ or 
$\zeta(\lambda) \geq 0$ we can choose a positive integer 
$r_{\lambda}$ such that
$$
\left\{\; m \in \Z \;\left| \; m + \frac{B(\lambda)}{2\pi}
+ \frac{1}{2r}\zeta(\lambda) \geq 0 \;\right. \right\} 
= \left\{ \;m \in \Z \;\left| \; m 
+ \frac{B(\lambda)}{2\pi} \geq 0 \;\right.\right\}
$$
for all $r \geq r_{\lambda}$. If $B(\lambda) \in 2\pi\Z$ and 
$\zeta(\lambda) <0$ we can choose $r_{\lambda}$ such that
$$
\left\{ \;m \in \Z \;\left| \; m + \frac{B(\lambda)}{2\pi}
+ \frac{1}{2r}\zeta(\lambda) \geq 0 \;\right. \right\} 
= \left\{ \;m \in \Z \;\left| \; m 
+ \frac{B(\lambda)}{2\pi} > 0 \;\right.\right\}
$$
for all $r \geq r_{\lambda}$. Put 
$r_{0} = \max\{ r_{\lambda} \}_{\lambda \in \Lambda}$.
For $r \geq r_{0}$ we have
\begin{eqnarray*}
Z_{\pol}^{I}(\delta,\eta) &=& A_{1}(\delta,\eta) 
+ 2\pi i\sum_{l \in \Z} \sum_{\lambda \in \Lambda}
g(\lambda) e^{2\pi ir\left( \Bzeta \right)l}(-1)^{kl} \\
&& \times \sum_{\stackrel{m \in \Z}{m +\frac{B}{2\pi} \geq 0}}
\Res_{z = i\eta} \left\{ e^{-\pi\delta(z+l)^{2}}
\frac{q_{1}(z+l)e^{2\pi ir\nu(m,\lambda)z}}
{\sin^{k}(\pi(z-i\eta))} \right\},
\end{eqnarray*}
where
$$
A_{1}(\delta,\eta) = -2\pi i\sum_{l \in \Z} 
\sum_{\stackrel{\lambda \in \Lambda}{B \in 2\pi\Z;\;\zeta <0}}
g(\lambda) e^{\pi i\zeta l}(-1)^{kl}
\Res_{z = i\eta} \left\{ e^{-\pi\delta(z+l)^{2}}
\frac{q_{1}(z+l)e^{\pi i\zeta z}}
{\sin^{k}(\pi(z-i\eta))} \right\}
$$
We therefore have
\begin{eqnarray*}
Z_{\pol}^{I}(\delta,\eta) &=& A_{1}(\delta,\eta) +
A_{2}(\delta,\eta) + 
2\pi i\sum_{l \in \Z} \sum_{\lambda \in \Lambda}
g(\lambda) e^{2\pi ir\left( \Bzeta \right)l}(-1)^{kl} \\
&& \times \sum_{\stackrel{m \in \Z}{m +\frac{B}{2\pi} \geq 0}}
\frac{1}{\Sym_{\pm}\left( m + \frac{B}{2\pi} \right)}
\Res_{z = i\eta} \left\{ e^{-\pi\delta(z+l)^{2}}
\frac{q_{1}(z+l)e^{2\pi ir\nu(m,\lambda)z}}
{\sin^{k}(\pi(z-i\eta))} \right\},
\end{eqnarray*}
where
$$
A_{2}(\delta,\eta) = \pi i\sum_{l \in \Z} 
\sum_{\stackrel{\lambda \in \Lambda}{B \in 2\pi\Z}}
g(\lambda) e^{\pi i\zeta l}(-1)^{kl}
\Res_{z = i\eta} \left\{ e^{-\pi\delta(z+l)^{2}}
\frac{q_{1}(z+l)e^{\pi i\zeta z}}
{\sin^{k}(\pi(z-i\eta))} \right\}.
$$
We note that
$$
A_{1}(\delta,\eta) + A_{2}(\delta,\eta) 
=A_{+}(\delta,\eta) + A_{0}(\delta,\eta) - A_{-}(\delta,\eta),
$$
where
$$
A_{j}(\delta,\eta) = \pi i\sum_{l \in \Z} 
\sum_{\lambda \in \Lambda_{j}}
g(\lambda) e^{\pi i\zeta l}(-1)^{kl}
\Res_{z = i\eta} \left\{ e^{-\pi\delta(z+l)^{2}}
\frac{q_{1}(z+l)e^{\pi i\zeta z}}
{\sin^{k}(\pi(z-i\eta))} \right\}
$$
for $j=+,0,-$, where 
$\Lambda_{j} = \{ \lambda \in \Lambda\;|\; B \in 2\pi\Z,\;
j\zeta(\lambda) > 0 \}$
for $j=+,-$ and
$\Lambda_{0} = \{ \lambda \in \Lambda\;|\; B \in 2\pi\Z,\;
\zeta(\lambda) = 0 \}$.
Similarly, we find that
\begin{eqnarray*}
Z_{\pol}^{\emptyset}(\delta,\eta) &=& A_{+}(\delta,\eta) -
A_{0}(\delta,\eta) - A_{-}(\delta,\eta)
+ 2\pi i\sum_{l \in \Z} \sum_{\lambda \in \Lambda}
g(\lambda) e^{2\pi ir\left( \Bzeta \right)l}(-1)^{kl} \\
&& \times \sum_{\stackrel{m \in \Z}{m +\frac{B}{2\pi} \geq 0}}
\frac{1}{\Sym_{\pm}\left( m + \frac{B}{2\pi} \right)}
\Res_{z = i\eta} \left\{ e^{-\pi\delta(z+l)^{2}}
\frac{q_{1}(z+l)e^{2\pi ir\nu(m,\lambda)z}}
{\sin^{k}(\pi(z-i\eta))} \right\}.
\end{eqnarray*}
By letting $I_{+}=I$ and $I_{-}=\emptyset$ we therefore have

\begin{thm}\label{thm:A03}
There exists a family of positive integers 
$\{ M_{\lambda} \}_{\lambda \in \Lambda}$ such that
\begin{eqnarray*}
Z_{\pol}^{I_{\tau}}(r) &=& r\lidel \sdel 
\lieta Z_{\pol}^{I_{\tau}}(\delta,\eta) \\
&=& r\left(A_{+} + \tau A_{0}- A_{-}\right) 
+ 2\pi ir\sum_{\lambda \in \Lambda}
\frac{g(\lambda)}{M_{\lambda}} 
\sum_{l \in \Z \bmod{M_{\lambda}}}
e^{2\pi ir\left( \Bzeta \right)l} (-1)^{kl} \\
&& \times \left( \frac{i}{2}\Res_{z=0}\left\{  
\frac{q_{1}(z+l)e^{2\pi ir \left( \Bzeta \right)z}}
{\sin^{k}(\pi z)} \cot(\pi r z) \right\} \right. \\
&& \left. + \sum_{\stackrel{m \in \Z}
{0 \leq m \leq \left| \frac{B}{2\pi} \right|}}
\frac{\sign(B) \Res_{z=0} \left\{ 
\frac{q_{1}(z+l)e^{2\pi ir\left( \sign(B)m + \Bzeta \right)z}}
{\sin^{k}(\pi z)} \right\}}
{\Sym_{\pm}(m)
\Sym_{\pm}\left( m - \left| \frac{B}{2\pi} \right| \right)}
\right)
\end{eqnarray*}
for all $r \geq r_{0}$ and $\tau =+,-$, where 
$A_{j} = \lidel\sdel\lieta A_{j}(\delta,\eta)$ is given by
$$
A_{j} = \pi i\sum_{\lambda \in \Lambda_{j}}
\frac{g(\lambda)}{M_{\lambda}} 
\sum_{l \in \Z \bmod{M_{\lambda}}}
e^{\pi i\zeta l} (-1)^{kl}\Res_{z=0} \left\{ 
\frac{q_{1}(z+l)e^{\pi i \zeta z}}{\sin^{k}(\pi z)} \right\}
$$
for $j=+,0,-$. If $a_{0}'-k\pi<0$, then 
$Z(r)=Z_{\pol}^{I_{\tau}}(r)$, $\tau =+,-$. In particular, 
$A_{0}=0$ in this case.
\end{thm}

If $q_{1}=q$ (so $\zeta(\lambda)=0$ for all 
$\lambda \in \Lambda$) and $a_{0}-k\pi <0$, then $A_{0}=0$ 
and $A_{-}=A_{+}=0$ since $\Lambda_{j}=\emptyset$, $j=-,+$.
Note that \refthm{thm:A03} is a corollary of \refthm{thm:A02}
in this case. In fact $(m,\lambda) \mapsto \lambda$ defines a 
bijection $I \to \Lambda_0$ and 
$A_{0} = A_{+}^{I}(r) = A_{-}^{\emptyset}(r)$ (which is 
independent of $r$ in this case). 

Let us finally take a closer look at the Seifert case, i.e.\ 
the case $Z(r)=Z(X;r)$. Then 
$q(z)=\exp\left(-i\pi \pmualphasum z\right)$ and
$B=-2\pi \nalphasuma$. By using \refthm{thm:A02} and 
\refthm{thm:A03} directly we would end up with a result similar
to \refcor{cor:polarSeifert}. As argued below that corollary we
obtain in this way a sum with pairs of terms with the same
Chern--Simons value. We will therefore follow a strategy 
similar to the one followed after \refcor{cor:polarSeifert}.

Recall that we always assume that $a_{\ep}g$ is even. 
We will futher in what follows assume that 
$n+a_{\ep}g-2>0$. This is not a serious assumption. Since we 
treat the case $n=0$ as the case $n=1$ with 
$(\alpha_1,\beta_1)=(1,0)$ we have $n+a_{\ep}g-2>0$ except in
case $g=0$ (hence $\ep=\os$) and $n=1,2$. Since we
futher assume that $E=0$, we find that $X=S^{2} \times S^{1}$ 
in these cases (compare with 
\cite[Sect.~4 pp.~28--31]{JankinsNeumann} and the discussion
following \refthm{thm:asymptotics-E0}).

Let us first put $q_{0}(z)=q(z)$, so 
$\zeta = \zeta(\umu) = - \mualphasum$. Since $q_{1}=1$ (so 
$a_{0}'=0$) we get from the general argument that
$$
Z(X;r) = r\lidel \sdel \lieta Z_{\pol}^{J}(\delta,\eta),
$$
where
\begin{eqnarray*}
Z_{\pol}^{J}(\delta,\eta) &=& 2\pi i\sum_{l \in \Z} 
\musum \sum_{\un \in \mS}g(\umu,\un) \\
&& \times \sum_{\stackrel{m \in \Z}{\sign(\nu(m,\umu,\un))=1}}
\Res_{z=l+i\eta} \left\{ e^{-\pi\delta z^{2}}
\frac{e^{2\pi i r \nu(m,\umu,\un)z}}
{\sin^{n+a_{\ep}g-2}(\pi(z-i\eta))} \right\},
\end{eqnarray*} 
where 
$\nu(m,\umu,\un) = m - \nalphasuma - \frac{1}{2r}\mualphasum$ 
and $g(\umu,\un)$ is given by (\ref{eq:Seifert-case}). Here
$\sign(\nu(m,\umu,\un)=1$ if $(m,\umu,\un) \in J$ and
$\sign(\nu(m,\umu,\un)=-1$ if $(m,\umu,\un) \in I \sm J$
(so here the terms with $\nu(m,\umu,\un)=0$ do not contribute 
to $Z(X;r)$). Let us specialize to the case $J=I$ and let us 
write $Z_{\pol}$ for $Z_{\pol}^{I}$. Like in the $E \neq 0$ 
case we first rewrite the above expression by introducing 
$n_{j}'=n_{j}+\frac{1}{2}\beta_{j}l$. This leads to
\begin{eqnarray*}
Z_{\pol}(\delta,\eta) &=& A_{+}^{I}(\delta,\eta)
+ 2\pi i\left(\frac{i}{2}\right)^{n}\sum_{l \in \Z} 
(-1)^{nl} \sum_{\un' \in J_{l}} \exp(2\pi irq_{(l,\un')}) \\
&& \times \musum \muprod  \mgex \\
&& \times \sum_{\stackrel{m \in \Z}{\nu(m,\umu,\un')) \geq 0}}
\frac{1}{\Sym_{\pm} \left( \nu(m,\umu,\un') \right)}
\Res_{z=i\eta} \left\{ e^{-\pi\delta (z+l)^{2}}
\frac{e^{2\pi i r \nu(m,\umu,\un')z}}
{\sin^{n+a_{\ep}g-2}(\pi(z-i\eta))} \right\},
\end{eqnarray*} 
where $J_{l}$ is given by (\ref{eq:Jl}) and 
\begin{eqnarray*}
A_{+}^{I}(\delta,\eta) &=& \pi i\left(\frac{i}{2}\right)^{n}
\sum_{l \in \Z} 
(-1)^{nl} \sum_{\un' \in J_{l}} \exp(2\pi irq_{(l,\un')}) \\
&& \times \musum \muprod  \mgex \\
&& \times \sum_{\stackrel{m \in \Z}{\nu(m,\umu,\un')) = 0}}
\Res_{z=i\eta} \left\{ e^{-\pi\delta (z+l)^{2}}
\frac{1}{\sin^{n+a_{\ep}g-2}(\pi(z-i\eta))} \right\}.
\end{eqnarray*} 
Note that $\nmalphasum=\nalphasuma$ since $E=0$. Like in the 
case $E \neq 0$ we can change the sum $\sum_{\un' \in J_{l}}$ 
to
$$
\mumsum \sum_{\un' \in J_{l}'} \msyma
$$
if we at the same time substitute $\mu_{j}'n_{j}'$ for $n_{j}'$
everywhere. This leads to the identity
\begin{eqnarray*}
Z_{\pol}(\delta,\eta) &=& A_{+}^{I}(\delta,\eta) +
2\pi i\left(\frac{i}{2}\right)^{n}\sum_{l \in \Z} 
(-1)^{nl} \sum_{\un' \in J_{l}'} \msyma 
\exp(2\pi irq_{(l,\un')}) \\
&& \times \mumsum
\sum_{\stackrel{m \in \Z}{\nu(m,\umu,\un',\umu')) \geq 0}}
\musum \muprod 
\frac{1}{\Sym_{\pm}\left( \nu(m,\umu,\un',\umu') \right)} \\
&& \times \mmugex \\
&& \times \Res_{z=i\eta} \left\{ e^{-\pi\delta (z+l)^{2}}
\frac{e^{2\pi i r \nu(m,\umu,\un',\umu')z}}
{\sin^{n+a_{\ep}g-2}(\pi(z-i\eta))} \right\},
\end{eqnarray*} 
where $J_{l}'$ is given by (\ref{eq:Jl'}) and
$\nu(m,\umu,\un',\umu') = m - \mualphammsum -
\frac{1}{2r}\mualphasum$.
By an application of \reflem{lem:Lemma8.14} we thus find that
\begin{eqnarray*}
Z_{\pol}(\delta,\eta) &=& A_{+}^{I}(\delta,\eta)
+ 2\pi i\left(\frac{i}{2}\right)^{n}\sum_{l \in \Z} 
(-1)^{nl} \sum_{\un' \in J_{l}'} \msyma \\
&& \times \exp(2\pi irq_{(l,\un')})
\left( Z_{0}^{(l,\un')}(\delta,\eta) 
+ \tilde{Z}_{1}^{(l,\un')}(\delta,\eta) \right),
\end{eqnarray*}
where $Z_{0}^{(l,\un')}(\delta,\eta)$ is given by 
(\ref{eq:Z0lun'}) (with $E=0$) and
{\allowdisplaybreaks
\begin{eqnarray*}
\tilde{Z}_{1}^{(l,\un')}(\delta,\eta) &=& -\mumsum \musum 
\muprod \mmugex \\*
&& \times \sum_{\stackrel{m \in \Z}{0 \leq m \leq |a|} } 
\frac{\sign(a)}
{\Sym_{\pm} (m) \Sym_{\pm}\left( m - |a| \right)} \\
&& \times \Res_{z=i\eta} \left\{ e^{-\pi \delta (z+l)^{2}} 
\frac{\exp \left( 2\pi ir \sign(a)
\left( m - |a| \right) z \right)}
{\sin^{n+a_{\ep}g-2}(\pi(z-i\eta))} \right\}, \\
\end{eqnarray*}}\noindent
where $a=a(\umu,\umu',\un') = \mualphammsuma$. Since the terms
with $\nu(m,\umu,\un)=0$ do not contribute to $Z(X;r)$ we have
$\lidel\sdel\lieta A_{+}^{I}(\delta,\eta) = 0$.
The first part of \refthm{thm:asymptotics-E0} thus follows by 
calculating $\lidel\sdel\lieta Z_{\pol}(\delta,\eta)$ like in 
the case $E \neq 0$.

Like \refthm{thm:A02} the disadvantage with the above obtained
expression for $Z(X;r)$ is that the summation limit in 
$\tilde{Z}_{1}^{(l,\un')}(r)$, see \refthm{thm:asymptotics-E0},
depends on the level $r$. Let us proceed like in the proof of
\refthm{thm:A03} to obtain a slightly different expression for
$Z(X;r)$. In fact, there exists a positive integer $r_{0}$ 
(depending on $X$) such that
\begin{eqnarray*}
Z_{\pol}^{I_{\tau}}(\delta,\eta) &=& A_{+}(\delta,\eta) + 
\tau A_{0}(\delta,\eta) - A_{-}(\delta,\eta) \\
&& + 2\pi i\left(\frac{i}{2}\right)^{n}
\sum_{l \in \Z} (-1)^{nl} \sum_{\un' \in J_{l}'} 
\msyma \exp(2\pi irq_{(l,\un')}) \mumsum \\
&& \times \musum \muprod \msumc 
\frac{1}{\Sym_{\pm}\left( m - \mualphammsum \right)} \\
&& \times \mmugex \\
&& \times \Res_{z=i\eta} \left\{ e^{-\pi\delta (z+l)^{2}}
\frac{e^{2\pi i r \nu(m,\umu,\un',\umu')z}}
{\sin^{n+a_{\ep}g-2}(\pi(z-i\eta))} \right\}
\end{eqnarray*}
for all $r \geq r_{0}$, where 
$\nu(m,\umu,\un',\umu') = m - \mualphammsum -
\frac{1}{2r}\mualphasum$
as above and
\begin{eqnarray*}
A_{j}(\delta,\eta) &=& \pi i\left(\frac{i}{2}\right)^{n}
\sum_{l \in \Z} (-1)^{nl} \sum_{\un' \in J_{l}'} \msyma 
\exp(2\pi irq_{(l,\un')}) \\
&& \times \sum_{(\umu,\umu') \in \Lambda(\un',j)} \muprod 
\mmugex \\
&& \times \Res_{z=i\eta} \left\{ e^{-\pi\delta (z+l)^{2}}
\frac{\exp\left[-\pi i \pmualphasum z \right]}
{\sin^{n+a_{\ep}g-2}(\pi(z-i\eta))} \right\}
\end{eqnarray*}
for $j=+,0,-$, where
$$
\Lambda(\un',j) = \left\{ (\umu,\umu') \in \B
\times \B \; \left| \; -\sign\pmualphasum = j, 
\mualphammsum \in \Z \right\},\right.
$$
where $\sign(0)=0$. (Note the minus in front of 
$\sign\pmualphasum$ due to the fact that 
$\zeta=\zeta(\umu)=-\mualphasum$.) By \reflem{lem:Lemma8.14} we
then find for $r \geq r_{0}$ that
\begin{eqnarray*}
Z_{\pol}^{I_{\tau}}(\delta,\eta) &=& \pi i
\left(\frac{i}{2}\right)^{n}\sum_{l \in \Z} 
(-1)^{nl} \sum_{\un' \in J_{l}'} \msyma 
\exp(2\pi irq_{(l,\un')})
\left( 2Z_{0}^{(l,\un')}(\delta,\eta) \right.\\
&& \left. + 2Z_{1}^{(l,\un')}(\delta,\eta) + 
A_{+}^{(l,\un')}(\delta,\eta) +
\tau A_{0}^{(l,\un')}(\delta,\eta)
- A_{-}^{(l,\un')}(\delta,\eta) \right),
\end{eqnarray*}
where $Z_{0}^{(l,\un')}(\delta,\eta)$ and 
$Z_{1}^{(l,\un')}(\delta,\eta)$ are given as in the case 
$E \neq 0$, i.e.\ they are given by respectively 
(\ref{eq:Z0lun'}) and (\ref{eq:Z1lun'}) with $E=0$, and where
\begin{eqnarray*}
A_{j}^{(l,\un')}(\delta,\eta) &=&  
\sum_{(\umu,\umu') \in \Lambda(\un',j)} \muprod \mmugex \\
&& \times \Res_{z=i\eta} \left\{ e^{-\pi\delta (z+l)^{2}}
\frac{\exp\left[-\pi i \pmualphasum z \right]}
{\sin^{n+a_{\ep}g-2}(\pi(z-i\eta))} \right\}
\end{eqnarray*}
for $j=+,0,-$. This expression leads to
\begin{equation}\label{eq:ZE0}
Z(X;r) = Z_{\pol}(X;r) + r Z_{\spec}(X;r)
\end{equation}
for $r \geq r_{0}$, where $Z_{\pol}(X;r)$ is as in 
\refthm{thm:asymptotics-Enot0} (with $E=0$) and
$$
Z_{\spec}(X;r)= \sum_{(l,\un') \in \mI_{2}^{a}} b_{(l,\un')}
\exp(2\pi irq_{(l,\un')})
\left( A_{+}^{(l,\un')} - A_{-}^{(l,\un')} \right),
$$
where $\mI_{2}$ and $b_{(l,\un')}$ are as in 
\refthm{thm:asymptotics-Enot0} and
\begin{eqnarray*}
A_{j}^{(l,\un')} &=& \frac{\pi i}{2} 
\left(\frac{i}{4}\right)^{n}
\sum_{(\umu,\umu') \in \Lambda(\un',j)} \muprod
\exp\left( 2\pi i\sum_{j=1}^{n} 
\mu_{j}\frac{\rho_{j}}{\alpha_{j}}\mu_{j}'n_{j}'\right) \\
&& \times \Res_{z=0} 
\left\{ \frac{\exp\left[-\pi i \pmualphasum z \right]}
{\sin^{n+a_{\ep}g-2}(\pi z)} \right\}
\end{eqnarray*}
for $j=+,-$. (Note that we necessarily have 
$\lidel\sdel\lieta A_{0}(\delta,\eta)=0$.) 

Now put $q_{1}(z)=q(z)$ (so $q_{0}(z)=1$, $\zeta=0$). In this 
case (\ref{eq:qfunction-symmetry1}) is valid with
$a_{0}=\pi \sum_{j=1}^{n} \frac{1}{\alpha_{j}}$ (and 
$b_{0}=1$). Assume in the following that 
$a_{0} < \pi(n+a_{\ep}g-2)$. Then we have from the general 
case, see \refthm{thm:A02}. that
$$
Z(X;r)=r\lidel\sdel\lieta Z_{\pol}^{J}(\delta,\eta),
$$
where
\begin{eqnarray*}
Z_{\pol}^{J}(\delta,\eta) &=& 2\pi i
\left(\frac{i}{2}\right)^{n}\sum_{l \in \Z} 
(-1)^{nl} \sum_{\un' \in J_{l}} \exp(2\pi irq_{(l,\un')}) \\
&& \times \musum \muprod  \mgex \\
&& \times \sum_{\stackrel{m \in \Z}{\sign(\nu(m,\un'))=1}}
\Res_{z=i\eta} \left\{ e^{-\pi\delta (z+l)^{2}}
\frac{q(z)e^{2\pi i r \nu(m,\un')z}}
{\sin^{n+a_{\ep}g-2}(\pi(z-i\eta))} \right\},
\end{eqnarray*} 
where $\nu(m,\un')=m-\nmalphasum$ and $\sign(\nu(m,\un'))=1$ 
if $(m,\un') \in J$ and $\sign(\nu(m,\un'))=-1$ if 
$(m,\un') \in I \sm J$, where $J$ is any subset of $I$. By the
same proceedure as used above we immediately get that
$Z(X;r) = Z_{\pol}(X;r)$, where $Z_{\pol}(X;r)$ is as in 
\refthm{thm:asymptotics-Enot0} (with $E=0$). A comparison with 
(\ref{eq:ZE0}) reveals that 
$Z_{\spec}(X;r)=0$ for $r \geq r_{0}$ in case 
$a_{0} < \pi(n+a_{\ep}g-2)$, and since the $q_{(l,\un')}$'s are
rational, we in fact have that $Z_{\spec}(X;r)=0$ for all
levels $r$ if $a_{0} < \pi(n+a_{\ep}g-2)$.

Let $(l,\un') \in \mI_{2}$ be fixed. If $\umu' \in \B$ 
satisfies that $\mualphammsum \in \Z$, then the same is true 
for $-\umu'$. Moreover, if $\umu \in \B$ satisfies that 
$\mualphasum<0$, then $-\umu$ satisfies the opposite 
inequality. These observations show that
$Z_{\spec}(X;r)$ is given by the formula in
\refthm{thm:asymptotics-E0}.

The level $r_{0}$ is chosen so that we for all $r \geq r_{0}$,
$l \in \Z$ and 
$(\umu,\umu',\un') \in \B \times \B \times J_{l}'$ have
$$
\left\{\; m \in \Z \;\left| \; m - \mualphammsum
- \frac{1}{2r}\mualphasum \geq 0 \;\right. \right\} 
= \left\{ \;m \in \Z \;\left| \; m 
-\mualphammsum \geq 0 \;\right.\right\}
$$
in case $\mualphasum \notin \Z$ or $\mualphasum \leq 0$, and 
$$
\left\{ \;m \in \Z \;\left| \; m - \mualphammsum
- \frac{1}{2r}\mualphasum \geq 0 \;\right. \right\} 
= \left\{ \;m \in \Z \;\left| \; m 
- \mualphammsum > 0 \;\right.\right\}
$$
if $\mualphasum \in \Z$ and $\mualphasum > 0$. In particular we
can put $r_{0}=2$ (the minimum level) for $n=1$. In general we 
see that $r_{0}$ only depends on the pairs of Seifert 
invariants 
$(\alpha_{1},\beta_{1}),\ldots,(\alpha_{n},\beta_{n})$ 
(so is independent of the data for the base space, i.e.\ the 
genus $g$ and the invariant $\ep$).

\section{Appendices}\label{sec-Appendices}

\noindent We have in the following appendices collected 
material of a technical nature.

\subsection{Appendix A. The periodicity lemma}

The following lemma is proven in bigger generality than needed 
in this paper. However, to make parts of the ideas in this 
paper useable to other situations, in particular to 
multi-dimensional cases, we give it in this generality. One can
find an even more general version of the periodicity lemma in
\cite[Appendix]{HansenTakata}, where it is used in a 
multi-dimensional setting.

\begin{lem}\label{lem:periodicity}
Let $F : \Z ^{n} \ra \C$ be a function periodic in all 
variables with a period of $N_{i}$ in the $i$th variabel. Let 
$P_{i}(x)=A_{i}x^{2}+B_{i}x+C_{i}$ be polynomials with complex
coefficients such that $\rea(A_{i})>0$ for all $i$. Then
\begin{eqnarray*}
&& \sum_{k_{1}=0}^{N_{1}-1} \ldots \sum_{k_{n}=0}^{N_{n}-1} 
F(k_{1},\ldots,k_{n}) \\
&=& \left( \prod_{i=1}^{n}\sqrt{A_{i}} N_{i} \right) \livep 
\vep^{n/2} \sum_{(k_{1},\ldots,k_{n}) \in \Z ^{n}} 
e^{-\pi \vep f(\vep )(P_{1}(k_{1})+\cdots+P_{n}(k_{n}))}
F(k_{1},\ldots,k_{n}) 
\end{eqnarray*}
for all functions $f:]0,\vep_{0}] \ra ]0,\infty[$ with 
$\livep f(\vep) =1$, where $\vep_{0}>0$ is arbitrary but fixed.
\end{lem}

\begin{prf}
The series 
$\sum_{(k_{1},\ldots,k_{n}) \in \Z ^{n}} 
e^{-\pi \vep f(\vep )(P_{1}(k_{1})+\cdots+P_{n}(k_{n}))}
F(k_{1},\ldots,k_{n})$
is absolutely convergent by periodicity of $F$, hence
\begin{eqnarray*}
&& \sum_{(k_{1},\ldots,k_{n}) \in \Z ^{n}} 
e^{-\pi \vep f(\vep )(P_{1}(k_{1})+\cdots+P_{n}(k_{n}))}
F(k_{1},\ldots,k_{n}) \\
&=& \sum_{j_{1}=0}^{N_{1}-1} \ldots 
\sum_{j_{n}=0}^{N_{n}-1} F(j_{1},\ldots,j_{n}) 
\sum_{(k_{1},\ldots,k_{n}) \in \Z ^{n}} 
e^{-\pi \vep f(\vep )(P_{1}(j_{1}+k_{1}N_{1}) + \cdots
+ P_{n}(j_{n}+k_{n}N_{n}))}.
\end{eqnarray*}
Thus we have to show that
\begin{equation}\label{eq:limit}
\livep \svep \sqrt{A} N Q(\vep) = 1,
\end{equation}
where
$Q(\vep) = \sum_{k \in \Z} e^{-\pi \vep f(\vep)P(j+kN)}$
for a polynomial $P(x)=Ax^{2}+Bx+C$ with $\rea(A)>0$. 
Let $j\in \{0,1,\ldots,N-1 \}$ be fixed. By Poisson's formula 
(\ref{eq:Poisson}) we have
$$
Q(\vep) =  \sum_{m \in \Z} \iny e^{2 \pi i mx} 
e^{-\pi \vep f(\vep)P(j+xN)} \dte x
$$
which is a sum of Gaussian integrals. Thus
$$
Q(\vep) = \frac{e^{-\pi k(\vep)C}}{N\sqrt{Ak(\vep)}} 
\exp\left( \frac{\pi}{4} \frac{B^{2}k(\vep)}{A} \right) 
\sum_{m \in \Z} e^{-\pi i \frac{2jA+B}{AN}m} 
\exp\left( -\frac{\pi}{Ak(\vep)N^{2}}m^{2} \right),
$$
where $k(\vep)=\vep f(\vep)$ from which it follows that 
$Q(\vep)$ is absolutely convergent. The term coming from $m=0$ 
gives the contribution
$$
\livep \frac{e^{-\pi k(\vep)C}}{\sqrt{f(\vep)}} 
\exp\left( \frac{\pi}{4} \frac{B^{2}k(\vep)}{A} \right) = 1
$$
to the left-hand side of (\ref{eq:limit}). To see that the 
remaining part is zero we simply note that
\begin{eqnarray*}
&& \left|  \sum_{m \in \Z \sm \{0\}}  
e^{-\pi i \frac{2jA+B}{AN}m} 
\exp\left( -\frac{\pi}{Ak(\vep)N^{2}}m^{2} \right)\right| \\
&\leq& \exp\left( -\frac{1}{2}\frac{\pi}{Ak(\vep)N^{2}}\right)
\sum_{m \in \Z \sm \{0\}} 
\exp\left( -\frac{1}{2} \frac{\pi}{Ak(\vep)N^{2}}m^{2} +
\frac{\pi\ima(B)}{AN}m \right).
\end{eqnarray*}
where the factor in front of the $m$-sum in the right-hand side
converges to zero as $\vep \ra 0_{+}$ and the $m$-sum in the
right-hand side is majorized by the convergent and 
$\vep$--independent series
$\sum_{m \in \Z \sm \{0\}} 
\exp\left( -\frac{1}{2} \frac{\pi}{2A\vep_{0}N^{2}}m^{2} +
\frac{\pi\ima(B)}{AN}m \right)$
on an interval on the form $]0,\vep_1]$,
where $\vep_1 \leq \vep_0$ is chosen such that $f(\vep)\leq 2$
for all $\vep \in ]0,\vep_{1}]$.
\end{prf}

\subsection{Appendix B. Proofs of 
\reflem{lem:Lemma8.7}--\reflem{lem:Lemma8.13}}

For the following proofs we need certain estimates. Let 
$a_{\eta}(t;m,\lambda) = 
| \sin( \pi (\gamma_{(m,\lambda)}(t) - i \eta)) |$,
where $\gamma_{(m,\lambda)}$ is the steepest descent contour, 
see (\ref{eq:sdparametrization}). For $u,v \in \R$ we have
\begin{equation}\label{eq:normsin}
| \sin(u+iv)|= \frac{1}{2} \sqrt{ e^{2v}+e^{-2v}-2\cos(2u)}
\end{equation}
and therefore
\begin{equation}\label{eq:normsin1}
a_{\eta}(t;m,\lambda) = \frac{1}{2}
\sqrt{ e^{ 2\pi ( \sgA \frac{t}{\sqrt{2}} - \eta)} + 
e^{ -2\pi ( \sgA \frac{t}{\sqrt{2}} - \eta)}
- 2 \cos \left(2 \pi \left[ z_{\st}(m,\lambda) + 
\frac{t}{\sqrt{2}} \right] \right) }. 
\end{equation}
Let $\eta_{0}>0$ be as in (\ref{eq:zstassumption}). 
Then $z_{\st}(m,\lambda) + \sgA \eta \in \R \sm \Z$ for all
$(m,\lambda) \in \Z \times \Lambda$ and all 
$\eta \in ]0,\eta_{0}]$. Since the function 
$x \mapsto e^{x} +e^{-x}$ is strictly decreasing on 
$]-\infty,0]$ and strictly increasing on $[0,\infty]$ and takes
the value $2$ for $x=0$ we thus have
\begin{equation}\label{eq:sinestimate1}
d_{\eta} := \inf_{t \in \R, m \in \Z, \lambda \in \Lambda} 
| a_{\eta}(t;m,\lambda)  | > 0
\end{equation}
for every fixed $\eta \in ]0,\eta_{0}]$. Moreover, letting
$\rho=\sqrt{2}\eta_{0}$, there exists a $\alpha \in [0,1[$ such
that
$\cos \left(2 \pi \left[z_{\st}(m,\lambda) + t/\sqrt{2} 
\right]\right) \leq \alpha$
for all $t \in [-\rho,\rho]$ and all $(m,\lambda) \in W$. 
(We use here that 
$\{ \; z_{\st}(m,\lambda) \pmod{\Z} \; | \;
(m,\lambda) \in \Z \times \Lambda \;\}$
is finite. If we only assume (\ref{eq:zstassumption}) then put
$\rho=\eta_{0}$.) Let $s=-\sgA \frac{t}{\sqrt{2}}$. Then
$a_{\eta}(t;m,\lambda) \geq \frac{1}{2} 
\sqrt{ e^{2\pi(s+\eta)} + e^{-2\pi(s+\eta)} -2 }.$
Letting $\eta_{1} =\frac{\rho}{2\sqrt{2}} < \eta_{0}$ we then 
get
$a_{\eta}(t;m,\lambda) \geq \frac{1}{2} 
\sqrt{ e^{\pi\eta_{1}} + e^{-\pi\eta_{1}} -2 }$
for $|t|>\rho$ and $\eta \in [0,\eta_{1}]$. We conclude that
\begin{equation}\label{eq:sinestimate2}
a_{\eta}(t;m,\lambda) \geq M : = 
\min \left\{ \frac{1}{\sqrt{2}}\sqrt{1-\alpha},
\frac{1}{2}\sqrt{ e^{\pi\eta_{1}} + e^{-\pi\eta_{1}} -2 } 
\right\} >0
\end{equation}
for all $t \in \R$, $\eta \in [0,\eta_{1}]$ and 
$(m,\lambda) \in W$.

\begin{prfa}{\bf of \reflem{lem:Lemma8.7}\hspace{.1in}}
Let $(m,\lambda) \in \Z \times \Lambda$ be fixed in the
following and let $C_{\sd}=C_{\sd}(m,\lambda)$, 
$\gamma=\gamma_{(m,\lambda)}$, $z_{\st}=z_{\st}(m,\lambda)$, 
and $f(z)=f(z;m,\lambda)$. Assume first that $k>0$. Observe 
that $\sgA ( \rea z_{l} - z_{\st} ) > \ima z_{l}$ is equivalent
to $\sgA ( l- z_{\st} ) > \eta$, explaining the summation range
in the sum of residues in (\ref{eq:Lemma8.7}). We will first
show that the integral $\int_{C_{\sd}} f(z) \dte z$ is 
convergent for all sufficiently small positive $\eta$. By 
(\ref{eq:ffunction}), (\ref{eq:sdintegral}) and 
(\ref{eq:sinestimate1}) we have
$$
\left| f(\gamma(t)) \right| \leq 
\frac{1}{d_{\eta}^{k}} |e^{-\pi\delta \gamma(t)^{2}}
q(\gamma(t))|\exp \left( -r|A| t^{2} \right)
$$
for every $\delta>0$ and $\eta \in ]0,\eta_{0}]$. Here 
$\left| e^{-\pi\delta \gamma(t)^{2}} \right|
= e^{-\pi \delta z_{\st}^{2}}e^{-\pi\sqrt{2}z_{\st}\delta t}$,
so by (\ref{eq:qfunction-symmetry1}) we get
$$
\left| f(\gamma(t)) \right| \leq 
\frac{b_0}{d_{\eta}^{k}} e^{-\pi \delta z_{\st}^{2}}
e^{-\pi\sqrt{2}z_{\st}\delta t}e^{a_{0}(|z_{\st}|+|t|)}
\exp \left( -r|A| t^{2} \right)
$$
where the right-hand side is integrable. Thus 
$\int_{C_{\sd}} f(z) \dte z$ is convergent and
\begin{eqnarray}\label{eq:ulighedLemma8.7}
\left| \int_{C_{\sd}} f(z) \dte z \right|
&\leq& \frac{b_{0}}{d_{\eta}^{k}} e^{-\pi \delta z_{\st}^{2}} 
e^{a_{0}|z_{\st}|} \iny e^{-\pi\sqrt{2}z_{\st}\delta t}
\left( e^{a_{0}t}+e^{-a_{0}t}\right) 
\exp \left( -r|A| t^{2} \right) \dte t \nonumber \\
&=& \frac{b_{0}}{d_{\eta}^{k}} e^{-\pi \delta z_{\st}^{2}} 
e^{a_{0}|z_{\st}|} \sqrt{\frac{\pi}{r|A|}} 
\left( \exp\left( \frac{a_{+}^{2}}{4r|A|}\right)
+ \exp\left(\frac{a_{-}^{2}}{4r|A|}\right) \right)
\end{eqnarray}
where $a_{\pm} = \pi \delta \sqrt{2} z_{\st} \pm a_{0}$.

Next we estimate the contributions coming from integrals along
curves connecting $C_{\sd}$ by the real axes. Let 
$R_{n} = n+ \frac{1}{2}$ for $n \in \Z$ and let 
$\tau_{n}(t)=R_{n}+it$, $t \in [a,b]$, where $a,b \in \R$. Then
\begin{eqnarray*}
\left| \int_{\tau_{n}} f(z) \dte z \right| &\leq& 
e^{-\pi \delta R_{n}^{2}} \int_{a}^{b}
\frac{ |q(\tau_{n}(t)| }
{ |\sin^{k}(\pi (\gamma(t) - i \eta)) | } e^{\pi \delta t^{2}}
|e^{rQ_{m}(\tau_{n}(t))}| \dte t \\
&\leq& b_{0}e^{a_{0}|R_{n}|}e^{-\pi \delta R_{n}^{2}}
\int_{a}^{b} \frac{ e^{a_{0}|t|} }
{ |\sin^{k}(\pi (\tau_{n}(t) - i \eta))|} e^{\pi \delta t^{2}} 
e^{ 2rA (z_{\st}-R)t} \dte t 
\end{eqnarray*}
by (\ref{eq:qfunction-symmetry1}). By (\ref{eq:normsin}) we 
have
$|\sin(\pi(\tau_{n}(t)-i \eta))| \geq 
\frac{1}{\sqrt{2}} \sqrt{1-\cos(2 \pi R_{n})} =1$
so
$$
\left| \int_{\tau_{n}} f(z) \dte z \right| \leq b_{0}
e^{a_{0}|R_{n}|}e^{- \pi \delta R_{n}^{2}}
\int_{a}^{b} \left(e^{a_{0}t}+e^{-a_{0}t} \right) 
e^{\pi \delta t^{2}} e^{2r A (z_{\st}-R_{n}) t } \dte t.
$$
In our case $\{a,b\}=\{0,\sgA (R_{n}-z_{\st}) \}$, where we 
assume that $a \leq b$. For $c \in \R$ we have
\begin{eqnarray*}
&& \int_{a}^{b} e^{ct} e^{\pi \delta t^{2}} 
e^{2r A (z_{\st}-R_{n}) t } \dte t \\
&=& \int_{a}^{b} e^{\pi \delta t^{2}} 
e^{-\sgA \left[ r |A| (R_{n}-z_{\st}) -\sgA c \right] t}
e^{-\sgA r |A| (R_{n}-z_{\st}) t} \dte t.
\end{eqnarray*}
We can choose a $N_{c} \in \N$ such that
$$
\sgA \left( r|A| (R_{n}-z_{\st}) - \sgA c \right) t \geq 0
$$
for all $t \in [a,b]$ and all $n \in \Z$ with $|n| > N_{c}$. 
Therefore
$$
\left| \int_{\tau_{n}} f(z) \dte z \right| \leq 2b_{0}
e^{a_{0}|R_{n}|}e^{-\pi \delta R_{n}^{2}}
\int_{a}^{b} e^{- \left[ r |A| \sgA (R_{n}-z_{\st}) - 
\pi \delta t \right] t} \dte t
$$
for all $n \in \Z$ with 
$|n| > N=\max\{ N_{-a_{0}},N_{a_{0}}\}$. For $|n| > N$ and 
$\delta>0$ with $r|A| \geq \pi \delta$ we get
\begin{eqnarray*}
\left| \int_{\tau_{n}} f(z) \dte z \right| &\leq& 2b_{0}
e^{a_{0}|R_{n}|}e^{-\pi \delta R_{n}^{2}} \int_{a}^{b} 
e^{- r |A| \left[ \sgA (R_{n}-z_{\st}) - t \right] t} \dte t \\
&\leq& 2b_{0}e^{a_{0}|R_{n}|}e^{-\pi \delta R_{n}^{2}}(b-a)
= 2|R_{n}-z_{\st}|b_{0}e^{a_{0}|R_{n}|}
e^{-\pi \delta R_{n}^{2}}
\end{eqnarray*}
which converges to zero as $|n| \longrightarrow \infty$. Now, 
since $\iny f(z) \dte z$ is convergent for all $\delta>0$ and 
all $\eta \in ]0,\eta_{0}]$, it follows by the residue theorem
and the above that
$$
\sum_{\stackrel{l \in \Z}{\sgA (l-z_{\st}) > \eta}} 
\Res _{z=z_{l}(\eta)} \left\{ f(z) \right\}
$$
is convergent and (\ref{eq:Lemma8.7}) is satisfied for all 
$\eta \in ]0,\eta_{0}]$ and all $\delta \in ] 0, r|A|/\pi ]$. 
We will actually see in what follows that this last sum is 
absolutely convergent. 

The case $k \leq 0$ is handled as above but is easier since the
$\sin$--factor is included in the function $q$ and $\eta$ is 
not present. Note also that in this case we do not need the 
assumption (\ref{eq:zstassumption}).
\end{prfa}

For the next proof and later we note the following elementary 
fact which is stated for the sake of easy reference. Let 
$a>0$, $\delta>0$ and $b \in \C$ be constants and let 
$j \geq 1$. Then
\begin{equation}\label{eq:P1}
\left( \frac{d}{dz} \right)^{j} e^{-a \delta (z+b)^{2}} = 
P_{j}(z+b,\delta) e^{-a \delta (z+b)^{2}},
\end{equation}
where $P_{j}(y,\delta)$ is a polynomial in $\delta$ and $y$. By
induction we find that
$P_{j}(y,\delta) = 
\sum_{k=0}^{j} a_{k}^{j}(\delta) (\delta y)^{k}$
where 
$a_{k}^{j}(\delta) = \sum_{s=0}^{m_{k}^{j}} b_{k}^{j}(s) 
\delta^{s}$ 
is a polynomial such that for every $j \geq 1$ either 
$a_{0}^{j}(\delta) = 0$ or else $b_{0}^{j}(0) = 0$ and 
$m_{k}^{j} \geq 1$. Of importance here is that all terms in a 
derivative of $e^{-a \delta (z+b)^{2}}$ are of the form 
$Az^{k} \delta^{j}e^{-a \delta (z+b)^{2}}$, $A \in \C$, 
$k \in \N$, $j >0$ (i.e.\ there are no terms with $j=0$). In 
what follows $P_{0}(y,\delta)=1$ by definition.

Throughtout we will use the following elementary fact: If $l$
is a nonnegative integer and $P(z)$ is a polynomial of degree 
$\geq 1$ with complex coefficients with leading coefficient 
having positive real part, then 
$\sum_{k =1}^{\infty} k^{l} e^{-P(k)}$ is absolutely convergent
and $\sum_{k \in \Z} k^{l} e^{-P(k)}$ is absolutely convergent
if $P$ is of even degree.

\begin{prfa}{\bf of \reflem{lem:Lemma8.12}\hspace{.1in}}
Let $\eta_{0}$ be as in (\ref{eq:zstassumption}) and let 
$\eta \in ]0,\eta_{0}]$ and $\lambda \in \Lambda$ be arbitrary 
but fixed in the following. Write $C_{\sd}(m)$, $f(z;m)$ and 
$z_{\st}(m)$ for $C_{\sd}(m,\lambda)$, $f(z;m,\lambda)$ and 
$z_{\st}(m,\lambda)$ respectively, $m \in \Z$. By 
(\ref{eq:ulighedLemma8.7}) we have
$$
\left| \int_{C_{\sd}(m)} f(z;m) \dte z \right| \leq C 
e^{-\pi \delta z_{\st}^{2}(m)} e^{a_{0}|z_{\st}(m)|}
\left( \exp\left(\frac{a_{+}^{2}}{4r|A|}\right) + 
\exp\left(\frac{a_{-}^{2}}{4r|A|}\right) \right)
$$
for all $\delta >0$, where $C$ and $a_{0}$ are positive 
constants independent of $m$, and 
$a_{\pm} = \pi \delta \sqrt{2} z_{\st}(m) \pm a_{0}$. We have
$$
e^{-\pi \delta z_{\st}^{2}(m)} 
\exp\left(\frac{a_{\pm}^{2}}{4r|A|}\right)
= e^{-\pi \delta d z_{\st}^{2}(m)}
\exp\left( \frac{1}{4r|A|} 
(a_{0}^{2} \pm 2\sqrt{2}\pi \delta a_{0} z_{\st}(m))\right),
$$
where $d=1-\frac{\pi \delta}{2r|A|}$. By this and 
(\ref{eq:stationarypoint}) we get
$$
\sum_{m \in \Z} \left| \int_{C_{\sd}(m)} f(z;m) \dte z \right| 
\leq C' \sum_{\mu = \pm 1} e^{\mu a } \sum_{m\in \Z}
e^{\frac{\pi a_{0}}{|A|}|m|} e^{(b+\mu c) m} 
\exp\left( - \frac{\pi^{3} \delta}{A^{2}} dm^{2} \right), 
$$
where $C'$, $a$, $b$ and $c$ are constants independent of $m$.
If we choose $\delta \in ]0, 2r|A|/\pi[$ we have $d>0$ and
$$
\sum_{m \in \Z} \left| \int_{C_{\sd}(m)} f(z;m) \dte z \right| 
< \infty.
$$
Next let us consider the sum $\Sigma_{2}$. Let $w=z-l$. Then
$f(z;m)=\beta(l) \phi(w;l,m)$, where $\beta(l)$ and $\phi$ are 
given by (\ref{eq:Lemma8.12}). Since $A$ and $B$ are assumed to
be real we get
$$
|\Sigma_{2}| \leq \sum_{m \in \Z} 
\sum_{\stackrel{l \in \Z}
{m +\frac{B}{2\pi} +\frac{A}{\pi} l > \frac{|A|}{\pi}\eta}}
\left| \Res _{z=i\eta} \left\{ \phi(w;l,m) \right\} \right|.
$$
Write $\phi(w;l,m)=\frac{F(w;l,m)}{\sin^{k}(\pi(w-i\eta))}$, 
where
$F(w;l,m) = e^{-\pi \delta (w+l)^{2}}q(w+l)
e^{ir\left(Aw^{2} + 
2\pi\left(m+\frac{B}{2\pi}+\frac{A}{\pi}l\right)w\right)}$.
For $|w-i\eta|<1$ we have
$\frac{1}{\sin^{k}(\pi(w-i\eta))} = 
\sum_{j=-k}^{\infty} c_{j} (w-i\eta)^{j}$,
where the $c_{j}$'s are independent of $l$ and $m$. Moreover,
$F(w;l,m) = \sum_{j=0}^{\infty} a_{j}(w-i\eta)^{j}$ for 
$w \in \C$. Thus 
$$
\Res _{w=i\eta} \left\{ \phi(w;l,m) \right\} = 
\sum_{j=0}^{k-1} a_{j}c_{-1-j}.
$$
We have 
$a_j = \frac{1}{j!} \frac{\dte ^{j} F(w;l,m)}{{\dte w}^{j}}$,
and thus by Leibnitz' formula and (\ref{eq:P1})
\begin{eqnarray*}
a_{j} &=& \frac{1}{j!} 
\sum_{\stackrel{k_{1},k_{2},k_{3},k_{4} \geq 0}
{k_{1}+k_{2}+k_{3}+k_{4}=j}} 
\frac{j!}{k_{1}!k_{2}!k_{3}!k_{4}!}
\left. \frac{\dte^{k_{1}} e^{-\pi \delta (w+l)^{2}}}
{{\dte w}^{k_{1}}} \right|_{w=i\eta} \\
&& \times \left. \frac{\dte ^{k_{2}} 
e^{2\pi i r \left(m + \frac{B}{2\pi} + \frac{A}{\pi} l \right)
w}}{{\dte w}^{k_{2}}} \right|_{w=i\eta}
\left. \frac{\dte^{k_{3}} q(w+l)}{{\dte w}^{k_{3}}} 
\right|_{w=i\eta} \left. \frac{\dte^{k_{4}} 
e^{irAw^{2}}}{{\dte w}^{k_{4}}} \right|_{w=i\eta} \\
&=& \sum_{\stackrel{k_{1},k_{2},k_{3},k_{4} \geq 0}
{k_{1}+k_{2}+k_{3}+k_{4}=j}} 
\frac{M_{k_{4}}(\eta)}{k_{1}!k_{2}!k_{3}!k_{4}!}
e^{-irA\eta^{2}}P_{k_{1}}(i\eta + l,\delta) 
e^{-\pi \delta (i\eta+l)^{2}} \\
&& \times \left(2 \pi i r \left(m + \frac{B}{2\pi} +
\frac{A}{\pi} l \right)\right)^{k_{2}}
e^{-2 \pi r\left(m +\frac{B}{2\pi} +\frac{A}{\pi} l \right) 
\eta} \left. \frac{\dte^{k_{3}} q(w+l)}{{\dte w}^{k_{3}}} 
\right|_{w=i\eta},
\end{eqnarray*}
where $P_{m}(i\eta + l,\delta)$, $m \in \{0,1,2,\ldots\}$,
is a polynomial in $l$ for fixed $\eta>0$ and $\delta>0$ by 
(\ref{eq:P1}), and $M_{k_{4}}(\eta)$ is a constant independent
of $m$ and $l$. We therefore get 
\begin{eqnarray*}
|\Sigma_{2}| &\leq& \sum_{j=0}^{k-1} |c_{-1-j}| 
\sum_{\stackrel{k_{1},k_{2},k_{3},k_{4} \geq 0}
{k_{1}+k_{2}+k_{3}+k_{4}=j}}
\frac{|M_{k_{4}}(\eta)|}{k_{1}!k_{2}!k_{3}!k_{4}!} 
\sum_{m \in \Z}
\sum_{\stackrel{l \in \Z}
{m +\frac{B}{2\pi} +\frac{A}{\pi} l > \frac{|A|}{\pi}\eta}} 
| P_{k_{1}}(i\eta + l,\delta)| \\
&& \times  \left| q^{(k_{3})}(i\eta+l)\right|
e^{-\pi \delta (l^{2}-\eta^{2})} 
\left(2 \pi r \left(m + \frac{B}{2\pi} + 
\frac{A}{\pi} l \right)\right)^{k_{2}}
e^{-2 \pi r\left(m + \frac{B}{2\pi} +
\frac{A}{\pi} l \right) \eta}.
\end{eqnarray*}
By (\ref{eq:qfunction-symmetry1}) and 
(\ref{eq:qfunction-symmetry2}) we have
$\left| q^{(k_{3})}(i\eta+l)\right| \leq 
b_{k_3}e^{a_{k_3}(\eta + m_{q})}$, and for $l \in \Z$ fixed we
have that
$$
\sum_{\stackrel{m \in \Z}
{m +\frac{B}{2\pi} +\frac{A}{\pi} l > \frac{|A|}{\pi}\eta}}
\left(2 \pi r \left(m + \frac{B}{2\pi} +
\frac{A}{\pi} l \right)\right)^{k_{2}}
e^{-2 \pi r\left(m +\frac{B}{2\pi} +\frac{A}{\pi} l \right) 
\eta} = \sum_{\stackrel{m' \in \Z +\frac{B}{2\pi} +
\frac{A}{\pi} l}{m' > \frac{|A|}{\pi}\eta}} 
(2 \pi r m' )^{k_{2}} e^{-2 \pi r m' \eta}.
$$
is bounded from above by
$$
\beta_{k_{2}}(\eta) = \max_{x \in [0,\infty[ } f(x) + 
\int_{0}^{\infty} f(x) \dte x = 
\left( \frac{k_{2}}{\eta} \right)^{k_{2}} e^{-k_{2}} + 
\frac{k_{2}!}{2\pi r \eta^{k_{2}+1}},
$$
where $f(x)=(2 \pi r x )^{k_{2}} e^{-2 \pi r \eta x}$. But then
\begin{eqnarray*}
|\Sigma_{2}| &\leq& \sum_{j=0}^{k-1} |c_{-1-j}| 
\sum_{\stackrel{k_{1},k_{2},k_{3},k_{4} \geq 0}
{k_{1}+k_{2}+k_{3}+k_{4}=j}}
\frac{|M_{k_{4}}(\eta)|\beta_{k_{2}}(\eta)b_{k_{3}}}
{k_{1}!k_{2}!k_{3}!k_{4}!} e^{a_{k_{3}}(\eta+m_{q})} \\
&& \times \sum_{l \in \Z} | P_{k_{1}}(i\eta + l,\delta)| 
      e^{-\pi \delta (l^{2}-\eta^{2})},
\end{eqnarray*}
where the right-hand side is convergent.
\end{prfa}

\begin{prfa}{\bf of \reflem{lem:Lemma8.13}\hspace{.1in}}
Let $\lambda \in \Lambda$, $l \in \Z$, $\eta \in ]0,1[$, and 
$\delta>0$ be fixed in the following and let 
$\F(y)=\F(y;l,\lambda,\delta,\eta)$ and let 
$C(\rho)=C(\lambda,\rho)$ for all $\rho \in \R$. To show the 
first identity in \reflem{lem:Lemma8.13}, first note that
\begin{eqnarray}\label{eq:H6}
\int_{C(\eta)} \F^{-}(y) \dte y &=& \frac{1}{2} 
\left( \int_{C(\eta)} \F^{-}(y) \dte y - 
\int_{C(-\eta)} \F^{-}(y) \dte y \right) \nonumber \\
&=& \pi i \Res_{y=0} \F^{-}(y) = \pi i \Res_{y=0} \F(y),
\end{eqnarray}
where we in the last equality use that $\F^{+}(y)$ has zero 
residue in $y=0$ as an even function. The first equality simply
follows by the fact that
$\int_{C(-\rho)} f(-y) \dte y = \int_{C(\rho)} f(y) \dte y$
for all $\rho >0$ and every complex function such that 
$\int_{C(\rho)} f(y) \dte y$ exists. By this we have
\begin{eqnarray*}
- \frac{1}{2} \int_{C(-\eta)} \F^{-}(y) \dte y &=& \frac{1}{4}
\left( \int_{C(-\eta)} \F(-y) \dte y  - 
\int_{C(-\eta)} \F(y) \dte y \right) \\
&=& \frac{1}{4} \left( \int_{C(\eta)} \F(y) \dte y  - 
\int_{C(\eta)} \F(-y) \dte y \right) = \frac{1}{2} 
\int_{C(\eta)} \F^{-}(y) \dte y.
\end{eqnarray*}
The second equality in (\ref{eq:H6}) follows by the residue 
theorem and the following estimation (since $y=0$ is the only 
pole for $F^{-}(y)$ between the contours $C(-\eta)$ and 
$C(\eta)$ since $\eta \in ]0,1[$): Let $\tau_{R}(t)=iR+t$, 
$t \in [a,b]$, $R \in \R$, where $a= \sgA R -\eta$ and 
$b= \sgA R +\eta$ (so $\tau_{R}(t)$ is a contour parallel with
the real axes and connecting $C(-\eta)$ and $C(\eta)$). Since
$$
\int_{\tau_{R}} \F(-y) \dte y = \int_{\tau_{-R}} \F(y) \dte y
$$
we only have to show that $|\int_{\tau_{R}} \F(y) \dte y|$ 
converges to zero as $|R|$ converges to infinity. By 
(\ref{eq:normsin}) one finds
$|\sin(\pi\tau_{R}(t))| \geq \frac{1}{4}e^{\pi |R|}$
for $|R|$ sufficiently large. By this and 
(\ref{eq:qfunction-symmetry1}) we therefore get
$$
\left| \int_{\tau_{R}} \F(y) \dte y \right| \leq b_{0} 4^{k} 
e^{-\pi k |R|} e^{a_{0}(|R|+l+\eta)}
\int_{a}^{b} \left| e^{-\pi \delta (iR+t+l+i\eta)^{2}} \right| 
e^{-2rA(R+\eta)t} \dte t
$$
for $|R|$ sufficiently large. Now, for $|R|$ large, $a$ and $b$
have the same sign and
$$
\sup_{t \in [a,b]} e^{-2rA(R+\eta)t} \leq 
e^{-2rA(R+\eta)a} + e^{-2rA(R+\eta)b}
\leq 2 e^{ -2r|A| R^{2}} e^{2 r|A| |R| \eta} 
e^{2r|A| \eta ^{2}}
$$
so
\begin{eqnarray*}
\left| \int_{\tau_{R}} \F(y) \dte y \right| &\leq& 2b_{0} 4^{k}
e^{-\pi k |R|} e^{a_{0}(|R|+l+\eta)} e^{ -2r|A| R^{2}} \\
&& \times \exp ( 2 r|A| |R| \eta ) \exp ( 2r|A| \eta ^{2} )
\int_{a}^{b} \left| e^{-\pi \delta (iR+t+l+i\eta)^{2}} \right| 
\dte t.
\end{eqnarray*}
Here
$\int_{a}^{b} \left| e^{-\pi \delta (iR+t+l+i\eta)^{2}} \right|
\dte t \leq e^{\pi \delta (R+\eta)^{2}} 
\iny e^{-\pi \delta t^{2}} \dte t
= \frac{1}{\sqrt{\delta}} e^{\pi \delta (R+\eta)^{2}}$
so
$$
\left| \int_{\tau_{R}} \F(y) \dte y \right| \leq 
    c e^{d|R|} e^{-(2r|A|-\pi\delta)R^{2}},
$$
where $c$ and $d$ are constants independent of $R$, so the 
right-hand side converges to zero as $|R|$ converges to 
infinity if $\delta \in ]0,2r|A|/\pi[$.

The last claim about the absolute convergency is proved by the
same technique as used in the proof of the absolute convergency
of $\Sigma_{2}$ in \reflem{lem:Lemma8.12}. However, the proof 
is shorter here since we only have a single sum to consider.
\end{prfa}

\subsection{Appendix C. Proofs of \reflem{lem:Lemma8.16} and 
\refprop{prop:Lemma8.17}}

\begin{prfa}{\bf of \reflem{lem:Lemma8.16}\hspace{.1in}}
We will follow the same strategy as in the last part of the
proof of \reflem{lem:Lemma8.12}. In particular we let $c_j$,
$j=-k,-k+1,\ldots$ be as in that proof. Let
\begin{eqnarray*}
\tau(z) &=& \cot(\pi rz) \sum_{\nu \in \Gamma(l)} 
\K_{l}(z;\nu)e^{-2\pi irVz}, \\
\theta(z) &=& \sum_{\nu \in \Gamma(l)} 
\sum_{\stackrel{m\in\Z}{0 \leq m \leq |V|}} \K_{l}(z;\nu)
e^{2\pi ir \sgV (m-|V|)z)}.
\end{eqnarray*}
By (\ref{eq:zero-3}), $\tau$ is analytic in a neighborhood
of $z=0$, in fact on $D(0,1/r)$. For $\eta \in ]0,1/r[$ we have
\begin{eqnarray*}
e^{-\pi \delta (z+l)^{2}} \tau(z) &=& \sum_{j=0}^{\infty} 
a_{j}(l,\eta)(z-i \eta)^{j}, \\
e^{-\pi \delta (z+l)^{2}} \theta(z) &=& \sum_{j=0}^{\infty} 
b_{j}(l,\eta)(z-i \eta)^{j}
\end{eqnarray*}
for $z$ in a neighborhood of $i\eta$. We thus have
\begin{eqnarray*}
&&\left|\Res_{z=i\eta} \left\{ e^{-\pi \delta (z+l)^{2}} 
\frac{\tau(z)}{\sin^{k}(\pi(z-i \eta))}\right\}\right|
\leq \sum_{j=0}^{k-1} |c_{-1-j}||a_{j}(l,\eta)| \\
&\leq& \sum_{j=0}^{k-1} |c_{-1-j}| \frac{1}{j!} \sum_{n=0}^{j}
\binom{j}{n} \left| P_{n}(l+i\eta,\delta) \right|
e^{-\pi \delta (l^{2}-\eta^{2}) } \left| \left. 
\left( \frac{\dte}{\dte z} \right)^{j-n} \tau(z) 
\right|_{z=i\eta} \right|
\end{eqnarray*}
and
\begin{eqnarray*}
&&\left|\Res_{z=i\eta} \left\{ e^{-\pi \delta (z+l)^{2}} 
\frac{\theta(z)}{\sin^{k}(\pi(z-i \eta))}\right\}\right|
\leq \sum_{j=0}^{k-1} |c_{-1-j}||b_{j}(l,\eta)| \\
&\leq& \sum_{j=0}^{k-1} |c_{-1-j}| \frac{1}{j!} \sum_{n=0}^{j}
\binom{j}{n} \left| P_{n}(l+i\eta,\delta) \right|
e^{-\pi \delta (l^{2}-\eta^{2}) } \left| \left. 
\left( \frac{\dte}{\dte z} \right)^{j-n} \theta(z) 
\right|_{z=i\eta} \right|.
\end{eqnarray*}
Now let $\eta_{2} \in ]0,1/r[$ be fixed and observe that
$A_{j} = \sup_{l \in \Z, \eta \in [0,\eta_{2}]} \left| \left. 
\left( \frac{\dte}{\dte z} \right)^{j} \tau(z) 
\right|_{z=i\eta} \right|$
and
$B_{j} = \sup_{l \in \Z, \eta \in [0,\eta_{2}]} \left| \left. 
\left( \frac{\dte}{\dte z} \right)^{j} \theta(z) 
\right|_{z=i\eta} \right|$
are both finite since both $\tau$ and $\theta$ are periodic in
$l$ with a period of $L$ by the assumptions on $\Gamma(l)$ and 
$K_{l}$. By the above we get
$$
|Z^{l}(\delta,\eta)| \leq \pi e^{\pi\delta\eta_{2}^{2}}
\sum_{j=0}^{k-1} |c_{-1-j}| \frac{1}{j!} \sum_{n=0}^{j} 
\binom{j}{n} (A_{j-n}+2B_{j-n})
e^{-\pi \delta l^{2} }\left| P_{n}(l+i\eta,\delta) \right|.
$$
By using the explicit expression for the polynomials $P_n$, 
see (\ref{eq:P1}), we can now appeal to Weierstrass' test for
uniform convergence. Note finally that
\begin{eqnarray*}
&& \lieta \Res_{z=i\eta} \left\{ e^{-\pi \delta (z+l)^{2}} 
\frac{\tau(z)}{\sin^{k}(\pi(z-i \eta))}\right\}
= \sum_{j=0}^{k-1} c_{-1-j}\lieta a_{j}(l,\eta) \\
&=& \sum_{j=0}^{k-1} c_{-1-j}a_{j}(l,0)
= \Res_{z=0} \left\{ e^{-\pi \delta (z+l)^{2}} 
\frac{\tau(z)}{\sin^{k}(\pi z)}\right\}
\end{eqnarray*}
since 
$a_{j}(l,\eta) = \left.\frac{1}{j!} 
\left( \frac{\dte}{\dte z} \right)^{j}
e^{-\pi \delta (z+l)^{2}} \tau(z) \right|_{z=0}$,
and similarly with the other residue term in 
$Z^{l}(\delta,\eta)$.
\end{prfa}

\noindent For the proof of \refprop{prop:Lemma8.17} we need the
following small lemma which we state for the sake of easy 
reference.

\begin{lem}\label{lem:Lemma8.11}
Let $P(z)=Az^{2}+Bz+C$ be a polynomial with complex 
coefficients such that $\rea(A)>0$. Then
$$
\lidel \sdel \sum_{k \in \Z} |k\delta|^{l} e^{-\delta P(k)} = 0
$$
for all $l \in \{1,2,\ldots\}$ and
$$
\lidel \sdel \sum_{k \in \Z} \delta e^{-\delta P(k)} = 0.
$$
\end{lem}

\noindent The first statement follows by comparing the
sums by integrals of the form 
$\int_{0}^{\infty} x^{l} e^{-\delta a(x+d)^{2}} \dte x$.
The second statement is a simple corollary of the periodicity 
lemma, \reflem{lem:periodicity}.

\begin{prfa}{\bf of \refprop{prop:Lemma8.17}\hspace{.1in}}
Let us use the notation from the proof of 
\reflem{lem:Lemma8.16}. Thus we have
$Z_{\pol}(r) = r\lidel \sdel \sum_{l\in \Z} Z^{l}(\delta)$,
where
$$
Z^{l}(\delta) = 2\pi i e^{-\pi \delta l^{2}} \sum_{j=0}^{k-1} 
c_{-1-j}\frac{1}{j!} \sum_{n=0}^{j} \binom{j}{n} 
P_{n}(l,\delta) \left. \left( \frac{\dte}{\dte z} \right)^{j-n}
\left(\frac{i}{2} \tau(z) + \theta(z) \right)\right|_{z=0}.
$$
For $j \in \{0,1,\ldots,k-1\}$ and $n \in \{0,1,\ldots,j\}$ we 
let
$$
C(j,n,\delta) = 2\pi i c_{-1-j}\frac{1}{j!}\binom{j}{n} 
\sum_{l \in\Z} P_{n}(l,\delta)e^{-\pi \delta l^{2}}
\left. \left( \frac{\dte}{\dte z} \right)^{j-n} 
\left(\frac{i}{2} \tau(z) + \theta(z) \right)\right|_{z=0}.
$$
Then
$$
|C(j,n,\delta)| \leq 2\pi |c_{-1-j}|\frac{1}{j!}\binom{j}{n} 
\left( \frac{1}{2} A_{j} + B_{j} \right)\sum_{l \in\Z} 
|P_{n}(l,\delta)|e^{-\pi \delta l^{2}}.
$$
By \reflem{lem:Lemma8.11} and (\ref{eq:P1}) we therefore get
\begin{eqnarray*}
Z_{\pol}(r) &=& r\lidel \sdel \sum_{j=0}^{k-1} \sum_{n=0}^{j} 
C(j,n,\delta) =  r\lidel \sdel \sum_{j=0}^{k-1} 
C(j,0,\delta) \\
&=& 2\pi ir\lidel \sdel \sum_{l \in\Z} \e^{-\pi \delta l^{2}}
\sum_{j=0}^{k-1} c_{-1-j}\frac{1}{j!} \left. 
\left( \frac{d}{dz} \right)^{j} 
\left(\frac{i}{2} \tau(z) + \theta(z) \right)\right|_{z=0} \\
&=& \frac{2\pi ir}{L} \sum_{l =0}^{L-1} \sum_{j=0}^{k-1} 
c_{-1-j}\frac{1}{j!} \left. \left( \frac{d}{dz} \right)^{j} 
\left(\frac{i}{2} \tau(z) + \theta(z) \right)\right|_{z=0},
\end{eqnarray*}
where the last identity follows by \reflem{lem:periodicity} and
the fact that $\tau$ and $\theta$ are periodic in $l$ with a 
period of $L$.
\end{prfa}

\subsection{Appendix D. Proofs of \reflem{lem:Lemma8.18}, 
\reflem{lem:Lemma8.21} and \refprop{prop:Lemma8.22}}

\begin{prfa}{\bf of \reflem{lem:Lemma8.18}\hspace{.1in}}
By following the first part of the proof of 
\reflem{lem:Lemma8.7} leading to (\ref{eq:ulighedLemma8.7}) and
by the remarks before the proof of \reflem{lem:Lemma8.7} we 
find that
\begin{eqnarray*}
\left| \int_{C_{\sd}(m,\lambda)} f(z;m,\lambda) \dte z \right|
&\leq& \frac{b_{0}}{M^{k}} \sqrt{\frac{\pi}{r|A|}} 
e^{-\pi \delta z_{\st}^{2}(m,\lambda)} 
e^{a_{0}|z_{\st}(m,\lambda)|} \\
&& \times 
\left( \exp\left( \frac{a_{+}^{2}(m,\lambda)}{4r|A|}\right)
+ \exp\left(\frac{a_{-}^{2}(m,\lambda)}{4r|A|}\right) \right)
\end{eqnarray*}
for all $\eta \in [0,\eta_{1}]$ and all $(m,\lambda) \in W$.
Here $a_{0}$ and $b_{0}$ are the constants from 
(\ref{eq:qfunction-symmetry1}) and
$a_{\pm}(m,\lambda) = \pi \delta \sqrt{2} z_{\st}(m,\lambda) 
\pm a_{0}$. 
By following the first part of the proof of 
\reflem{lem:Lemma8.12} we find that
$$
\sum_{(m,\lambda) \in W} \frac{|g(\lambda)|}{\sqrt{|A|}} 
e^{-\pi \delta z_{\st}^{2}(m,\lambda)}
e^{a_{0}|z_{\st}(m,\lambda)|}
\left( \exp\left( \frac{a_{+}^{2}(m,\lambda)}{4r|A|}\right) 
+ \exp\left(\frac{a_{-}^{2}(m,\lambda)}{4r|A|}\right) \right)
$$
is convergent for $\delta \in ]0,2rA_{0}/\pi[$, where $A_{0}$ 
is given by (\ref{eq:A0}).
\end{prfa}

\begin{prfa}{\bf of \reflem{lem:Lemma8.21}\hspace{.1in}}
Keep $N \in \Z_{\geq 0}$ fixed in the following. For 
$\lambda \in \Lambda$, let 
$W_{\lambda} = \{ m \in \Z \;|\; (m,\lambda) \in W \}$. We write 
$z_{\st}$ for $z_{\st}(m,\lambda)$. Let us start by showing 
that $Z_{\inte,1}(N;\delta)$ and $\Sigma^{1}(N;\delta)$ are 
absolutely convergent. By Cauchy's integral formula for the
derivatives we have
$$
a_{j}(\delta) = \frac{1}{j!} \kappa^{(j)}(z_{\st};\delta) 
= \frac{1}{2\pi i} \int_{|z-z_{\st}|=3\rho} 
\frac{\kappa(z;\delta)}{(z-z_{\st})^{j+1}} \dte z.
$$
Let
$L(\rho,\delta) = \sup_{|z-z_{\st}|=3\rho} |\kappa(z;\delta)|$.
Then
\begin{equation}\label{eq:adeltaestimate}
|a_{j}(\delta)| \leq (3\rho)^{-j} L(\rho,\delta).
\end{equation}
By (\ref{eq:kappafunction}) and (\ref{eq:qestimation}) we get
$L(\rho,\delta) \leq C_{\rho} \sup_{|z|=3\rho} 
\left| e^{-\pi\delta (z+z_{\st})^{2}}
\sin^{-k}(\pi (z+z_{\st})) \right|$.
Since $z_{\st} \in \R$ we get by (\ref{eq:normsin}) that
\begin{eqnarray*}
|\sin(\pi (z+z_{\st}))| &=& \frac{1}{2}
\sqrt{ e^{2\pi\ima (z)} + e^{-2\pi\ima(z)} - 
2\cos(2\pi(z_{\st}+\rea(z)))} \\
&\geq& \frac{1}{\sqrt{2}} 
\sqrt{ 1 - \cos(2\pi(z_{\st}+\rea(z)))}.
\end{eqnarray*}
By our choice of $\rho$ we have that
$z_{\st}(m,\lambda) + \rea(z) \notin \Z + [-\rho,\rho]$
for all $(m,\lambda) \in W$ and all $z \in \C$ with 
$|z|=3\rho$, so
$|\sin(\pi (z+z_{\st}(m,\lambda)))| \geq M_{\rho}
:= \sqrt{ 1 - \cos(2\pi\rho)} >0$
for all $(m,\lambda) \in W$ and all $z \in \C$ with 
$|z|=3\rho$. Therefore
$L(\rho,\delta) \leq C_{\rho}M_{\rho}^{-k} 
\sup_{|z|=3\rho} \left| e^{-\pi\delta (z+z_{\st})^{2}} 
\right|$.
Now use that
\begin{eqnarray*}
\sup_{|z|=3\rho} \left| e^{-\pi\delta (z+z_{\st})^{2}} \right| 
&=& e^{-\pi\delta z_{\st}^{2}} \sup_{|z|=3\rho} 
\left| e^{-\pi\delta (z^{2}+2z_{\st}z)} \right| \\
&\leq& e^{-\pi\delta z_{\st}^{2}} 
\left( \sup_{|z|=3\rho}e^{-\pi\delta \rea(z^{2})} \right)
\left( \sup_{|z|=3\rho}e^{-2\pi\delta z_{\st}\rea(z)}\right) \\
&=& e^{-\pi\delta z_{\st}^{2}} e^{9\pi\delta\rho^{2}} 
e^{6\pi\delta |z_{\st}| \rho}.
\end{eqnarray*}
Letting 
$C_{\rho}' = C_{\rho}M_{\rho}^{-k} e^{9\pi\delta\rho^{2}}$ 
and $d=6\pi\rho$ we see that
\begin{equation}\label{eq:Lestimate}
L(\rho,\delta) \leq C_{\rho}' e^{d\delta|z_{\st}|} 
e^{-\pi\delta z_{\st}^{2}}.
\end{equation}
By (\ref{eq:adeltaestimate}) and (\ref{eq:Lestimate}) we 
therefore get
\begin{eqnarray*}
\sum_{m \in W_{\lambda}}
\sum_{\stackrel{0 \leq j \leq N}{j \in 2\Z}} \left| 
a_{j}(\delta) \left( \frac{i}{rA} \right)^{\frac{j+1}{2}} 
\Gamma\left( \frac{j+1}{2} \right) \right|
&\leq& C_{N}^{(0)}\Sigma(\lambda;\delta), \\
\sum_{ m \in W_{\lambda}}
\sum_{\stackrel{0 \leq j \leq N}{j \in 2\Z}} \left|
a_{j}(\delta) \left( \frac{i}{rA} \right)^{\frac{j+1}{2}}
\Gamma\left( \frac{j+1}{2},4|A|\rho^{2} r \right) \right|
&\leq& C_{N}^{(1)}\Sigma(\lambda;\delta),
\end{eqnarray*}
where
\begin{equation}\label{eq:lambdasum}
\Sigma(\lambda;\delta) = 
\sum_{m \in \Z} e^{-\pi\delta z_{\st}(m)^{2}}
\left(e^{d\delta z_{\st}(m)} + e^{-d\delta z_{\st}(m)}\right) 
\end{equation}
with $z_{\st}(m)=z_{\st}(m,\lambda)$, and where
\begin{eqnarray}\label{eq:CN1}
C_{N}^{(0)} &=& C_{\rho}' 
\sum_{\stackrel{0 \leq j \leq N}{j \in 2\Z}}
\left( rA_{0} \right)^{-\frac{j+1}{2}} (3\rho)^{-j} 
\Gamma\left( \frac{j+1}{2} \right), \nonumber \\
C_{N}^{(1)} &=& C_{N}^{(1)}(r) = C_{\rho}' 
\sum_{\stackrel{0 \leq j \leq N}{j \in 2\Z}}
\left( rA_{0} \right)^{-\frac{j+1}{2}} (3\rho)^{-j} 
\Gamma\left( \frac{j+1}{2},4A_{0}\rho^{2}r \right)
\end{eqnarray}
are independent of $\lambda$. Since $z_{\st}$ is a real affine
expression in $m$ for $\lambda \in \Lambda$ fixed, the sum
$\Sigma(\lambda;\delta)$ is convergent for all $\delta>0$.

Next we show that $\Sigma^{2}_{\nu}(N;\delta)$ are absolutely 
convergent. By (\ref{eq:kappapowerseries}) and 
(\ref{eq:kappaR}) we have
$$
R_{N}(z;\delta) = \sum_{j=N+1}^{\infty} a_{j}(\delta)
(z-z_{st})^{j}
$$
for $z \in D(z_{\st},4\rho)$. By (\ref{eq:adeltaestimate}) we 
have for $|z-z_{\st}| \leq 2\rho$ that
$$
|R_{N}(z;\delta)| \leq 3L(\rho,\delta) 
\left( \frac{|z-z_{\st}|}{3\rho} \right)^{N+1}.
$$
Therefore
\begin{eqnarray*}
|E_{\nu}(N;\delta)| &\leq& \int_{0}^{2\rho} \left| 
R_{N}\left(\gamma\left( (-1)^{\nu}t\right);\delta\right)\right|
\exp\left(-r|A|t^{2}\right)\dte t \\
&\leq& 3L(\rho,\delta)(3\rho)^{-N-1} \int_{0}^{\infty} 
t^{N+1} \exp\left(-r|A|t^{2}\right) \dte t \\
&=& \frac{3}{2} \Gamma\left( \frac{N+2}{2} \right) 
(3\rho)^{-N-1}(r|A|)^{-(N+2)/2} L(\rho,\delta).
\end{eqnarray*}
Note that this estimate is independent of $\nu \in \{ 0,1 \}$ 
and that the only quantity depending on $m$ is 
$L(\rho,\delta)$. Therefore we have for $\nu \in \{ 0,1 \}$ 
that
$$
\sum_{ m \in W_{\lambda} } \left|E_{\nu}(N;\delta)\right|
\leq C_{N}' r^{-(N+2)/2} \sum_{ m \in W_{\lambda} } 
L(\rho,\delta),
$$
where
$C_{N}' = \frac{1}{2A_{0}^{(N+2)/2}}(3\rho)^{-N-1}
\Gamma\left( \frac{N+2}{2} \right)$.
By (\ref{eq:Lestimate}) we therefore get
\begin{equation}\label{eq:errorestimate2}
\sum_{ m \in W_{\lambda} } \left|E_{\nu}(N;\delta)\right| 
\leq C_{N}^{(2)} r^{-(N+2)/2} \Sigma(\lambda;\delta)
\end{equation}
for $\nu \in \{ 0,1 \}$, where $\Sigma(\lambda;\delta)$ is the 
convergent sum in (\ref{eq:lambdasum}) and
\begin{equation}\label{eq:CN2}
C_{N}^{(2)} = C_{\rho}'C_{N}' = 
\frac{C_{\rho}}{2M_{\rho}^{k}A_{0}^{(N+2)/2}} 
e^{9\pi\delta\rho^{2}}(3\rho)^{-N-1}
\Gamma\left( \frac{N+2}{2} \right).
\end{equation}
Note that the right-hand side of (\ref{eq:errorestimate2}) is
independent of $\nu$ and $C_{N}^{(2)}$ is independent of 
$\lambda$.

Let us finally consider the sums $\Sigma_{\nu}^{3}$. We begin
by estimating the integrals
$$
J_{\nu}(\delta) = \int_{2\rho}^{\infty} 
\psi\left( (-1)^{\nu}t;\delta \right) 
\exp\left(-r|A|t^{2}\right) \dte t.
$$
By the remarks above the proof of \reflem{lem:Lemma8.7} we have
$$
|\psi(t;\delta)| \leq \frac{1}{M^{k}} 
\left| e^{-\pi\delta \gamma(t)^{2}} q(\gamma(t)) \right|
= \frac{1}{M^{k}} e^{-\pi\delta z_{\st}^{2}} 
e^{-\sqrt{2}\pi\delta z_{\st} t} |q(\gamma(t))|.
$$
By our periodicity assumptions on $q$ and $z_{\st}$ we have
$$
q(z+z_{\st}(m+m_{q}H,\lambda);\lambda) =
q(z+z_{\st}(m,\lambda);\lambda)
$$
for all $(m,\lambda) \in \Z \times \Lambda$ and all $z \in \C$.
Combining this with (\ref{eq:qfunction-symmetry1}) we get
$$
\left|q\left(\gamma\left((-1)^{\nu}t \right)\right)\right| 
\leq b_{0} e^{a_{0}K} e^{a_{0}|t|}
$$
for all $(m,\lambda) \in \Z \times \Lambda$ and $t \in \R$, 
where
$$
K = \max \left\{ \; |z_{\st}(m,\lambda)| \; \left| \; 
m=0,1,\ldots,m_{q}H-1, \hspace{.1in}\lambda \in \Lambda \; 
\right\} \right..
$$ 
We therefore have the estimate
$$
|J_{\nu}(\delta)| \leq K_{1} e^{-\pi\delta z_{\st}^{2}} 
\int_{2\rho}^{\infty} 
e^{(a_{0}-(-1)^{\nu}\sqrt{2}\pi\delta z_{\st})t} 
e^{-r|A|t^{2}} \dte t,
$$
where $K_{1}=\frac{b_{0}}{M^{k}} e^{a_{0}K}$. Here
$e^{-(-1)^{\nu}\sqrt{2}\pi\delta z_{\st}t} = 
e^{(-1)^{\nu}\pi \frac{B}{\sqrt{2}A} \delta t} 
e^{(-1)^{\nu} \frac{\pi^{2}}{\sqrt{2}A} \delta m t}$.
Choose $r_{\lambda} \in \Z_{\geq 2}$ such that
$e^{a_{0}t} \leq e^{\frac{1}{2}r|A|t^{2}}$
for all $t \geq 2\rho$ and all $r \geq r_{\lambda}$, and choose
$\delta_{\lambda}>0$ such that
$e^{(-1)^{\nu}\pi \frac{B}{\sqrt{2}A} \delta t} \leq 
e^{\frac{1}{4}r|A|t^{2}}$
for all $t \geq 2\rho$, $r \geq r_{\lambda}$ and
$\delta \in ]0,\delta_{\lambda}]$. Then
$$
|J_{\nu}(\delta)| \leq K_{1} e^{-\pi\delta z_{\st}^{2}} 
\int_{2\rho}^{\infty} e^{c_{\delta} |m| t} e^{-at^{2}} \dte t 
$$
for $r \geq r_{\lambda}$ and 
$\delta \in ]0,\delta_{\lambda}]$, where 
$c_{\delta}=\frac{\pi^{2}}{\sqrt{2}|A|} \delta$ and 
$a=\frac{1}{4}r|A|$. Now use that 
$e^{c_{\delta}|m|t} \leq e^{\frac{a}{2} t^{2}}$ if and only if 
$2c_{\delta}|m|/a \leq t$. If $2\rho \geq 2c_{\delta}|m|/a$ we
therefore get
$$
\int_{2\rho}^{\infty} e^{c_{\delta} |m| t} e^{-at^{2}} \dte t
\leq \int_{2\rho}^{\infty} e^{-\frac{a}{2}t^{2}} \dte t
= \sqrt{\frac{2}{|A|}} r^{-1/2} 
\Gamma\left( \frac{1}{2}, \frac{|A|\rho^{2}}{2}r \right).
$$
For $2\rho < 2c_{\delta}|m|/a$ we have
$$
\int_{2\rho}^{\infty} e^{c_{\delta} |m| t} e^{-at^{2}} \dte t
\leq \int_{2\rho}^{2c_{\delta}|m|/a} e^{c_{\delta} |m| t} 
e^{-at^{2}} \dte t + \int_{2c_{\delta}|m|/a}^{\infty} 
e^{-\frac{a}{2}t^{2}} \dte t.
$$
Here
$$
\int_{2c_{\delta}|m|/a}^{\infty} e^{-\frac{a}{2}t^{2}} \dte t
\leq \int_{2\rho}^{\infty} e^{-\frac{a}{2}t^{2}} \dte t
= \sqrt{\frac{2}{|A|}} r^{-1/2} 
\Gamma\left( \frac{1}{2}, \frac{|A|\rho^{2}}{2}r \right),
$$
and
\begin{eqnarray*}
&&\int_{2\rho}^{2c_{\delta}|m|/a} e^{c_{\delta} |m| t} 
e^{-at^{2}} \dte t \leq e^{-4a\rho^{2}} 
\int_{2\rho}^{2c_{\delta}|m|/a} e^{c_{\delta} |m| t} \dte t \\
&\leq& e^{-r|A|\rho^{2}} \int_{0}^{2c_{\delta}|m|/a} 
e^{c_{\delta} |m| t} \dte t
\leq \frac{2c_{\delta}|m|}{a} 
e^{\frac{2}{a} c_{\delta}^{2} m^{2}} e^{-r|A|\rho^{2}}.
\end{eqnarray*}
We therefore get
$$
\sum_{ m \in W_{\lambda} } \left|J_{\nu}(\delta) \right|
\leq K_{1} \sum_{m \in W_{\lambda}} e^{-\pi\delta z_{\st}^{2}} 
\left( \sqrt{\frac{2}{r|A|}} 
\Gamma\left( \frac{1}{2}, \frac{|A|\rho^{2}}{2}r \right) 
+ \frac{2c_{\delta}|m|}{a} e^{\frac{2}{a} c_{\delta}^{2}m^{2}} 
e^{-r|A|\rho^{2}} \right)
$$
for $r \geq r_{\lambda}$ and $\delta \in ]0,\delta_{\lambda}]$.
Here
$$
\frac{2c_{\delta}|m|}{a} e^{\frac{2}{a} c_{\delta}^{2}m^{2}} 
e^{-\pi\delta z_{\st}^{2}}
= \frac{8\pi^{2}}{\sqrt{2}A^{2}} 
e^{-\pi\delta \left( \frac{B}{2A} \right)^{2}} \frac{1}{r} 
\delta |m|
e^{-\pi\delta \left( \frac{\pi}{A} \right)^{2}
\left( 1 - \frac{4\pi}{r|A|} \delta\right) m^{2}}
e^{-\delta \left( \frac{\pi}{A}\right)^{2} Bm}.
$$
Finally choose $\delta_{\lambda}>0$ so small that
$\frac{4\pi}{rA} \delta \leq \frac{1}{2}$
for all $r \in \Z_{\geq 2}$ and all 
$\delta \in ]0,\delta_{\lambda}]$. Then
$$
\sum_{m \in \Z} \frac{2c_{\delta}|m|}{a} 
e^{\frac{2}{a} c_{\delta}^{2}m^{2}} e^{-\pi\delta z_{\st}^{2}}
\leq K_{2}\frac{1}{r} \sum_{m \in \Z} \delta |m| 
e^{-\delta P(m;\lambda)},
$$
where $K_{2}=\frac{8\pi^{2}}{\sqrt{2}A_{0}^{2}}$ and
$P(m;\lambda) = \frac{\pi^{3}}{2A^{2}} m^{2} + 
\left( \frac{\pi}{A} \right)^{2} B m$.
We therefore get that
$\sum_{ m \in W_{\lambda} } \left| J_{\nu}(\delta) \right| 
\leq A_{3}(\lambda;\delta)$
for $\nu \in \{0,1\}$, $r \geq r_{\lambda}$ and 
$\delta \in ]0,\delta_{\lambda}]$, where
\begin{equation}\label{eq:A3term}
A_{3}(\lambda;\delta) = K_{1}
\left(\sqrt{\frac{2}{A_{0}r}} 
\Gamma\left( \frac{1}{2}, \frac{A_{0}\rho^{2}}{2}r \right) 
\sum_{m \in \Z} e^{-\pi\delta z_{st}^{2}}
+ \frac{K_{2}}{r} e^{-rA_{0}\rho^{2}} \sum_{m \in \Z} 
\delta |m| e^{-\delta P(m;\lambda)}\right)
\end{equation}
is convergent. Finally put
$r_{0}=\max\{ \; r_{\lambda} \; | \; \lambda \in \Lambda \;\}$ 
and
$\delta_{0} = \min\{ \; \delta_{\lambda} \; | \; \lambda \in 
\Lambda \;\}$. 
\end{prfa}

\begin{prfa}{\bf of \refprop{prop:Lemma8.22}\hspace{.1in}}
The idea of the proof is to use the periodicity lemma,
\reflem{lem:periodicity}. To this end we have to show that the 
summand in the infinite sum over $m$ is periodic in $m$ (for 
each fixed $\lambda \in \Lambda$). Keep $\lambda \in \Lambda$ 
fixed. First note that
$$
irAz_{\st}(m)^{2}= \pi i r \frac{P}{H} m^{2} + 
Bir \frac{P}{H}m + c,
$$
where $c$ is a constant independent of $m$. By (\ref{eq:B}) is
follows that $\exp\left(-irAz_{\st}(m)^{2}\right)$ is periodic 
in $m$, lets say with a period of $N_{\lambda}$. By 
(\ref{eq:kappapowerseries}) and (\ref{eq:P1}) we get
$$
a_{j}(\delta) = \frac{1}{j!} \left. 
\frac{\dte^{(j)}\kappa(z;\delta)}{\dte z^{j}} 
\right|_{z=z_{\st}}
= \frac{1}{j!} \sum_{n=0}^{j} 
\left( \begin{array}{c} j \\ n \end{array} \right) 
P_{n}(z_{\st},\delta) e^{-\pi\delta z_{\st}^{2}} \left. 
\frac{\dte^{(j-n)}h(z)}{\dte z^{j-n}} \right|_{z=z_{\st}},
$$
where $h(z)=h(z;\lambda)=q(z;\lambda)/\sin^{k}(\pi z)$. Put
$$
C_{j,n,\lambda}(\delta) = \msumd g(\lambda)
\exp\left(-irAz_{\st}^{2}\right)P_{n}(z_{\st},\delta)
e^{-\pi\delta z_{\st}^{2}} \left. 
\frac{\dte^{(j-n)}h(z)}{\dte z^{j-n}} \right|_{z=z_{\st}}.
$$
By (\ref{eq:qfunction-symmetry2}) $h$ is periodic in $z$ with 
a period of $2m_{q}$, so the derivatives of $h$ in $z_{\st}$ 
are periodic in $m$ with a period of $2m_{q}H$. We conclude 
that the functions $m \mapsto G_{j}(m,\lambda)$, 
$j \in \Z_{\geq 0}$, in (\ref{eq:Gjfunctions}) are periodic 
with a period of $2m_{q}H N_{\lambda}$.
Moreover, there exists for each $\lambda \in \Lambda$ a 
constant $K_{N}(\lambda)$ such that
$$
\sup_{m \in \Z} \left| \left. 
\frac{\dte^{(n)}h(z;\lambda)}{\dte z^{n}} 
\right|_{z=z_{\st}(\lambda,m)} \right| \leq K_{N}(\lambda)
$$
for all $n=0,1,\ldots,N$. But then
$\lidel \sdel |C_{j,n,\lambda}(\delta)|=0$ since
$$
\lidel \sdel \sum_{m \in \Z} |P_{n}(z_{\st},\delta)|
e^{-\pi\delta z_{\st}^{2}} =0
$$
for all $j$ and $n$ with $1 \leq n \leq j \leq N$ by 
\reflem{lem:Lemma8.11}. Note here that
$$
\sum_{m \in \Z} |P_{n}(z_{\st},\delta)|
e^{-\pi\delta z_{\st}^{2}}
\leq \sum_{s=0}^{n} \left(\frac{\pi}{|A|}\right)^{s} 
\sum_{m \in \Z} |a_{s}^{n}(\delta)|\delta^{s}
\left| m + \frac{B}{2\pi} \right|^{s} e^{-\delta P(m)},
$$
where the $a_{s}^{n}(\delta)$ are described below (\ref{eq:P1})
and
$P(m) = \pi z_{\st}(m)^{2} = 
\pi\left(\frac{\pi}{A}m+\frac{B}{2A}\right)^{2}$. Thus
\begin{eqnarray*}
r\lidel\sdel Z_{\inte,1}(N;\delta) &=& r 
\sum_{\stackrel{0 \leq j \leq N}{j \in 2\Z}} 
\frac{1}{j!}\Gamma \left( \frac{j+1}{2} \right)
\left(\frac{i}{r} \right)^{\frac{j+1}{2}} 
\sum_{\lambda \in \Lambda} A^{-\frac{j+1}{2}} g(\lambda) \\
&& \times \lidel\sdel \msumd e^{-\pi\delta z_{\st}^{2}}
G_{j}(m,\lambda)
\end{eqnarray*}
by \reflem{lem:Lemma8.21}.
The identity (\ref{eq:Zintmain1}) now follows by the 
periodicity lemma and the fact that
$\Gamma(m+1/2)/(2m)!=\sqrt{\pi}/(4^{m}m!)$.

Let us next show the result about the remainder term 
$R_{1}(N;\delta)$. We stress again that we do not show that 
$\lidel \sdel R_{1}(N;\delta)$ exists, see the discussion 
below (\ref{eq:Zpolar-lterm1}). By the proof of
\reflem{lem:Lemma8.21} we can let
$$
\vep_{1}(N,\lambda;\delta)
= \left( C_{N}^{(1)} + 2C_{N}^{(2)}r^{-\frac{N+2}{2}} 
\right) \Sigma(\lambda;\delta) + 2 A_{3}(\lambda;\delta),
$$
where $\Sigma(\lambda;\delta)$ is given by 
(\ref{eq:lambdasum}), $A_{3}(\lambda;\delta)$ is given by 
(\ref{eq:A3term}) and $C_{N}^{(1)}$ and $C_{N}^{(2)}$ are given
by respectively (\ref{eq:CN1}) and (\ref{eq:CN2}). We note that
$C_{N}^{(1)}$ and $C_{N}^{(2)}$ are continuous in $\delta$ on 
$[0,\delta_{0}]$. Moreover, since the polynomial $z_{\st}(m)$ 
in $m$ has $(\pi/A)^{2}$ as leading coefficient, we get 
immediately from \reflem{lem:periodicity} that
$$
\lidel\sdel \Sigma(\lambda;\delta) = \frac{A}{\pi}.
$$
By (\ref{eq:A3term}), \reflem{lem:periodicity} and 
\reflem{lem:Lemma8.11} we get
$$
\lidel\sdel A_{3}(\lambda;\delta) = \frac{K_{1}A}{\pi}
\sqrt{\frac{2}{A_{0}r}} 
\Gamma\left( \frac{1}{2}, \frac{A_{0}\rho^{2}}{2}r \right).
$$
We note that $C_{N}^{(2)}$ is independent of $r$.
For $\alpha$ and $u_0$ positive we have the following estimate
$$
\Gamma(\alpha,u_0) = \int_{u_0}^{\infty} u^{\alpha -1}
e^{-u} \dte u \leq e^{-\frac{u_0}{2}} 
\int_{0}^{\infty} u^{\alpha-1} e^{-\frac{u}{2}} \dte u
= 2^{\alpha} \Gamma(\alpha) e^{-\frac{u_{0}}{2}}
$$
finalizing the proof.
\end{prfa}

\subsection{Appendix E. Proofs of \reflem{lem:Lemma8.25} and 
\refprop{prop:Lemma8.26}}

\begin{prfa}{\bf of \reflem{lem:Lemma8.25}\hspace{.1in}}
For the following proof, recall that 
$\beta(l)=\beta(l,\lambda)$ has norm $1$ for all 
$(l,\lambda) \in \Z \times \Lambda$, see (\ref{eq:Lemma8.12}).
We parametrize $C_{\eta}$ by 
$C_{\eta}(t)= \eksgA t -i\eta$, $t \in \R$. Let
$$
\tilde{C}_{\eta}(t) = e^{-\frac{\pi i}{4}\sgA} C_{\eta}(t)
= t - \frac{1}{\sqrt{2}}\sgA \eta - \frac{i}{\sqrt{2}}\eta,
$$
$t \in \R$, and get
$$
I(N;l,\delta,\eta) = 
\sum_{\stackrel{\nu=0}{\nu -k \in 2 \Z}}^{N}
d_{\nu}(l,\delta,\eta)\left( \eksgA \right)^{\nu-k+1}
\int_{\tilde{C}_{\eta}} (z^{2})^{\frac{\nu -k}{2}} 
e^{-r|A|z^{2}} \dte z.
$$
For $s \in \C$ and $z \in \C$ we let 
$z^{s}=\exp(s\log(z))$, where 
$\log(re^{i\theta}) = \log(r) + i\theta$, $r>0$, 
$\theta \in [0,2\pi[$, with $\log(r) \in \R$ being the usual
(principal) logarithm of a positive real number. 
In particular, $r^{s}>0$ for $r \in ]0,\infty[$ and 
$s \in \Z$. Let us also in the following text reserve the 
symbol $\sqrt{\cdot}$ for the continuous square root on 
$\C \sm ]-\infty,0]$, being positive on the positive real axes.
Using these conventions we have
$$
\int_{\tilde{C}_{\eta}} (z^{2})^{\frac{\nu -k}{2}}
e^{-r|A|z^{2}} \dte z
= - \frac{1}{2} \left(\frac{1}{r|A|}\right)^{\frac{\nu -k}{2}}
\frac{1}{\sqrt{r|A|}}
\int_{\chi_{\eta}} z^{\frac{\nu -k}{2} -\frac{1}{2}} 
e^{-z} \dte z,
$$
where $\chi_{\eta}(t) = r|A|\tilde{C}_{\eta}(t)^{2}$
is shown on Figure~\ref{fig-1}.

\begin{figure}

\begin{center}
\begin{texdraw}
\drawdim{cm}
\linewd 0.01

\move(-4 0) \avec(6 0)
\move(0 -4) \avec(0 4)
\linewd 0.02
\move(0 0) \larc r:.5 sd:100 ed:260
\move(-0.09 0.49) \lvec(5.2 1.0)
\move(-0.09 -0.49) \lvec(5.2 -1.0)
\move(2 0.69) \lvec(2.4 0.9)
\move(2 0.69) \lvec(2.4 0.5)
\move(2 -0.69) \lvec(1.6 -0.45)
\move(2 -0.69) \lvec(1.6 -0.85)
\move(1.5 -1.0) \textref h:C v:C \htext{$\chi_{\eta}$}
\linewd 0.01
\move(-0.5 0.3) \lvec(-0.5 -0.3)
\move(-0.7 -0.8) \htext{$-\frac{r|A|\eta^{2}}{2}$}

\end{texdraw}
\end{center}

\caption{}\label{fig-1}
\end{figure}

We have the following integral representation for the gamma 
function due to Hankel valid for all $s \in \C \sm \Z$:
$$
\Gamma(s+1) = \frac{e^{-i\pi s}}{2i\sin(\pi s)}
\int_{\infty}^{(0+)} z^{s} e^{-z} \dte z.
$$
Here the symbol $\int_{\infty}^{(0+)}$ means that we integrate
along a contour which originates at $+\infty$, runs in toward 
the origin just above the real axis, circles the origin once 
counterclockwise, and then returns to $+\infty$ just below the
real axes, see the contour $\Gamma$ on Figure~\ref{fig-2}.

\begin{figure}
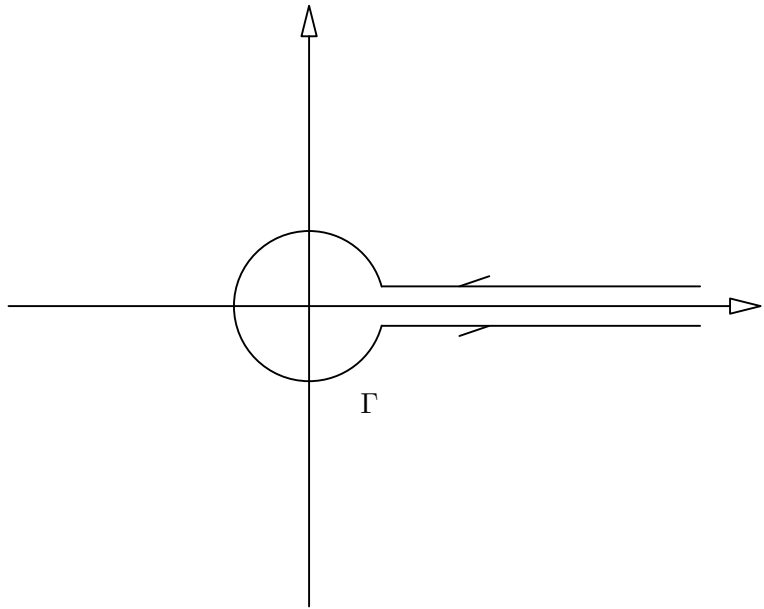


\begin{center}
\begin{texdraw}
\drawdim{cm}
\move(-4 0) \avec(6 0)
\move(0 -4) \avec(0 4)
\move(0 0) \larc r:1 sd:15 ed:345
\move(0.97 0.26) \lvec (5.2 0.26)
\move(0.97 -0.26) \lvec (5.2 -0.26)
\move(2 0.26) \lvec (2.4 0.4)
\move(2.4 -0.26) \lvec (2.0 -0.4)
\move(0.8 -1.3) \textref h:C v:C \htext{$\Gamma$}

\end{texdraw}
\end{center}

\caption{Contour of integration for representation of the gamma
function.}\label{fig-2}
\end{figure}

We have
$$
\Gamma(s+1) = \frac{e^{-i\pi s}}{2i\sin(\pi s)}
\int_{\chi_{\eta}} z^{s} e^{-z} \dte z
$$
for all $s \in \R \sm \Z$ and all fixed $\eta >0$. To see this,
note that the integrand is holomorphic on $\C \sm [0,\infty[$.
By Cauchy's theorem we therefore only have to show that
$\int_{\varsigma_{R}^{\pm}} z^{s} e^{-z} \dte z$ converges to 
zero as $R$ converges to infinite, where
$\varsigma_{R}^{\pm}(t) = R + it$, 
$t \in [0 , t_{\eta}^{\pm}(R)]$, $R>0$, and
$$
t_{\eta}^{\pm}(R)= \frac{1}{\sqrt{2}} \left( 
\eta \sgA \pm \sqrt{\eta^{2} + \frac{2R}{r|A|}} \; \right)
$$
are the two $t$--values for which the line $R+it$ intersects 
the contour $\chi_{\eta}$. But this follows by the estmate
$$
\left| \int_{\varsigma_{R}^{\pm}} z^{s} e^{-z} dz \right| \leq
 e^{-R} \left| \int_{0}^{t_{\eta}^{\pm}(R)} |R+it|^{s} \dte t 
\right| \leq \left| t_{\eta}^{\pm}(R) \right|
\left( R^{s}+(R+|t_{\eta}^{\pm}(R)|)^{s} \right) e^{-R}.
$$ 
Therefore
\begin{equation}\label{eq:Int2}
I(N;l,\delta,\eta) = 
\frac{e^{\frac{i\pi}{4}\sgA}}{\sqrt{r|A|}}
\sum_{\stackrel{\nu=0}{\nu -k \in 2 \Z}}^{N} 
d_{\nu}(l,\delta,\eta)
\left(\frac{i}{rA}\right)^{\frac{\nu-k}{2}}
\Gamma\left( \frac{\nu -k +1}{2} \right).
\end{equation}
To show that $Z_{\inte,2}(N;\delta,\eta)$ is absolutely 
convergent it is enough to show that
$\sum_{l \in \Z} \left| d_{\nu}(l,\delta,\eta) \right|$ is
convergent for all $\nu \in \Z_{\geq 0}$ and all
$\lambda \in \Lambda$ (where $\lambda$ is suppressed from the 
notation). By (\ref{eq:dnu}) and (\ref{eq:P1}) we have
\begin{equation}\label{eq:dnu-1}
d_{\nu}(l,\delta,\eta) = \frac{1}{\pi^{k}} 
\sum_{j=0}^{\nu} \frac{1}{j!(\nu-j)!} 
P_{j}(i\eta +l,\delta) e^{-\pi \delta (i\eta +l)^{2}}
\partial_{y}^{(\nu -j)} h_{\eta}(y) |_{y=0},
\end{equation}
where
$h_{\eta}(y) = q(y+i\eta+l)e^{-2rA\eta y}
\left( \frac{\pi y}{\sin(\pi y)} \right)^{k}$.
Let $\eta_{3}>0$. By the periodicity assumption on $q$ (and 
continuity) we have that
$$
M_{\nu}(\eta_{3}):=\max\left\{ \, \left| \partial_{y}^{(j)} 
h_{\eta}(y) |_{y=0} \right| \, : \, 
j \in \{0,1,\ldots,\nu\}, \lambda \in \Lambda, l \in \Z,
\eta \in [0,\eta_{3}] \, \right\}
$$
is finite for $\nu \in \Z_{\geq 0}$. Using notation from 
(\ref{eq:P1}) we therefore get
\begin{equation}\label{eq:dnuestimate}
|d_{\nu}(l,\delta,\eta)| \leq 
\frac{M_{\nu}(\eta_{3})}{\pi^{k}} e^{\pi\delta \eta_{3}^{2}}
\sum_{j=0}^{\nu} \frac{1}{j!(\nu-j)!}
\sum_{n=0}^{j} |a_{n}^{j}(\delta)|(\delta(\eta_{3} + |l|))^{n}
e^{-\pi \delta l^{2}},
\end{equation}
where the $a_{n}^{j}(\delta)$'s are polynomials in $\delta$
independent of $l$. It follows that 
$\sum_{l \in \Z} \left| d_{\nu}(l,\delta,\eta) \right|$
is convergent. In fact, we see that this series is uniformly
convergent w.r.t.\ $\eta$ on $[0,\eta_{3}]$.

Let us next show that $\Sigma^{5}(N;\delta,\eta)$ is absolutely
convergent. To this end we consider the series
$\sum_{l \in \Z} \left| J_{2}(N;l,\delta,\eta) \right|$ for
$\lambda \in \Lambda$ arbitrary but fixed. Here
$$
\left| J_{2}(N;l,\delta,\eta) \right| \leq 
\sum_{\stackrel{\nu=0}{\nu -k \in 2 \Z}}^{N}
\left| d_{\nu}(l,\delta,\eta) \right|
\int_{|t| \geq \frac{1}{2}} |C_{\eta}(t)|^{\nu -k} 
\left|e^{irAC_{\eta}(t)^{2}} \right|\dte t.
$$
We have
$\left|e^{irAC_{\eta}(t)^{2}} \right| = 
e^{-r|A| t^{2}} e^{\sqrt{2} r|A| \eta t}$.
Moreover, $|C_{\eta}(t)| \leq |t| + \eta$ and
$|C_{\eta}(t)| \geq |t|/\sqrt{2}$, so
$$
\left| J_{2}(N;l,\delta,\eta) \right| \leq (2\sqrt{2})^{k} 
\sum_{\stackrel{\nu=0}{\nu -k \in 2 \Z}}^{N} 2^{\nu}
\left| d_{\nu}(l,\delta,\eta) \right|
\int_{|t| \geq \frac{1}{2}} |t|^{\nu}
e^{\sqrt{2} r|A| \eta t} e^{-r|A| t^{2}} \dte t
$$
for $\eta \leq 1/2$. For $t \geq 2\sqrt{2}\eta$ we have
$\exp\left(\sqrt{2} r|A| \eta t \right) \leq 
\exp\left( \frac{1}{2}r|A| t^{2} \right)$, so
$$
\left| J_{2}(N;l,\delta,\eta) \right| \leq (2\sqrt{2})^{k}
\sum_{\stackrel{\nu=0}{\nu -k \in 2 \Z}}^{N} 2^{\nu}
\left| d_{\nu}(l,\delta,\eta) \right|
\int_{|t| \geq \frac{1}{2}} |t|^{\nu}
e^{-\frac{1}{2}r|A| t^{2}} \dte t
$$
for $\eta \in ]0,\eta']$, where $\eta'=\frac{1}{4\sqrt{2}}$.
Therefore
$$
\left| J_{2}(N;l,\delta,\eta) \right| \leq 
\frac{1}{4}(2\sqrt{2})^{k}
\sum_{\stackrel{\nu=0}{\nu -k \in 2 \Z}}^{N} 
\left( \frac{8}{r|A|} \right)^{\frac{\nu}{2}+1} 
\left| d_{\nu}(l,\delta,\eta) \right|
\Gamma\left( \frac{\nu+1}{2}, \frac{1}{8} |A|r \right).
$$
This implies that
$\sum_{l \in \Z} \left| J_{2}(N;l,\delta,\eta) \right| 
\leq A_{5}(N,\lambda;\delta,\eta)$, 
where
\begin{equation}\label{eq:A5}
A_{5}(N,\lambda;\delta,\eta) = \frac{1}{4}
\left(2\sqrt{2}\right)^{k}
\sum_{\stackrel{\nu=0}{\nu -k \in 2 \Z}}^{N}
\left( \frac{8}{r|A|} \right)^{\frac{\nu}{2}+1}
\Gamma\left( \frac{\nu+1}{2}, \frac{1}{8} |A| r \right)
\sum_{l \in \Z} \left| d_{\nu}(l;\delta,\eta) \right|.
\end{equation}
Thus $\Sigma^{5}(N;\delta,\eta)$ is absolutely convergent.

Next we show that $\Sigma^{4}(N;\delta,\eta)$ is absolutely 
convergent. It is enough to prove that the series
$\sum_{l \in \Z} \left| J_{1}(N;l,\delta,\eta) \right|$ is
convergent. Thus we have to estimate 
$R_{N}(y;l,\delta,\eta) =
\sum_{\stackrel{\nu = N+1}{\nu-k \in 2\Z}}^{\infty}
d_{\nu}(l,\delta,\eta)$ for $y=C_{\eta}(t)$, $|t| \leq 1/2$.
By (\ref{eq:dnu}) and Cauchy's formula we have
$$
d_{\nu}(l,\delta,\eta) = \frac{1}{\pi^{k}}\frac{1}{2\pi i}
\int_{|z| = \rho} \frac{\Psi(z;l,\delta,\eta)}{z^{\nu +1}} 
\dte z
$$
for $\nu \in \Z_{\geq 0}$ and $\rho \in ]0,1[$. Therefore
$\left| d_{\nu}(l,\delta,\eta) \right| \leq \pi^{-k}\rho^{-\nu}
L(\rho,l,\delta,\eta)$,
where
$L(\rho,l,\delta,\eta) = 
\max_{|z|=\rho} \left| \Psi(z;l,\delta,\eta) \right|$.
Let $\nu_{1} \in \{N+1,N+2\}$ such that $\nu_{1} - k \in 2\Z$.
Assume first that $\nu_{1} \geq k$. Then
$$
\left| R_{N}(y;l,\delta,\eta) \right| \leq \frac{1}{\pi^{k}}
L(\rho,l,\delta,\eta) 
\sum_{\stackrel{\nu=\nu_{1}}{\nu -k \in 2\Z}}^{\infty}
\frac{1}{\rho^{\nu}} |y|^{\nu -k}
\leq \frac{1}{\pi^{k}} L(\rho,l,\delta,\eta) 
\frac{1}{\rho^{\nu_{1}}} |y|^{\nu_{1} -k} \sum_{\nu=0}^{\infty}
\left( \frac{|y|}{\rho} \right)^{\nu}
$$
for $\rho \in ]0,1[$ and $|y| < \rho$. Now let $\rho =3/4$ and
let $L(l,\delta,\eta)=L(3/4,l,\delta,\eta)$. Then
$$
\left| R_{N}(y;l,\delta,\eta) \right| \leq \frac{9}{\pi^{k}} 
\left(\frac{4}{3} \right)^{\nu_{1}} L(l,\delta,\eta) 
|y|^{\nu_{1} -k}
$$
for $|y| \leq 2/3$. For $\eta \leq 1/6$ we have 
$|C_{\eta}(t)| \leq 2/3$ for all $t \in [-1/2,1/2]$ so
$$
\left| J_{1}(N;l,\delta,\eta) \right| \leq \frac{18}{\pi^{k}}
\left(\frac{4}{3} \right)^{\nu_{1}}
\exp\left( \frac{r|A|\eta}{\sqrt{2}} \right) L(l,\delta,\eta) 
\int_{0}^{1/2} (t+\eta)^{\nu_{1}-k} e^{-r|A|t^{2} }\dte t.
$$
By (\ref{eq:Psi}) and the assumptions on $q$, see
(\ref{eq:qfunction-symmetry1}) and 
(\ref{eq:qfunction-symmetry2}), we have
$$
L(\rho,l,\delta,\eta) \leq M_{1}(\rho)^{k} b_{0}
e^{a_{0}(\rho+\eta+m_{q})} e^{2r|A|\eta\rho}
\max_{|z|=\rho} \left| e^{-\pi\delta(z+i\eta+l)^{2}} \right|,
$$
where
$M_{1}(\rho) = \max_{|z|=\rho} 
\left| (\pi z)/\sin(\pi z) \right|$.
Here
$\max_{|z|=\rho} \left| e^{-\pi\delta(z+i\eta+l)^{2}} \right|
\leq e^{-\pi\delta l^{2} + 2\pi\delta\rho|l|} 
e^{\pi \delta (\eta+\rho)^{2}}$
so
$$
L(l,\delta,\eta)=L(3/4,l,\delta,\eta) = C(\delta,\eta) 
e^{\frac{3}{2} |A| r \eta}
\left( e^{-\frac{3}{2}\pi\delta l} + 
e^{\frac{3}{2}\pi\delta l}\right)
e^{-\pi\delta l^{2}},
$$
where
$C(\delta,\eta) = M_{1}(3/4)^{k} b_{0}e^{a_{0}(3/4+\eta+m_{q})}
e^{\pi\delta(\eta+3/4)^{2}}$.
Therefore $\sum_{l \in \Z} |J_{1}(N;l,\delta,\eta)|$
is bounded from above by 
\begin{equation}\label{eq:A4}
A_{4}^{0}(N,\lambda;\delta,\eta) := C_{1}(\lambda;\delta,\eta)
\int_{0}^{\infty} (t+\eta)^{\nu_{1}-k} e^{-r|A|t^{2}} \dte t
\sum_{l \in \Z} \left( e^{-\frac{3}{2}\pi\delta l}
+ e^{\frac{3}{2}\pi\delta l}\right) e^{-\pi\delta l^{2}},
\end{equation}
where
$$
C_{1}(\lambda;\delta,\eta) = \frac{18C(\delta,\eta)}{\pi^{k}} 
\left( \frac{4}{3} \right)^{\nu_{1}}
\exp\left(\left( \frac{3}{2} + \frac{1}{\sqrt{2}}\right) 
|A|r\eta\right).
$$
We note that $A_{4}^{0}(N,\lambda;\delta,\eta)$ is convergent, 
hence $\Sigma_{4}(N;\delta,\eta)$ is absolutely convergent.

Assume next that $\nu_{1} < k$ and write
$$
R_{N}(y;l,\delta,\eta) = 
\sum_{\stackrel{\nu=\nu_{1}}{\nu -k \in 2\Z}}^{k-2}
d_{\nu}(l,\delta,\eta) y^{\nu -k} + R(y;l,\delta,\eta).
$$
Then
$J_{1}(N;l,\delta,\eta) = I_{3}(N;l,\delta,\eta) + 
J_{3}(l,\delta,\eta)$,
where
\begin{eqnarray*}
I_{3}(N;l,\delta,\eta) &=& 
\sum_{\stackrel{\nu=\nu_{1}}{\nu -k \in 2\Z}}^{k-2}
d_{\nu}(l,\delta,\eta) \int_{C_{\eta}^{0}} y^{\nu -k} 
e^{irAy^{2}} \dte y, \\
J_{3}(l,\delta,\eta) &=& \int_{C_{\eta}^{0}} R(y;l,\delta,\eta)
e^{irAy^{2}} \dte y.
\end{eqnarray*}
Let
$$
\Sigma^{7}(\delta,\eta) = \sum_{\lambda \in \Lambda} g(\lambda)
e^{-irA\eta^{2}} \sum_{l \in I_{\lambda}} \beta(l)
J_{3}(l,\delta,\eta).
$$
According to the above we have that
$\sum_{l \in \Z} |J_{3}(l,\delta,\eta)|$
is bounded from above by an expression equal to the right-hand
side of (\ref{eq:A4}) with $\nu_{1}$ replaced by $k$, i.e.\ by
\begin{equation}\label{eq:A4-1}
A_{4}^{1}(\lambda;\delta,\eta) := C_{1}(\lambda;\delta,\eta)) 
\int_{0}^{\infty} e^{-r |A| t^{2}} \dte t
\sum_{l \in \Z} \left( e^{-\frac{3}{2}\pi\delta l} 
+ e^{\frac{3}{2}\pi\delta l}\right) e^{-\pi\delta l^{2}},
\end{equation}
where $C_{1}(\lambda,\delta,\eta)$ is as in (\ref{eq:A4}) with
$\nu_{1}$ replaced by $k$. The series $A_{4}^{1}(\delta,\eta)$ 
is convergent, hence $\Sigma^{7}(\delta,\eta)$ is absolutely 
convergent. We further partition
$I_{3}(N;l,\delta,\eta) = I_{4}(N;l,\delta,\eta) + 
J_{4}(N;l,\delta,\eta)$
with
\begin{eqnarray*}
I_{4}(N;l,\delta,\eta) &=& 
\sum_{\stackrel{\nu=\nu_{1}}{\nu -k \in 2\Z}}^{k-2}
d_{\nu}(l,\delta,\eta) \int_{C_{\eta}} y^{\nu -k} 
e^{irAy^{2}} \dte y, \\
J_{4}(N;l,\delta,\eta) &=& 
-\sum_{\stackrel{\nu=\nu_{1}}{\nu -k \in 2\Z}}^{k-2} 
d_{\nu}(l,\delta,\eta)\int_{C_{\eta}^{\infty}} y^{\nu -k} 
e^{irAy^{2}} \dte y.
\end{eqnarray*}
Put
\begin{eqnarray*}
\Sigma^{8}(N;\delta,\eta) = \sum_{\lambda \in \Lambda} 
g(\lambda)e^{-irA\eta^{2}} \sum_{l \in I_{\lambda}} 
\beta(l)I_{4}(N;l,\delta,\eta), \\
\Sigma^{9}(N;\delta,\eta) = \sum_{\lambda \in \Lambda} 
g(\lambda)e^{-irA\eta^{2}} \sum_{l \in I_{\lambda}}
\beta(l)J_{4}(N;l,\delta,\eta).
\end{eqnarray*}
As in the case $Z_{\inte,2}(N;\delta,\eta)$ we get that
$\sum_{l \in \Z} |I_{4}(N;l,\delta,\eta)|$ is bounded from 
above by
\begin{equation}\label{eq:A4-2}
A_{4}^{2}(N,\lambda;\delta,\eta) :=
\sum_{\stackrel{\nu=\nu_{1}}{\nu -k \in 2 \Z}}^{k-2} 
\left( \frac{1}{r|A|} \right)^{\frac{\nu-k+1}{2}}
\left|\Gamma\left( \frac{\nu -k +1}{2} \right)\right| 
\sum_{l \in \Z} \left| d_{\nu}(l,\delta,\eta)\right|,
\end{equation}
thus $\Sigma^{8}(N;\delta,\eta)$ is absolutely convergent.
Moreover, similar to the case $\Sigma^{5}(N;\delta,\eta)$ we 
get that $\sum_{l \in \Z} |J_{4}(N;l,\delta,\eta)|$ is bounded
from above by
\begin{equation}\label{eq:A4-3}
A_{4}^{3}(N,\lambda;\delta,\eta) := \frac{1}{4}(2\sqrt{2})^{k} 
\sum_{\stackrel{\nu=\nu_{1}}{\nu -k \in 2 \Z}}^{k-2}
\left( \frac{8}{r|A|} \right)^{\frac{\nu}{2}+1}
\Gamma\left( \frac{\nu+1}{2}, \frac{1}{8} |A|r \right)
\sum_{l \in \Z} \left| d_{\nu}(l;\delta,\eta) \right|
\end{equation}
showing that $\Sigma^{9}(N;\delta,\eta)$ is absolutely
convergent.

Finally let us consider $\Sigma^{6}(\delta,\eta)$. By
(\ref{eq:normsin}) we have
$$
\left| \sin\left(\pi C_{\eta}(t)\right) \right| \geq 
\frac{1}{2} \sqrt{ e^{2v_{\eta}(t)} + e^{-2v_{\eta}(t)} -2},
$$
where 
$v_{\eta}(t)=\pi\left( \frac{1}{\sqrt{2}}\sgA t -\eta \right)$.
Let $\eta_{2} < \frac{1}{2\sqrt{2}}$ and let
$M_{2} = \inf_{\eta \in [0,\eta_{2}], |t| \geq \frac{1}{2}} 
\left| \sin\left(\pi C_{\eta}(t)\right) \right| >0$.
Then
$$
\sum_{l \in \Z} |I_{\infty}(l,\delta,\eta)| \leq 
\frac{1}{M_{2}^{k}} \sum_{l \in \Z} 
\int_{|t| \geq \frac{1}{2}} f(t;l,\delta,\eta) \dte t
$$
for $\eta \in ]0,\eta_{2}]$ and $\delta>0$, where
$$
f(t;l,\delta,\eta) = \left( 
|\tilde{G}(C_{\eta}(t);l,\delta,\eta)| + 
|\tilde{G}(-C_{\eta}(t);l,\delta,\eta)| \right)
\left|e^{irAC_{\eta}(t)^{2}}\right|,
$$
where
$\tilde{G}(y;l,\delta,\eta) = e^{-\pi\delta(y+i\eta+l)^{2}} 
q(y+i\eta+l)e^{-2rA\eta y}$.
By (\ref{eq:qfunction-symmetry1}),
(\ref{eq:qfunction-symmetry2}) and the fact that
$|C_{\eta}(t)| \leq |t|+\eta$ we get
$$
\left| q\left( \pm C_{\eta}(t) +i\eta+l\right) \right| \leq 
b_{0} e^{a_{0}(2\eta+m_{q}+|t|)}.
$$
Moreover,
$\left|e^{-2rA\eta C_{\eta}(t)} \right| \leq 
e^{\sqrt{2} r|A|\eta t}$,
and
$
\left|e^{irAC_{\eta}(t)^{2}}\right| \leq 
e^{-r|A|t^{2}} e^{\sqrt{2}r|A|\eta |t|}$
so we have
$$
f(t;l,\delta,\eta) \leq b_{0} e^{a_{0}(2\eta+m_{q})} 
e^{a_{0}|t|} e^{2\sqrt{2}r|A|\eta |t| } e^{-r|A|t^{2}}
\sum_{\mu \in \{\pm 1\}} \left|
e^{-\pi\delta\left(\mu C_{\eta}(t)+i\eta+l\right)^{2}}\right|.
$$
Here
$\left|e^{-\pi\delta\left(-C_{\eta}(t)+i\eta+l\right)^{2}}
\right| = g(t;l)e^{4\pi\delta\eta^{2}}
e^{-2\sqrt{2}\pi\delta \eta\sgA t}$
and
$\left|e^{-\pi\delta\left(C_{\eta}(t)+i\eta+l\right)^{2}}
\right| = g(t;l)$,
where $g(t;l)=e^{-\pi\delta l^{2}}e^{-\sqrt{2}\pi\delta lt}$ so
$$ 
f(t;l,\delta,\eta) \leq C_{2}(\delta,\eta) e^{-\pi\delta l^{2}}
e^{\sqrt{2}\pi\delta |l||t|}
e^{D(\delta,\eta)|t|}e^{2\sqrt{2}r|A|\eta |t| }
\exp\left(-r|A|t^{2}\right),
$$
where
$C_{2}(\delta,\eta) = b_{0}e^{a_{0}(2\eta+m_{q})} 
e^{4\pi\delta\eta^{2}}$
and $D(\delta,\eta)= a_{0}+2\sqrt{2}\pi\delta\eta$. Choose for
each $\lambda \in \Lambda$ a $r_{\lambda} \in \Z_{\geq 2}$ so
$$
e^{(D(\delta,\eta)+2\sqrt{2}r|A|\eta)t} \leq 
\exp\left(\frac{1}{2} r|A| t^{2}\right)
$$
for all $t \geq 1/2$, $\delta, \eta \in [0,1/2]$ and 
$r \geq r_{\lambda}$.
For $\delta \in [0,1/2]$, $\eta \in [0, \eta_{2}]$, and 
$r \geq r_{\lambda}$ we then get
$$
\int_{\frac{1}{2}}^{\infty} e^{\sqrt{2}\pi\delta |l|t}
e^{D(\delta,\eta)t}e^{2\sqrt{2}r|A|\eta t}e^{-r|A|t^{2}} \dte t
\leq \int_{\frac{1}{2}}^{\infty} e^{\sqrt{2}\pi\delta |l|t}
e^{-\frac{1}{2}r|A|t^{2}} \dte t.
$$
We therefore have
$$
\sum_{l \in \Z} |I_{\infty}(l,\delta,\eta)| \leq 
\frac{2C_{2}(\delta,\eta)}{M_{2}^{k}} \sum_{l \in \Z} 
e^{-\pi \delta l^{2}} \int_{\frac{1}{2}}^{\infty}
e^{\sqrt{2}\pi\delta |l|t}e^{-\frac{1}{2}r|A|t^{2}} \dte t
$$
for $\delta \in [0,1/2]$, $\eta \in [0, \eta_{2}]$ and 
$r \geq r_{\lambda}$. For $t>0$ we have that
$e^{\sqrt{2}\pi\delta |l|t} \leq e^{\frac{r}{4}|A|t^{2}}$ if 
and only if $t \geq c\delta |l|$, where 
$c=4\sqrt{2}\pi/(r|A|)$. If $c\delta |l| \leq \frac{1}{2}$ we 
therefore get
$$
\int_{\frac{1}{2}}^{\infty} e^{\sqrt{2}\pi\delta |l|t}
e^{-\frac{1}{2}r|A|t^{2}} \dte t \leq 
\int_{\frac{1}{2}}^{\infty} e^{-\frac{1}{4}r|A|t^{2}} \dte t
= \frac{1}{\sqrt{r|A|}} \int_{\frac{r|A|}{16}}^{\infty} 
u^{\frac{1}{2}-1} e^{-u} \dte u
= \frac{1}{\sqrt{r|A|}}
\Gamma\left( \frac{1}{2},\frac{1}{16}r|A| \right).
$$
For $\frac{1}{2} < c\delta |l|$ we have
$$
\int_{\frac{1}{2}}^{\infty} e^{\sqrt{2}\pi\delta |l|t}
e^{-\frac{1}{2}r|A|t^{2}} \dte t \leq 
\int_{\frac{1}{2}}^{c\delta |l|} e^{\sqrt{2}\pi\delta |l|t}
e^{-\frac{1}{2}r|A|t^{2}} \dte t + 
\int_{c\delta |l|}^{\infty} e^{-\frac{1}{4}r|A| t^{2}} \dte t.
$$
Here
$\int_{c\delta |l|}^{\infty} e^{-\frac{1}{4}r|A| t^{2}} \dte t
\leq \int_{\frac{1}{2}}^{\infty} e^{-\frac{1}{4}r|A| t^{2}} 
\dte t = \frac{1}{\sqrt{r|A|}}
\Gamma\left( \frac{1}{2},\frac{1}{16}r|A| \right)$
and
$$
\int_{\frac{1}{2}}^{c\delta |l|} e^{\sqrt{2}\pi\delta |l|t}
e^{-\frac{1}{2}r|A|t^{2}} \dte t \leq e^{-\frac{1}{16}r|A|} 
\int_{0}^{c\delta |l|} e^{\sqrt{2}\pi\delta |l|t} \dte t
\leq e^{-\frac{1}{16}r|A|} c\delta |l|
e^{\sqrt{2}\pi\delta c \delta^{2} l^{2}}.
$$
By putting everything together we get that
$\sum_{l \in \Z} |I_{\infty}(l,\delta,\eta)|$
is bounded from above by
$$ 
\frac{2C_{2}(\delta,\eta)}{M_{2}^{k}}
\sum_{l \in \Z} e^{-\pi\delta l^{2}}
\left( \frac{1}{\sqrt{r|A|}} 
\Gamma\left( \frac{1}{2},\frac{1}{16}r|A| \right)
+ \frac{4\sqrt{2}\pi e^{-\frac{1}{16}r|A|}}{r|A|} \delta |l| 
e^{\frac{8\pi\delta}{r|A|} \pi\delta l^{2}} \right)
$$
for $\delta \in ]0,1/2]$, $\eta \in ]0,\eta_{2}]$ and 
$r \in \Z_{\geq r_{\lambda}}$. By choosing 
$\delta_{0} \in ]0,1/2]$ so small that 
$\frac{8\pi\delta_{0}}{A_{0}} \leq \frac{1}{2}$ we get
$\sum_{l \in \Z} |I_{\infty}(l,\delta,\eta)| \leq 
A_{6}(\lambda;\delta,\eta)$, 
where
\begin{equation}\label{eq:A6}
A_{6}(\lambda;\delta,\eta) = C_{3}(\delta,\eta)
\left( \frac{1}{\sqrt{r|A|}} 
\Gamma\left( \frac{1}{2},\frac{1}{16}|A|r \right) 
\sum_{l \in \Z} e^{-\pi\delta l^{2}} + C_{4}(r)\sum_{l \in \Z}
\delta |l| e^{-\frac{1}{2}\pi\delta l^{2}} \right)
\end{equation}
for all $\eta \in [0,\eta_{2}]$, $\delta \in ]0,\delta_{0}]$
and $r \in \Z_{\geq r_{\lambda}}$. Here
$C_{3}(\delta,\eta)=2C_{2}(\delta,\eta)/M_{2}^{k}$ is 
independent of $r$ and 
$C_{4}(r) = \frac{4\sqrt{2}\pi}{|A|r}e^{- \frac{1}{16}|A|r}$.
Finally put 
$r_{0}=\max_{\lambda \in \Lambda}\{ r_{\lambda} \}$. The claims
about uniform convergence immediately follows from the above.
\end{prfa}

\begin{prfa}{\bf of \refprop{prop:Lemma8.26}\hspace{.1in}}
By the remarks following \reflem{lem:Lemma8.25} we have
$$
Z_{\inte,2}(N;\delta) = \frac{1}{\pi^{k}}
\sum_{\lambda \in \Lambda} g(\lambda)
\frac{e^{\frac{i\pi}{4}\sgA}}{\sqrt{r|A|}}
\left( D(N,\lambda;\delta) + E(N,\lambda;\delta) \right),
$$
where
{\allowdisplaybreaks
\begin{eqnarray*}
D(N,\lambda;\delta) &=& \sum_{l \in I_{\lambda}} 
\beta(l,\lambda) \sum_{\stackrel{\nu=0}{\nu -k \in 2 \Z}}^{N}
\left(\frac{i}{rA}\right)^{\frac{\nu-k}{2}}
\frac{\Gamma\left( \frac{\nu -k +1}{2} \right)}{\nu !}
e^{-\pi\delta l^{2}} \partial_{y}^{(\nu)} h(y)|_{y=0}, \\ 
E(N,\lambda;\delta) &=& \sum_{l \in I_{\lambda}}
\beta(l,\lambda) \sum_{\stackrel{\nu=0}{\nu -k \in 2 \Z}}^{N}
\left(\frac{i}{rA}\right)^{\frac{\nu-k}{2}}
\Gamma\left( \frac{\nu -k +1}{2} \right) \\
&& \times \sum_{j=1}^{\nu} \frac{1}{j!(\nu-j)!}
P_{j}(l,\delta) e^{-\pi\delta l^{2}}
\partial_{y}^{(\nu-j)} h(y;l)|_{y=0},
\end{eqnarray*}}\noindent
where 
$h(y;l)=h_{0}(y) = q(y+l)
\left( \frac{\pi y}{\sin(\pi y)} \right)^{k}$. 
By \reflem{lem:Lemma8.11} together with the remarks about the 
polynomials $P_{j}(l,\delta)$ in connection to (\ref{eq:P1}) 
we have
$$
\lidel \sdel E(N,\lambda;\delta) = 0.
$$
By (\ref{eq:qfunction-symmetry2}) the derivatives of 
$y \mapsto h(y;l)$ in $y=0$ are periodic in $l$ with a period
of $m_{q}$. Moreover, the function $\beta(l,\lambda)$ is 
periodic in $l$ for each fixed $\lambda$ with a period equal to
the least common multiplum of $2$, $P(\lambda)$ and 
$b(\lambda)$, where $P(\lambda)$ and $b(\lambda)$ are as in 
\refprop{prop:Lemma8.26}. Finally we observe that if
$l_{0} \in I_{\lambda}$ then 
$l_{0} + P(\lambda)\Z \subseteq I_{\lambda}$. Thus 
(\ref{eq:Zintmain2}) follows by \reflem{lem:periodicity}.

Next let us estimate the remainder term $R_{2}(N;\delta)$. By 
the proof of \reflem{lem:Lemma8.25} and the remarks after that
lemma we have
$$
\ep(N,\lambda;\delta) \leq A_{4}(N,\lambda;\delta) +
A_{5}(N,\lambda;\delta) + A_{6}(\lambda;\delta),
$$
where 
$A_{\nu}(N,\lambda;\delta) = \lieta 
A_{\nu}(N,\lambda;\delta,\eta) = 
A_{\nu}(N,\lambda;\delta,0)$ for $\nu=4,5$, and 
$A_{6}(\lambda;\delta)$ is equal to 
$\lieta A_{6}(\lambda;\delta,\eta) = A_{6}(\lambda;\delta,0)$, where 
$A_{5}(N,\lambda;\delta,\eta)$ and $A_{6}(\lambda;\delta,\eta)$
are given by respectively (\ref{eq:A5}) and (\ref{eq:A6}), and
where 
$A_{4}(N,\lambda;\delta,\eta) = 
A_{4}^{0}(N,\lambda;\delta,\eta)$ in case $\nu_{1} \geq k$ and
\begin{equation}\label{eq:A4case2}
A_{4}(N,\lambda;\delta,\eta) = A_{4}^{1}(\lambda;\delta,\eta) 
+ A_{4}^{2}(N,\lambda;\delta,\eta) + 
A_{4}^{3}(N,\lambda;\delta,\eta)
\end{equation}
in case $\nu_{1} < k$, where 
$A_{4}^{0}(N,\lambda;\delta,\eta)$, 
$A_{4}^{1}(\lambda;\delta,\eta)$, 
$A_{4}^{2}(N,\lambda;\delta,\eta)$ and 
$A_{4}^{3}(N;\lambda;\delta,\eta)$ are given by respectively 
(\ref{eq:A4}), (\ref{eq:A4-1}), (\ref{eq:A4-2}) and
(\ref{eq:A4-3}). Recall here that $\nu_{1} \in \{N+1,N+2\}$ 
such that $\nu_{1} - k \in 2\Z$.  

As in the calculation of $Z_{\inte,2}(r;N)$ above we get
$$
A_{5}(N,\lambda;\delta) \leq 
\tilde{A}_{5}(N,\lambda;\delta) := \frac{1}{4}
\left(\frac{2\sqrt{2}}{\pi}\right)^{k}
\left(D(N,\lambda;\delta) + E(N,\lambda;\delta)\right),
$$
with
$$
D(N,\lambda;\delta) = 
\sum_{\stackrel{\nu=0}{\nu -k \in 2 \Z}}^{N}
\left( \frac{8}{r|A|} \right)^{\frac{\nu}{2} +1}
\frac{1}{\nu !}
\Gamma\left( \frac{\nu+1}{2}, \frac{1}{8}|A|r \right) \\
\sum_{l \in \Z} e^{-\pi \delta l^{2}}
\left| \partial_{y}^{(\nu)} h(y;l)|_{y=0} \right|
$$
and with $\lidel \sdel E(N,\lambda;\delta)=0$. We let
\begin{eqnarray*}
A_{5}(r;N,\lambda) &=& r\lidel\sdel 
\tilde{A}_{5}(N,\lambda;\delta)
= \frac{r}{4}
\left(\frac{2\sqrt{2}}{\pi}\right)^{k}
\lidel \sdel D(N,\lambda;\delta) \\
&=& \sum_{\stackrel{\nu=0}{\nu -k \in 2 \Z}}^{N}
c_{\nu} r^{-\frac{\nu}{2}}
\Gamma\left( \frac{\nu+1}{2}, \frac{1}{8}|A| r \right),
\end{eqnarray*}
where
$c_{\nu} = \frac{1}{4m_{q}\nu !}
\left(\frac{2\sqrt{2}}{\pi}\right)^{k}
\left( \frac{8}{|A|} \right)^{\frac{\nu}{2} +1}
\sum_{l=0}^{m_{q}-1} 
\left| \partial_{y}^{(\nu)} h(y;l)|_{y=0} \right|$ 
are independent of $r$. 

Next we get
$$
A_{6}(r;\lambda) := r\lidel \sdel A_{6}(\lambda;\delta)
= C \sqrt{\frac{r}{|A|}}
\Gamma\left( \frac{1}{2}, \frac{1}{16}|A|r\right),
$$
where $C =C_{3}(0,0)$ is independent of $r$.

Let us finally consider $A_{4}(N,\lambda;\delta)$. We 
first consider the case $\nu_{1} \geq k$. In this case
$A_{4}(N,\lambda;\delta)$ is equal to 
$A_{4}^{0}(N,\lambda;\delta,\eta)$ in (\ref{eq:A4}) with 
$\eta=0$. We have
$$
A_{4}^{0}(r;N,\lambda) := r\lidel \sdel A_{4}(N,\lambda;\delta)
= C r^{\frac{k-\nu_{1}+1}{2}},
$$
where 
$C= \frac{1}{2} C_{1}(\lambda;0,0) 
\Gamma\left(\frac{\nu_{1}-k+1}{2}\right)
|A|^{\frac{k-\nu_{1}-1}{2}}$ 
is independent of $r$. (Actually $C_{1}(\lambda;0,0)$ is
independent of $\lambda$.)
 
Next assume that $\nu_{1} <k$. Here 
$A_{4}(N,\lambda;\delta)$ is given by the right-hand side of
(\ref{eq:A4case2}) with $\eta=0$. Recall here that 
$A_{4}^{1}(\lambda;\delta,\eta)$ is simply equal to 
$A_{4}^{0}(N,\lambda;\delta,\eta)$ with $\nu_1$ replaced by 
$k$, hence
$$
A_{4}^{1}(r;\lambda) := 
r\lidel\sdel A_{4}^{1}(\lambda;\delta,0) = C\sqrt{r},
$$
where $C$ is a constant independent of $r$.

The analysis of $A_{4}^{2}(N,\lambda;\delta,0)$ and
$A_{4}^{3}(N,\lambda;\delta,0)$ are similar to the analysis 
of $A_{5}(N,\lambda;\delta)$, the main point being to 
analyse the series $\sum_{l \in \Z} |d_{\nu}(l;\delta,0)|$.
Thus we find that $A_{4}^{2}(N,\lambda;\delta,0)$ is bounded 
from above by an expression 
$\tilde{A}_{4}^{2}(N,\lambda;\delta)$ for which
$A_{4}^{2}(r;N,\lambda) := r\lidel 
\tilde{A}_{4}^{2}(N,\lambda;\delta)$ 
exists and is equal to
$$
A_{4}^{2}(r;N,\lambda) = 
\sum_{\stackrel{\nu=\nu_{1}}{\nu -k \in 2\Z}}^{k-2}
c_{\nu} r^{\frac{k-\nu +1}{2}},
$$
where
$$
c_{\nu} = \frac{1}{m_{q}\pi^{k}}
\left| \Gamma\left( \frac{\nu-k+1}{2} \right)\right|
|A|^{-\frac{\nu-k+1}{2}}
\sum_{l=0}^{m_{q}-1}
\left| \partial_{y}^{(\nu)} h(y;l)|_{y=0} \right|
$$
are independent of $r$. Finally we find that 
$A_{4}^{3}(N,\lambda;\delta,0)$ is bounded from above by an 
expression $\tilde{A}_{4}^{3}(N,\lambda;\delta)$ for which 
$A_{4}^{3}(r;N,\lambda) := r\lidel\sdel\lieta 
\tilde{A}_{4}^{3}(N,\lambda;\delta)$ 
exists and is equal to
$$
A_{4}^{3}(r;N,\lambda) =
\sum_{\stackrel{\nu=\nu_{1}}{\nu -k \in 2\Z}}^{k-2}
c_{\nu} r^{-\frac{\nu}{2}}
\Gamma\left( \frac{\nu+1}{2},\frac{1}{8}A_{0}r \right),
$$
where
$$
c_{\nu} = \frac{1}{4m_{q}}
\left(\frac{2\sqrt{2}}{\pi}\right)^{k}
\left( \frac{8}{|A|} \right)^{\frac{\nu}{2}+1}
\sum_{l=0}^{m_{q}-1} 
\left| \partial_{y}^{(\nu)} h(y;l)|_{y=0} \right|
$$
are independent of $r$. Thus we can put
$$
\vep_{2}(N,\lambda;\delta) = 
\tilde{A}_{4}(N,\lambda;\delta) +
\tilde{A}_{5}(N,\lambda;\delta) + A_{6}(\lambda;\delta),
$$
where 
$\tilde{A}_{4}(N,\lambda;\delta) = 
A_{4}^{0}(N,\lambda;\delta,0)$
if $\nu_{1} \geq k$ and
$\tilde{A}_{4}(N,\lambda;\delta) = 
A_{4}^{1}(\lambda;\delta,0) + 
\tilde{A}_{4}^{2}(N,\lambda;\delta) +
\tilde{A}_{4}^{3}(N,\lambda;\delta)$
if $\nu_{1} < k$. By the above we have
$$
r\lidel\sdel\vep_{2}(N,\lambda;\delta) =
A_{4}(r;N,\lambda) + A_{5}(r;N,\lambda) + A_{6}(r;\lambda),
$$
where $A_{4}(r;N,\lambda) = A_{4}^{0}(r;N,\lambda)$ for 
$\nu_{1} \geq k$ and
$A_{4}(r;N,\lambda) = A_{4}^{1}(r;\lambda)
+ A_{4}^{2}(r;N,\lambda) + A_{4}^{3}(r;N,\lambda)$
for $\nu_{1} <k$. By the remarks about the incomplete gamma
function at the very end of the proof of 
\refprop{prop:Lemma8.22} we see that this limit is
$O\left(r^{\frac{k-\nu_{1}+1}{2}}\right)$. 
\end{prfa}

\subsection{Appendix F. Proof of \reflem{lem:A01}}

\begin{prfa}{\bf of \reflem{lem:A01}\hspace{.1in}}
Let $(m,\lambda) \in \Z \times \Lambda$ be arbitrary and put
$\nu=\nu(m,\lambda)$. Then
$$
\int_{\kappa_{\nu}^{\rho}} f(z;m,\lambda) \dte z = 
= e^{\pi \delta \rho^{2}} e^{-2\pi r |\nu|\rho}
\iny e^{-\pi\delta t^{2}} e^{-2\pi\delta i\sign(\nu)\rho t}
\frac{q_{1}(t+i\sign(\nu)\rho)e^{2\pi ir\nu t}}
{\sin^{k}(\pi(t+i\sign(\nu)\rho-i\eta))} \dte t.
$$
But then
$$
\left| \int_{\kappa_{\nu}^{\rho}} f(z;m,\lambda) \dte z \right|
\leq b_{0}'e^{a_{0}'(m_{q}'+\rho)}e^{\pi \delta \rho^{2}} 
e^{-2\pi r |\nu|\rho} \iny e^{-\pi\delta t^{2}}
\frac{1}{\left|\sin(\pi(t+i\sign(\nu)\rho-i\eta))\right|^{k}}
\dte t.
$$
Here
$$
\left|\sin(\pi(t+i\sign(\nu)\rho-i\eta))\right|
= \frac{1}{2}\sqrt{ e^{2a} + e^{-2a} -2\cos(2\pi t)}
\leq \frac{1}{2}\sqrt{ e^{2a} + e^{-2a} -2},
$$
where $a=\pi(\sign(\nu)\rho -\eta)$. For 
$\eta \in ]0,\eta_{0}]$ we have 
$|a| \geq \pi(\rho-\eta_{0}) >0$. If we put
$c_{\rho} = \sqrt{ e^{2\pi(\rho-\eta_{0})} + 
e^{-2\pi(\rho-\eta_{0})} -2}$
we get
\begin{eqnarray}\label{eq:kappaintegral}
\left| \int_{\kappa_{\nu}^{\rho}} f(z;m,\lambda) \dte z \right|
&\leq& b_{0}'e^{a_{0}'(m_{q}'+\rho)} 
\left(\frac{2}{c_{\rho}}\right)^{k}
e^{\pi \delta \rho^{2}} e^{-2\pi r |\nu|\rho} 
\iny e^{-\pi\delta t^{2}} \dte t \nonumber\\
&=& b_{0}'e^{a_{0}'(m_{q}'+\rho)} \frac{1}{\sqrt{\delta}} 
\left(\frac{2}{c_{\rho}}\right)^{k}
e^{\pi \delta \rho^{2}} e^{-2\pi r |\nu|\rho}.
\end{eqnarray}
Therefore
$$
\sum_{m \in \Z}
\left| \int_{\kappa_{\nu}^{\rho}} f(z;m,\lambda) \dte z \right|
\leq \frac{1}{\sqrt{\delta}} b_{0}'e^{a_{0}'(m_{q}'+\rho)}
\left(\frac{2}{c_{\rho}}\right)^{k} e^{\pi \delta \rho^{2}} 
\sum_{m \in \Z} e^{-2\pi r \rho |\nu(m,\lambda)|},
$$
which is convergent since $\nu(m,\lambda)$ is an affine 
function of $m$, showing that $Z_{\inte}^{J}(\delta,\eta)$ is
absolutely convergent.

For the sum $Z_{\pol}^{J}(\delta,\eta)$ we note that
\begin{eqnarray*}
&&\Res_{z=l+i\eta} \left\{ e^{-\pi\delta z^{2}}
\frac{q_{1}(z)e^{2\pi ir \nu z}}{\sin^{k}(\pi(z-i\eta))} 
\right\} \\
&=& (-1)^{kl} e^{2\pi ir\nu l}
\Res_{z=i\eta} \left\{ e^{-\pi\delta (z+l)^{2}}
\frac{q_{1}(z+l)e^{2\pi ir \nu z}}{\sin^{k}(\pi(z-i\eta))} 
\right\}.
\end{eqnarray*}
We have to show that
$$
\Sigma(\lambda) := \sum_{l \in \Z}
\sum_{m \in \Z : \sign(\nu(m,\lambda))=1} 
\left|\Res_{z=i\eta} \left\{ e^{-\pi\delta (z+l)^{2}}
\frac{q_{1}(z+l)e^{2\pi ir \nu z}}{\sin^{k}(\pi(z-i\eta))} 
\right\} \right|
$$
is convergent for each $\lambda \in \Lambda$. But this follows 
exactly as in the 2nd part of the proof of 
\reflem{lem:Lemma8.12}. In fact we get by using notation from 
that proof that
$$
\Sigma(\lambda) \leq \sum_{j=0}^{k-1} \sum_{l \in \Z} 
\sum_{m \in \Z : \nu(m,\lambda)) \geq 0} |c_{-1-j}| |a_{j}|,
$$
where
\begin{eqnarray*}
a_{j} &=& \left.\frac{1}{j!} \frac{\dte^{j} 
e^{-\pi\delta(z+l)^{2}}q_{1}(z+l) 
e^{2\pi ir\nu z}}{\dte z^{j}}\right|_{z=i\eta} \\
&=& \sum_{\stackrel{k_{1},k_{2},k_{3} \geq 0}
{k_{1}+k_{2}+k_{3}=j}}
\frac{1}{k_{1}!k_{2}!k_{3}!} P_{k_{1}}(i\eta+l,\delta) 
e^{-\pi\delta(i\eta+l)^{2}} (2\pi ir\nu)^{k_{2}} 
e^{2\pi ir \nu(i\eta +l)} q_{1}^{(k_{3})}(i\eta +l).
\end{eqnarray*}
We therefore get
\begin{eqnarray*}
|a_{j}| &\leq& \sum_{\stackrel{k_{1},k_{2},k_{3} \geq 0}
{k_{1}+k_{2}+k_{3}=j}}
\frac{1}{k_{1}!k_{2}!k_{3}!} |P_{k_{1}}(i\eta+l,\delta)| 
e^{\pi\delta\eta^{2}}e^{-\pi\delta l^{2}} (2\pi r\nu)^{k_{2}} 
e^{-2\pi r \nu \eta} b_{k_{3}}'e^{a_{k_{3}}'(\eta +m_{q}')}
\end{eqnarray*}
showing the convergence of $\Sigma(\lambda)$.
\end{prfa}

\end{document}